\documentclass{amsart}

\usepackage[T1]{fontenc}
\usepackage[utf8]{inputenc}
\usepackage[english]{babel}

\usepackage{amsthm,amssymb,verbatim}
\usepackage{etoolbox}
\usepackage{mathtools}
\usepackage{mathrsfs}
\usepackage[shortlabels]{enumitem}
\usepackage{graphicx}
\usepackage{float}
\usepackage[font=scriptsize]{caption}
\usepackage{subcaption}
\usepackage{xcolor}
\usepackage{tikz}
\usetikzlibrary{arrows.meta,calc}
\usepackage{pgfplots}
\pgfplotsset{compat=1.16}
\usepackage{chngcntr}
\usepackage{apptools}
\usepackage{ragged2e}
\usepackage{amscd}
\usepackage{hyperref}
\usepackage[msc-links,backrefs]{amsrefs}
\usepackage[foot]{amsaddr}
\usepackage{lmodern}  
\usepackage{microtype}
\makeatletter
\g@addto@macro\normalsize{%
	\setlength\abovedisplayskip{5pt plus 2pt minus 3pt}%
	\setlength\belowdisplayskip{5pt plus 2pt minus 3pt}%
	\setlength\abovedisplayshortskip{2pt plus 2pt minus 2pt}%
	\setlength\belowdisplayshortskip{3pt plus 2pt minus 2pt}%
}
\makeatother
\setlist{topsep=3pt,itemsep=2pt,parsep=0pt}
\makeatletter
\def\thm@space@setup{\thm@preskip=7pt plus 2pt minus 3pt
	\thm@postskip=\thm@preskip}
\makeatother

\usepackage{placeins}  
\makeatletter
\let\oldincludegraphics\includegraphics
\renewcommand{\includegraphics}[2][]{%
	\IfFileExists{#2}{\oldincludegraphics[#1]{#2}}{%
		\fbox{\parbox[c][0.28\textwidth][c]{0.68\textwidth}{\centering Missing figure file:\par\texttt{\detokenize{#2}}}}%
	}%
}
\makeatother

\newcommand{\rr}{{\mathbb R}}

\newcommand{\Sp}{{\mathbb S}}

\newcommand{\eps}{\varepsilon}

\newcommand{\set}[1]{\left\{#1\right\}}
\newcommand{\pare}[1]{\left(#1\right)}

\newcommand{\pI}[1]{\left <{#1}\right >}

\newcommand\restr[2]{{%
		\left.\kern-\nulldelimiterspace
		#1
		\vphantom{\big|}
		\right|_{#2}%
}}

\theoremstyle{plain}
\newtheorem{theorem}{Theorem}[section]
\newtheorem{lemma}[theorem]{Lemma}
\newtheorem{proposition}[theorem]{Proposition}
\newtheorem{corollary}[theorem]{Corollary}

\newtheorem{definition}[theorem]{Definition}  

\newtheoremstyle{mainintrothm}
{\topsep}
{\topsep}
{\itshape}
{}
{\bfseries}
{.}
{.5em}
{#1}
\theoremstyle{mainintrothm}
\newtheorem{maintheorema}{Theorem A}

\newtheorem{maintheoremb}{Theorem B}

\newtheorem{maintheoremc}{Theorem C}

\newtheorem{maintheoremd}{Theorem D}

\theoremstyle{definition}

\newtheorem{notation}[theorem]{Notation}
\newtheoremstyle{remarkstyle}
	{\topsep}{\topsep}{\small}{}{\bfseries}{.}{.5em}{}
\theoremstyle{remarkstyle}
\newtheorem{remark}[theorem]{Remark}
\AtEndEnvironment{remark}{\hfill\ensuremath{\lozenge}}
\theoremstyle{definition}
\newtheorem{op}[theorem]{Open Problem}

\title[Rotationally symmetric translating solitons]%
{Rotationally symmetric translating solitons of fully nonlinear extrinsic geometric flows: Classification and Applications}
\author{Jose Torres Santaella}
\keywords{Translating solitons; Fully nonlinear extrinsic geometric flows; Bowl soliton; Catenoidal translators; Asymptotic expansions.}
\makeatletter\@namedef{subjclassname@2020}{\textup{2020} Mathematics Subject Classification}\makeatother
\subjclass[2020]{53A10, 53C21, 53C42, 53E10, 35J60, 58J70}

\begin{document}
	
	\begin{abstract}
	We study rotationally symmetric translating solitons of fully nonlinear extrinsic geometric flows in $\rr^{n+1}$ driven by positive $\alpha$-homogeneous curvature functions, building on Rengaswami's classification of bowl-type translators. For positive ends at infinity and every $\alpha>0$, we derive two correction terms beyond the leading order in the nondegenerate case and the growth rate in each degenerate regime. Imposing the translator equation across the neck on a signed extension of the curvature function, we construct catenoidal translators and classify their continuations according to a branch of the rotational level set. These yield complete embedded translators or maximal pieces, for which we determine the regularity at each transition and the shrinking-neck asymptotics; each alternative is realized in a model family of curvature functions.

	As applications of the tangency principle, we obtain rigidity and nonexistence results governed by the global geometry of the bowl-type solutions. When the bowl is entire, it is unique up to vertical translation among convex translators with a single graphical end asymptotic to it (also under a weighted condition in the strictly convex case), and its existence rules out both convex translating graphs over bounded domains that blow up at the boundary and entire ones growing uniformly faster or slower than the bowl. When the bowl is not entire, it may have finite height for $\alpha>1$; in this setting, any strictly convex translating graph over a ball that is continuous up to the boundary and has vertical tangent planes there is necessarily a bowl of finite height.
\end{abstract}
	
	\maketitle
	
\section{Introduction}\label{sec:introduction}

Extrinsic geometric flows deform an immersed hypersurface of $\rr^{n+1}$ by prescribing its normal velocity in terms of its principal curvatures. The model case is the mean curvature flow, whose normal velocity is the mean curvature $H$, the sum of the principal curvatures. We refer the reader to \cite{mantegazza2011lecture} for an introduction to the mean curvature flow and to \cite{andrews2020extrinsic} for the broader theory of fully nonlinear extrinsic flows.

In the fully nonlinear setting, we fix a symmetric open cone $\Gamma\subset\rr^n$ containing the positive cone $\Gamma_+=\{\lambda\in\rr^n:\lambda_i>0\ \text{for all }i\}$, in which the principal curvatures of the evolving hypersurface are required to lie. We then replace $H$ by a symmetric function $\gamma:\Gamma\to[0,\infty)$ of class $\mathcal C^1$ that is positive, $\alpha$-homogeneous and strictly increasing in each entry (Definition \ref{def:alpha-positive-curvature-function}). More precisely, let $F_0:\Sigma^n\to\rr^{n+1}$ be an oriented immersion. The hypersurface $\Sigma_0=F_0(\Sigma)$ evolves by the $\gamma$-flow if there exists a one-parameter family of immersions $F:\Sigma\times[0,T)\to\rr^{n+1}$ satisfying
\begin{equation}\label{eq:flow}
	\begin{cases}
		\langle \partial_tF,\vec{\nu}\rangle = \gamma(\lambda)
		& \text{in }\Sigma\times(0,T),\\
		F(\cdot,0)=F_0
		& \text{on }\Sigma,
	\end{cases}
\end{equation}
where $\vec\nu$ denotes the unit normal along $\Sigma_t=F(\Sigma,t)$, and $\lambda\in\Gamma$ is the principal curvature vector with respect to $\vec\nu$, computed via the Weingarten operator $\mathcal W(\vec X)=-\nabla_{\vec X}\vec\nu$. With this convention, a round sphere has positive principal curvatures with respect to its inner normal. The strict monotonicity of $\gamma$ makes the flow \eqref{eq:flow} parabolic while the principal curvatures remain in $\Gamma$ \cite[Section 18.4]{andrews2020extrinsic}.

Typical examples of such curvature functions include the elementary symmetric polynomials $S_k(\lambda)=\sum_{1\leq i_1<\cdots<i_k\leq n}\lambda_{i_1}\cdots\lambda_{i_k}$, with $1\leq k\leq n$, defined on the G\aa rding cones $\Gamma_k=\{\lambda\in\rr^{n}:S_j(\lambda)>0\text{ for }j=1,\ldots,k\}$, where $\Gamma_n=\Gamma_+$. This family contains the mean curvature $H=S_1$ and the Gauss--Kronecker curvature $K=S_n$, the two classical examples, whose study draws on the theories of quasilinear equations and of equations of Monge--Amp\`ere type, respectively. Other standard choices are the Hessian quotients $S_k/S_j$ with $1\leq j<k\leq n$, also defined on $\Gamma_k$, which we normalize as $\mathcal Q_{k,j}$ in Subsection \ref{subsec:model-examples-branches}, and the positive powers of these functions, such as the powers of $K$, whose flows are studied, for instance, in \cite{andrews2016flow}.

The main objects of study in this article are translating solitons, or $\gamma$-translators for short. These are solutions to \eqref{eq:flow} that evolve by translation with a constant velocity, which, up to a rotation and, by homogeneity, a dilation of the ambient space, we may assume to be $\vec e_{n+1}$. Up to tangential diffeomorphisms, the evolution is given by $F(\cdot,t)=F_0+t\,\vec e_{n+1}$, and each time slice satisfies the translator equation
\begin{equation}\label{eq:translator}
\gamma(\lambda)=\langle\vec\nu,\vec e_{n+1}\rangle .
\end{equation}

Complete $\gamma$-translators without boundary are noncompact \cite[Theorem 1.1]{Yo}, and they play a central role in the analysis of the flow for three reasons. First, they are eternal solutions, defined for all times $t\in(-\infty,\infty)$. Second, under suitable assumptions on the curvature function, they arise as blow-up limits of type II singularities, and classification results such as the rotational symmetry of two-convex limits have been extended to a class of fully nonlinear translator equations \cite[Theorem 1.3]{bourni-langford2016typeII}. Third, they enter the classification of ancient solutions, those defined for all times $t\in(-\infty,T)$ with $T\leq\infty$ \cite{wang2011convex,langford-lynch2020sharp}, where they serve as barriers in \cite{lynch2022uniqueness}.

In this article we develop a structural theory for the rotationally symmetric $\gamma$-translators, in which the mean curvature and the Gauss--Kronecker curvature lie on the two sides of a single dichotomy. The theory has three parts. The first determines the asymptotics of positive ends for every homogeneity $\alpha>0$, in particular those of bowl-type translators, under regularity and decay hypotheses on the positive branch. The second constructs and classifies the catenoidal $\gamma$-translators. The third proves the uniqueness of the bowl within its asymptotic class and uses an entire bowl as a barrier to exclude convex translating graphs over bounded domains that blow up at the boundary.

The foundation for this article is Rengaswami's classification of rotational bowl-type translators for general fully nonlinear curvature functions \cite[Theorems 1.1--1.4]{rengaswami2024classification}. Our approach uses equations and methods from Rengaswami's work and from the rotational analysis of Clutterbuck--Schn\"urer--Schulze for $\gamma=H$ \cite{clutterbuck2007stability}, in particular the reformulation of the rotational translator equation as an implicit positive-branch problem, combined with barrier arguments.
Part of this program was initiated in the unpublished preprint \cite{RengaswamiTorresSantaella2023Bowl} with S. Rengaswami, as described in Subsection \ref{subsec:intro-previous}.

\subsection{Background and classical examples}\label{subsec:intro-background}
Before stating our results, we recall the known rotationally symmetric translators. In dimension one, they are the complete convex translating curves of the $\kappa^\alpha$-flow, where $\kappa$ is the curvature of the curve. They are entire graphs when $\alpha\leq\frac{1}{2}$ and graphs over bounded intervals when $\alpha>\frac{1}{2}$ \cite[Proposition 2.4]{choi2021convergence}. For $\alpha=1$ the curve is the grim reaper $u(x)=\log\sec x$ on $\left(-\frac{\pi}{2},\frac{\pi}{2}\right)$. If $e_1=(1,0,\ldots,0)\in\Gamma$, the product with $\rr^{n-1}$ of a dilation of the curve with the same exponent $\alpha$ is a $\gamma$-translator, the grim reaper cylinder of $\gamma$ (Remark \ref{rem:langford-rengaswami-condition}).

In dimension $n\geq2$, the first examples come from the mean curvature flow, whose entire rotationally symmetric translating graph is the bowl soliton. It was obtained by Altschuler--Wu in $\rr^3$, who showed that it is asymptotic to a paraboloid \cite{altschuler1994translating}, and Clutterbuck--Schn\"urer--Schulze extended their argument to every dimension and refined this expansion \cite[Lemma 2.2]{clutterbuck2007stability}. For fully nonlinear curvature functions, the most studied translators in the literature are those of the powers $K^p$ of the Gauss--Kronecker curvature. Urbas proved that they exhibit a dichotomy in the power $p$ \cite[Theorems 2 and 3 and Section 4]{urbas1998complete}. More precisely, he proved that the complete convex $K^p$-translators lie over bounded convex domains when $p>\frac{1}{2}$, whereas for $p\leq\frac{1}{2}$ these domains are unbounded, and in that case he constructed an entire rotationally symmetric translating graph, now called the bowl-type solution.

Rengaswami formalized the rotational theory of translating graphs for general curvature functions \cite[Section 2]{rengaswami2024classification}, and we follow his formulation to state our results. A rotational hypersurface has the principal curvature vector $(x,y,\ldots,y)$, where $x$ is the meridional curvature and $y$, of multiplicity $n-1$, is the rotational one, so $\gamma$ enters the translator equation \eqref{eq:translator} only through its rotational restriction $\widetilde\gamma(x,y)=\gamma(x,y,\ldots,y)$. Since $\widetilde\gamma$ is $\alpha$-homogeneous, wherever the angle function $z=\pI{\vec{\nu},\vec{e}_{n+1}}$ is positive, the equation $\widetilde\gamma(x,y)=z$ is equivalent to $\widetilde\gamma\bigl(z^{-\frac{1}{\alpha}}x,z^{-\frac{1}{\alpha}}y\bigr)=1$, so the rescaled pair lies on the level set $\widetilde\gamma=1$. Rengaswami solved this level set for the meridional curvature by the implicit function theorem \cite[Definition 6]{rengaswami2024classification}, obtaining the positive branch $x=g_+(y)$ for $y$ in an interval $I_+$ (Proposition \ref{prop:positive-branch}). On vertical graphs $u(r)$, this gives a first-order equation for the slope $v=u'$, whose solutions we call \emph{positive ends}. 

The bowl-type solutions are the rotational graphs $B$ of functions $u_B$ with $u_B'(0)=0$ on a maximal interval $[0,R_\gamma)$, and away from the axis their slopes are positive ends. They exist, are strictly convex and are unique among rotational graphs up to vertical translation \cite[Theorem 1.1 and Sections 3, 4 and 8]{rengaswami2024classification}. Rengaswami also gave criteria to decide whether the bowl is \emph{entire} ($R_\gamma=\infty$) or defined over a ball ($R_\gamma<\infty$) \cite[Theorems 1.1--1.4]{rengaswami2024classification}, and in the second case the bowl is called \emph{cylindrical} when $u_B(r)\to\infty$ as $r\to R_\gamma^-$ \cite[Definition 5]{rengaswami2024classification}.

In Rengaswami's terminology, $\gamma$ is degenerate or nondegenerate according to whether $\gamma(0,1,\ldots,1)$ vanishes, where this value is defined as $\lim_{s\to0^+}\gamma(s,1,\ldots,1)$ when $(0,1,\ldots,1)\notin\Gamma$ \cite[Definition 3]{rengaswami2024classification}. This means that the round cylinder $\rr\times\Sp^{n-1}$ is stationary in the first case and collapses to its axis in the second. In our setting, $\gamma$ is \emph{normalized $1$-nondegenerate} when $(0,1)\in\widetilde\Gamma$ and $\widetilde\gamma(0,1)=1$, and then $I_+$ is bounded with supremum $y=1$ and the expansion of $g_+$ near $y=1$ determines the coefficients in Theorem \ref{thm:main-a}. In contrast, $\gamma$ is \emph{$1$-degenerate} when $(0,1)\in\partial\widetilde\Gamma$ and $\widetilde\gamma(t,1)\to0$ as $t\to0^+$, and then $I_+$ is unbounded and the decay of $g_+$ at infinity determines the growth of the end (Definition \ref{def:one-nondegeneracy}). The mean curvature and its powers are $1$-nondegenerate, whereas the Gauss--Kronecker curvature and its roots are $1$-degenerate. Finally, Remark \ref{rem:langford-rengaswami-condition} compares this notion with the nondegeneracy condition of Langford--Rengaswami \cite{rengaswami-langford2023pancake}.

Catenoidal translators, one of the topics of this article, go back to Clutterbuck--Schn\"urer--Schulze, who constructed them for the mean curvature flow and called them \emph{winglike} \cite[Subsection 2.2]{clutterbuck2007stability}. For this flow, they are complete embedded rotationally symmetric translators with the topology of an annulus $\Sp^{n-1}\times\rr$, formed by two graphs over the exterior of a ball $B_R\subset\rr^n$ that meet with vertical tangent planes along a horizontal $(n-1)$-sphere of radius $R$, called the \emph{neck}. Analogous solutions have also been constructed for the powers $H^p$ with $p>0$ by Sheng--Wu \cite[Corollary 5.1]{sheng2010rotationally}, and, when $k$ is odd, for the elementary symmetric polynomials $S_k$ in $\rr^n\times\rr$ and $\mathbb H^n\times\rr$ by de Lima--Pipoli \cite[Theorem 1(ii)]{lima2025translators} and for $S_k$ and $S_k^{1/k}$ in warped products by de Lira--de Farias \cite[Theorem 4(c)]{lira2026translating}.

We notice that, for a general curvature function, the construction of catenoidal translators requires a signed extension of $\widetilde\gamma$ beyond $\widetilde\Gamma$. Indeed, the angle function vanishes at the neck, and so the translator equation requires $\widetilde\gamma(x,y)=0$ there, which cannot hold in $\widetilde\Gamma$, where $\widetilde\gamma$ is positive. In addition, the angle function is negative on the lower side of the generating curve, below the neck, so the extension must take negative values there. Since reversing the unit normal changes the sign of both the principal curvatures and the angle function, the translator equation is independent of the choice of unit normal whenever the extension is odd, as happens for $S_k$ with $k$ odd, which includes $H$ \cite[Remark 2]{lima2025translators}, \cite[Remark 2]{lira2026translating}. For even $k$, de Lima--Pipoli obtain catenoid-type translators with boundary \cite[Theorem 1(iii)]{lima2025translators}. Accordingly, we require the extension to be odd on the sheets selected at the neck (Definition \ref{def:selected-signed-neck-datum}).

\subsection{Results}\label{subsec:intro-results}
Our first result describes positive ends at infinity, and in particular the bowl-type solutions.

\begin{maintheorema}
	\label{thm:main-a}
	Let $\alpha>0$ and let $v=u'$ be a positive end, the slope of a rotational graphical $\gamma$-translator. In the $1$-degenerate case, assume that $g_+(y)=c_\gamma y^{-k_\gamma}(1+o(1))$ as $y\to\infty$, for some constants $c_\gamma>0$ and $k_\gamma>0$.
	\begin{enumerate}[label=\textup{(\roman*)}]
		\item \emph{Normalized $1$-nondegenerate case.} If $\widetilde\gamma$ is of class $\mathcal C^3$ near $(0,1)$, then
		\begin{equation*}
			v(r)
			=
			r^\alpha
			-
			\frac{a_\gamma}{r^\alpha}
			+
			\frac{b_\gamma}{r^{3\alpha}}
			+
			o(r^{-4\alpha}) \quad \text{as } r\to\infty,
		\end{equation*}
		where the constants $a_\gamma$ and $b_\gamma$ depend only on $\alpha$ and on the asymptotic expansion of $g_+$ at the endpoint $y=1$.
		
		\item \emph{Power-growth regime.} If $k_\gamma>2\alpha-1$, then
		\begin{equation*}
			v(r)
			=
			A_\gamma r^{d_\gamma}(1+o(1)) \quad \text{as } r\to\infty,
		\end{equation*}
		where $d_\gamma=\frac{\alpha(k_\gamma+1)}{k_\gamma+1-2\alpha}$ and $A_\gamma>0$ satisfies $A_\gamma^{\frac{k_\gamma+1-2\alpha}{\alpha}}=\frac{c_\gamma}{d_\gamma}$.
		
		\item \emph{Critical regime.} If $k_\gamma=2\alpha-1$ (so that $\alpha>\frac{1}{2}$), then
		\begin{equation*}
			\log v(r)
			=
			\frac{c_\gamma}{2\alpha}r^{2\alpha}
			+
			o(r^{2\alpha}) \quad \text{as } r\to\infty.
		\end{equation*}
	\end{enumerate}
\end{maintheorema}

The full statement is Theorem \ref{thm:positive-bowl-asymptotics}, which adds the second-order correction in the power-growth regime and is followed by the literature on each term. The leading term $r^\alpha$ is due to Rengaswami \cite[Theorem 1.1 and Proposition 7.2]{rengaswami2024classification}, and the coefficient $b_\gamma$ and the asymptotics in the degenerate regimes are new. Under the hypotheses of each item, every positive end is defined up to $r=\infty$ (Propositions \ref{prop:positive-end-maximal-interval} and \ref{prop:positive-end-at-infinity}). For $\alpha\leq\frac{1}{2}$, every $k_\gamma>0$ lies in the power-growth regime, in agreement with \cite[Theorem 1.2]{rengaswami2024classification}. In contrast, for $k_\gamma<2\alpha-1$, which requires $\alpha>\frac{1}{2}$, every positive end ends at a finite radius, so the bowl-type solution is not entire in that range (Lemma \ref{lem:subcritical-degenerate-blow-up}).

Mart\'in--Savas-Halilaj--Smoczyk proved that the bowl soliton is the only complete embedded translator of the mean curvature flow with a single end smoothly asymptotic to the translating paraboloid \cite[Theorem A]{Paco_2014}. Our second result is its counterpart for an entire bowl and for a bowl of finite height. For this theorem and for Theorem \ref{thm:main-d} we require the curvature function to be smooth on $\Gamma_+$, as the tangency principle of \cite[Theorem 1.2]{Yo} does, and we assume that the principal curvature vector of a $\gamma$-translator lies in $\Gamma$, where the translator equation is elliptic. The latter holds automatically for graphs when $\gamma$ extends continuously by zero to $\partial\Gamma$.

\begin{maintheoremb}
	\label{thm:main-b}
	Let $B$ be the bowl-type solution for the $\gamma$-flow. In parts \textup{(a)} and \textup{(b)}, assume that $B$ is entire, and let $M$ be a connected, properly immersed $\gamma$-translator without boundary. Assume that, for some $R_0>0$, the part of $M$ over $\{|x|<R_0\}$ is bounded and the part over $\{|x|\geq R_0\}$ is contained in the graph of a continuous function. Then $M$ is the graph of a function $u\in\mathcal C^2(\rr^n)$. Let $\varphi(x)=u(x)-u_B(|x|)-c$ for some $c\in\rr$.
	\begin{enumerate}[label=\textup{(\alph*)}]
		\item If $M$ is convex and $\varphi(x)\to0$ as $|x|\to\infty$, then $M=B+c\,\vec e_{n+1}$.
		\item Assume that one of the following conditions holds.
		\begin{enumerate}[label=\textup{(\roman*)}]
			\item $\gamma$ is normalized $1$-nondegenerate and $\alpha\geq1$.
			\item The positive end of $B$ lies in the power-growth regime of Theorem \ref{thm:main-a}\textup{(ii)}, with $d_\gamma\geq1$.
			\item The positive end of $B$ lies in the critical regime of Theorem \ref{thm:main-a}\textup{(iii)}.
		\end{enumerate}
		If $M$ is strictly convex and $|\varphi(x)|+|D\varphi(x)|=o\bigl(|x|^{-1}u_B'(|x|)\bigr)$ as $|x|\to\infty$, with $|x|^{-1}$ replaced by $|x|^{1-2\alpha}$ in case \textup{(iii)} when $\alpha>1$, then $M$ is a vertical translate of $B$.
		\item Let $M$ be a strictly convex $\gamma$-translator given by the graph of a function $u\in\mathcal C^2(B_R)\cap\mathcal C^0(\overline B_R)$ over a ball $B_R\subset\rr^n$ of radius $R$, and assume that $|Du(x)|\to\infty$ as $x\to\partial B_R$. Then, $B$ is defined over a ball of radius $R$ and has finite height. Furthermore, $M$ is a translate of $B$.
	\end{enumerate}
\end{maintheoremb}

We prove parts \textup{(a)} and \textup{(b)} of Theorem \ref{thm:main-b} in Theorem \ref{thm:bowl-asymptotic-uniqueness}, and part \textup{(c)} in Theorem \ref{thm:finite-height-bowl-uniqueness}. As in \cite[Theorem A]{Paco_2014}, part \textup{(b)} follows from the Alexandrov reflection method, where the growth of the bowl in Theorem \ref{thm:main-a} gives the separation at infinity, while part \textup{(a)} follows by sliding $B$ vertically. In particular, the asymptotic condition of part \textup{(b)} allows $\varphi$ to be unbounded when $u_B'$ grows faster than linearly. In our setting, the embeddedness assumed by Mart\'in--Savas-Halilaj--Smoczyk is replaced by properness, since $\lambda\in\Gamma$ gives a positive angle function $\gamma(\lambda)$, so that the vertical projection of $M$ is a proper local diffeomorphism onto $\rr^n$, and hence a diffeomorphism (Lemma \ref{lem:one-end-entire-graph}). For cylindrical bowls, \cite[Theorem 1.4]{Yo} proves the rotational symmetry of strictly convex graphs over a ball asymptotic to a vertical cylinder, and hence their uniqueness. Part \textup{(c)} is the counterpart of this result for bowls of finite height, which occur for $\bigl(\frac{S_n}{S_{n-1}}\bigr)^\alpha$ with $\alpha>1$ and for the powers $K^p$ with $p>1$ (Proposition \ref{prop:finite-height-ends}). Its proof is again a reflection argument, where the infinite slope at the boundary of the ball keeps the reflected boundary above the graph. The ball cannot be replaced by another convex domain, since for $K^p$ with $p>1$ the translators of Urbas over nonround ellipsoids satisfy the remaining hypotheses of part \textup{(c)} but are not rotationally symmetric \cite[Theorem 2 and Section 4]{urbas1998complete} (Remark \ref{rem:finite-height-bowl}).

Without an asymptotic condition, Wang proved for the mean curvature flow that the bowl soliton is the only entire convex translating graph in $\rr^3$ \cite[Theorem 1.1]{wang2011convex}, whereas this rigidity fails in higher dimensions \cite[Theorem 1.2]{wang2011convex}. The bowl is then characterized under curvature assumptions, such as the noncollapsing and uniform two-convexity of Haslhofer \cite[Theorem 1.2]{haslhofer2015uniqueness}, or the cylindrical and gradient estimates under which Bourni--Langford proved rotational symmetry for $n\geq3$ \cite[Theorem 1.1]{bourni-langford2016typeII}. We collect the corresponding questions for $\gamma$-translators in Open Problem \ref{op:bowl-entire}.

The catenoidal $\gamma$-translators are generated by a curve with a neck of radius $R>0$, a strict minimum of the distance to the axis $\rr\vec e_{n+1}$. For them, the translator equation \eqref{eq:translator} is written with a signed rotational level profile $\widehat\gamma$ extending $\widetilde\gamma$ (Definition \ref{def:signed-rotational-extension}), whose zero level is the neck.

The upper side of the curve enters the positive branch at a finite radius and continues as a positive end. The lower side, after a reorientation, follows a reflected branch $(I_-,g_-)$ of the level set of $\widehat\gamma$, uniquely determined by the angular neck data (Definitions \ref{def:selected-signed-neck-datum} and \ref{def:neck-accessible-branch}). The behavior of $g_-$ at the right endpoint of $I_-$ sorts the lower sides into four types (Definition \ref{def:lower-continuation-outcomes}).
\begin{itemize}
	\item The branch is axis-crossing if $I_-$ extends to $y=0$ and $g_-$ has a positive limit there.
	\item It has finite-angle zero output if $I_-$ extends to $y=0$ and $g_-$ vanishes there with $y^{-1}g_-(y)$ converging to a finite limit.
	\item It has infinite-angle zero output if $g_-$ vanishes at $y=0$ while $y^{-1}g_-(y)\to-\infty$, and fractional zero output when in addition $g_-$ behaves like a fractional power of $|y|$.
	\item It is regularity-limited if $I_-$ ends strictly before $y=0$.
\end{itemize} A regularity-limited branch may extend along its component of the level set up to $y=0$, either as a curvature-blow-up continuation, where the meridional value diverges at $y=0$, or as a finite-angle continuation, with a finite angular limit (Definition \ref{def:lower-continuation}). Moreover, the regular and the singular positive-matching conditions of Definitions \ref{def:regular-positive-matching} and \ref{def:singular-positive-matching}, for axis-crossing branches and for curvature-blow-up continuations, respectively, allow the lower side to reach a second positive end.

\begin{maintheoremc}
	\label{thm:main-c}
	Let $\gamma$ be a normalized $1$-nondegenerate $\alpha$-positive curvature function with $\alpha>\frac{1}{2}$, equipped with selected angular neck data satisfying conditions \ref{cat:all-radii} through \ref{cat:lower-branch} of Section \ref{sec:catenoidal-translators}, and let $(I_-,g_-)$ be the maximal regular continuation of the neck-accessible branch. For every $R>0$, there is a unique catenoidal neck of radius $R$ on the selected sheet. This neck generates a rotational hypersurface $W_R$ whose upper side becomes a complete positive end beyond a finite entry radius, and whose lower side continues according to the specific type of $(I_-,g_-)$:
	\begin{enumerate}[label=\textup{(\roman*)}]
		\item \emph{Axis-crossing with regular positive matching.} The lower side turns at a bottom circle and rises into a second complete positive end. Thus, $W_R$ is a complete embedded two-ended catenoidal translator.
		
		\item \emph{Finite-angle zero output, or regularity-limited with a finite-angle continuation, with angular limit $\ell<0$.} The lower side never turns, and its slope is $-r^{\ell+o(1)}$ as $r\to\infty$. Consequently, $W_R$ is a complete embedded catenoidal translator with a single positive end.
		
		\item \emph{Fractional zero output.} The maximal continuation of the lower side ends at a finite radius on a circle of horizontal tangency, beyond which we make no further assertions.
		
		\item \emph{Regularity-limited.} The selected continuation ends at a finite radius. However, under a curvature-blow-up continuation of order $q>0$ and the singular positive-matching condition, the lower side extends past this radius and crosses a subsequent circle of radius $r_R^{\mathrm{sing}}$. Along this circle, the two pieces match with $\mathcal C^{1,1/(q+1)}$ regularity, sharing a continuous unit normal despite a divergent meridional curvature. The surface then rises into a second complete positive end, and $W_R$ is a complete embedded two-ended catenoidal translator.
	\end{enumerate}
\end{maintheoremc}

Theorem \ref{thm:main-c} summarizes Theorem \ref{thm:rotational-translators}, which also gives the regularity at the neck circle, the asymptotic expansion of the positive ends, and the dependence on $R$. For an axis-crossing branch without regular positive matching, and for infinite-angle zero output that is not fractional, Theorem \ref{thm:rotational-translators} gives only the maximal continuation along the selected branch.

\begin{figure}[tb]
	\centering
	\captionsetup[sub]{labelfont=up,justification=raggedright,singlelinecheck=false}
	\begin{subfigure}[t]{0.32\textwidth}
		\centering
		\begin{tikzpicture}\node[anchor=south west,inner sep=0] (im) {\includegraphics[width=\textwidth,trim=156bp 3.4bp 6.72bp 4.2bp,clip]{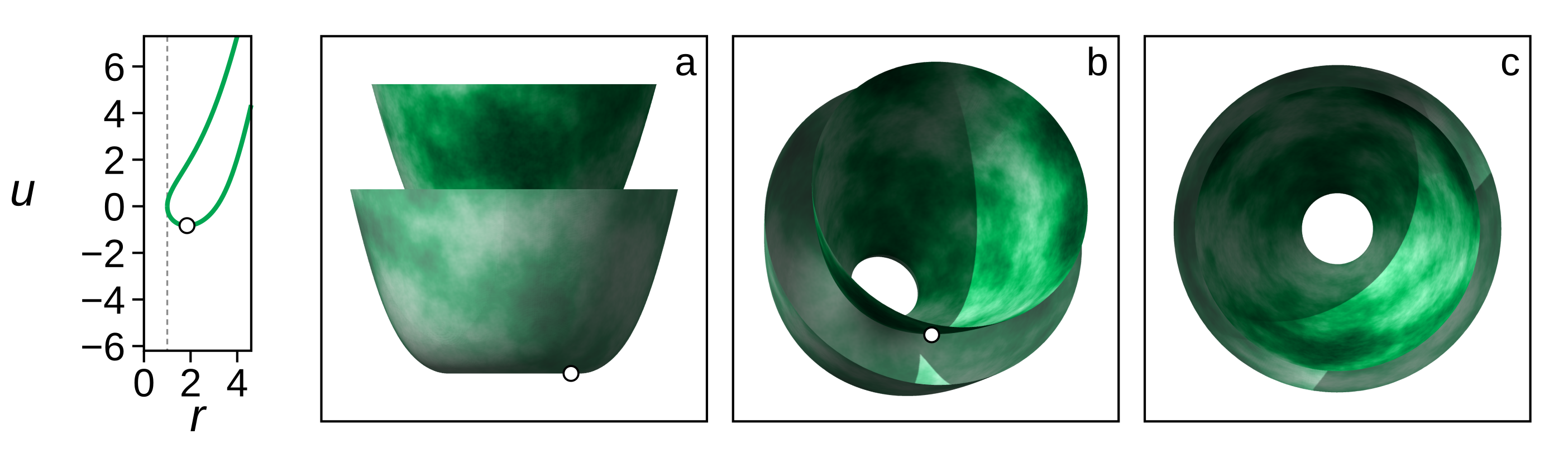}};\begin{scope}[x={(im.south east)},y={(im.north west)}]\fill[white] (0.4236,0.8416) rectangle (0.4821,0.9472);\fill[white] (0.9311,0.8416) rectangle (0.9897,0.9472);\end{scope}\end{tikzpicture}
		\caption{\scriptsize $\frac{H}{n-1}$, $n=2$, $R=1$: case \textup{(i)}.}
		\label{fig:catenoid-axis-crossing}
	\end{subfigure}\hfill
	\begin{subfigure}[t]{0.32\textwidth}
		\centering
		\begin{tikzpicture}\node[anchor=south west,inner sep=0] (im) {\includegraphics[width=\textwidth,trim=156bp 3.4bp 6.72bp 4.2bp,clip]{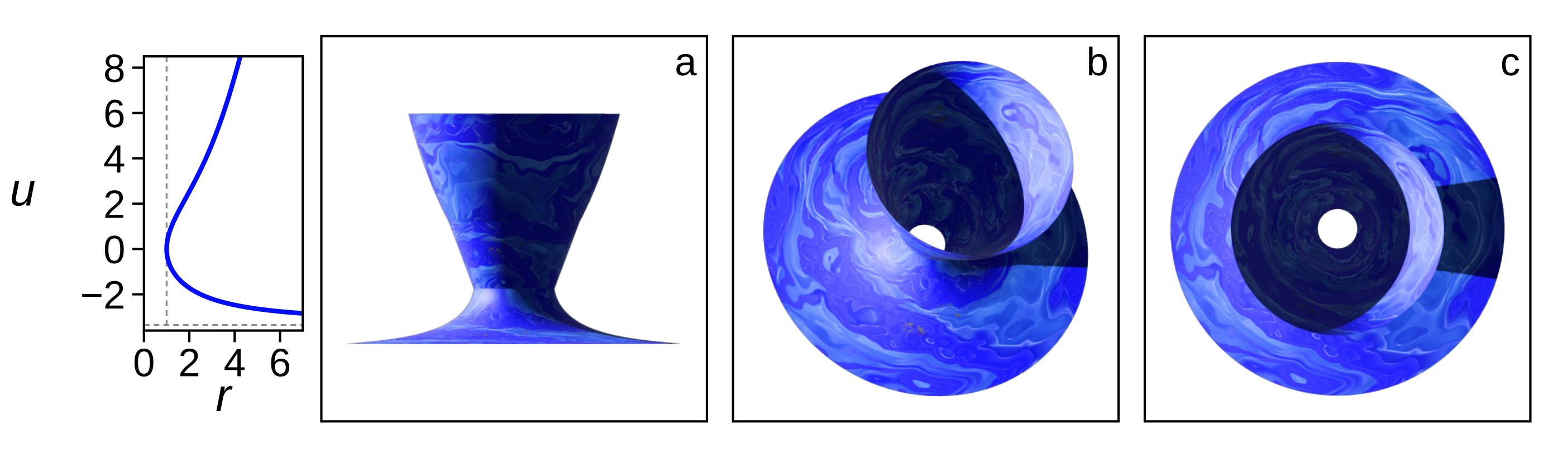}};\begin{scope}[x={(im.south east)},y={(im.north west)}]\fill[white] (0.4236,0.8416) rectangle (0.4821,0.9472);\fill[white] (0.9311,0.8416) rectangle (0.9897,0.9472);\end{scope}\end{tikzpicture}
		\caption{\scriptsize $\mathcal Q_{2,1}$, $n=3$, $R=1$: case \textup{(ii)}, $\ell=-2$.}
		\label{fig:catenoid-flat}
	\end{subfigure}\hfill
	\begin{subfigure}[t]{0.32\textwidth}
		\centering
		\begin{tikzpicture}\node[anchor=south west,inner sep=0] (im) {\includegraphics[width=\textwidth,trim=156bp 3.4bp 6.72bp 4.2bp,clip]{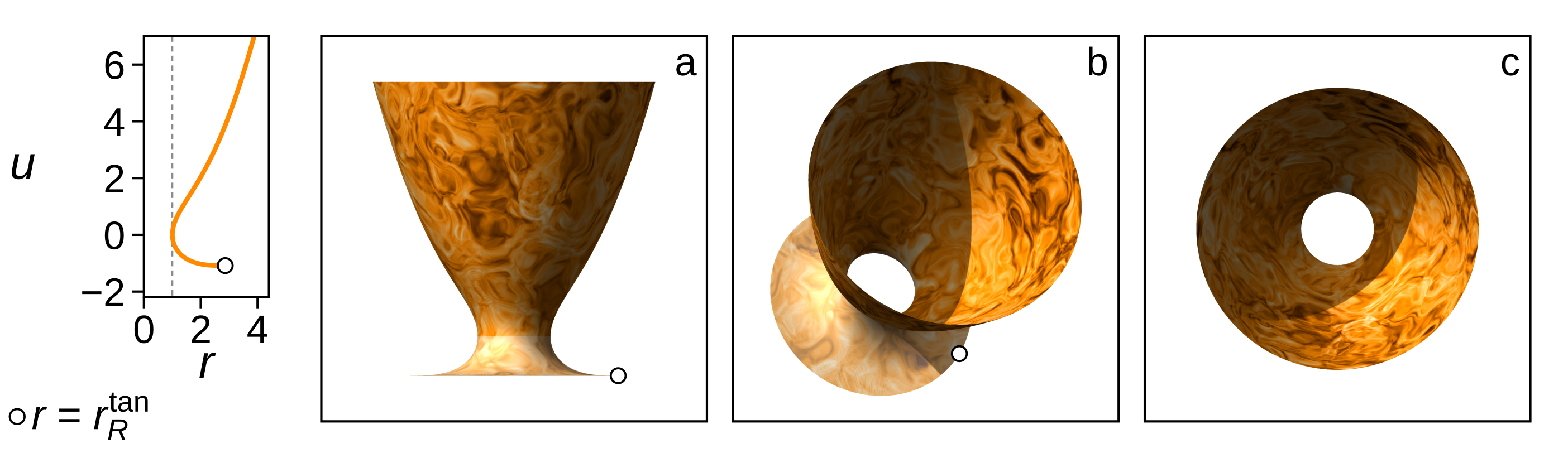}};\begin{scope}[x={(im.south east)},y={(im.north west)}]\fill[white] (0.4236,0.8416) rectangle (0.4821,0.9472);\fill[white] (0.9311,0.8416) rectangle (0.9897,0.9472);\end{scope}\end{tikzpicture}
		\caption{\scriptsize $\gamma_\eps$, $\eps=1$, $n=3$, $R=1$: case \textup{(iii)}.}
		\label{fig:catenoid-fractional}
	\end{subfigure}

	\medskip

	\hspace*{\fill}	\begin{subfigure}[t]{0.32\textwidth}
		\centering
		\begin{tikzpicture}\node[anchor=south west,inner sep=0] (im) {\includegraphics[width=\textwidth,trim=156bp 3.4bp 6.72bp 4.2bp,clip]{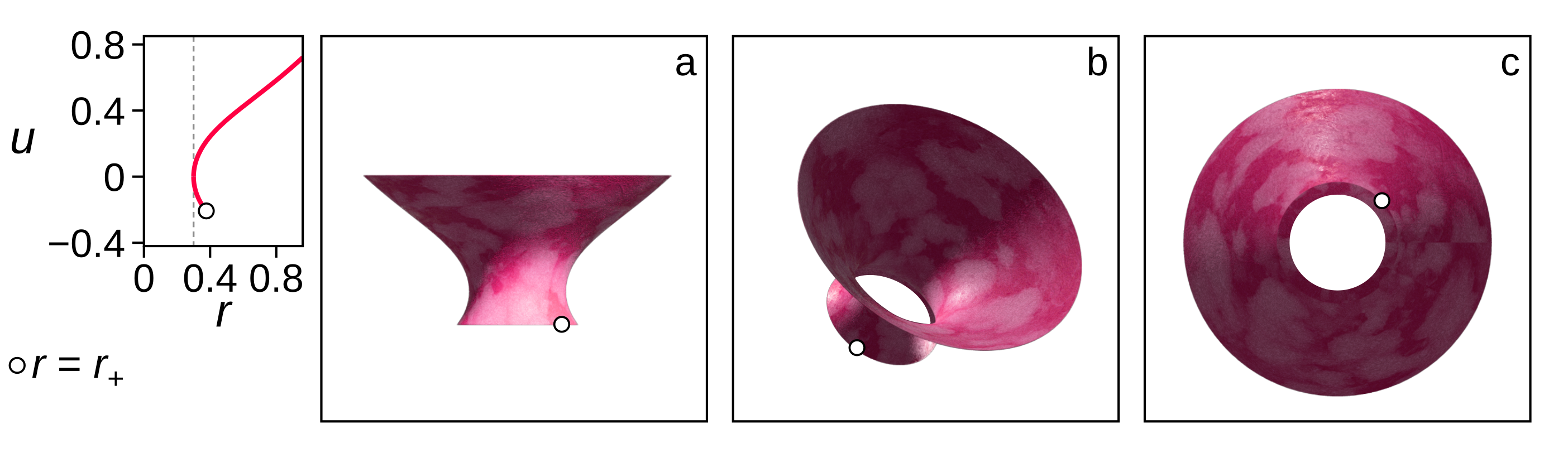}};\begin{scope}[x={(im.south east)},y={(im.north west)}]\fill[white] (0.4236,0.8416) rectangle (0.4821,0.9472);\fill[white] (0.9311,0.8416) rectangle (0.9897,0.9472);\end{scope}\end{tikzpicture}
		\caption{\scriptsize $\gamma_{\mathrm{tri}}$, $n=3$, $R=\frac{3}{10}$: case \textup{(iv)}, no continuation.}
		\label{fig:catenoid-ellipticity-limited}
	\end{subfigure}\hspace{0.04\textwidth}
	\begin{subfigure}[t]{0.32\textwidth}
		\centering
		\begin{tikzpicture}\node[anchor=south west,inner sep=0] (im) {\includegraphics[width=\textwidth,trim=156bp 3.4bp 6.72bp 4.2bp,clip]{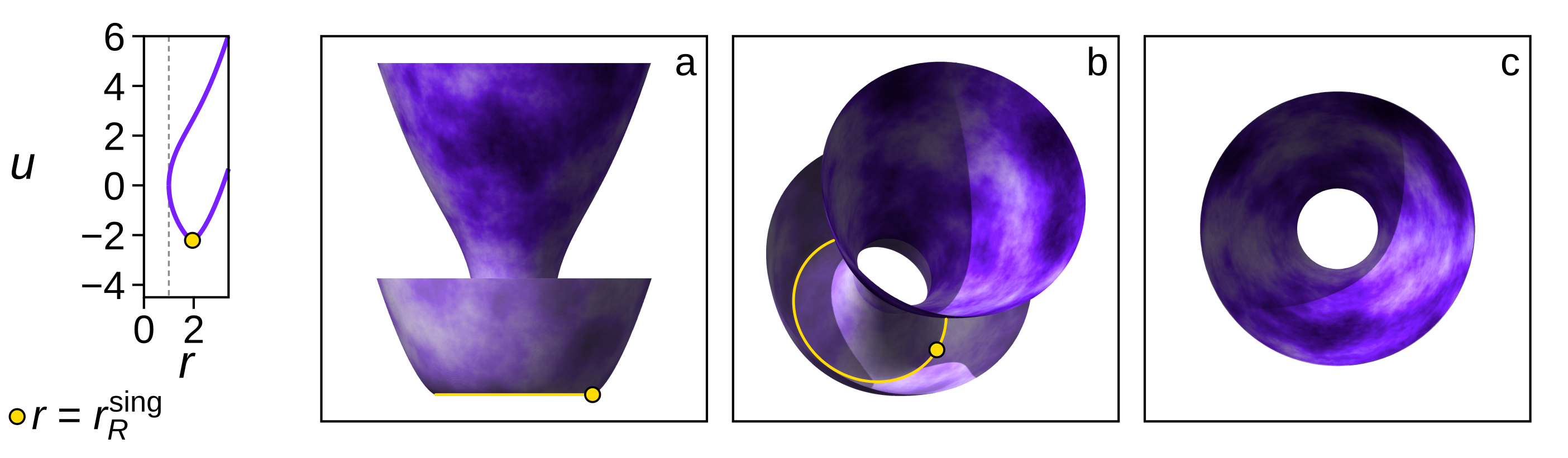}};\begin{scope}[x={(im.south east)},y={(im.north west)}]\fill[white] (0.4236,0.8416) rectangle (0.4821,0.9472);\fill[white] (0.9311,0.8416) rectangle (0.9897,0.9472);\end{scope}\end{tikzpicture}
		\caption{\scriptsize normalized $S_3^{1/3}$, $n=4$, $R=1$: case \textup{(iv)}, $q=2$.}
		\label{fig:catenoid-singular}
	\end{subfigure}\hspace*{\fill}
	\caption{Rotational models of catenoidal translators for the cases of Theorem \ref{thm:main-c}, computed from the graphical equations of Section \ref{sec:barrier-method} for the models of Subsection \ref{subsec:model-examples-branches}. Each panel shows an oblique view and a view along the axis. The dots mark the circle where the lower side has a horizontal tangent plane in \textup{(a)}, and the circles $r=r_R^{\mathrm{tan}}$ in \textup{(c)}, $r=r_+$, where the maximal piece ends, in \textup{(d)}, and $r=r_R^{\mathrm{sing}}$, the singular circle, highlighted in yellow, in \textup{(e)}. The generating curves are shown in Figure \ref{fig:catenoid-profiles}. Images courtesy of Ignacio McManus.}
	\label{fig:catenoid-cases}
\end{figure}

The model families of Subsection \ref{subsec:model-examples-branches} realize every alternative of Theorem \ref{thm:main-c} (Figure \ref{fig:catenoid-cases}). For the mean curvature and the other curvature functions defined on the half-space $\{H>0\}$, the lower side returns to the cone once the normal is reversed. For curvature functions on a G{\aa}rding cone $\Gamma_k$ with $k\geq2$, such as $S_k$ and $S_k^{1/k}$, it lies outside the cone for either normal, so it is a translator only for $\widehat\gamma$ (Lemma \ref{lem:lower-side-outside-garding}). The classification recovers the catenoids recalled in Subsection \ref{subsec:intro-background}, those of Clutterbuck--Schn\"urer--Schulze \cite[Subsection 2.2]{clutterbuck2007stability} and Sheng--Wu \cite[Corollary 5.1]{sheng2010rotationally} in case \textup{(i)}, and those of de Lima--Pipoli \cite[Theorem 1(ii)]{lima2025translators} and de Lira--de Farias \cite[Theorem 4(c)]{lira2026translating} in case \textup{(iv)} for $k$ odd.

For the mean curvature flow, Shahriyari in $\rr^3$ and M{\o}ller in every dimension excluded complete translating graphs from vertical cylinders over bounded sets \cite{shahriyari2015translating,moller2014nonexistence}, and P\'erez-Garc\'ia ruled out noncompact embedded translators with compact boundary from any round cylinder \cite[Theorem 2.2]{perez-garcia2016translating}. Rotational translators are natural barriers for results of this kind, and in \cite[Corollary 1.3]{Yo} we used the bowl to prove the nonexistence of entire convex $\gamma$-translators when the bowl is defined only on a ball. When the bowl is entire, we obtain the following counterpart.

\begin{maintheoremd}
	\label{thm:main-d}
	Let $\gamma$ be an $\alpha$-positive curvature function whose bowl-type solution is entire. Then there exists no convex $\gamma$-translator given by the graph of a function $u\in\mathcal C^2(\Omega)$ over a bounded domain $\Omega\subset\rr^n$ satisfying $u(x)\to+\infty$ as $x\to\partial\Omega$.
\end{maintheoremd}

Theorem \ref{thm:main-d} corresponds to Theorem \ref{thm:vertical-cylinder-graphs}. Together with Corollary \ref{cor:all-round-cylinders}, it is in the spirit of the results of Shahriyari and P\'erez-Garc\'ia for the mean curvature flow, with convexity in place of the ellipticity of the equation of $H$ at every curvature vector (Remark \ref{rem:vertical-cylinders-literature}). By \cite[Theorems 1.1, 1.2 and 1.4]{rengaswami2024classification}, this result applies when $\gamma$ is $1$-nondegenerate, when $\alpha\leq\frac{1}{2}$, and when $\gamma$ is $1$-degenerate with $\alpha>\frac{1}{2}$ and $g_+(y)=O\bigl(y^{-(2\alpha-1)}\bigr)$ as $y\to\infty$, since in these cases the bowl is entire, an assumption that cannot be dropped (Remark \ref{rem:vertical-cylinders-literature}).

The same idea of proof allows us to conclude that an entire convex $\gamma$-translating graph cannot grow uniformly faster or slower than the bowl-type solution at infinity (Corollary \ref{cor:growth-entire-graphs}), as Mart\'in--Savas-Halilaj--Smoczyk observed for the mean curvature flow \cite[Remark 3.1]{Paco_2014}.

The homogeneity enters these results in different ways. Theorem \ref{thm:main-a} holds for every $\alpha>0$, and the threshold $\alpha>\frac{1}{2}$ of Rengaswami's classification \cite[Theorems 1.2--1.4]{rengaswami2024classification} appears in it only through the critical regime. Theorem \ref{thm:main-c} assumes $\alpha>\frac{1}{2}$, whereas parts \textup{(a)} and \textup{(b)} of Theorem \ref{thm:main-b} and Theorem \ref{thm:main-d} require no bound on $\alpha$ beyond the bowl being entire. In contrast, the hypotheses of part \textup{(c)} of Theorem \ref{thm:main-b} force $\alpha>1$ (Remark \ref{rem:finite-height-bowl}). The places where the threshold is used are recorded in Subsection \ref{subsec:graphical-equation}.

\subsection{Related results}\label{subsec:intro-related}
\label{par:flying-wings}The mean curvature flow also has complete convex translating graphs without rotational symmetry. Wang proved that a complete convex translating graph is either entire or defined over a slab \cite[Corollaries 2.1 and 2.2]{wang2011convex}. The strictly convex ones over slabs of width greater than $\pi$ are the ``flying wings'', whose existence was predicted by a conjecture attributed to Hamilton and which Ilmanen called $\Delta$-wings \cite[Problem 22]{ilmanen2003problems}. Wang first proved their existence in every dimension \cite[Theorem 1.2]{wang2011convex}.

They were then constructed independently and at about the same time by Bourni--Langford--Tinaglia, whose preprint was posted first, and by Hoffman--Ilmanen--Mart\'in--White. Bourni--Langford--Tinaglia obtained them in every dimension and over every slab of width greater than $\pi$, and determined their asymptotics and symmetry \cite[Theorems 1 and 3 and Corollaries 2 and 4]{bourni2020existence}. Hoffman--Ilmanen--Mart\'in--White constructed them in $\rr^3$ over prescribed strips and proved that the flying wing over each strip is unique up to translation \cite[Theorems 4.1 and 5.7]{hoffman2019graphical}. In higher dimensions, Hoffman--Ilmanen--Mart\'in--White also constructed complete translating graphs over slabs \cite[Theorem 11.1]{hoffman2019graphical}, and the convexity of the rotationally symmetric ones follows from the convexity theorem of Bourni--Langford--Tinaglia for translators with at most two distinct principal curvatures \cite[Theorem 3.2]{bourni2020existence}. Moreover, the barriers of Bourni--Langford--Tinaglia show that their rotationally symmetric examples satisfy the slope bound of the uniqueness theorem of Hoffman--Ilmanen--Mart\'in--White \cite[Theorem 11.3]{hoffman2019graphical}. In $\rr^3$, Spruck--Xiao proved that every complete immersed mean convex translator is convex \cite[Theorem 1.1]{spruck-xiao}, and this applies to every translating graph, whose mean curvature equals its angle function and is therefore positive. Together with the constructions above, this completes the classification in $\rr^3$, where the complete translating graphs are, up to isometries, the grim reaper cylinder, the tilted grim reaper cylinders, the flying wings and the bowl soliton \cite[Theorem 7.1]{hoffman2019graphical}. While the general form of Hamilton's conjecture remains open, Bourni--Langford--Tinaglia subsequently constructed convex translators from unbounded circumscribed regular polytopes and simplices, which decompose at infinity into grim reaper cylinders and flying wings \cite[Theorems 4.2 and 4.4, Proposition 4.6 and Remark 4.8]{bourni2022ancient}. Whether the flying wings have analogues for curvature functions with $\alpha>\frac{1}{2}$ is the question of Open Problem \ref{op:flying-wings}.

\subsection{Relation to previous work}\label{subsec:intro-previous}
The present article supersedes the unpublished preprint with S. Rengaswami \cite{RengaswamiTorresSantaella2023Bowl}, which concerns $1$-homogeneous curvature functions. The catenoids for $S_k^{1/k}$ with $k$ odd have since been obtained by de Lira--de Farias \cite[Theorem 4(c)]{lira2026translating}. Every statement from the preprint that we use is proved here again under the hypotheses of this paper, as the $1$-homogeneous or leading-order cases of Theorems \ref{thm:main-a} and \ref{thm:main-b}, of the Gauss-root example in Subsection \ref{subsubsec:gauss-root-bowl}, of cases \textup{(ii)} and \textup{(iv)} of Theorem \ref{thm:main-c}, and of Corollary \ref{cor:growth-entire-graphs}. The even-root construction announced in the preprint requires a signed real root and is omitted.

\subsection{Organization}\label{subsec:intro-organization}
Section \ref{sec:preliminaries} fixes the notation. Section \ref{sec:curvature-functions-branches} introduces the curvature-function classes, the branches, the signed neck data, and the model families. Section \ref{sec:barrier-method} develops the comparison tools for the graphical equation. Section \ref{sec:positive-branch-asymptotics} proves Theorem \ref{thm:main-a}, Section \ref{sec:catenoidal-translators} proves Theorem \ref{thm:main-c}, and Section \ref{sec:uniqueness-catenoidal-barriers} proves Theorems \ref{thm:main-b} and \ref{thm:main-d}. Section \ref{sec:conclusions} collects five open problems, and Appendix \ref{app:q-eta-properties} records the properties of the phase map.
\section{Geometric preliminaries and rotational ansatz}
\label{sec:preliminaries}

\subsection{Background in extrinsic geometry}\label{subsec:extrinsic-notation}
Throughout this article $n\geq2$. Let $F_0:\Sigma^n\to\rr^{n+1}$ be an oriented immersion of a smooth manifold $\Sigma$, write $\Sigma_0=F_0(\Sigma)$, and let $\vec\nu$ be a chosen unit normal along $\Sigma_0$ and $\nabla$ the Euclidean connection. At each point $p\in\Sigma_0$, the Weingarten operator with respect to $\vec\nu$ is given by $\mathcal W(\vec X)=-\nabla_{\vec X}\vec\nu$, which is a self-adjoint endomorphism of $T_p\Sigma_0$. Its eigenvalues are the principal curvatures $\lambda(p)=(\lambda_1(p),\ldots,\lambda_n(p))$, and the second fundamental form is $\pI{\mathcal W\vec X,\vec Y}$. With this convention, the round sphere of radius $\rho$ with the inward normal has $\lambda=(\rho^{-1},\ldots,\rho^{-1})$, so the flow \eqref{eq:flow} moves it inward and it contracts to its center in finite time. 

Our hypersurfaces are noncompact, and \eqref{eq:flow} is not invariant under $\vec\nu\mapsto-\vec\nu$ unless $\gamma$ is odd, so the orientation must be fixed. We choose it so that the principal curvatures of a convex hypersurface are nonnegative.

\subsection{Rotationally symmetric hypersurfaces}\label{subsec:rotational-generating-curves}
Let $\mathtt c(s)=(r(s),u(s))$ be a regular curve in $\rr^2$ with $r(s)>0$, parametrized by arclength on $(s_-,s_+)$, so that $\dot r(s)^2+\dot u(s)^2=1$ with $\dot{(\cdot)}=\frac{d}{ds}$. The rotational hypersurface generated by $\mathtt c$ around the vertical axis is
\begin{equation*}
	\Sigma(\mathtt c)=\set{(r(s)\omega,u(s)):s\in(s_-,s_+),\ \omega\in\Sp^{n-1}},
\end{equation*}
where $\Sp^{n-1}\subset\rr^{n}$ denotes the unit sphere.

To express the curvature in scalar form, we define the tangent angle $\theta$, measured with respect to the horizontal axis, by $\dot r=\cos\theta$ and $\dot u=\sin\theta$. We choose along $\mathtt c$ the unit normal $\vec N=(-\dot u,\dot r)=(-\sin\theta,\cos\theta)$, so that $\Sigma(\mathtt c)$ carries the unit normal $\vec\nu(s,\omega)=(-\sin\theta(s)\,\omega,\cos\theta(s))$ and its angle function is $\pI{\vec\nu,\vec e_{n+1}}=\cos\theta$. If $\mathtt c$ is of class $\mathcal C^2$, its signed curvature is $\dot\theta=\ddot u\dot r-\dot u\ddot r$. Then, the principal curvatures of $\Sigma(\mathtt c)$ are the meridional curvature $\lambda_1=\dot\theta$ and the rotational curvature $\lambda_2=\cdots=\lambda_n=r^{-1}\sin\theta$, with multiplicity $n-1$.

Consequently, if the principal curvature vector $\lambda=\pare{\dot\theta,r^{-1}\sin\theta,\ldots,r^{-1}\sin\theta}$ lies in $\Gamma$, then $\Sigma(\mathtt c)$ satisfies the translator equation \eqref{eq:translator} precisely when the profile curve satisfies
\begin{equation}
	\label{eq:rotational-translator-arclength}
	\widetilde\gamma\left(\dot\theta, r^{-1}\sin\theta\right)=\cos\theta,
\end{equation}	
where $\widetilde\gamma(x,y)=\gamma(x,y,\ldots,y)$. We call $\widetilde\Gamma=\set{(x,y)\in\rr^2:(x,y,\ldots,y)\in\Gamma}$ the cone of rotational curvature pairs.

Now, on an interval where $\dot r>0$ the profile curve is the vertical graph $\mathtt c(r)=(r,u(r))$ with the upward normal, and with $v=u'$ the principal curvatures are
\begin{equation*}
	\lambda_1=\frac{v'}{(1+v^2)^{\frac{3}{2}}},
	\quad
	\lambda_2=\cdots=\lambda_n
	=
	\frac{v}{r\sqrt{1+v^2}}.
\end{equation*}
Since $\pI{\vec\nu,\vec e_{n+1}}=(1+v^2)^{-\frac{1}{2}}$ on this graph, the $\alpha$-homogeneity of $\gamma$ turns \eqref{eq:rotational-translator-arclength} into
\begin{equation}
	\label{eq:vertical-graph-translator}
	\widetilde\gamma
	\left(\frac{v'}{(1+v^2)^{1+\beta}},\frac{v}{r(1+v^2)^\beta}\right)=1,\quad\text{with}\quad\beta=\frac{\alpha-1}{2\alpha}.
\end{equation}
\section{Curvature functions and phase diagrams}\label{sec:curvature-functions-branches}

The class of curvature functions introduced in Section \ref{sec:introduction} is studied in this section for its later use throughout the article. In particular, we discuss the phase diagram associated with the rotational theory, together with the signed extension given in Subsection \ref{subsec:signed-level-sheets}.

\subsection{Positive curvature functions and the positive branch}\label{subsec:positive-curvature-functions}
\begin{definition}\label{def:alpha-positive-curvature-function}
	Let $\Gamma\subset\rr^n$ be an open symmetric cone containing the positive cone $\Gamma_+=\{\lambda\in\rr^n:\lambda_i>0\text{ for every }i\}$. We say that a function $\gamma:\Gamma\to[0,\infty)$ is an \emph{$\alpha$-positive curvature function} if:
	\begin{enumerate}[label=\textup{(\alph*)}]
		\item $\gamma$ is symmetric, i.e., invariant under permutations of its entries;
		\item $\gamma\in\mathcal C^1(\Gamma)$, $\gamma>0$ on $\Gamma$, and $\gamma$ is elliptic, i.e., $\partial_i\gamma>0$ for every $i$;
		\item $\gamma$ is positively $\alpha$-homogeneous of some degree $\alpha>0$, i.e., $\gamma(c\lambda)=c^\alpha\gamma(\lambda)$ for every $c>0$ and every $\lambda\in\Gamma$.
	\end{enumerate}
\end{definition}

Recall that, by \eqref{eq:rotational-translator-arclength}, a rotational profile curve generates a $\gamma$-translator when its meridional curvature $x=\dot\theta$, its rotational curvature $y=r^{-1}\sin\theta$ and its angle function $z=\cos\theta$ satisfy
\begin{equation}\label{eq:level-set-equation}
	\widetilde\gamma(x,y)=z\qquad\text{with }(x,y)\in\widetilde\Gamma.
\end{equation}
Particularly, the study of rotational $\gamma$-translators reduces to the study of the level sets of $\widetilde\gamma$. Moreover, as we mentioned in the introduction, the behavior of these level sets depends on the behavior of $\widetilde\gamma$ at the cylindrical direction $(0,1)$, which belongs to $\widetilde\Gamma$ or to its boundary. We notice that, by ellipticity, the map $t\mapsto\widetilde\gamma(t,1)$ is positive and strictly increasing on $(0,\infty)$, and therefore the limit $\lim_{t\to0^+}\widetilde\gamma(t,1)$ exists. This motivates the following definition.

\begin{definition}	\label{def:one-nondegeneracy}
	Let $\gamma:\Gamma\to[0,\infty)$ be an $\alpha$-positive curvature function.
	\begin{enumerate}[label=\textup{(\roman*)}]
		\item We say that $\gamma$ is \emph{$1$-nondegenerate} if $(0,1)\in\widetilde\Gamma$. It is \emph{normalized $1$-nondegenerate} if, in addition, $\widetilde\gamma(0,1)=1$\footnote{Every $1$-nondegenerate $\gamma$ becomes normalized $1$-nondegenerate after replacing $\gamma$ by $\widetilde\gamma(0,1)^{-1}\gamma$, which does not change $\widetilde\Gamma$.}.
		
		\item We say that $\gamma$ is $1$-degenerate if $(0,1)\in\partial\widetilde\Gamma$ and $\lim_{t\to0^+}\widetilde\gamma(t,1)=0$.
	\end{enumerate}
\end{definition}

We note that this is the cylindrical nondegeneracy condition of \cite[Definition 3]{rengaswami2024classification} when $(0,1)\in\widetilde\Gamma$, which concerns the curvature direction $(0,1,\ldots,1)$, that is, up to scale, the principal curvature vector of the round cylinder $\rr\times\Sp^{n-1}$. However, the two cases of Definition \ref{def:one-nondegeneracy} do not exhaust the $\alpha$-positive curvature functions. We say that $\gamma$ is \emph{boundary $1$-nondegenerate} when $(0,1)\in\partial\widetilde\Gamma$ and $\lim_{t\to0^+}\widetilde\gamma(t,1)>0$, as for $H+S_n^{1/n}$ on $\Gamma_+$. These curvature functions are nondegenerate in the sense of \cite[Definition 3]{rengaswami2024classification}, but the positive branch reaches the direction $(0,1)$ on the boundary of $\widetilde\Gamma$. For $H+S_n^{1/n}$, the branch vanishes there to order $n$. We do not consider this case in the present article, and we return to it in Open Problem \ref{op:boundary-nondegenerate}.

\begin{definition}
	\label{def:bowl-type-solution}
	Let $\gamma:\Gamma\to[0,\infty)$ be an $\alpha$-positive curvature function. A \emph{bowl-type solution} for the $\gamma$-flow is the rotational graph of a function $u:[0,R_\gamma)\to\rr$, with $u'(0)=0$ and $R_\gamma\in(0,\infty]$ maximal, whose vertical translates $u(|x|)+t$ solve \eqref{eq:flow}. Equivalently, its slope $v=u'$ solves \eqref{eq:vertical-graph-translator} with $v(0)=0$. It is \emph{entire} when $R_\gamma=\infty$, and \emph{cylindrical} when $R_\gamma<\infty$ and $\lim_{r\to R_\gamma^-}u(r)=\infty$.
\end{definition}

We note that every $\alpha$-positive curvature function of Definition \ref{def:alpha-positive-curvature-function} meets the hypotheses on the curvature function in \cite{rengaswami2024classification}. Therefore, by \cite[Theorem 1.1 and Sections 3, 4 and 8]{rengaswami2024classification}, the bowl-type solution of each such $\gamma$ exists and is unique. Moreover, its profile is of class $\mathcal C^2([0,R_\gamma))$ and it is strictly convex, that is, its principal curvature vector lies in $\Gamma_+$.

\begin{remark}\label{rem:langford-rengaswami-condition}
	Recall that the degeneracy condition of Definition \ref{def:one-nondegeneracy} concerns the behavior of the cylinders $\rr\times\Sp^{n-1}$ along the $\gamma$-flow. Indeed, the radius of such a cylinder satisfies
	\begin{equation*}
		\dot\rho=-\rho^{-\alpha}\,\widetilde\gamma(0,1),
	\end{equation*}
	see \cite[Chapter 18]{andrews2020extrinsic}. Therefore, if $\gamma$ is $1$-nondegenerate, the cylinder collapses in finite time, with $\rho(t)=\left(\rho(0)^{\alpha+1}-(\alpha+1)\widetilde\gamma(0,1)\,t\right)^{\frac{1}{\alpha+1}}$, whereas if $\gamma$ is $1$-degenerate, where $\widetilde\gamma(0,1)$ denotes the limit in Definition \ref{def:one-nondegeneracy}\textup{(ii)}, the cylinder is a static solution of the flow. In particular, it is invariant under vertical translations, which makes it a $\gamma$-translator. The index $1$ counts the Euclidean factors of the cylinder, as in \cite[Section 4]{lynch2022uniqueness}, where the analogous behavior of the cylinders $\rr^m\times\Sp^{n-m}$ is considered.

	In contrast, Langford--Rengaswami use in \cite{rengaswami-langford2023pancake} a nondegeneracy condition at the direction $(1,0,\ldots,0)$, which corresponds to the cylinders $\rr^{n-1}\times\Sp^{1}$. On the one hand, they require $\widetilde\gamma(1,0)>0$, and therefore, under the normalization $\widetilde\gamma(1,0)=1$, the cylinder over a solution of the curve shortening flow is a solution of the flow. On the other hand, they require the map $s\mapsto\widetilde\gamma(1,s)$ to be of class $\mathcal C^2$ up to $s=0$, which they use to expand the curvature function to second order at this direction in their area and displacement estimates \cite[Section 5]{rengaswami-langford2023pancake}. We emphasize that the bowl-type solutions and the catenoidal translators of the present article require only $\mathcal C^1$ regularity. Higher regularity enters only in the expansions of Section \ref{sec:positive-branch-asymptotics}, where it is required near $(0,1)$, and in the results of Section \ref{sec:uniqueness-catenoidal-barriers}, where the curvature function is smooth on $\Gamma_+$.

	Apart from this condition, their curvature functions are $1$-positive curvature functions on $\Gamma_+$ in the sense of Definition \ref{def:alpha-positive-curvature-function}, which in addition are smooth and inverse-concave, where the inverse-concavity allows them to use the differential Harnack inequality of Andrews. Under these conditions they construct ancient pancake solutions, analogous to those of \cite{bourni2021collapsing} for the mean curvature flow. We remark that the two conditions are not equivalent when $n\geq3$.

	The rank-one direction is the one that determines the grim reaper cylinders of Subsection \ref{subsec:intro-background}. A cylinder $\mathcal C\times\rr^{n-1}$ over a planar curve $\mathcal C$, parametrized by arclength with tangent angle $\theta$ as in Subsection \ref{subsec:rotational-generating-curves}, has principal curvatures $(\dot\theta,0,\ldots,0)$ and angle function $\cos\theta$. If $e_1\in\Gamma$, then $c=\gamma(e_1)>0$, and by $\alpha$-homogeneity the translator equation \eqref{eq:translator} reduces to $c\,\dot\theta^{\alpha}=\cos\theta$. On a G{\aa}rding cone $\Gamma_k$ with $k\geq2$, the vector $e_1$ is not in $\Gamma_k$, and there is no grim reaper cylinder. For $\alpha\leq\frac{1}{2}$ the solution of this equation is an entire graph, and for $\alpha>\frac{1}{2}$ it is a graph over a bounded interval \cite[Proposition 2.4]{choi2021convergence}. In the second case we denote by $w$ the width of the slab that contains the grim reaper cylinder. The graph is complete when $\alpha\leq1$. When $\alpha>1$ it ends at a finite height, and the complete translator then contains two vertical half-lines and is $\mathcal C^1$ but not smooth \cite[Remark 1.3]{choi2021convergence}.
\end{remark}

In the next proposition we describe the maximal positive branch $x=g_+(y,z)$ of $\widetilde\gamma(x,y)=z$. At the level $z=1$, this corresponds to the branch of \cite[Definition 6 and Section 2.3.2]{rengaswami2024classification} under the same hypotheses on $\gamma$, and this is the only level that we use after this section.

\begin{proposition}\label{prop:positive-branch}
	Let $\gamma:\Gamma\to[0,\infty)$ be a normalized $1$-nondegenerate or $1$-degenerate $\alpha$-positive curvature function. Define
	\begin{align*}
		\mathcal D_+=
		\begin{cases}
			\left\{(y,z)\in(0,\infty)^2:
			\frac{z}{\widetilde\gamma(1,1)}<y^\alpha<z\right\},
			& \text{if $\gamma$ is $1$-nondegenerate},\\[6pt]
			\left\{(y,z)\in(0,\infty)^2:
			\frac{z}{\widetilde\gamma(1,1)}<y^\alpha\right\},
			& \text{if $\gamma$ is $1$-degenerate}.
		\end{cases}
	\end{align*}
	Then there exists a unique function $g_+\in\mathcal C^1(\mathcal D_+)$ such that
	\begin{equation*}
		0<g_+(y,z)<y,\qquad\widetilde\gamma(g_+(y,z),y)=z
	\end{equation*}
	for every $(y,z)\in\mathcal D_+$. Moreover,
	\begin{align*}
		\partial_y g_+(y,z)&=-\partial_y\widetilde\gamma(g_+(y,z),y)\,\partial_x\widetilde\gamma(g_+(y,z),y)^{-1}<0,
		\\
		\partial_z g_+(y,z)&=\partial_x\widetilde\gamma(g_+(y,z),y)^{-1}>0.
	\end{align*}
	At the level $z=1$, writing $g_+(y)=g_+(y,1)$ and
	\begin{align*}
		I_+=
		\begin{cases}
			(y_*,1), & \text{if $\gamma$ is $1$-nondegenerate},
			\\
			(y_*,\infty), & \text{if $\gamma$ is $1$-degenerate},
		\end{cases}
		\quad\text{with}\quad y_*=\widetilde\gamma(1,1)^{-\frac{1}{\alpha}},
	\end{align*}
	the function $g_+$ extends continuously to the left endpoint $y_*$ of $I_+$, with $g_+(y_*)=y_*$.
\end{proposition}

\begin{proof}
	Firstly, we determine the domain $\mathcal D_+$. We fix $y>0$, and we note that by ellipticity the map $x\mapsto\widetilde\gamma(x,y)$ is strictly increasing on the positive sheet. On the one hand, let us assume that $\gamma$ is $1$-nondegenerate. Then, by the normalization $\widetilde\gamma(0,1)=1$ and homogeneity, we have $\widetilde\gamma(0,y)=y^\alpha$ and $\widetilde\gamma(y,y)=\widetilde\gamma(1,1)y^\alpha$. Therefore, the level $z$ is crossed by a unique $x\in(0,y)$ precisely when $y^\alpha<z<\widetilde\gamma(1,1)y^\alpha$, which is the stated domain $\mathcal D_+$ in the $1$-nondegenerate case.
	
	On the other hand, let us assume that $\gamma$ is $1$-degenerate. Then the map $x\mapsto\widetilde\gamma(x,y)$ is strictly increasing on $(0,y]$, and by homogeneity we have
	\begin{equation*}
		\lim_{x\to0^+}\widetilde\gamma(x,y)=y^\alpha\lim_{t\to0^+}\widetilde\gamma(t,1)=0,\qquad\widetilde\gamma(y,y)=\widetilde\gamma(1,1)y^\alpha.
	\end{equation*}
	Therefore, the level $z>0$ is crossed by a unique $x\in(0,y)$ precisely when $0<z<\widetilde\gamma(1,1)y^\alpha$, which is the stated domain $\mathcal D_+$ in the $1$-degenerate case.
	
	Next, we prove the regularity, which follows from the implicit function theorem. Indeed, at any point $(y_0,z_0)\in\mathcal D_+$, from the inequality $\partial_x\widetilde\gamma(x_0,y_0)>0$ we obtain a local $\mathcal C^1$ branch $x=g_+(y,z)$. Moreover, by uniqueness these local branches agree on overlaps, and they define $g_+\in\mathcal C^1(\mathcal D_+)$ with $0<g_+(y,z)<y$ and $\widetilde\gamma(g_+(y,z),y)=z$. Then, by differentiating this identity with respect to $y$ and $z$, we obtain the formulas for $\partial_yg_+$ and $\partial_zg_+$ of the statement, whose signs are determined by $\partial_x\widetilde\gamma>0$ and $\partial_y\widetilde\gamma>0$ on the positive sheet.
	
	Now, we consider the level $z=1$. The condition $(y,1)\in\mathcal D_+$ reads $\widetilde\gamma(1,1)^{-1}<y^\alpha<1$ in the $1$-nondegenerate case and $\widetilde\gamma(1,1)^{-1}<y^\alpha$ in the $1$-degenerate case, which is the stated interval $I_+$.
	
	Finally, we prove the continuous extension to $y_*$ at the level $z=1$. Since $g_+'<0$ and $0<g_+(y)<y$ on $I_+$, the function $g_+(y)$ increases as $y$ decreases to $y_*$, and therefore the limit $x_0=\lim_{y\to y_*^+}g_+(y)$ exists and satisfies $0<x_0\leq y_*$. In particular, the point $(x_0,y_*)$ has positive coordinates, and consequently it belongs to $\widetilde\Gamma$, where $\widetilde\gamma$ is continuous. Then we have $1=\lim_{y\to y_*^+}\widetilde\gamma(g_+(y),y)=\widetilde\gamma(x_0,y_*)$. Moreover, since $x\mapsto\widetilde\gamma(x,y_*)$ is strictly increasing on $(0,y_*]$ with $\widetilde\gamma(y_*,y_*)=1$, we conclude that $x_0=y_*=g_+(y_*)$.
\end{proof}

\begin{lemma}\label{lem:positive-endpoint-regularity}
	Let $\gamma$ be a normalized $1$-nondegenerate $\alpha$-positive curvature function, and assume that $\widetilde\gamma$ is of class $\mathcal C^m$ near $(0,1)$, with $m\geq1$. Then there are $\delta>0$ and a function $\phi$ of class $\mathcal C^m$ on $(1-\delta,1+\delta)$ such that
	\begin{equation*}
		\widetilde\gamma(\phi(y),y)=1,\qquad\phi(1)=0,\qquad\phi>0\text{ on }(1-\delta,1),\qquad\phi<0\text{ on }(1,1+\delta),
	\end{equation*}
	and $\phi=g_+$ on $(1-\delta,1)$. In particular, $g_+$ extends to $y=1$ with $g_+(1)=0$.
\end{lemma}

\begin{proof}
	Recall that the point $(0,1)$ lies in the open set $\widetilde\Gamma$, where $\widetilde\gamma(0,1)=1$ and $\partial_x\widetilde\gamma(0,1)>0$. Then, by the implicit function theorem, there are $\delta>0$ and a function $\phi$ of class $\mathcal C^m$ on $(1-\delta,1+\delta)$, with $\phi(1)=0$, whose graph is the level-one set near $(0,1)$. Moreover, since $\partial_y\widetilde\gamma(0,1)>0$, we have $\phi'(1)<0$. Therefore, after decreasing $\delta$, the function $\phi$ is positive on $(1-\delta,1)$ and negative on $(1,1+\delta)$, and $\phi(y)<y$. In particular, on $(1-\delta,1)$ the value $\phi(y)$ solves the level-one equation in $(0,y)$, and by the uniqueness in Proposition \ref{prop:positive-branch} we conclude that $\phi=g_+$ there.
\end{proof}

Accordingly, we call $\phi$ the level-one chart at the cylindrical direction. In addition, we also write $g_+$ for $\phi$ when its argument ranges over a two-sided neighborhood of $y=1$, as in the expansions of Section \ref{sec:positive-branch-asymptotics}.

\begin{notation}\label{not:asymptotic-notation}
	We use the big $O$-notation and the little $o$-notation along the text to describe the asymptotic behavior of the functions under consideration. To be more precise, given real-valued functions $f$ and $h$ defined near a limit point $\ell$, finite or infinite, we say that $f=O(h)$ if $f(x)h(x)^{-1}$ is bounded near $\ell$, and that $f=o(h)$ if $\lim_{x\to\ell}f(x)h(x)^{-1}=0$. In addition, the asymptotic symbols always refer to the limit stated in the corresponding result, and the limits at an endpoint are taken from within the domain.
\end{notation}

Recall that in the $1$-degenerate case the interval $I_+$ is unbounded. In Section \ref{sec:positive-branch-asymptotics} we prove that the growth of a positive end is then determined by the decay of $g_+(y)$ as $y\to\infty$, for which we require the power-type behavior of the following definition. In particular, the balanced case of item \textup{(iii)} is the one where the second-order term of $g_+$ and the correction produced by the prefactor $(1+v^2)^{\beta+1}$ in \eqref{eq:vertical-graph-translator} are of the same order.

\begin{definition}\label{def:degenerate-positive-branch-asymptotics}
	Let $\gamma$ be a $1$-degenerate $\alpha$-positive curvature function, and let $g_+(y)=g_+(y,1)$ be the positive branch.
	\begin{enumerate}[label=\textup{(\roman*)}]
		\item The function $\gamma$ has a \emph{power-decay positive branch} if there are constants $c_\gamma>0$ and $k_\gamma>0$ such that $g_+(y)=c_\gamma y^{-k_\gamma}(1+o(1))$ as $y\to\infty$.
		
		\item The function $\gamma$ has a \emph{second-order power-decay positive branch} if it has a power-decay positive branch and there are constants $c_{\gamma,1}\in\rr$ and $\sigma_\gamma>0$ such that $g_+(y)=c_\gamma y^{-k_\gamma}+c_{\gamma,1}y^{-k_\gamma-\sigma_\gamma}+o(y^{-k_\gamma-\sigma_\gamma})$ as $y\to\infty$.
		
		\item A second-order power-decay positive branch is \emph{balanced} if $\sigma_\gamma=k_\gamma+1$.
	\end{enumerate}
\end{definition}

\subsection{Signed extensions and maximal implicit branches}\label{subsec:signed-level-sheets}
To construct catenoidal translators for curvature functions other than $H^p$, $S_k$ and $S_k^{1/k}$, for which they are known (Subsection \ref{subsec:intro-background}), we have to deal with two difficulties. The first one is the neck, where the generating curve has a vertical tangent and attains a minimum of the distance to the axis. Indeed, near a neck of radius $R$ the curve is a graph over the height, $\mathtt c(u)=(r(u),u)$ with $r(0)=R$, whose arclength parametrization satisfies $(r(0),u(0),\theta(0))=\left(R,0,\frac{\pi}{2}\right)$. Therefore, \eqref{eq:rotational-translator-arclength} becomes
\begin{equation}\label{eq:neck-equation}
	\widetilde\gamma\left(\dot\theta(0),R^{-1}\right)=0,
\end{equation}
that is, the neck lies on the zero level of $\widetilde\gamma$. The second one is the branch below the neck. Along it the angle function is negative, and therefore the curvature function must take negative values. For $S_k$ with $k$ odd, the curvature function is a polynomial, defined for every curvature vector and odd, and both difficulties are solved by a change of orientation \cite[Remark 2]{lima2025translators}, \cite[Remark 2]{lira2026translating}. We describe these examples from our point of view in Subsubsection \ref{subsubsec:ellipticity-limited-models}. Since our aim is to cover the largest possible class of curvature functions, which are defined only on a cone, we require instead the extension of the following definition.

\begin{definition}\label{def:signed-rotational-extension}
	Let $\gamma:\Gamma\to[0,\infty)$ be an $\alpha$-positive curvature function. A \emph{signed rotational level profile} for $\gamma$ is a real-valued function $\widehat\gamma:\widehat\Gamma\to\rr$, defined on an open set $\widehat\Gamma\subset\rr^2\setminus\{(0,0)\}$ containing $\widetilde\Gamma$, such that $\widehat\gamma=\widetilde\gamma$ on $\widetilde\Gamma$.
\end{definition}

Given such an extension, we read the neck equation \eqref{eq:neck-equation} with $\widehat\gamma$ in place of $\widetilde\gamma$. Moreover, by differentiating $\dot r=\cos\theta$ we have $\ddot r(0)=-\dot\theta(0)$, and therefore the profile has a strict neck at $s=0$ precisely when $\dot\theta(0)<0$. When $\widehat\gamma$ is homogeneous, both conditions depend only on the direction $a=R\,\dot\theta(0)$ of the curvature pair, which leads to the following definition.

\begin{definition}\label{def:selected-signed-neck-datum}
	Let $\widehat\gamma:\widehat\Gamma\to\rr$ be a signed rotational level profile for $\gamma$. \emph{Selected angular neck data} consist of an open interval $A_g\subset\rr$ and a number $\rho_0>0$ such that $-\rho_0\in A_g$ and $\widehat\gamma(-\rho_0,1)=0$. After setting $a=xy^{-1}$ for the angular phase variable, the \emph{conical sheet} selected by $A_g$ is $U_g=\{(ay,y):a\in A_g,\ y>0\}$.
	
	The selected data are required to satisfy:
	\begin{enumerate}[label=\textup{(H\arabic*)}, ref=H\arabic*]
		\item \label{hyp:signed-domain}
		$U_g\subset\widehat\Gamma$ and $-U_g\subset\widehat\Gamma$.
		
		\item \label{hyp:signed-homogeneity}
		The function $\widehat\gamma$ is positively $\alpha$-homogeneous on $U_g\cup(-U_g)$, meaning $\widehat\gamma(tx,ty)=t^\alpha\widehat\gamma(x,y)$ for every $t>0$ and every $(x,y)\in U_g\cup(-U_g)$.
		
		\item \label{hyp:signed-oddness}
		The function $\widehat\gamma$ is odd under central reflection between the two selected sheets, meaning $\widehat\gamma(-x,-y)=-\widehat\gamma(x,y)$ for every $(x,y)\in U_g$.
		
		\item \label{hyp:signed-angular-monotonicity}
		The angular profile $\Psi_g(a)=\widehat\gamma(a,1)$ with $a\in A_g$ is continuous and strictly increasing on $A_g$.
	\end{enumerate}
\end{definition}

By \ref{hyp:signed-homogeneity} and \ref{hyp:signed-angular-monotonicity}, a neck of radius $R$ on the selected sheet has $R\,\dot\theta(0)=-\rho_0$, the unique zero of $\Psi_g$ in $A_g$, so it is strict because $\rho_0>0$. The requirement $\rho_0>0$ excludes the Gauss curvature and its powers. For $\gamma=S_n^{\alpha}$, the rotational restriction is $\widetilde\gamma(x,y)=\left(xy^{n-1}\right)^{\alpha}$, so the zero level with $y>0$ is $x=0$ and $\rho_0$ vanishes, whereas the polynomial models of Subsection \ref{subsec:model-examples-branches} have $\rho_k=\frac{n-k}{k}>0$ for $k\leq n-1$. When $\rho_0=0$ the zero direction is horizontal, so $\dot\theta$ vanishes at a vertical point of a rotational profile. In this case no neck is strict, and we lose the transversality on which the local construction of Section \ref{sec:catenoidal-translators} rests. However, rotational profiles of this excluded type exist. Examples are the catenoidal $K^p$-translators of $\rr^3$ classified by Aydin--L\'opez \cite[Theorem 1]{aydin2023translators} and the $n$-grim reapers of de Lima--Pipoli for $\gamma=S_n$ \cite[Theorems 4(iii) and 5(iii)]{lima2025translators}, which for $n$ odd are $\mathcal C^2$-singular along their sphere of minimal height.

\begin{definition}\label{def:regular-signed-neck-datum}
	The selected angular neck data from Definition \ref{def:selected-signed-neck-datum} are called \emph{regular} when there exists an open set $\Omega_g\subset\widehat\Gamma$ such that:
	\begin{enumerate}[label=\textup{(R\arabic*)}, ref=R\arabic*]
		\item \label{hyp:regular-neighborhood}
		$U_g\subset\Omega_g$ and $\widehat\gamma$ is of class $\mathcal C^1$ on $\Omega_g$.
		
		\item \label{hyp:regular-ellipticity}
		The ellipticity inequalities $\partial_x\widehat\gamma>0$ and $\partial_y\widehat\gamma>0$ hold on $U_g$.
	\end{enumerate}
\end{definition}

\begin{lemma}\label{lem:reduced-ellipticity}
	Let the selected angular neck data satisfy condition \ref{hyp:regular-neighborhood}. Then $\Psi_g\in\mathcal C^1(A_g)$ and, for every $a\in A_g$ and $y>0$,
	\begin{equation}
		\label{eq:reduced-ellipticity-identities}
		\partial_x\widehat\gamma(ay,y)=y^{\alpha-1}\,\Psi_g'(a),\qquad\partial_y\widehat\gamma(ay,y)=y^{\alpha-1}\left(\alpha\,\Psi_g(a)-a\,\Psi_g'(a)\right).
	\end{equation}
	In particular, condition \ref{hyp:regular-ellipticity} is equivalent to the \emph{reduced ellipticity inequalities}
	\begin{equation*}
		\Psi_g'(a)>0 \quad \text{and} \quad \alpha\,\Psi_g(a)-a\,\Psi_g'(a)>0 \quad \text{for every } a\in A_g.
	\end{equation*}
\end{lemma}

\begin{proof}
	We fix $a\in A_g$ and $y>0$. Since $(a,1)\in U_g$, by homogeneity \ref{hyp:signed-homogeneity} (applied with $t=y$) we have
	\begin{equation*}
		\widehat\gamma(ay,y)=y^\alpha\,\widehat\gamma(a,1)=y^\alpha\,\Psi_g(a).
	\end{equation*}
	Recall that $\Psi_g$ is of class $\mathcal C^1$ on $A_g$ by condition \ref{hyp:regular-neighborhood}. Then, by differentiating in $a$ and dividing by $y$, we obtain the first identity in \eqref{eq:reduced-ellipticity-identities}. Moreover, by differentiating in $y$, we have $a\,\partial_x\widehat\gamma(ay,y)+\partial_y\widehat\gamma(ay,y)=\alpha\,y^{\alpha-1}\,\Psi_g(a)$, and the second identity follows from the first one. Finally, since $y^{\alpha-1}>0$, condition \ref{hyp:regular-ellipticity} is equivalent to the positivity of both angular quantities on $A_g$.
\end{proof}

Regularity imposes ellipticity on the selected sheet, which makes the implicit graph $x=g(y,z)$ of the next proposition differentiable.

\begin{proposition}
	\label{prop:conical-sheet-maximal-chart}
	Let $(A_g,\rho_0)$ be selected angular neck data, and define
	\begin{equation*}
		\mathcal D_g=\left\{(y,z)\in(0,\infty)\times\rr: zy^{-\alpha}\in\Psi_g(A_g)\right\}.
	\end{equation*}
	The selected signed sheet is represented on $\mathcal D_g$ by the unique function
	\begin{equation*}
		g(y,z)=y\,\Psi_g^{-1}\!\left(zy^{-\alpha}\right).
	\end{equation*}
	Equivalently, $(g(y,z),y)\in U_g$ and $\widehat\gamma(g(y,z),y)=z$ for every $(y,z)\in\mathcal D_g$. The map $g$ is continuous, and $\mathcal D_g$ serves as the maximal domain on which the selected sheet can be expressed as a graph over $(y,z)$. Moreover, for $k\geq1$, when condition \ref{hyp:regular-ellipticity} holds and $\widehat\gamma\in\mathcal C^k(\Omega_g)$ for some open $\Omega_g\supset U_g$, the function $g$ is of class $\mathcal C^k$ on $\mathcal D_g$.
\end{proposition}

\begin{proof}
	We fix $y>0$. Firstly, we reduce the equation to the angular variable. The points of the selected sheet $U_g$ are those of the form $(ay,y)$ with $a\in A_g$, and by condition \ref{hyp:signed-homogeneity} we have $\widehat\gamma(ay,y)=y^\alpha\widehat\gamma(a,1)=y^\alpha\Psi_g(a)$. Therefore, solving $\widehat\gamma(x,y)=z$ in $U_g$ is equivalent to solving $\Psi_g(a)=zy^{-\alpha}$ for $a\in A_g$.
	
	Next, we solve the angular equation. By condition \ref{hyp:signed-angular-monotonicity}, the function $\Psi_g$ is continuous and strictly increasing on the interval $A_g$. Therefore, the angular equation has at most one solution, and it has one if and only if $zy^{-\alpha}\in\Psi_g(A_g)$. This determines the domain $\mathcal D_g$ and the angular solution $a(y,z)\in A_g$, and we set $g(y,z)=y\,a(y,z)$.
	
	Now, we prove the continuity. Let $(y_j,z_j)\to(y_0,z_0)$ in $\mathcal D_g$, and let $a_j,a_0\in A_g$ denote the corresponding angular solutions for $q_j=z_jy_j^{-\alpha}\to q_0=z_0y_0^{-\alpha}$. For every $\eps>0$ such that $a_0\pm\eps\in A_g$, by strict monotonicity we have $\Psi_g(a_0-\eps)<q_0<\Psi_g(a_0+\eps)$, and therefore $a_0-\eps<a_j<a_0+\eps$ for all sufficiently large $j$. Consequently, $a_j\to a_0$, and we obtain $g(y_j,z_j)=y_ja_j\to y_0a_0=g(y_0,z_0)$.
	
	Finally, let us assume that the selected angular neck data are regular. Then, by condition \ref{hyp:regular-neighborhood}, $\widehat\gamma$ is of class $\mathcal C^1$ on an open neighborhood of $U_g$. Let $\mathcal G(a,y,z)=\Psi_g(a)-zy^{-\alpha}$ for $a\in A_g$ and $y>0$. By condition \ref{hyp:regular-ellipticity} and Lemma \ref{lem:reduced-ellipticity}, we have $\partial_a\mathcal G(a,y,z)=\Psi_g'(a)=\partial_x\widehat\gamma(a,1)>0$. Therefore, by the implicit function theorem, the angular solution $a=a(y,z)$ is locally of class $\mathcal C^1$, and so is $g(y,z)=y\,a(y,z)$. Moreover, if $\widehat\gamma\in\mathcal C^k(\Omega_g)$, the same argument implies the local $\mathcal C^k$ regularity. We note that these local solutions agree on overlaps by the uniqueness of the angular solution.
\end{proof}

From now on, $g:\mathcal D\to\rr$, with $\mathcal D=\mathcal D_g$, denotes the selected signed neck chart of Proposition \ref{prop:conical-sheet-maximal-chart}. Since $-\rho_0\in A_g$, $\Psi_g(-\rho_0)=0$ and $\Psi_g$ is strictly increasing, the zero level of this chart is the ray
\begin{equation}
	\label{eq:zero-trace-conical}
	g(y,0)=-\rho_0 y,\quad y>0.
\end{equation}
In particular, the neck curvature found after Definition \ref{def:selected-signed-neck-datum} is $\dot\theta(0)=g\left(R^{-1},0\right)$.

The selected zero level places the curvature vector of the neck circle outside the open cone $\Gamma$. In the following lemma we decide whether it lies on $\partial\Gamma$ or outside $\overline\Gamma$.

\begin{lemma}
	\label{lem:neck-cone-boundary}
	Let $(A_g,\rho_0)$ be selected angular neck data, assume that $\Gamma$ is convex, and let $(a_-,\infty)$ be the angular slice $\left\{a\in\rr:(a,1)\in\widetilde\Gamma\right\}$. Then $-\rho_0\leq a_-\leq0$. If $a_-=-\rho_0$, the curvature vector $R^{-1}(-\rho_0,1,\ldots,1)$ of a neck of radius $R>0$ lies on $\partial\Gamma$, and if $a_->-\rho_0$, it lies outside $\overline\Gamma$.
\end{lemma}

\begin{proof}
	Firstly, we prove that no $a\in A_g$ with $a<-\rho_0$ satisfies $(a,1)\in\overline{\widetilde\Gamma}$. Let us assume by contradiction that such a point exists. Since $U_g$ is open and contains $(a,1)$, this point is a limit of points $p_i\in\widetilde\Gamma\cap U_g$, where $\widehat\gamma(p_i)=\widetilde\gamma(p_i)>0$. Then, by the continuity of $\widehat\gamma$ on $U_g$, we have $\Psi_g(a)=\widehat\gamma(a,1)=\lim_{i\to\infty}\widehat\gamma(p_i)\geq0$. However, by condition \ref{hyp:signed-angular-monotonicity} and $\Psi_g(-\rho_0)=0$, we have $\Psi_g(a)<0$ for every $a<-\rho_0$ in $A_g$, leading to the desired contradiction.
	
	Next, we describe the angular slice of $\widetilde\Gamma$, that is, its intersection with the line $y=1$. Since $\Gamma$ is convex, so is $\widetilde\Gamma$, being its preimage under the linear map $(x,y)\mapsto(x,y,\ldots,y)$, and therefore the angular slice is an interval. Moreover, it is open because $\widetilde\Gamma$ is open, and it contains $(0,\infty)$ because $\Gamma_+\subset\Gamma$. Consequently, it is of the form $(a_-,\infty)$ with $a_-\leq0$.
	
	Now, we prove that $a_-\geq-\rho_0$. Recall that the interval $A_g$ is open and contains $-\rho_0$. If $a_-<-\rho_0$, there would be $a\in A_g$ with $a_-<a<-\rho_0$ and $(a,1)\in\widetilde\Gamma$, which is excluded by the first step. Finally, we distinguish the two cases. On the one hand, if $a_-=-\rho_0$, the slice contains $(-\rho_0,0)$, and $(-\rho_0,1)$ is a limit of points of $\widetilde\Gamma$. Therefore, the curvature vector lies on $\partial\Gamma$. On the other hand, if $a_->-\rho_0$, then $(-\rho_0,1)\notin\overline{\widetilde\Gamma}$, since for a cone $(a,1)\in\overline{\widetilde\Gamma}$ precisely when $a\geq a_-$.
\end{proof}

In Section \ref{sec:catenoidal-translators} we construct the neck by integrating the rotational system near the zero level, where the radius is the unknown and the angle function $z=\cos\theta$ is the independent variable. To do so, we need the chart $g$ to be Lipschitz in the rotational variable $y$, since this gives the existence and the uniqueness of the neck. We also need a H\"older bound for $g$ in $z$, because its exponent determines the regularity of the generating curve across the neck circle.

\begin{definition}\label{def:neck-regular-radius}
	We say that a radius $R>0$ is \emph{neck-admissible} for the selected signed chart if there exist an open rectangle $\mathcal V=J\times K\subset\mathcal D$ containing $(R^{-1},0)$ and a constant $L>0$ such that
	\begin{equation*}
		\left|g(y_1,z)-g(y_2,z)\right|
		\leq
		L\left|y_1-y_2\right|
		\qquad
		\text{for all }(y_1,z),(y_2,z)\in\mathcal V .
	\end{equation*}
	In addition, we say that $R$ is \emph{$\varkappa$-regular at the neck}, for $\varkappa\in(0,1]$, if there exists a constant $C>0$ such that $\left|g(y,z_1)-g(y,z_2)\right|\leq C\left|z_1-z_2\right|^\varkappa$ for $(y,z_1),(y,z_2)\in\mathcal V$.
\end{definition}

By Proposition \ref{prop:conical-sheet-maximal-chart}, the chart is continuous on $\mathcal D$ for all selected angular neck data, but neither condition of Definition \ref{def:neck-regular-radius} follows from continuity alone. When the data are regular, $g$ is of class $\mathcal C^1$ on $\mathcal D$, and therefore every radius $R>0$ is neck-admissible and $1$-regular at the neck. Consequently, the definition is a restriction only for data that are not regular, and in the next lemma we determine it from the order of vanishing of $\Psi_g$ at the zero ray.

\begin{lemma}\label{lem:neck-regularity-vanishing-order}
	Let $(A_g,\rho_0)$ be selected angular neck data, and suppose that near $a=-\rho_0$
	\begin{equation*}
		\Psi_g(a)=\operatorname{sgn}(a+\rho_0)\,|a+\rho_0|^{p}\,\phi(a),\qquad p>0,\qquad\phi\in\mathcal C^1,\qquad\phi(-\rho_0)>0.
	\end{equation*}
	Then there exists a continuous positive function $\psi$, defined near the origin and of class $\mathcal C^1$ when $p\leq1$, with $\psi(0)=\phi(-\rho_0)^{-1/p}$, fulfilling
	\begin{equation*}
		\Psi_g^{-1}(w)=-\rho_0+\operatorname{sgn}(w)|w|^{\frac1p}\,\psi(w),\quad g(y,z)=-\rho_0y+y^{1-\frac{\alpha}{p}}\operatorname{sgn}(z)|z|^{\frac1p}\,\psi\!\left(zy^{-\alpha}\right).
	\end{equation*}
	Moreover, every radius $R>0$ is $\frac{1}{p}$-regular at the neck when $p>1$ and $1$-regular at the neck when $p\leq1$, and the selected data are regular only when $p=1$.
\end{lemma}

\begin{proof}
	Firstly, we prove the two formulas. Near the zero ray, by the positivity of $\phi$ we have $\operatorname{sgn}\Psi_g(a)=\operatorname{sgn}(a+\rho_0)$ and $\left|\Psi_g(a)\right|^{\frac1p}=|a+\rho_0|\,\phi(a)^{\frac1p}$, and $\Psi_g$ is invertible near $a=-\rho_0$ by condition \ref{hyp:signed-angular-monotonicity}. Then, with $w=\Psi_g(a)$ and $\psi(w)=\phi\!\left(\Psi_g^{-1}(w)\right)^{-\frac1p}$ we obtain the first formula, and the second one follows from the chart $g(y,z)=y\,\Psi_g^{-1}\!\left(zy^{-\alpha}\right)$ of Proposition \ref{prop:conical-sheet-maximal-chart}. Moreover, by continuity, $\psi(0)=\phi(-\rho_0)^{-\frac1p}$. Now, by differentiating the hypothesis, we have $\Psi_g'(a)=p\,|a+\rho_0|^{p-1}\phi(a)+O\bigl(|a+\rho_0|^{p}\bigr)$ as $a\to-\rho_0$, which tends to $\infty$ when $p<1$ and to $\phi(-\rho_0)$ when $p=1$. In both cases, $(\Psi_g^{-1})'(w)=\Psi_g'(\Psi_g^{-1}(w))^{-1}$ has a limit as $w\to0$. Therefore, $\Psi_g^{-1}$ is of class $\mathcal C^1$ near the origin by the mean value theorem, and $\psi$ inherits the regularity of $\phi$.
	
	Next, we examine the three regimes. For $p=1$, the profile is differentiable at the zero ray with $\Psi_g'(-\rho_0)=\phi(-\rho_0)>0$, and the chart is of class $\mathcal C^1$ near $z=0$. For $p>1$, the derivative exists but $\Psi_g'(-\rho_0)=0$, and therefore condition \ref{hyp:regular-ellipticity} fails at the zero ray only. In this case the chart is $\frac1p$-H\"older in $z$ near $z=0$, and no better, because $|z|^{\frac1p}$ is multiplied by a factor bounded away from zero. Finally, we notice that condition \ref{hyp:regular-neighborhood} fails for $p<1$, since $\Psi_g'$ blows up at the zero ray. However, since $\psi$ is of class $\mathcal C^1$ by the first step, we have that the chart is still of class $\mathcal C^1$ near $z=0$, with $\partial_zg=0$ there.
	
	Now, it remains to verify the Lipschitz bound in the rotational variable $y$, which turns out to be the same in the three regimes. For this purpose, we set $w=zy^{-\alpha}$, so that $g(y,z)=y\,\Psi_g^{-1}(w)$, and then
	\begin{equation*}
		\partial_yg=\Psi_g^{-1}(w)-\alpha\,w\left(\Psi_g^{-1}\right)'(w).
	\end{equation*}
	Furthermore, since $\Psi_g$ is invertible, the angular variable $a=\Psi_g^{-1}(w)$ satisfies $w\left(\Psi_g^{-1}\right)'(w)=\Psi_g(a)\Psi_g'(a)^{-1}$. Then, by the hypothesis, we notice that this quotient satisfies $\Psi_g(a)\Psi_g'(a)^{-1}=O\!\left(|a+\rho_0|\right)$ as $a\to-\rho_0$, which implies that it vanishes at the zero ray. Particularly, since $\partial_yg$ is bounded near $(R^{-1},0)$, we obtain the neck admissibility of every $R>0$, together with the H\"older exponents in $z$ given by the three regimes. Consequently, we conclude the proof by noticing that the data can be regular only when $p=1$, since condition \ref{hyp:regular-ellipticity} fails for $p>1$ and condition \ref{hyp:regular-neighborhood} fails for $p<1$.
\end{proof}

The powers $\gamma=\left(\frac{H}{n-1}\right)^p$ illustrate the lemma. Their level-one diagram is the line $x+(n-1)y=n-1$ for every $p$, while $\Psi_g(a)=\left(\frac{a+n-1}{n-1}\right)^p$ vanishes to order $p$ at the zero ray, so the three regimes occur within a single family.

\subsection{Reflection and lower level-one branches}\label{subsec:reflected-lower-components}
To study the lower side of a catenoidal profile $\mathtt c(s)=(r(s),u(s))$ with the barrier method of Section \ref{sec:barrier-method}, we consider $\overline{\mathtt c}(s)=\mathtt c(-s)$ with $\overline\theta(s)=\theta(-s)+\pi$, as in \cite[Subsection 4.2.2]{martinez2022equilibrium} and \cite[Remark 2]{lira2026translating}. This reparametrization does not change the curve and replaces the unit normal $\vec N$ of Subsection \ref{subsec:rotational-generating-curves} by $-\vec N$, so that the phase variables $(x,y,z)=(\dot\theta,r^{-1}\sin\theta,\cos\theta)$ change to $(-x,-y,-z)$. By the oddness condition \ref{hyp:signed-oddness}, the reparametrized curve solves the same signed equation $\widehat\gamma(x,y)=z$, now with $\overline z>0$, and by positive homogeneity we normalize it at the level $\overline z=1$.

In terms of the chart $g$, we describe the reparametrized side by the reflected chart
\begin{equation*}
	g^\sharp(y,z)=-g(-y,-z) \quad \text{over} \quad \mathcal D^\sharp=\{(y,z):(-y,-z)\in\mathcal D\}.
\end{equation*}
Indeed, by the oddness condition \ref{hyp:signed-oddness}, for every $(y,z)\in\mathcal D^\sharp$ we have $\widehat\gamma(g^\sharp(y,z),y)=z$ and $(g^\sharp(y,z),y)\in-U_g$. In particular, the reflected lower graphical branch at output level $1$ is $g^\sharp(y,1)=-g(-y,-1)$.

We emphasize that the level-one branch is not evaluated at the neck curvature $-R^{-1}$. Indeed, by homogeneity, a point $(\overline x,\overline y)\in-U_g$ at output level $\overline z=\widehat\gamma(\overline x,\overline y)>0$ satisfies $\widehat\gamma(\overline z^{-1/\alpha}\overline x,\overline z^{-1/\alpha}\overline y)=\overline z^{-1}\widehat\gamma(\overline x,\overline y)=1$, and therefore its level-one representative $\overline z^{-1/\alpha}(\overline x,\overline y)$ is the point $(g^\sharp(t,1),t)$ with $t=\overline z^{-1/\alpha}\overline y$. Since near a neck of radius $R$ we have $\overline y\to-R^{-1}$ and $\overline z\to0^+$, it follows that $t\to-\infty$ for every $R>0$. Consequently, all necks lie at the end $t\to-\infty$ of the same branch $g^\sharp(\cdot,1)$.

\begin{definition}\label{def:selected-lower-branch}
	A \emph{selected fixed-level lower branch} is a pair formed by an interval $J\subset(-\infty,0)$ and a function $h:J\to\rr$ such that $\widehat\gamma(h(y),y)=1$ for every $y\in J$. In addition, we say that the branch is \emph{regular} if there exists an open set $\Omega_-\subset\widehat\Gamma$ containing $\{(h(y),y):y\in J\}$ such that $\widehat\gamma\in\mathcal C^1(\Omega_-)$ and both partial derivatives $\partial_x\widehat\gamma$ and $\partial_y\widehat\gamma$ are positive along this branch.
\end{definition}

We notice that a regular selected lower branch is strictly decreasing on $J$. Indeed, by differentiating $\widehat\gamma(h(y),y)=1$, we obtain
\begin{equation}\label{eq:lower-branch-monotonicity}
	h'(y)=-\partial_y\widehat\gamma(h(y),y)\,\partial_x\widehat\gamma(h(y),y)^{-1}<0.
\end{equation}
Moreover, on the selected sheet the level-one branch is $y\mapsto g^\sharp(y,1)$, and by homogeneity we have
\begin{equation}\label{eq:reflected-lower-angular-relation}
	y^{-1}g^\sharp(y,1)=\Psi_g^{-1}\left(-(-y)^{-\alpha}\right).
\end{equation}
As $\Psi_g$ is strictly increasing with $\Psi_g(-\rho_0)=0$ and the reflected level $-(-y)^{-\alpha}$ is negative, we obtain
\begin{equation}\label{eq:lower-branch-positivity}
	y^{-1}g^\sharp(y,1)<-\rho_0, \quad \text{or equivalently} \quad g^\sharp(y,1)>-\rho_0y>0,
\end{equation}
wherever the branch is defined. The reflected level tends to $0^-$ as $y\to-\infty$ and therefore $\lim_{y\to-\infty}y^{-1}g^\sharp(y,1)=-\rho_0$. In addition, it is strictly decreasing on $(-\infty,0)$, and $\Psi_g(A_g)$ is an interval containing $0$. Consequently, the set of $y<0$ with $(y,1)\in\mathcal D^\sharp$ has the form $(-\infty,y_{\mathrm{dom}})$ with $y_{\mathrm{dom}}\leq0$.

\begin{definition}\label{def:neck-accessible-branch}
	Let $(A_g,\rho_0)$ be selected angular neck data, and let $y\mapsto g^\sharp(y,1)$ on $(-\infty,y_{\mathrm{dom}})$ be the \emph{neck-accessible reflected branch}. We say that this branch \emph{admits a maximal regular continuation} if it is regular in the sense of Definition \ref{def:selected-lower-branch} on some interval $(-\infty,y')$. In this case, the \emph{maximal regular continuation} is the pair $(I_-,g_-)$, where $I_-=(-\infty,y_{\mathrm{reg}})$ with
	\begin{equation*}
		y_{\mathrm{reg}}=\sup\left\{y'\leq y_{\mathrm{dom}}:\ \text{the branch is regular on }(-\infty,y')\right\},
	\end{equation*}
	and $g_-$ is the restriction of the branch to $I_-$.
\end{definition}

Since the regularity of Definition \ref{def:selected-lower-branch} is a local condition, the branch is regular on $I_-$. Moreover, regular selected angular neck data make the whole neck-accessible branch regular. Indeed, the set $-U_g$ is open, and by the oddness condition \ref{hyp:signed-oddness} we have $\widehat\gamma(p)=-\widehat\gamma(-p)$ for $p\in-U_g$. Therefore, by conditions \ref{hyp:regular-neighborhood} and \ref{hyp:regular-ellipticity}, $\widehat\gamma$ is of class $\mathcal C^1$ on $-U_g$ with $\partial_x\widehat\gamma(p)=\partial_x\widehat\gamma(-p)>0$ and $\partial_y\widehat\gamma(p)=\partial_y\widehat\gamma(-p)>0$. Since the branch lies in $-U_g$, it is regular with $\Omega_-=-U_g$ and $y_{\mathrm{reg}}=y_{\mathrm{dom}}$. The same argument applies when these conditions fail only on the zero ray, as for the powers of the mean curvature. The branch never meets the reflected zero ray, because its angular phase $a=y^{-1}g^\sharp(y,1)$ satisfies $\Psi_g(a)=-(-y)^{-\alpha}<0$ by \eqref{eq:reflected-lower-angular-relation}.

\begin{definition}\label{def:lower-continuation-outcomes}
	Let $(I_-,g_-)$ be the maximal regular continuation of a neck-accessible reflected branch, with $I_-=(-\infty,y_{\mathrm{reg}})$. We distinguish the following four alternatives.
	\begin{enumerate}[label=\textup{(\roman*)}]
		\item The branch is \emph{axis-crossing} if $y_{\mathrm{reg}}=0$ and, for some $x_0>0$, $g_-(y)\to x_0$ as $y\to0^-$.
		
		\item The branch has \emph{finite-angle zero output} if $y_{\mathrm{reg}}=0$ and $y^{-1}g_-(y)\to\ell$ as $y\to0^-$ for some $\ell\in\rr$. It is called \emph{flat} if, in addition, there is a constant $\mu>0$ such that $g_-(y)=\ell y+O(|y|^{1+\mu})$ as $y\to0^-$.
		
		\item The branch has \emph{infinite-angle zero output} if $y_{\mathrm{reg}}=0$, $g_-(y)\to0$, and $y^{-1}g_-(y)\to-\infty$ as $y\to0^-$. It has \emph{fractional zero output} if there are constants $d>0$ and $\sigma\in(0,1)$ such that $g_-(y)=d(-y)^\sigma(1+o(1))$ as $y\to0^-$.
		
		\item The branch is \emph{regularity-limited} if $y_{\mathrm{reg}}<0$, and we write $x_{\mathrm{reg}}=\lim_{y\to y_{\mathrm{reg}}^-}g_-(y)$ for its endpoint limit. It is called \emph{ellipticity-limited} if $\widehat\gamma$ is of class $\mathcal C^1$ on a neighborhood of $(x_{\mathrm{reg}},y_{\mathrm{reg}})$ in $\widehat\Gamma$, with $\partial_x\widehat\gamma>0$ and $\partial_y\widehat\gamma=0$ at that point.
	\end{enumerate}
\end{definition}

In all our regularity-limited models $y_{\mathrm{reg}}=y_{\mathrm{dom}}$, because the left endpoint of $A_g$ is the direction where the ellipticity condition \ref{hyp:regular-ellipticity} fails. The branch then reaches the boundary of the selected sheet, and the ellipticity-limited condition describes this boundary point.

\begin{proposition}\label{prop:general-lower-phase-alternative}
	The maximal regular continuation of a neck-accessible reflected branch belongs to exactly one of the four alternatives in Definition \ref{def:lower-continuation-outcomes}. If $y_{\mathrm{reg}}=0$, then $y^{-1}g_-(y)\to\inf A_g$ as $y\to0^-$. In particular, when $y_{\mathrm{reg}}=0$ the branch has finite-angle zero output exactly when $\inf A_g>-\infty$, in which case $\ell=\inf A_g<-\rho_0$. In the regularity-limited case the endpoint limit satisfies $x_{\mathrm{reg}}\geq\rho_0|y_{\mathrm{reg}}|>0$.
\end{proposition}

\begin{proof}
	Firstly, we note that the branch is positive by \eqref{eq:lower-branch-positivity}, and that $g_-'<0$ on $I_-$ by regularity and \eqref{eq:lower-branch-monotonicity}. Therefore, $g_-$ has a finite nonnegative limit at $y_{\mathrm{reg}}$.
	
	Now, we distinguish two cases. On the one hand, if $y_{\mathrm{reg}}<0$, the branch is regularity-limited, and from \eqref{eq:lower-branch-positivity} we have $g_-(y)>-\rho_0y>\rho_0|y_{\mathrm{reg}}|$. Consequently, $x_{\mathrm{reg}}\geq\rho_0|y_{\mathrm{reg}}|>0$.
	
	On the other hand, if $y_{\mathrm{reg}}=0$, then $y_{\mathrm{dom}}=0$ because $y_{\mathrm{reg}}\leq y_{\mathrm{dom}}\leq0$, and therefore $\Psi_g(A_g)$ contains every negative value. Moreover, as $y\to0^-$ the reflected level tends to $-\infty$, and by \eqref{eq:reflected-lower-angular-relation} we have $y^{-1}g_-(y)\to\inf A_g$. If $\inf A_g$ is finite, it is strictly less than $-\rho_0$ because $-\rho_0$ is an interior point of $A_g$. In this case $g_-(y)\to0$, and the branch has finite-angle zero output with $\ell=\inf A_g$. If $\inf A_g=-\infty$, then $y^{-1}g_-(y)\to-\infty$, and the branch is axis-crossing when the endpoint limit is positive and has infinite-angle zero output when it is zero.
	
	Finally, the four alternatives are mutually exclusive. Indeed, items \textup{(i)}--\textup{(iii)} require $y_{\mathrm{reg}}=0$ and item \textup{(iv)} requires $y_{\mathrm{reg}}<0$. Moreover, item \textup{(i)} is the only one of the first three with a positive endpoint limit, and \textup{(ii)} differs from \textup{(iii)} by a finite angular limit.
\end{proof}

\begin{proposition}\label{prop:angular-growth-lower-branch}
	Assume that the maximal regular continuation of the neck-accessible reflected branch satisfies $y_{\mathrm{reg}}=0$. Assume also that $A_g$ is unbounded below and that $\Psi_g(a)=-c_\infty(-a)^{p_\infty}(1+o(1))$ as $a\to-\infty$ for some constants $c_\infty>0$ and $p_\infty>0$. Then $p_\infty\geq\alpha$ and
	\begin{equation*}
		g_-(y)=c_\infty^{-1/p_\infty}(-y)^{1-\alpha/p_\infty}(1+o(1))\qquad\text{as }y\to0^-.
	\end{equation*}
	Moreover, in the notation of Definition \ref{def:lower-continuation-outcomes},
	\begin{enumerate}[label=\textup{(\roman*)}]
		\item if $p_\infty=\alpha$, the branch is axis-crossing with $x_0=c_\infty^{-1/\alpha}$;
		
		\item if $p_\infty>\alpha$, the branch has fractional zero output with $\sigma=1-\alpha/p_\infty$ and $d=c_\infty^{-1/p_\infty}$.
	\end{enumerate}
\end{proposition}

\begin{proof}
	Recall that $y^{-1}g_-(y)\to\inf A_g=-\infty$ as $y\to0^-$ by Proposition \ref{prop:general-lower-phase-alternative}. Then, by applying $\Psi_g$ to \eqref{eq:reflected-lower-angular-relation} and by using the hypothesis on $\Psi_g$, we obtain
	\begin{equation*}
		c_\infty\left(-y^{-1}g_-(y)\right)^{p_\infty}(1+o(1))=(-y)^{-\alpha},
	\end{equation*}
	and therefore $g_-(y)=c_\infty^{-1/p_\infty}(-y)^{1-\alpha/p_\infty}(1+o(1))$.
	
	Moreover, by \eqref{eq:lower-branch-monotonicity} and \eqref{eq:lower-branch-positivity} the branch is positive and decreasing on $(-\infty,0)$, and consequently it has a finite limit as $y\to0^-$. This implies that $1-\alpha/p_\infty\geq0$, that is, $p_\infty\geq\alpha$. On the one hand, if $p_\infty=\alpha$, then $g_-(y)\to c_\infty^{-1/\alpha}>0$ and the branch is axis-crossing. On the other hand, if $p_\infty>\alpha$, then $\sigma=1-\alpha/p_\infty\in(0,1)$ and the branch has fractional zero output with $d=c_\infty^{-1/p_\infty}$.
\end{proof}

We remark that the four alternatives describe only the maximal regular part of the branch. For instance, an axis-crossing branch approaches the point $(x_0,0)$ of the level-one set, but beyond this point the level-one component could leave the domain of the signed profile, lose ellipticity, cease to be a graph over the rotational variable, or turn back before reaching the diagonal. To exclude these possibilities, we require the following condition.

\begin{definition}\label{def:regular-positive-matching}
	We say that an axis-crossing lower branch satisfies the \emph{regular positive-matching condition} if $\widehat\gamma$ is of class $\mathcal C^1$ on a neighborhood of $(x_0,0)$ and the connected component of the level set $\{\widehat\gamma=1\}$ containing its graph extends through $(x_0,0)$ as a $\mathcal C^1$ graph $x=g_{\mathrm{ax}}(y)$ on $0\leq y\leq y_*$, with $y_*=\widetilde\gamma(1,1)^{-\frac{1}{\alpha}}$, such that
	\begin{equation*}
		g_{\mathrm{ax}}(0)=x_0,\quad g_{\mathrm{ax}}(y_*)=y_*,\quad \partial_x\widehat\gamma(g_{\mathrm{ax}}(y),y)>0\ \text{ and }\ g_{\mathrm{ax}}(y)>y\ \text{ for } 0\leq y<y_*.
	\end{equation*}
\end{definition}

We observe that the endpoint $y_*$ is determined by homogeneity. Indeed, by the identity $\widetilde\gamma(y,y)=y^\alpha\widetilde\gamma(1,1)$, the diagonal meets the level-one set only at $y=y_*$.

\begin{lemma}\label{lem:matching-c1-junctions}
	The matching graphs join the neighboring branches as follows.
	\begin{enumerate}[label=\textup{(\roman*)}]
		\item If $(I_-,g_-)$ is an axis-crossing branch that satisfies the regular positive-matching condition, the function equal to $g_-$ on $I_-$ and to $g_{\mathrm{ax}}$ on $[0,y_*]$ is of class $\mathcal C^1$ across $y=0$.
		\item Let $h$ be a continuous function on a left neighborhood of $y_*$ such that $\widetilde\gamma(h(y),y)=1$ there and $h(y)\to y_*$ as $y\to y_*^-$. Then the function equal to $h$ on the left of $y_*$ and to $g_+$ on the right is of class $\mathcal C^1$ across $y_*$. This applies to $h=g_{\mathrm{ax}}$.
	\end{enumerate}
\end{lemma}

\begin{proof}
	\textup{(i)} By Definition \ref{def:regular-positive-matching}, $\widehat\gamma$ is of class $\mathcal C^1$ near $(x_0,0)$ with $\partial_x\widehat\gamma(x_0,0)>0$. Then, by the implicit function theorem, the level-one set near this point is a single $\mathcal C^1$ graph over a two-sided neighborhood of $y=0$. Moreover, since $g_-(y)\to x_0$ as $y\to0^-$ and $g_{\mathrm{ax}}(0)=x_0$, this graph coincides with $g_-$ on the left of $0$ and with $g_{\mathrm{ax}}$ on the right.
	
	\textup{(ii)} We note that the point $(y_*,y_*)$ lies in $\widetilde\Gamma$, where $\widetilde\gamma$ is of class $\mathcal C^1$ and elliptic, and that $\widetilde\gamma(y_*,y_*)=1$ by homogeneity. Therefore, by the same argument, there is a single $\mathcal C^1$ graph through $(y_*,y_*)$. This graph coincides with $h$ on the left of $y_*$ and with $g_+$ on the right, because $g_+(y)\to y_*$ as $y\to y_*^+$ by Proposition \ref{prop:positive-branch}.
\end{proof}

A regularity-limited branch may continue past its endpoint even though the selected ellipticity condition fails there. At an ellipticity-limited endpoint, the branch stops only because the strict inequality $\partial_y\widehat\gamma>0$ of Definition \ref{def:selected-lower-branch} fails.

\begin{lemma}\label{lem:ellipticity-limited-endpoint}
	Let $(I_-,g_-)$ be an ellipticity-limited branch with endpoint $(x_{\mathrm{reg}},y_{\mathrm{reg}})$. Then, near this endpoint, the connected component of the level set $\{\widehat\gamma=1\}$ containing the graph of $g_-$ coincides with the graph of a $\mathcal C^1$ function $G$ defined on a two-sided neighborhood of $y_{\mathrm{reg}}$. Furthermore, $G$ agrees with $g_-$ on a left neighborhood of $y_{\mathrm{reg}}$ and satisfies $G'(y_{\mathrm{reg}})=0$.
\end{lemma}

\begin{proof}
	By Definition \ref{def:lower-continuation-outcomes}\textup{(iv)}, $\widehat\gamma$ is of class $\mathcal C^1$ near $(x_{\mathrm{reg}},y_{\mathrm{reg}})$ with $\partial_x\widehat\gamma>0$ there. Then, by the implicit function theorem, there is a function $G$ of class $\mathcal C^1$ on a two-sided neighborhood of $y_{\mathrm{reg}}$ whose graph coincides with the level set in a box around this point. Moreover, since $g_-(y)\to x_{\mathrm{reg}}$ as $y\to y_{\mathrm{reg}}^-$, the graph of $g_-$ enters this box, and therefore $G$ agrees with $g_-$ on a left neighborhood of $y_{\mathrm{reg}}$. Finally, by differentiating the level identity $\widehat\gamma(G(y),y)=1$, we obtain
	\begin{equation*}
		G'(y)=-\partial_y\widehat\gamma(G(y),y)\,\partial_x\widehat\gamma(G(y),y)^{-1}
	\end{equation*}
	as in \eqref{eq:lower-branch-monotonicity}, and from $\partial_y\widehat\gamma(x_{\mathrm{reg}},y_{\mathrm{reg}})=0$ we conclude that $G'(y_{\mathrm{reg}})=0$.
\end{proof}

By Lemma \ref{lem:ellipticity-limited-endpoint}, the component containing the graph of $g_-$ extends beyond $y_{\mathrm{reg}}$ as a graph, which is increasing exactly where $\partial_y\widehat\gamma<0$, as long as $\partial_x\widehat\gamma>0$. Among these continuations, we are interested in those that remain graphs over the whole negative rotational axis, which leads to the following definition.

\begin{definition}\label{def:lower-continuation}
	Let $(I_-,g_-)$ be a regularity-limited branch. An \emph{extension} of $g_-$ is a function $g_-^{\mathrm{ext}}$ of class $\mathcal C^1$ on $(-\infty,0)$ whose graph is the part of the connected component of the level set $\{\widehat\gamma=1\}$ containing the graph of $g_-$ that lies in $\{y<0\}$, such that $g_-^{\mathrm{ext}}=g_-$ on $I_-$ and, for every $y<0$,
	\begin{equation*}
		\partial_x\widehat\gamma\left(g_-^{\mathrm{ext}}(y),y\right)>0,\qquad g_-^{\mathrm{ext}}(y)>0.
	\end{equation*}
	No sign condition on $\partial_y\widehat\gamma$ is imposed outside $I_-$. We say that $g_-$ admits
	\begin{enumerate}[label=\textup{(\roman*)},ref=\thedefinition\,\textup{(\roman*)}]
		\item \label{def:curvature-blow-up-continuation}
		a \emph{curvature-blow-up continuation} if it has an extension with $g_-^{\mathrm{ext}}(y)\to+\infty$ as $y\to0^-$. This continuation has \emph{order $q>0$} if $g_-^{\mathrm{ext}}(y)=c_-(-y)^{-q}(1+o(1))$ as $y\to0^-$ for some $c_->0$.
		
		\item \label{def:finite-angle-continuation}
		a \emph{finite-angle continuation} with angular limit $\ell<0$ if it has an extension which is nonincreasing on some interval $(-\delta,0)$ and satisfies $y^{-1}g_-^{\mathrm{ext}}(y)\to\ell$ as $y\to0^-$. This continuation is called \emph{flat} if, in addition, there exists a constant $\mu>0$ such that $g_-^{\mathrm{ext}}(y)=\ell y+O(|y|^{1+\mu})$ as $y\to0^-$.
	\end{enumerate}
\end{definition}

In a finite-angle continuation, $y^{-1}g_-^{\mathrm{ext}}(y)\to\ell<0$ gives $c|y|\leq g_-^{\mathrm{ext}}(y)\leq C|y|$ near $y=0^-$, so the extension reaches the phase origin in the same manner as a regular branch with finite-angle zero output. In the Hessian quotients of Subsubsection \ref{subsubsec:finite-angle-models} the continuation leaves the selected sheet at $y_{\mathrm{reg}}$, which is the first case of the following lemma.

\begin{lemma}\label{lem:finite-angle-continuation-angle}
	Let $(I_-,g_-)$ be a regularity-limited branch with a finite-angle continuation $g_-^{\mathrm{ext}}$ of angular limit $\ell$, and assume that $0\in A_g$. Then $\ell\leq\inf A_g<-\rho_0$. Moreover, if $y_{\mathrm{dom}}=0$, then $g_-^{\mathrm{ext}}=g^\sharp(\cdot,1)$ on $(-\infty,0)$.
\end{lemma}

\begin{proof}
	Firstly, we note that, since $A_g$ is an open interval containing $-\rho_0$ and $0$, it contains $[-\rho_0,0]$ and $\inf A_g<-\rho_0$. Let $a(y)=y^{-1}g_-^{\mathrm{ext}}(y)$, which is negative because $g_-^{\mathrm{ext}}>0$ and $y<0$.
	
	Now, let us assume that $y_{\mathrm{dom}}<0$, and let $y\in(y_{\mathrm{dom}},0)$. If $(g_-^{\mathrm{ext}}(y),y)$ belonged to $-U_g$, then by \ref{hyp:signed-homogeneity} and \ref{hyp:signed-oddness} we would have $\Psi_g(a(y))=-(-y)^{-\alpha}$. This would imply that $(y,1)\in\mathcal D^\sharp$ and $y<y_{\mathrm{dom}}$, which is impossible. Therefore, $a(y)\notin A_g$, and since $a(y)<0$ and $[-\rho_0,0]\subset A_g$, we obtain $a(y)\leq\inf A_g$. Then, by letting $y\to0^-$, we conclude that $\ell\leq\inf A_g$.
	
	Next, let us assume that $y_{\mathrm{dom}}=0$. Then $g^\sharp(\cdot,1)$ is defined on $(-\infty,0)$, and it is continuous because $g^\sharp(y,z)=-g(-y,-z)$ and $g$ is continuous on $\mathcal D$ by Proposition \ref{prop:conical-sheet-maximal-chart}. Moreover, its graph lies in $-U_g\subset\widehat\Gamma$ by \ref{hyp:signed-domain} and satisfies $\widehat\gamma=1$, and therefore it is a connected subset of $\{\widehat\gamma=1\}$. We note that this graph contains the graph of $g_-$, since the two functions coincide on $I_-$. Consequently, it lies in the component containing that graph, whose part in $\{y<0\}$ is the graph of $g_-^{\mathrm{ext}}$ by Definition \ref{def:lower-continuation}. Since both functions are defined on $(-\infty,0)$, and a graph has exactly one point over each ordinate, we have $g_-^{\mathrm{ext}}=g^\sharp(\cdot,1)$. Finally, since $\Psi_g(A_g)\supset(-\infty,0)$, by \eqref{eq:reflected-lower-angular-relation} we obtain $a(y)\to\inf A_g$ as $y\to0^-$, and we conclude that $\ell=\inf A_g$.
\end{proof}

The cubic model of Subsubsection \ref{subsubsec:ellipticity-limited-models} shows that an ellipticity-limited branch may admit neither of these two continuations. In that model the component turns at a point where $\partial_x\widehat\gamma$ vanishes, and therefore it has a vertical tangent there and ceases to be a graph over the rotational variable.

For a curvature-blow-up continuation, the challenge is to reach the positive branch. Indeed, since the meridional curvature diverges as $y\to0^-$, the continuation cannot be joined to $g_+$ through $y=0$. The junction is made instead with a second branch of the level-one set, defined on $0<y<y_*$, which has the same blow-up at $y=0^+$ and reaches the diagonal at $y_*$, as in the following condition.

\begin{definition}\label{def:singular-positive-matching}
	We say that a curvature-blow-up continuation of order $q$ satisfies the \emph{singular positive-matching condition} if there exists a function $g_{\mathrm{sing}}\in\mathcal C^1((0,y_*))$ such that, for every $0<y<y_*$,
	\begin{equation*}
		\widetilde\gamma(g_{\mathrm{sing}}(y),y)=1,\qquad\partial_x\widetilde\gamma(g_{\mathrm{sing}}(y),y)>0,\qquad g_{\mathrm{sing}}(y)>y.
	\end{equation*}
	In addition, we require $g_{\mathrm{sing}}(y)=c_+y^{-q}(1+o(1))$ as $y\to0^+$ for some $c_+>0$, and $g_{\mathrm{sing}}(y)\to y_*$ as $y\to y_*^-$.
\end{definition}

We emphasize that the two matching conditions apply to different profiles. On the one hand, the singular profile is defined on $0<y<y_*$, where the phase coordinates are positive and $\widehat\gamma=\widetilde\gamma$. On the other hand, the regular profile passes through $(x_0,0)$, which need not belong to $\overline{\widetilde\Gamma}$ and is reached only through the signed extension. Moreover, the order $q$ is the same on both sides of $y=0$ and measures how fast the meridional curvature diverges as $y\to0$. Finally, by Lemma \ref{lem:matching-c1-junctions}\textup{(ii)}, $g_{\mathrm{sing}}$ joins $g_+$ in a $\mathcal C^1$ manner at $y_*$.

\subsection{Model phase diagrams}\label{subsec:model-examples-branches}
We treat each model in the subsubsection of the lower-branch type that it realizes. In particular, none of the alternatives of Definition \ref{def:lower-continuation-outcomes} is empty. Within each example, $\Psi_g$, $g$ and $g_-$ denote the angular profile, the selected chart and the neck-accessible branch of the model under consideration. Throughout, $A_i=\binom{n-1}{i}$ and $\rho_i=\frac{A_i}{A_{i-1}}=\frac{n-i}{i}$, which is strictly decreasing in $i$. The models are the mean curvature $\frac{H}{n-1}$ and its powers $\left(\frac{H}{n-1}\right)^p$ with $p>1$, the Hessian quotients $\mathcal Q_{k,j}$ with $k-j$ odd, the polynomial models $\frac{S_k}{A_k}$ with $k$ odd and their real $k$-th roots, and, for $n\geq3$, the models $\gamma_{\mathrm{tri}}$ and $\gamma_\eps$.

Here $1\leq j<k\leq n-1$, and we take $k\leq n-1$ so that the cylindrical direction stays inside $\Gamma_k$. Moreover,
\begin{align*}
	\gamma_{\mathrm{tri}}&=\frac{2\left(S_1^3+S_1S_2\right)}{(n-1)^2(3n-4)},\qquad\mathcal Q_{k,j}=\frac{A_j}{A_k}\,\frac{S_k}{S_j},
	\\
	\gamma_\eps&=\frac{1}{c_\eps}\left(H-\frac{\eps|A|^2}{H+|A|}\right),\qquad c_\eps=n-1-\frac{\eps(n-1)}{n-1+\sqrt{n-1}},
\end{align*}
with $0<\eps\leq1$, and all of them are normalized so that $\gamma(0,1,\ldots,1)=1$.

Except for $\gamma_{\mathrm{tri}}$ and $\gamma_\eps$, these are the curvature functions of the literature recalled in Subsection \ref{subsec:intro-background}, although, as far as we know, the catenoidal translators of the Hessian quotients have not been considered before. The models $\gamma_{\mathrm{tri}}$ and $\gamma_\eps$ are ours, and we use them to realize the behaviors of the lower branch that the previous models do not exhibit. Let us check that each model is a $1$-nondegenerate $\alpha$-positive curvature function on its cone. For the mean curvature and its powers the cone is $\{H>0\}$, where the conditions of Definition \ref{def:alpha-positive-curvature-function} are immediate. The models built from $S_k$ are defined on the G{\aa}rding cone $\Gamma_k$, where the claim follows from the identity $S_i(x,y,\ldots,y)=A_{i-1}y^{i-1}(x+\rho_iy)$ and Maclaurin's inequality on $\overline{\Gamma_k}$ \cite[Corollary 2.8 and Lemma 2.10]{spruck2005geometric}, together with the ellipticity of $\mathcal Q_{k,j}$ on $\Gamma_k$ \cite[Theorem 2.16]{spruck2005geometric}. We also note that Maclaurin's inequality extends $\mathcal Q_{k,j}$ continuously by zero to the points of $\partial\Gamma_k$ where $S_j$ vanishes. Finally, $\gamma_{\mathrm{tri}}$ and $\gamma_\eps$ are defined on $\Gamma_2$ and $\{H>|A|\}$, respectively, and we verify the claim in their subsubsections by using $|A|^2=S_1^2-2S_2$. In contrast, $S_n/S_j$ tends to zero at the direction $(0,1)$, so it is $1$-degenerate.

The parity restrictions on $k-j$ and $k$ come from the oddness condition \ref{hyp:signed-oddness}. Indeed, since $S_i(-\lambda)=(-1)^iS_i(\lambda)$, this condition holds for $\mathcal Q_{k,j}$ precisely when $k-j$ is odd, and for $\frac{S_k}{A_k}$ precisely when $k$ is odd. In the even cases the function is invariant under central reflection, and therefore the signed construction does not apply to them.

Finally, for the models on G{\aa}rding cones the lower side of a neck cannot be described with $\gamma$ itself, which is the reason for the signed extension. This is the content of the following lemma.
\begin{lemma}\label{lem:lower-side-outside-garding}
	Let $2\leq k\leq n-1$, let $y<0$ and let $a=xy^{-1}$. Then $(x,y,\ldots,y)\notin\Gamma_k$, whereas $-(x,y,\ldots,y)\in\Gamma_k$ precisely when $a>-\rho_k$. In particular, if $(A_g,\rho_0)$ are selected angular neck data for a curvature function on $\Gamma_k$ and $a\leq-\rho_0$, then neither direction belongs to $\Gamma_k$.
\end{lemma}

\begin{proof}
	Firstly, we prove that $(x,y,\ldots,y)\notin\Gamma_k$. By the identity $S_i(x,y,\ldots,y)=A_{i-1}y^{i-1}(x+\rho_iy)$ and the condition $y<0$, the inequality $S_1>0$ reads $x>\rho_1|y|$ and $S_2>0$ reads $x<\rho_2|y|$. These inequalities are incompatible because $\rho_1>\rho_2$. Therefore, since $\Gamma_k\subseteq\Gamma_2$, the direction $(x,y,\ldots,y)$ does not belong to $\Gamma_k$.
	
	Next, we consider the opposite direction. By the same identity, we have $S_i(-(x,y,\ldots,y))=A_{i-1}|y|^{i-1}\left(\rho_i|y|-x\right)$, which is positive for $i=1,\ldots,k$ precisely when $x<\rho_i|y|$ for all such $i$, or equivalently when $x<\rho_k|y|$, because $\rho_i$ decreases in $i$. Then, by substituting $x=ay$ and $|y|=-y$, this condition becomes $a>-\rho_k$.
	
	Finally, we prove the last assertion. Since $S_i(a,1,\ldots,1)=A_{i-1}(a+\rho_i)$, the angular slice of $\widetilde\Gamma_k$ is $(-\rho_k,\infty)$, and by Lemma \ref{lem:neck-cone-boundary} we have $-\rho_0\leq-\rho_k$. Consequently, $a\leq-\rho_0$ implies $a\leq-\rho_k$.
\end{proof}

\subsubsection{Axis-crossing models}\label{subsubsec:axis-crossing-models}
The mean curvature and its powers realize the first alternative, with homogeneities $1$ and $p$. We extend their rotational restrictions $\frac{x+(n-1)y}{n-1}$ and its $p$-th power to the signed sheet as the odd functions
\begin{equation*}
	\widehat\gamma(x,y)=\operatorname{sgn}(x+(n-1)y)\left|\frac{x+(n-1)y}{n-1}\right|^{p},\qquad
	\Psi_g(a)=\operatorname{sgn}(a+n-1)\frac{|a+n-1|^{p}}{(n-1)^{p}},
\end{equation*}
with $p=1$ for the mean curvature itself. Both have $A_g=\rr$, the single zero direction $-\rho_0=1-n$, the chart $g(y,z)=(n-1)\left(\operatorname{sgn}(z)|z|^{\frac1p}-y\right)$ and the level-one line $x+(n-1)y=n-1$ for every $p$. This line contains the positive branch $g_+(y)=(n-1)(1-y)$ on $I_+=\left(\frac{n-1}{n},1\right)$, its continuation $g(y,1)$ on $(1,\infty)$ used by the upper side, and the reflected branch $g_-(y)=(n-1)(1-y)$ on $I_-=(-\infty,0)$, with $y_{\mathrm{reg}}=0$ and $x_0=n-1$.

The signed profile of the mean curvature is linear, and therefore elliptic on $\rr^2\setminus\{(0,0)\}$. Consequently, the selected data are regular with $0\in A_g$. Moreover, the cone is $\Gamma=\{H>0\}$, whose slice is $(1-n,\infty)$, and the neck circle lies on $\partial\Gamma$ by Lemma \ref{lem:neck-cone-boundary}.

Since the reflected branch tends to $n-1$ as $y\to0^-$, it is axis-crossing. Moreover, since the signed profile is linear, the same line realizes the regular positive-matching condition of Definition \ref{def:regular-positive-matching}. It continues through $(n-1,0)$, stays above the diagonal and reaches it at $y_*=\frac{n-1}{n}$, where it agrees with $g_+$. This is the phase diagram of the smooth translating catenoid of Clutterbuck--Schn\"urer--Schulze \cite{clutterbuck2007stability}.

When $p>1$, the cone, the position of the neck circle and the matching are those of the mean curvature, since the signed profile is of class $\mathcal C^1$ near $(n-1,0)$, where $x+(n-1)y\neq0$. The exponent $p$ only changes the regularity at the zero ray. Indeed, we have $\Psi_g'(1-n)=0$, and therefore condition \ref{hyp:regular-ellipticity} fails on the zero ray. Moreover, we notice that the zero ray is the only direction where it fails, since away from it we have $\partial_x\widehat\gamma=\frac{p|x+(n-1)y|^{p-1}}{(n-1)^p}>0$. Then, since the profile vanishes to order $p$ at the zero ray, Lemma \ref{lem:neck-regularity-vanishing-order} tells us that every radius is neck-admissible and is $\frac1p$-regular at the neck. We notice that the graphical construction of Sheng--Wu \cite[Corollary 5.1]{sheng2010rotationally} corresponds precisely to this odd extension.

\subsubsection{Models with finite-angle zero output}
\label{subsubsec:finite-angle-models}

We observe that the Hessian quotients $\mathcal Q_{k,j}$, with $k-j$ odd, are the models that realize the second alternative of Theorem \ref{thm:rotational-translators}, where the lower phase reaches the origin at a finite angle. We notice that their degree of homogeneity is $\alpha=k-j$, and we write $\kappa=\frac{\rho_j}{\rho_k}>1$. In addition, depending on the dimension, which enters through $\kappa=\frac{k(n-j)}{j(n-k)}$, these examples present the following dichotomy. For $\alpha=1$, and for $\alpha\geq2$ with $\kappa<\left(\frac{\alpha+1}{\alpha-1}\right)^2$, the reflected branch has finite-angle zero output, whereas otherwise it is ellipticity-limited and reaches the origin through a finite-angle continuation. By the identity $S_i(x,y,\ldots,y)=A_{i-1}y^{i-1}(x+\rho_iy)$, we have
\begin{equation*}
	\widehat\gamma(x,y)=\kappa\,y^{\alpha}\,\frac{x+\rho_ky}{x+\rho_jy},\qquad
	\Psi_g(a)=\kappa\,\frac{a+\rho_k}{a+\rho_j},\qquad
	\Psi_g'(a)=\frac{\kappa(\rho_j-\rho_k)}{(a+\rho_j)^{2}},
\end{equation*}
and in particular $\Psi_g'>0$ away from $-\rho_j$.

\paragraph*{Profile geometry and exact charts:} We notice that the angular profile has a simple pole at $-\rho_j$ and a unique zero at $-\rho_k$, and that it is real analytic and strictly increasing on $(-\rho_j,\infty)$ with range $(-\infty,\kappa)$. Consequently, the zero direction is $-\rho_0=-\rho_k$, whereas the selected interval $A_g$ depends on the regime and is specified below. Moreover, by Lemma \ref{lem:lower-side-outside-garding}, the slice of $\widetilde\Gamma_k$ is $(-\rho_k,\infty)$, and therefore the neck circle lies on $\partial\Gamma_k$ by Lemma \ref{lem:neck-cone-boundary}.

By solving the angular equation $\Psi_g(a)=zy^{-\alpha}$, we obtain the chart and the reflected level-one curve in closed form,
\begin{equation}\label{eq:quotient-chart-closed-form}
	g(y,z)=\rho_j\,y\,\frac{z-y^{\alpha}}{\kappa y^{\alpha}-z},\qquad g_-(y)=-\rho_j\,y\,\frac{1+(-y)^{\alpha}}{1+\kappa(-y)^{\alpha}},\quad y<0,
\end{equation}
which is consistent with \eqref{eq:zero-trace-conical}. The chart is real analytic on $\{(y,z):y>0,\ zy^{-\alpha}<\kappa\}$, where its denominator does not vanish, and therefore every radius is neck-admissible and $1$-regular at the neck. The curve \eqref{eq:quotient-chart-closed-form} is real analytic on $(-\infty,0)$. By \eqref{eq:reflected-lower-angular-relation} its angular value tends to $-\rho_j$ as $y\to0^-$, with $g_-(y)=-\rho_jy+O\left(|y|^{1+\alpha}\right)$. In other words, this curve reaches the origin flatly with $\ell=-\rho_j$ and $\mu=\alpha$.

\paragraph*{The ellipticity threshold and two regimes:} The quotient $\mathcal Q_{k,j}$ is not always elliptic on the whole sheet $(-\rho_j,\infty)$. Indeed, $\Psi_g'>0$ there, but the second reduced ellipticity expression of Lemma \ref{lem:reduced-ellipticity} is
\begin{align*}
	\alpha\Psi_g(a)-a\Psi_g'(a)&=\frac{\kappa}{(a+\rho_j)^{2}}\,N_2(a),\\
	N_2(a)&=\alpha a^2+\left[\alpha(\rho_j+\rho_k)-(\rho_j-\rho_k)\right]a+\alpha\rho_j\rho_k.
\end{align*}
In particular, its sign is that of $N_2$, which may change. In fact, the discriminant of $N_2$ is $\rho_k^2(\kappa-1)\left[\kappa(\alpha-1)^2-(\alpha+1)^2\right]$, and therefore $N_2$ is strictly positive on $\rr$ if and only if
\begin{equation}\label{eq:quotient-ellipticity-threshold}
	\kappa=\frac{\rho_j}{\rho_k}=\frac{k(n-j)}{j(n-k)}<\left(\frac{\alpha+1}{\alpha-1}\right)^2=\left(\frac{k-j+1}{k-j-1}\right)^2.
\end{equation}
This condition depends on the dimension through $\kappa$. We call the case where \eqref{eq:quotient-ellipticity-threshold} holds the regular regime, and the case where it fails the ellipticity-limited regime.
\begin{enumerate}[label=\textup{(\roman*)}]
	\item \textbf{Regular regime.} Let us assume that \eqref{eq:quotient-ellipticity-threshold} holds, which is always the case when $\alpha=1$. Then the selected interval is $A_g=(-\rho_j,\infty)$ and the data are regular. Moreover, since $\Psi_g(A_g)$ is unbounded below, we have $y_{\mathrm{dom}}=y_{\mathrm{reg}}=0$, and the branch itself reaches the phase origin. In particular, this regime covers every adjacent quotient $\mathcal Q_{k,k-1}$ with $2\leq k\leq n-1$, for which $\kappa=\frac{1}{\tau_k}$ with $\tau_k=\frac{A_kA_{k-2}}{A_{k-1}^{2}}$, and \eqref{eq:quotient-chart-closed-form} reduces to $g_-(y)=\frac{\rho_ky(1-y)}{y-\tau_k}$.
	
	\item \textbf{Ellipticity-limited regime.} Let us assume that \eqref{eq:quotient-ellipticity-threshold} fails. Then $N_2$ has roots $a_1\leq a_2$ in $(-\rho_j,-\rho_k)$, because it is positive at both endpoints and its vertex lies between them. Consequently, we have $\partial_y\widehat\gamma\leq0$ for $a_1\leq a\leq a_2$, while $\partial_x\widehat\gamma>0$ on $(-\rho_j,\infty)$. In this regime the selected interval is $A_g=(a_2,\infty)$, and the branch is ellipticity-limited at $y_{\mathrm{reg}}=y_{\mathrm{dom}}=-\left(-\Psi_g(a_2)\right)^{-\frac1\alpha}$. The rational curve \eqref{eq:quotient-chart-closed-form} is a finite-angle continuation in the sense of Definition \ref{def:finite-angle-continuation}.
\end{enumerate}

In the ellipticity-limited regime, with $s=(-y)^\alpha$ we have $\left(g_-^{\mathrm{ext}}\right)'=-\rho_j N(s)(1+\kappa s)^{-2}$ with $N(s)=\kappa s^2+\left[1+\kappa-\alpha(\kappa-1)\right]s+1$, whose discriminant is $\rho_k^{-2}$ times that of $N_2$. Therefore, $N$ has positive roots $s_1\leq s_2$ and $y_{\mathrm{reg}}=-s_2^{\frac1\alpha}$. The continuation increases on $\left(-s_2^{\frac1\alpha},-s_1^{\frac1\alpha}\right)$ and is nonincreasing on $\left(-s_1^{\frac1\alpha},0\right)$, an interval that becomes all of $(-\infty,0)$ in the equality case $a_1=a_2$.

\paragraph*{Dimensional dependence:} Because $\kappa$ decreases to $\frac{k}{j}$ as $n\to\infty$, the quotient is regular in every sufficiently large dimension. It is ellipticity-limited in finitely many low dimensions when $\frac{k}{j}<\left(\frac{k-j+1}{k-j-1}\right)^2$, and ellipticity-limited in every dimension otherwise (as is the case for $\mathcal Q_{k,1}$ with $k\geq4$ even). For instance, for $\mathcal Q_{5,2}$ (where $\alpha=3$ and $\kappa=\frac{5(n-2)}{2(n-5)}$), the regularity threshold $\kappa<4$ holds for dimensions $n\geq11$, becomes an equality for $n=10$, and fails for $6\leq n\leq9$.

\subsubsection{A model with infinite-angle zero output}\label{subsubsec:fractional-model}
We realize the third alternative by perturbing the mean curvature with the quotient $\frac{|A|^2}{H+|A|}$. The cone of $\gamma_\eps$ is $\Gamma_\star=\{\lambda:H(\lambda)>|A|(\lambda)\}$, which contains $\Gamma_+$ and the cylindrical direction. The bounds $H-\frac{\eps|A|^2}{H+|A|}\geq\frac{H}{2}$ and $\partial_i\left(H-\frac{\eps|A|^2}{H+|A|}\right)\geq1-\frac{\eps}{2}$ hold strictly on $\Gamma_\star$. Therefore, $\gamma_\eps$ is $1$-positive on $\Gamma_\star$.

We define its signed profile on $y>0$ and extend it to the reflected sheet by central oddness, so that conditions \ref{hyp:signed-domain}--\ref{hyp:signed-oddness} hold. On $y>0$ the signed profile is
\begin{equation*}
	\widehat\gamma(x,y)=\frac{1}{c_\eps}\left[x+(n-1)y-\frac{\eps\left(x^2+(n-1)y^2\right)}{x+(n-1)y+\sqrt{x^2+(n-1)y^2}}\right],
\end{equation*}
and, writing $V=\sqrt{a^2+n-1}$, $T=a+n-1+V$ and $E=3a^2+2a(n-1)+(n-1)^2-2a+(2a+2n-3)V$,
\begin{align*}
	\Psi_g(a)&=\frac{1}{c_\eps}\left(a+n-1-\frac{\eps(a^2+n-1)}{T}\right),\\
	\Psi_g'(a)&=\frac{1}{c_\eps}\left[1-\frac{\eps\left(2aT-V(V+a)\right)}{T^{2}}\right],\\
	\Psi_g(a)-a\Psi_g'(a)&=\frac{n-1}{c_\eps T^{2}}\left[T^2-\eps\left(V+n-1+2a-a^2\right)\right]\\
	&=\frac{n-1}{c_\eps T^{2}}\left[(1-\eps)T^2+\eps E\right].
\end{align*}
The chart $g$ is not elementary.

\paragraph*{Profile geometry and ellipticity:} The quantity $T$ is positive for every $a$. Indeed, this is clear when $a+n-1\geq0$, and otherwise we have $V^2-(a+n-1)^2=(n-1)(2-n-2a)>0$, which implies $V>-(a+n-1)$. Consequently, $\widehat\gamma$ is smooth on $\{y\neq0\}$. It remains to check both reduced ellipticity inequalities of Lemma \ref{lem:reduced-ellipticity} on $\rr$. For $a>-\frac{n-2}{2}$ the identity $H^2-|A|^2=(n-1)(2a+n-2)$ places $(a,1,\ldots,1)$ in $\Gamma_\star$, and the inequalities follow from the ellipticity of $\gamma_\eps$. For $a<0$ we use the formulas above. Since $2aT<0$ and $V+a>0$, we have $\Psi_g'\geq c_\eps^{-1}$. When $2a+2n-3\geq0$ we have $E\geq2a^2+(a+n-1)^2-2a>0$. When $2a+2n-3<0$, the bound $0<V\leq-a+\sqrt{n-1}$ bounds $E$ from below by a quadratic in $a$ with discriminant $4n-3-4(n-1)^2-8(n-1)\sqrt{n-1}<0$ for $n\geq3$.

\paragraph*{Neck data and asymptotic branch:} As $a\to-\infty$, by the expansion $V=-a-\frac{n-1}{2a}+O\left(|a|^{-3}\right)$ we have $T\to n-1$ and $\Psi_g(a)=-\frac{\eps}{(n-1)c_\eps}(-a)^2(1+o(1))$. Therefore, $\Psi_g$ increases strictly from $-\infty$ on $A_g=\rr$, with a unique zero direction $-\rho_0$, and the selected data are regular with $0\in A_g$. The slice of $\widetilde\Gamma_\star$ is $\left(-\frac{n-2}{2},\infty\right)$, and at its endpoint we have $H=|A|=\frac n2$, which implies $\Psi_g\left(-\frac{n-2}{2}\right)=\frac{1}{c_\eps}\left(1-\frac{\eps}{2}\right)\frac n2>0$. Consequently, $-\rho_0<-\frac{n-2}{2}$, and the neck circle lies outside $\overline{\Gamma_\star}$, which is the second case of Lemma \ref{lem:neck-cone-boundary}.

Since the data are regular and $\Psi_g(A_g)$ is unbounded below, we have $y_{\mathrm{reg}}=y_{\mathrm{dom}}=0$. Then, by Proposition \ref{prop:angular-growth-lower-branch} with $\alpha=1$ and $p_\infty=2$, we obtain
\begin{equation*}
	g_-(y)=\sqrt{\frac{(n-1)c_\eps}{\eps}}\,(-y)^{\frac{1}{2}}(1+o(1)) \qquad \text{as }y\to0^-,
\end{equation*}
and therefore the branch has fractional zero output with $\sigma=\frac{1}{2}$.

\subsubsection{Ellipticity-limited models}\label{subsubsec:ellipticity-limited-models}
Three models realize the fourth alternative. The polynomial model and its odd root continue with a curvature blow-up, whereas $\gamma_{\mathrm{tri}}$ admits no continuation. We assume $1<k\leq n-1$ is odd.

\paragraph*{The polynomial quotient and its root:} For $\gamma_k=\frac{S_k}{A_k}$, by the identity for $S_k$ we have
\begin{equation*}
	\widehat\gamma(x,y)=\frac{y^{k-1}(x+\rho_ky)}{\rho_k},\qquad
	\Psi_g(a)=\frac{a+\rho_k}{\rho_k},\qquad
	g(y,z)=\rho_k\left(zy^{1-k}-y\right),
\end{equation*}
and for the odd root $\left(\frac{S_k}{A_k}\right)^{\frac1k}$ the real $k$-th roots of the first two, with chart $g(y,z)=\rho_k\left(z^ky^{1-k}-y\right)$. Both charts are consistent with \eqref{eq:zero-trace-conical}, the zero direction is $-\rho_0=-\rho_k$ in both cases, and both have the reflected branch $g_-(y)=\rho_k\left((-y)^{1-k}-y\right)$.

The data of $\gamma_k=\frac{S_k}{A_k}$ are not regular on all of $\rr$. While $\partial_x\widehat\gamma(ay,y)=\frac{y^{k-1}}{\rho_k}>0$ for $y>0$, the derivative $\partial_y\widehat\gamma(ay,y)=\frac{y^{k-1}\bigl((k-1)a+k\rho_k\bigr)}{\rho_k}$ is positive precisely for $a>-\frac{k\rho_k}{k-1}$, which defines the left endpoint of $A_g$. As for $\mathcal Q_{k,j}$, the neck circle lies on $\partial\Gamma_k$.

The branch is halted by the finiteness of $A_g$ on the left. The range $\Psi_g(A_g)=\left(-\frac{1}{k-1},\infty\right)$ contains the reflected level $-(-y)^{-k}$ precisely when $(-y)^k>k-1$. Since the data are regular there, the maximal interval is $I_-=\left(-\infty,-(k-1)^{\frac1k}\right)$, at whose endpoint $g_-$ takes the finite value $\rho_k\left((k-1)^{\frac{1-k}{k}}+(k-1)^{\frac1k}\right)$. At this endpoint, the derivative $g_-'$ vanishes, and $\partial_y\widehat\gamma=0$ whereas $\partial_x\widehat\gamma>0$. Therefore, the branch is ellipticity-limited.

Beyond the regular endpoint, the same formula $g_-^{\mathrm{ext}}(y)=\rho_k\left((-y)^{1-k}-y\right)$ on all of $(-\infty,0)$ has $\partial_x\widehat\gamma>0$ and $\partial_y\widehat\gamma<0$. Since $g_-^{\mathrm{ext}}(y)=\rho_k(-y)^{-(k-1)}(1+o(1))$ as $y\to0^-$, it constitutes a curvature-blow-up continuation of order $q=k-1$. On the positive side, $g_{\mathrm{sing}}(y)=\rho_k\left(y^{1-k}-y\right)$ defined on $0<y<y_*$ (with $y_*=\left(\frac{\rho_k}{1+\rho_k}\right)^{\frac1k}=\left(\frac{n-k}{n}\right)^{\frac1k}$) stays above the diagonal, has the same blow-up order at $y=0^+$, and reaches $g_{\mathrm{sing}}(y_*)=y_*$. Therefore, the singular positive-matching condition holds, which allows the continued phase to join $g_+$.

After undoing the normalization via a spatial rescaling, this phase diagram underlies the odd-$k$ translating catenoids constructed by de Lima--Pipoli \cite[Theorem 1(ii) and Remark 4]{lima2025translators} and, for $S_k$ and $S_k^{1/k}$ in warped products, by de Lira--de Farias \cite[Theorem 4(c)]{lira2026translating}. The even-$k$ families of \cite[Theorem 1(iii)]{lima2025translators} and \cite[Theorem 4(b) and (d)]{lira2026translating} belong to a different phase configuration.

\paragraph*{The rooted model:} Since $t=1$ is a fixed point of the real root, the odd root $\left(\frac{S_k}{A_k}\right)^{\frac1k}$ has the same level-one set as $\gamma_k$, and along it its partial derivatives are $\frac1k$ times those of $\gamma_k$. With the same interval $A_g=\left(-\frac{k\rho_k}{k-1},\infty\right)$, it therefore has the same maximal regular interval, ellipticity-limited endpoint, curvature-blow-up continuation and singular positive matching. For $k>1$ its profile vanishes to order $p=\frac1k$ at the zero ray $-\rho_k$ in the sense of Lemma \ref{lem:neck-regularity-vanishing-order}, and therefore the selected data are not regular. Nevertheless, its chart is smooth near $(R^{-1},0)$, and every radius is neck-admissible and $1$-regular at the neck.

\paragraph*{The symmetric cubic model $\gamma_{\mathrm{tri}}$:} The function $\gamma_{\mathrm{tri}}$ is a positive multiple of $S_1^3+S_1S_2$, which is elliptic on $\Gamma_2$. Indeed, we have $S_1^2=|A|^2+2S_2>|A|^2\geq\lambda_i^2$ there, which implies $|\lambda_i|<S_1$ and $\partial_i\left(S_1^3+S_1S_2\right)=4S_1^2+S_2-S_1\lambda_i>3S_1^2+S_2>0$. Writing $\nu_\pm=\frac{1}{2}\left(-3(n-1)\pm\sqrt{(n-1)(3n-1)}\right)$, we have
\begin{align*}
	\widehat\gamma(x,y)&=\frac{(x+(n-1)y)\left(2x^2+6(n-1)xy+(n-1)(3n-4)y^2\right)}{(n-1)^2(3n-4)},\\
	\Psi_g(a)&=\frac{(a+n-1)\left(2a^2+6(n-1)a+(n-1)(3n-4)\right)}{(n-1)^2(3n-4)},\\
	3\Psi_g(a)-a\Psi_g'(a)&=\frac{(4a+3n-4)(2a+3n-3)}{(n-1)(3n-4)},
\end{align*}
with the three zero directions $\nu_-<1-n<\nu_+<0$. Since the level-one equation is a cubic in $x$, the functions $g$ and $g_-$ are not elementary.

The factorization of $3\Psi_g-a\Psi_g'$ is positive exactly for $a>-\frac{3n-4}{4}$. The two critical points of $\Psi_g$ lie below this value, as shown below for the larger one, and therefore $\Psi_g'>0$ there. Consequently, we select $A_g=\left(-\frac{3n-4}{4},\infty\right)$, which contains $0$. Only $\nu_+$, among the three roots, lies in $A_g$, since $4\nu_++3n-4>0$ reduces to $n(3n-4)>0$ while $\nu_-<1-n<-\frac{3n-4}{4}$. The slice of $\widetilde\Gamma_2$ is $\left(-\frac{n-2}{2},\infty\right)$, and $-\nu_+>\frac{n-2}{2}$ reduces to $n^2>0$. Consequently, the neck circle lies outside $\overline\Gamma_2$ by Lemma \ref{lem:neck-cone-boundary}.

Along the neck-accessible branch, \eqref{eq:reflected-lower-angular-relation} drives the angular value from $\nu_+$ down to the endpoint $-\frac{3n-4}{4}$ of $A_g$. This endpoint is reached at the level $\Psi_g\left(-\frac{3n-4}{4}\right)=-\frac{n^2}{32(n-1)^2}$, that is, at $y_{\mathrm{reg}}=y_{\mathrm{dom}}=-\left(\frac{32(n-1)^2}{n^2}\right)^{\frac{1}{3}}$. There, we have $\partial_y\widehat\gamma=0$ and $\partial_x\widehat\gamma>0$, and therefore the branch is ellipticity-limited.

Beyond that endpoint, the level-one component continues as a curve on which the meridional value increases, but its angular value decreases only until it reaches the critical point $a_{\mathrm{crit}}=\frac{1}{6}\left(-8(n-1)+\sqrt{2(n-1)(5n-2)}\right)$ of $\Psi_g$. This point lies in $\left(1-n,-\frac{3n-4}{4}\right)$, where $\Psi_g'=0$ and hence $\partial_x\widehat\gamma=0$ by \eqref{eq:reduced-ellipticity-identities}. There, the component attains its largest rotational value, acquires a vertical tangent, and turns back toward $y\to-\infty$. Consequently, $\gamma_{\mathrm{tri}}$ admits neither a curvature-blow-up nor a finite-angle continuation.

Table \ref{tab:model-summary} collects the data of the seven models. Its last column records the type of the neck-accessible branch and its continuation, from which the case of Theorem \ref{thm:rotational-translators} realized by each model is read.

\begin{table}[htbp]
	\centering
	\scriptsize
	\setlength{\tabcolsep}{3pt}
	\renewcommand{\arraystretch}{1.3}
	\begin{tabular}{@{}lclllp{3.2cm}@{}}
		\hline
		$\gamma$ & $\alpha$ & $A_g$ & $-\rho_0$ & $y_{\mathrm{reg}}$ & neck-accessible branch\\
		\hline
		$\frac{H}{n-1}$ & $1$ & $\rr$ & $1-n$ & $0$ & axis-crossing; regular positive matching\\
		$\left(\frac{H}{n-1}\right)^p$, $p>1$ & $p$ & $\rr$ & $1-n$ & $0$ & axis-crossing; regular positive matching; $\Psi_g'=0$ at the zero ray\\
		$\mathcal Q_{k,j}$, \eqref{eq:quotient-ellipticity-threshold} holds & $k-j$ & $(-\rho_j,\infty)$ & $-\rho_k$ & $0$ & finite-angle, $\ell=-\rho_j$, through a regular branch\\
		$\mathcal Q_{k,j}$, \eqref{eq:quotient-ellipticity-threshold} fails & $k-j$ & $(a_2,\infty)$ & $-\rho_k$ & $-\left(-\Psi_g(a_2)\right)^{-\frac1\alpha}$ & ellipticity-limited; finite-angle continuation, $\ell=-\rho_j$\\
		$\gamma_\eps$ & $1$ & $\rr$ & $-\rho_0<-\frac{n-2}{2}$ & $0$ & fractional zero output, $\sigma=\frac{1}{2}$\\
		$\frac{S_k}{A_k}$, $k$ odd & $k$ & $\left(-\frac{k\rho_k}{k-1},\infty\right)$ & $-\rho_k$ & $-(k-1)^{\frac1k}$ & ellipticity-limited; curvature blow-up of order $k-1$, singular positive matching\\
		$\left(\frac{S_k}{A_k}\right)^{\frac1k}$, $k$ odd & $1$ & $\left(-\frac{k\rho_k}{k-1},\infty\right)$ & $-\rho_k$ & $-(k-1)^{\frac1k}$ & as above; $\Psi_g$ not $\mathcal C^1$ at the zero ray\\
		$\gamma_{\mathrm{tri}}$ & $3$ & $\left(-\frac{3n-4}{4},\infty\right)$ & $\nu_+$ & $-\left(\frac{32(n-1)^2}{n^2}\right)^{\frac{1}{3}}$ & ellipticity-limited; no continuation\\
		\hline
	\end{tabular}
	\caption{\scriptsize The models of Subsection \ref{subsec:model-examples-branches}, with $k-j$ odd for $\mathcal Q_{k,j}$ and $n\geq3$ for $\gamma_\eps$ and $\gamma_{\mathrm{tri}}$. The mean curvature and its powers fall in case \textup{(i)} of Theorem \ref{thm:rotational-translators}, $\mathcal Q_{k,j}$ in case \textup{(ii)}, $\gamma_\eps$ in case \textup{(iii)}, and the last three in case \textup{(iv)}.}
	\label{tab:model-summary}
\end{table}
\section{The first-order equation: barriers and phase dynamics}\label{sec:barrier-method}

The barrier method formalizes the study of the equation satisfied by the slope $v=u'$, where $u(r)$ is the height of the vertical graph $F(r,\omega)=(r\omega,u(r))$ of Section \ref{sec:preliminaries}. More precisely, in this section we prove the existence and uniqueness of its solutions along a branch of the level-one equation, a comparison principle together with the barriers that we use with it, and the behavior of the phase at the endpoints of the lower branch. We use these tools in Sections \ref{sec:positive-branch-asymptotics} and \ref{sec:catenoidal-translators}.

\subsection{The graphical equation}\label{subsec:graphical-equation}
Throughout this section $\alpha>0$, and we set $\beta=\frac{\alpha-1}{2\alpha}$, so that $1-2\beta=\alpha^{-1}>0$. Recall that the normalized translator equation \eqref{eq:vertical-graph-translator} is a first-order equation for $v$. At each radius the slope determines the phase point
\begin{equation*}
	y_v(r) = Q(r,v(r)) \quad \text{and} \quad x_v(r) = \frac{v'(r)}{(1+v(r)^2)^{\beta+1}},\quad\text{where}\quad Q(r,s) = \frac{s}{r(1+s^2)^\beta},
\end{equation*}
whose coordinates are the meridional and rotational principal curvatures of the graph multiplied by the common factor $(1+v^2)^{\frac{1}{2}-\beta}$. The variable $y_v$ is the one used by Rengaswami for the rotational translator equation \cite[Section 2.3]{rengaswami2024classification}, where level functions of this type also serve as barriers for the bowl-type solution \cite[Propositions 3.1--3.2, Lemma 6.1 and Section 7.1]{rengaswami2024classification}. Since we use both objects throughout the article, we collect the properties of $Q$ and of the level functions in Appendix \ref{app:q-eta-properties}, in our notation and with their sources.

The comparison and barrier arguments of this section use only $\alpha>0$, through the monotonicity of the rotational phase in the slope. The same holds for the expansions of Theorem \ref{thm:positive-bowl-asymptotics}, where the threshold $\alpha>\frac{1}{2}$ of Rengaswami's classification \cite[Theorems 1.2--1.4]{rengaswami2024classification} enters only through the critical regime, since $k_\gamma=2\alpha-1>0$ there. In contrast, for the height constants of Theorem \ref{thm:rotational-translators} we need $4\alpha>1$, or $3\alpha>1$ when $\widetilde\gamma$ is only of class $\mathcal C^2$ near $(0,1)$, and in Section \ref{sec:catenoidal-translators} we assume $\alpha>\frac{1}{2}$. Finally, for Theorem \ref{thm:bowl-asymptotic-uniqueness}\textup{(a)} and Theorem \ref{thm:vertical-cylinder-graphs} we do not need any bound on $\alpha$ beyond the bowl being entire, while for Theorem \ref{thm:bowl-asymptotic-uniqueness}\textup{(b)} we need in addition the ranges stated in it.

\begin{definition}\label{def:level-one-branch}
	A \emph{branch of the level-one equation} is a pair $(I,g)$, where $I\subset\rr$ is an open interval and $g\in\mathcal C^1(I)$.
\end{definition}

Given a branch $(I,g)$, a slope $v$ whose phase lies on the branch satisfies \eqref{eq:vertical-graph-translator} if and only if
\begin{equation}\label{eq:selected-branch-graphical-equation}
	v' = (1+v^2)^{\beta+1}g(Q(r,v)), \qquad Q(r,v)\in I,
\end{equation}
that is, if and only if $x_v=g(y_v)$ at every radius. For a reflected branch, $\widetilde\gamma$ is replaced by $\widehat\gamma$ in \eqref{eq:vertical-graph-translator}. We call \eqref{eq:selected-branch-graphical-equation} the graphical equation of $(I,g)$, and for $(I_-,g_-)$ and $(I_+,g_+)$ we call it the lower and the positive graphical equation, respectively. In addition, we denote its right-hand side by $\mathsf F(r,s)=(1+s^2)^{\beta+1}g(Q(r,s))$, which is defined on the open set $U=\{(r,s):r>0,\ Q(r,s)\in I\}$.

\begin{definition}\label{def:branch-admissible-functions}
	Let $J\subset(0,\infty)$ be an interval and let $v\in\mathcal C^1(J)$. The function $v$ is said to be \emph{branch-admissible} for $(I,g)$ on $J$ if $Q(r,v(r))\in I$ for every $r\in J$. Equivalently, a branch-admissible solution of \eqref{eq:selected-branch-graphical-equation} is a solution of $v'=\mathsf F(r,v)$ whose graph is contained in $U$.
\end{definition}

In what follows, we use the positive branch $(I_+,g_+)$ of Proposition \ref{prop:positive-branch} and the maximal regular continuation $(I_-,g_-)$ of Definition \ref{def:neck-accessible-branch}, with $I_-=(-\infty,y_{\mathrm{reg}})$.

\begin{proposition}\label{prop:branch-ivp}
	Let $(I,g)$ be a branch of the level-one equation, and let $r_0>0$ and $v_0\in\rr$ satisfy $Q(r_0,v_0)\in I$.
	\begin{enumerate}[label=\textup{(\roman*)}]
		\item There exists a unique solution $v\in\mathcal C^2((r_-,r_+))$ to \eqref{eq:selected-branch-graphical-equation} satisfying the initial condition $v(r_0)=v_0$, defined on a maximal branch-admissible interval $(r_-,r_+)$ with $0\leq r_-<r_0<r_+\leq\infty$.
		
		\item Suppose $g$ admits a continuous extension to a compact interval $K\subset\overline I$. If $r_+<\infty$ and $y_v([r_0,r_+))\subset K$, then $v$ extends to a function of class $\mathcal C^1([r_0,r_+])$ and the limit $y_v(r_+)$ is an endpoint of $I$. Analogously, if $r_->0$ and $y_v((r_-,r_0])\subset K$, the same conclusions hold at $r_-$.
	\end{enumerate}
\end{proposition}

\begin{proof}
Firstly, we note that, since $Q$ is smooth on $\{r>0\}$ and $g\in\mathcal C^1(I)$, the set $U$ is open and $\mathsf F\in\mathcal C^1(U)$. In particular, $\mathsf F$ is locally Lipschitz in $s$, uniformly in $r$, on compact subsets of $U$. Therefore, item \textup{(i)} follows from the Picard--Lindel\"of theorem and the existence of a maximal solution in $U$ \cite[Theorems 2.2 and 2.13]{teschl2012ordinary}. Moreover, from the equation $v'=\mathsf F(r,v)$ we obtain $v\in\mathcal C^2((r_-,r_+))$.

Next, we prove \textup{(ii)} at $r_+$, since the argument at $r_-$ is the same. Let $K=[k_1,k_2]$. Since $y_v\in K$ on $[r_0,r_+)$, by Lemma~\ref{lem:level-set-control}\textup{(ii)} we have $w_{k_1}\leq v\leq w_{k_2}$ there. Recall that both level functions are continuous on $[r_0,r_+]$. Therefore, $v$ is bounded, and $v'=(1+v^2)^{\beta+1}g(y_v)$ is bounded because $g$ is bounded on $K$. Consequently, $v$ is Lipschitz on $[r_0,r_+)$ and extends continuously to $r_+$. Moreover, by the equation, $v'$ extends as well, and we obtain $v\in\mathcal C^1([r_0,r_+])$ and $y_v(r_+)\in K$. Finally, let us assume by contradiction that $y_v(r_+)$ is an interior point of $I$. Then $(r_+,v(r_+))\in U$, and by the extensibility criterion \cite[Lemma 2.14]{teschl2012ordinary}, used similarly in \cite[Section 7.1]{rengaswami2024classification}, we can continue $v$ beyond $r_+$, which contradicts the maximality of $r_+$.
\end{proof}

\subsection{Comparison and level barriers}\label{subsec:residuals-comparison}\label{subsec:level-barriers}
We note that the comparison principle for \eqref{eq:selected-branch-graphical-equation} requires only that $g$ be nonincreasing along the phases being compared. In fact, this holds for the selected branch $g_+$ by Proposition \ref{prop:positive-branch} and for the selected branch $g_-$ by \eqref{eq:lower-branch-monotonicity}. In addition, it holds for a finite-angle continuation on its terminal interval, whatever its behavior elsewhere on $(-\infty,0)$.

\begin{definition}\label{def:residual-barriers}
	Let $J\subset(0,\infty)$ be an interval and let $w\in\mathcal C^1(J)$ be branch-admissible for $(I,g)$ on $J$. The \emph{residual} of $w$ is
	\begin{equation*}
		\mathcal R[w] = w' - (1+w^2)^{\beta+1}g(Q(r,w)).
	\end{equation*}
	The function $w$ is called a \emph{supersolution} or a \emph{subsolution} on $J$ depending on whether $\mathcal R[w]\geq0$ or $\mathcal R[w]\leq0$, respectively. The strict inequalities $\mathcal R[w]>0$ and $\mathcal R[w]<0$ define \emph{strict} supersolutions and \emph{strict} subsolutions. To make the underlying branch explicit, we write $\mathcal R_+[w]$ and $\mathcal R_-[w]$ for the residuals evaluated on $(I_+,g_+)$ and $(I_-,g_-)$, respectively.
\end{definition}

\begin{proposition}\label{prop:comparison-principle-ode}
	Let $J=[r_0,r_1)\subset(0,\infty)$, and suppose $v,w\in\mathcal C^1(J)$ are branch-admissible for the same branch $(I,g)$ of the level-one equation, where $g$ is nonincreasing on $I$. Assume that
	\begin{equation*}
		\mathcal R[v] \leq 0 \leq \mathcal R[w] \qquad \text{on } J,
	\end{equation*}
	and that $v(r_0)\leq w(r_0)$. Then $v\leq w$ on $J$. Furthermore, if $v(r_0)<w(r_0)$ and at least one of the residual inequalities is strict, then $v<w$ on $J$.
\end{proposition}

\begin{proof}
Let $G(s)=\int_0^s(1+t^2)^{-\beta-1}\,dt$ and $d=G(w)-G(v)$. Since $G'>0$, the function $d$ is of class $\mathcal C^1$ on $J$ and has the sign of $w-v$. Then, by the residual hypotheses, we have
	\begin{equation}
		\label{eq:comparison-derivative-bound}
		d'(r)
		=
		\frac{w'(r)}{(1+w(r)^2)^{\beta+1}} - \frac{v'(r)}{(1+v(r)^2)^{\beta+1}}
		\geq
		g(Q(r,w(r))) - g(Q(r,v(r))).
	\end{equation}
Firstly, we prove that $v\leq w$ on $J$. Let us assume by contradiction that $v(r_2)>w(r_2)$ for some $r_2\in J$. Since $d(r_0)\geq0>d(r_2)$, there is $r_3\in(r_0,r_2)$ with $d(r_3)<0$ and $d'(r_3)<0$. However, at $r_3$ we have $v>w$, and therefore $Q(r_3,v(r_3))>Q(r_3,w(r_3))$. Then, since $g$ is nonincreasing, the right-hand side of \eqref{eq:comparison-derivative-bound} is nonnegative, which contradicts $d'(r_3)<0$.

Next, we prove the strict inequality. Let us assume by contradiction that $v$ and $w$ meet in $J$, and let $r_4$ be the first radius where they do. Then $d>0$ on $[r_0,r_4)$ and $d(r_4)=0$, and consequently $d'(r_4)\leq0$. Moreover, at $r_4$ the right-hand side of \eqref{eq:comparison-derivative-bound} vanishes, and from the strict residual inequality we obtain $d'(r_4)>0$, leading to the desired contradiction.
\end{proof}

Recall that the level functions $w_m$ of Lemma \ref{lem:level-set-control}\textup{(ii)} are defined by $Q(r,w_m(r))=m$. We use them to turn slope bounds into phase bounds, since $w_{m_1}(r)\leq v(r)\leq w_{m_2}(r)$ if and only if $m_1\leq y_v(r)\leq m_2$. For $m\in I$, we notice that the function $w_m$ is branch-admissible for $(I,g)$ on $(0,\infty)$, and by \eqref{eq:wm-derivative-level} its residual is
\begin{equation}\label{eq:level-barrier-residual-section}
	\mathcal R[w_m]=(1+w_m^2)^{\beta+1}\left[\frac{m}{1+\alpha^{-1}w_m^2} - g(m)\right].
\end{equation}
Then, since $\alpha^{-1}=1-2\beta>0$, the sign of $\mathcal R[w_m]$ coincides with that of the bracket. Consequently, this allows us to enclose the phase of a solution between two levels. Indeed, if $v$ solves \eqref{eq:selected-branch-graphical-equation} on $J=[r_0,r_1)$, and $m_1<m_2$ in $I$ are such that $w_{m_1}$ is a subsolution and $w_{m_2}$ is a supersolution on $J$, with $m_1\leq y_v(r_0)\leq m_2$, then by Proposition \ref{prop:comparison-principle-ode} we obtain $m_1\leq y_v\leq m_2$ on all of $J$. When the level is a finite endpoint of $I$, we have to use the following variant of this argument.

\begin{proposition}\label{prop:endpoint-barrier}
	Let $b$ be a finite endpoint of $I$, and assume that $g$ extends continuously to $I\cup\{b\}$. Then, if $v\in\mathcal C^1([r_0,r_1])$ is a solution to \eqref{eq:selected-branch-graphical-equation} for the extended branch, whose phase level $y_v$ takes values in $I\cup\{b\}$, the following holds. If $b$ is the right endpoint, $y_v(r_0)<b$, and $\mathcal R[w_b]>0$ on $[r_0,r_1]$, then $y_v<b$ on $[r_0,r_1]$. Analogously, if $b$ is the left endpoint, $y_v(r_0)>b$, and $\mathcal R[w_b]<0$ on $[r_0,r_1]$, then $y_v>b$ on $[r_0,r_1]$.
\end{proposition}

\begin{proof}
We prove the statement for the right endpoint, since the argument for the left endpoint is analogous. Let us assume by contradiction that $y_v$ attains $b$, and let $\bar r\in(r_0,r_1]$ be the first radius with $y_v(\bar r)=b$. Then $v<w_b$ on $[r_0,\bar r)$ and $v(\bar r)=w_b(\bar r)$, and therefore $(v-w_b)'(\bar r)\geq0$. However, at $\bar r$, by the graphical equation and the definition of the residual, we have
	\begin{equation*}
		v'(\bar r) = (1+w_b(\bar r)^2)^{\beta+1}g(b),
		\qquad
		w_b'(\bar r) = \mathcal R[w_b](\bar r) + (1+w_b(\bar r)^2)^{\beta+1}g(b),
	\end{equation*}
and consequently $(v-w_b)'(\bar r)=-\mathcal R[w_b](\bar r)<0$, leading to the desired contradiction.
\end{proof}

In particular, the level functions at the finite endpoints of $I_+$ are barriers in this sense.

\begin{corollary}\label{cor:positive-branch-endpoint-barriers}
	Let $v\in\mathcal C^1([r_0,r_1])$ be a solution to the positive graphical equation for the extended branch, and assume its initial phase level $y_v(r_0)$ lies in the interior of $I_+$. Then $y_v>y_*$ on $[r_0,r_1]$. Furthermore, if $\gamma$ is normalized $1$-nondegenerate, then $y_v<1$ on $[r_0,r_1]$.
\end{corollary}

\begin{proof}
Firstly, we consider the left endpoint of $I_+$, which is $y_*=\widetilde\gamma(1,1)^{-\frac{1}{\alpha}}$. Recall that $g_+$ extends continuously to this endpoint with $g_+(y_*)=y_*$ by Proposition~\ref{prop:positive-branch}. Then, by \eqref{eq:level-barrier-residual-section}, we have
	\begin{equation*}
		\mathcal R_+[w_{y_*}]=(1+w_{y_*}^2)^{\beta+1}\,y_*\left[\frac{1}{1+\alpha^{-1}w_{y_*}^2}-1\right]<0.
	\end{equation*}
Therefore, by Proposition \ref{prop:endpoint-barrier}, we obtain $y_v>y_*$. Next, we consider the normalized $1$-nondegenerate case, where the right endpoint is $y=1$, since $\widetilde\gamma(0,1)=1$, and $g_+(1)=0$ by Lemma \ref{lem:positive-endpoint-regularity}. In this case we have
	\begin{equation*}
		\mathcal R_+[w_1]=\frac{(1+w_1^2)^{\beta+1}}{1+\alpha^{-1}w_1^2}>0,
	\end{equation*}
and Proposition \ref{prop:endpoint-barrier} implies that $y_v<1$.
\end{proof}

In contrast, on the lower branch the level functions are only one-sided barriers. Indeed, for $m\in I_-$ we have $m<0<g_-(m)$ by \eqref{eq:lower-branch-positivity}. Therefore, the bracket in \eqref{eq:level-barrier-residual-section} is negative, and $w_m$ is a strict subsolution of the lower graphical equation on $(0,\infty)$.

\subsection{Phase identity and endpoint dynamics}\label{subsec:level-barriers-negative-components}\label{subsec:regular-lower-endpoint-dynamics}
\begin{proposition}\label{prop:graphical-phase-evolution}
	Let $(I,g)$ be a branch of the level-one equation, and let $v$ be a branch-admissible solution to the associated graphical equation \eqref{eq:selected-branch-graphical-equation} on an interval $J\subset(0,\infty)$.
	\begin{enumerate}[label=\textup{(\roman*)},ref=\thetheorem\,\textup{(\roman*)}]
		\item On $J$, the phase evolution is governed by
		\begin{equation}\label{eq:graphical-phase-evolution}
			r y_v' = \left(1+\alpha^{-1}v^2\right)g(y_v)-y_v.
		\end{equation}
		
		\item Suppose $(I,g)$ is either the maximal regular lower branch $(I_-,g_-)$, or $I=(-\infty,0)$ and $g$ is an extension of $g_-$ in the sense of Definition \ref{def:lower-continuation}. Then $ry_v'>-y_v>0$ on $J$. In particular, the rotational phase $y_v$ is strictly increasing on $J$.
	\end{enumerate}
\end{proposition}

\begin{proof}
Firstly, by the chain rule, we have $y_v'=\partial_sQ(r,v)\,v'+\partial_rQ(r,v)$. Then, by substituting \eqref{eq:q-partial-s-section}, \eqref{eq:q-partial-r-section}, the identity $1-2\beta=\alpha^{-1}$ and the graphical equation, we obtain \eqref{eq:graphical-phase-evolution}. Next, under the hypotheses of \textup{(ii)} we have $y_v<0$ and $g(y_v)>0$, by \eqref{eq:lower-branch-positivity} on $I_-$ and by Definition \ref{def:lower-continuation} for an extension. Consequently, from \eqref{eq:graphical-phase-evolution} we obtain $ry_v'>-y_v>0$.
\end{proof}

In the finite-angle regime, we obtain the decay rate from the power barriers of Subsection \ref{subsec:power-barriers-finite-angle}.

\begin{lemma}\label{lem:continuation-phase-limits}
	Let $(I,g)$ be a branch of the level-one equation, and let $v$ be a branch-admissible solution of \eqref{eq:selected-branch-graphical-equation} on its maximal interval $[r_0,r_+)$.
	\begin{enumerate}[label=\textup{(\roman*)}]
		\item If $I=(a,b)$ with $b\leq0$ and $g>0$, then $y_v$ is strictly increasing and $y_v(r)\to b$ as $r\to r_+$. If moreover $g\geq c$ on $(b-\delta,b)$ for some $c,\delta>0$, then $r_+<\infty$.
		
		\item If $I=(a,b)$ with $a\geq1$ and $g<0$, then $y_v$ and $ry_v$ are strictly decreasing, so $y_v(r)<r_0y_v(r_0)r^{-1}$ for $r>r_0$. Moreover, $r_+<\infty$ and $y_v(r)\to a$ as $r\to r_+$.
		
		\item If $r_+<\infty$ and $y_v(r)\to0$ as $r\to r_+$, then $v(r)\to0$.
	\end{enumerate}
\end{lemma}

\begin{proof}
\textup{(i)} Since $y_v<b\leq0$ and $g(y_v)>0$, from \eqref{eq:graphical-phase-evolution} we have $ry_v'\geq g(y_v)-y_v>0$. Therefore, $y_v$ increases to a limit $L\leq b$. Let us assume by contradiction that $L<b$. Then $y_v([r_0,r_+))$ lies in the compact subinterval $[y_v(r_0),L]$ of $I$. On the one hand, if $r_+<\infty$, then by Proposition \ref{prop:branch-ivp}\textup{(ii)} the limit $L$ is an endpoint of $I$, which is impossible. On the other hand, if $r_+=\infty$, then $g\geq c_0>0$ on $[y_v(r_0),L]$. In this case we have $ry_v'\geq c_0$, and consequently $y_v(r)\geq y_v(r_0)+c_0\log\frac{r}{r_0}\to\infty$, which is also impossible. Therefore, $L=b$. Finally, if $g\geq c$ on $(b-\delta,b)$, the same logarithmic growth on the final part of the trajectory rules out $r_+=\infty$.

\textup{(ii)} Since $y_v>a\geq1$ and $g(y_v)<0$, from \eqref{eq:graphical-phase-evolution} we have $ry_v'<-y_v<0$. Therefore, $y_v$ is strictly decreasing and $(ry_v)'=ry_v'+y_v<0$, which implies the stated bound. Moreover, if $r_+=\infty$, then by this bound we would have $y_v<a$ at some radius. Consequently, $r_+<\infty$, and $y_v$ decreases to a limit $L\geq a$. Finally, if $L>a$, then $y_v([r_0,r_+))\subset[L,y_v(r_0)]$, and by Proposition \ref{prop:branch-ivp}\textup{(ii)} the limit $L$ would be an endpoint of $I$, which is impossible. Therefore, $L=a$.

\textup{(iii)} Since $r\leq r_+$, from the definition of $y_v$ we have $v(1+v^2)^{-\beta}=ry_v\to0$, and by Lemma \ref{lem:level-set-control}\textup{(i)} we conclude that $v(r)\to0$.
\end{proof}

\begin{proposition}\label{prop:zero-output-radius-criterion}
	Assume $y_{\mathrm{reg}}=0$, so that $I_-=(-\infty,0)$, and let $v$ be a branch-admissible solution to the lower graphical equation on its maximal interval $[r_0,r_+)$. Then $r_+<\infty$ if and only if
	\begin{equation*}
		\int_{-\delta}^{0}\frac{dy}{g_-(y)} < \infty
	\end{equation*}
	for some (equivalently, for every) $\delta>0$. If $r_+=\infty$, the slope satisfies $v(r)\to0^-$ as $r\to\infty$. In particular, the maximal radius $r_+$ is finite for axis-crossing and fractional zero-output branches, whereas $r_+=\infty$ for branches with finite-angle zero output.
\end{proposition}

\begin{proof}
Firstly, by \eqref{eq:lower-branch-positivity} and Lemma \ref{lem:continuation-phase-limits}\textup{(i)} with $b=0$, the phase $y=y_v$ increases to $0$ as $r\to r_+$. We note that, since $g_-$ is continuous and positive on $I_-$, the integrability of $1/g_-$ at $0^-$ does not depend on $\delta$. Moreover, on this branch \eqref{eq:graphical-phase-evolution} reads
	\begin{equation*}
		ry' = \left(1+\alpha^{-1}v^2\right)g_-(y)+|y|.
	\end{equation*}

Now, let us assume by contradiction that the integral is finite and $r_+=\infty$. Since $1+\alpha^{-1}v^2\geq1$, we have $ry'\geq g_-(y)$. Then, by integrating $y'/g_-(y)\geq r^{-1}$ from a radius $r_1$ with $y(r_1)>-\delta$, we obtain
	\begin{equation*}
		\int_{y(r_1)}^{y(r)}\frac{dy}{g_-(y)}\geq\log\frac{r}{r_1}\longrightarrow\infty\qquad\text{as } r\to\infty,
	\end{equation*}
a contradiction. Therefore, $r_+<\infty$.

Conversely, let us assume that $r_+<\infty$. Then $v(r)\to0$ by Lemma \ref{lem:continuation-phase-limits}\textup{(iii)}, and therefore $1+\alpha^{-1}v^2\leq2$ near $r_+$. Moreover, we have $|y|<\rho_0^{-1}g_-(y)$ by \eqref{eq:lower-branch-positivity}. Consequently, $ry'\leq\left(2+\rho_0^{-1}\right)g_-(y)$, that is,
	\begin{equation*}
		\frac{y'}{g_-(y)}\leq\frac{2+\rho_0^{-1}}{r}\leq\frac{2+\rho_0^{-1}}{r_0}.
	\end{equation*}
Then, by integrating over $[r_1,r)$ and letting $r\to r_+$, we obtain
	\begin{equation*}
		\int_{y(r_1)}^{0}\frac{dy}{g_-(y)}\leq\frac{\left(2+\rho_0^{-1}\right)(r_+-r_1)}{r_0}<\infty,
	\end{equation*}
and the integral is finite.

Next, we prove the limit of the slope when $r_+=\infty$. In this case, from the graphical equation we have $v'>0$, and therefore $v$ increases to a limit $v_\infty\leq0$. Let us assume by contradiction that $v_\infty<0$. Then
	\begin{equation*}
		y=Q(r,v) \leq Q(r,v_\infty) = -\frac{c}{r},\qquad\text{where } c=\frac{|v_\infty|}{(1+v_\infty^2)^\beta}>0,
	\end{equation*}
and by \eqref{eq:lower-branch-positivity} we have $g_-(y)>\rho_0|y|\geq\rho_0cr^{-1}$. This implies that $v'\geq\rho_0cr^{-1}$, whose integral diverges, which contradicts $v<0$. Consequently, $v(r)\to0^-$.

It remains to apply the criterion to the three types with $y_{\mathrm{reg}}=0$. Firstly, for an axis-crossing branch, $1/g_-$ extends continuously to $y=0$ because $g_-\to x_0>0$. Moreover, for a fractional branch, $g_-(y)^{-1}=d^{-1}(-y)^{-\sigma}(1+o(1))$ is integrable at $0^-$ because $\sigma<1$. In contrast, for a finite-angle branch, $g_-(y)^{-1}=|\ell|^{-1}|y|^{-1}(1+o(1))$ is not integrable at $0^-$.
\end{proof}

\begin{corollary}\label{cor:radius-criterion-on-continuations}
	Let $v$ be a branch-admissible solution of the lower graphical equation along a finite-angle continuation $g_-^{\mathrm{ext}}$. Then $r_+=\infty$ and $v(r)\to0^-$ as $r\to\infty$.
\end{corollary}

\begin{proof}
Since $y^{-1}g_-^{\mathrm{ext}}(y)\to\ell<0$, we have $c|y|\leq g_-^{\mathrm{ext}}(y)\leq C|y|$ near $0^-$. Therefore, the proof of Proposition \ref{prop:zero-output-radius-criterion} applies with the lower bound $c|y|\leq g_-^{\mathrm{ext}}(y)$ in place of \eqref{eq:lower-branch-positivity}.
\end{proof}

\begin{proposition}\label{prop:regular-lower-endpoint-dynamics}
	Let $v$ be a branch-admissible solution to the lower graphical equation \eqref{eq:selected-branch-graphical-equation} on its maximal interval $[r_0,r_+)$, and let $y=y_v$ denote its phase level.
	\begin{enumerate}[label=\textup{(\roman*)}]
		\item If $(I_-,g_-)$ is regularity-limited in the sense of Definition \ref{def:lower-continuation-outcomes}\textup{(iv)}, then $r_+<\infty$. Furthermore, both $v$ and $v'$ admit finite one-sided limits at $r_+$, and the phase point satisfies $(x_v(r),y_v(r))\to(x_{\mathrm{reg}},y_{\mathrm{reg}})$ as $r\to r_+^-$.
		
		\item If $(I_-,g_-)$ is axis-crossing, where $x_0>0$ denotes its endpoint limit from Definition \ref{def:lower-continuation-outcomes}\textup{(i)}, then $r_+<\infty$. Moreover, as $r\to r_+^-$, we have $y(r)\to0$, $v(r)\to0$, and $v'(r)\to x_0$.
		
		\item If $(I_-,g_-)$ has fractional zero output, where $d$ and $\sigma$ are the parameters from Definition \ref{def:lower-continuation-outcomes}\textup{(iii)}, then $r_+<\infty$ and the phase level satisfies the asymptotic expansion
		\begin{equation}\label{eq:fractional-phase-endpoint-asymptotics}
			-y(r)=\left(\frac{(1-\sigma)d}{r_+}(r_+-r)\right)^{\frac{1}{1-\sigma}}(1+o(1))\qquad \text{as } r\to r_+^-.
		\end{equation}
		In particular, $v(r)\to0$, $v'(r)\to0$, and $(x_v(r),y_v(r))\to(0,0)$ as $r\to r_+^-$.
	\end{enumerate}
\end{proposition}

\begin{proof}
Firstly, we prove a bound that we use in the three cases. By Proposition \ref{prop:graphical-phase-evolution}\textup{(ii)}, we have $y'>-yr^{-1}$, and therefore $w=-y>0$ satisfies $(\log w)'<-r^{-1}$. Consequently,
	\begin{equation}\label{eq:lower-phase-radial-bound}
		w(r) < \frac{r_0w(r_0)}{r}
	\end{equation}
as long as the phase stays on the lower branch.

\textup{(i)} Recall that branch admissibility requires $w>-y_{\mathrm{reg}}$. Since the right-hand side of \eqref{eq:lower-phase-radial-bound} falls below $-y_{\mathrm{reg}}$ at the radius $r_0w(r_0)/(-y_{\mathrm{reg}})$, we have $r_+<\infty$. Then, by Lemma \ref{lem:continuation-phase-limits}\textup{(i)} with $b=y_{\mathrm{reg}}$, we obtain $y(r)\to y_{\mathrm{reg}}$. Therefore, $x_v=g_-(y)\to x_{\mathrm{reg}}$, which is finite and positive by Proposition \ref{prop:general-lower-phase-alternative}. This means that $g_-$ extends continuously to $K=[y(r_0),y_{\mathrm{reg}}]$, and $K$ contains $y([r_0,r_+))$ because $y$ is increasing. Consequently, by Proposition \ref{prop:branch-ivp}\textup{(ii)}, we have $v\in\mathcal C^1([r_0,r_+])$, and from the graphical equation we obtain $v'(r)\to(1+v(r_+)^2)^{\beta+1}x_{\mathrm{reg}}$.

\textup{(ii)} Since $g_-(y)\to x_0>0$, we have $g_-\geq\frac{x_0}{2}$ near $y=0$. Then, by Lemma \ref{lem:continuation-phase-limits}\textup{(i)} with $b=0$, we obtain $r_+<\infty$ and $y\to0$. Moreover, $v\to0$ by Lemma \ref{lem:continuation-phase-limits}\textup{(iii)}, and from the graphical equation we conclude that $v'\to x_0$.

\textup{(iii)} By Lemma \ref{lem:continuation-phase-limits}\textup{(i)} with $b=0$, we have $w\to0$. We note that its last assertion does not apply, because $g_-$ vanishes at $0$. For this reason, we prove $r_+<\infty$ directly. For small $w$, by the expansion $g_-(-w)=dw^\sigma(1+o(1))$ and \eqref{eq:graphical-phase-evolution}, we have
	\begin{equation*}
		w'=-r^{-1}\left[\left(1+\alpha^{-1}v^2\right)dw^\sigma(1+o(1))+w\right]\leq-\frac{d}{2r}w^\sigma,
	\end{equation*}
since $1+\alpha^{-1}v^2\geq1$ and $w=o(w^\sigma)$ as $w\to0$. Therefore, $(w^{1-\sigma})'=(1-\sigma)w^{-\sigma}w'\leq-\frac{(1-\sigma)d}{2r}$. Since the integral of $r^{-1}$ diverges, $w^{1-\sigma}$ would become negative if $r_+=\infty$, and consequently $r_+<\infty$. We note that this is the separable computation of \cite[Lemma 2.5]{rengaswami2024classification}, applied to a vanishing rather than a diverging quantity. Now, by Lemma \ref{lem:continuation-phase-limits}\textup{(iii)}, we have $v\to0$, and therefore $1+\alpha^{-1}v^2\to1$. Then, by integrating up to $r_+$, we obtain
	\begin{align*}
		w'(r)&=-d\,r_+^{-1}w(r)^\sigma(1+o(1)),\qquad\text{hence}
		\\
		w(r)^{1-\sigma}&=(1-\sigma)d\,r_+^{-1}(r_+-r)(1+o(1)),
	\end{align*}
which is \eqref{eq:fractional-phase-endpoint-asymptotics}. Moreover, since $(1+v^2)^{-\beta}=1+o(1)$, from the definition of $y_v$ we have $v(r)=-r_+w(r)(1+o(1))$. Finally, $v'=(1+v^2)^{\beta+1}g_-(-w)=dw^\sigma(1+o(1))\to0$ and $x_v=g_-(-w)\to0$, and we conclude that $(x_v,y_v)\to(0,0)$.
\end{proof}

\subsection{Power barriers for finite-angle zero output}\label{subsec:power-barriers-finite-angle}
Recall that a level function $w_m$ keeps its phase at the level $m$, and therefore it cannot bound a phase that tends to $0^-$. For this reason, for branches with finite-angle zero output we use negative powers of the radius, whose phase also tends to $0^-$. We note that the argument uses the branch only near the phase origin, and consequently it also covers a finite-angle continuation. In what follows, we write $\overline g$ for $g_-$ or for $g_-^{\mathrm{ext}}$ on $(-\infty,0)$, so that $y^{-1}\overline g(y)\to\ell<0$ as $y\to0^-$.

\begin{proposition}\label{prop:power-barrier-flat-component}
	Let $\overline w(r)=-c_0r^\chi$ with constants $c_0>0$ and $\chi<0$, and let $t(r)=c_0r^{\chi-1}(1+c_0^2r^{2\chi})^{-\beta}$.
	\begin{enumerate}[label=\textup{(\roman*)},ref=\thetheorem\,\textup{(\roman*)}]
		\item The phase of $\overline w$ is $Q(r,\overline w(r))=-t(r)$, so $\overline w$ is branch-admissible for $\left((-\infty,0),\overline g\right)$ on $(0,\infty)$. Moreover, $t$ decreases to $0$, and the residual is
		\begin{equation}\label{eq:power-barrier-residual}
			\mathcal R[\overline w]=-c_0r^{\chi-1}\left[\chi-(1+\overline w^2)\frac{\overline g(-t(r))}{-t(r)}\right].
		\end{equation}
		
		\item Let $\eps,\delta>0$ satisfy $\ell+\eps<0$ and $\left|y^{-1}\overline g(y)-\ell\right|<\eps$ for $y\in(-\delta,0)$. If $\chi<\ell-\eps$, then $\overline w$ is a strict supersolution on the half-line where $t(r)<\delta$ and $\overline w(r)^2\leq\frac{\chi}{\ell-\eps}-1$. If $\ell+\eps<\chi<0$, it is a strict subsolution on the half-line where $t(r)<\delta$.
		
		\item For every $c_0>0$, the function $\overline w$ is a strict supersolution if $\chi<\ell$ and a strict subsolution if $\ell<\chi<0$, at all sufficiently large radii.
	\end{enumerate}
\end{proposition}

\begin{proof}
Firstly, we prove \textup{(i)}. The identity $Q(r,\overline w(r))=-t(r)$ is the definition of the phase. Then, by differentiating, we obtain
	\begin{equation*}
		\frac{t'(r)}{t(r)}=\frac{1}{r}\left[\chi\left(1-\frac{2\beta\,\overline w^2}{1+\overline w^2}\right)-1\right]<-\frac{1}{r},
	\end{equation*}
since $1-2\beta=\alpha^{-1}>0$ makes the factor in parentheses positive and $\chi<0$. Therefore, $t$ is strictly decreasing, and $t(r)\to0$ because $\chi<0$. Moreover, since $t>0$, the phase lies in $(-\infty,0)$. Finally, by substituting $\overline w'=-c_0\chi r^{\chi-1}$ into the residual, we obtain \eqref{eq:power-barrier-residual}.

Next, we prove \textup{(ii)}. Let $\Lambda(r)=\frac{\overline g(-t(r))}{-t(r)}$. We note that where $t(r)<\delta$ we have $\Lambda\in(\ell-\eps,\ell+\eps)$, an interval of negative numbers, and that the sign of $\mathcal R[\overline w]$ is opposite to that of the bracket in \eqref{eq:power-barrier-residual}. On the one hand, if $\ell+\eps<\chi<0$, then $(1+\overline w^2)\Lambda<\ell+\eps<\chi$. Therefore, the bracket is positive and $\overline w$ is a strict subsolution. On the other hand, if $\chi<\ell-\eps$, then $(1+\overline w^2)\Lambda>(1+\overline w^2)(\ell-\eps)$. In this case the bracket is negative as soon as $(1+\overline w^2)(\ell-\eps)\geq\chi$, that is, $\overline w^2\leq\frac{\chi}{\ell-\eps}-1$. We note that this bound is positive because $\chi<\ell-\eps<0$, and $\overline w$ is then a strict supersolution. Moreover, since $t$ and $\overline w^2=c_0^2r^{2\chi}$ decrease to $0$, each set of conditions holds on a half-line.

Finally, item \textup{(iii)} follows from \textup{(ii)} with $\eps<|\chi-\ell|$ and $\ell+\eps<0$.
\end{proof}

\subsection{Continuations past the regular endpoint}\label{subsec:post-elliptic-phase-dynamics}
Recall that an axis-crossing branch continues through the axis under the regular positive-matching condition. In contrast, past a regularity-limited endpoint the extended branch need not be monotone. In this case the comparison principle no longer applies, but the phase identity \eqref{eq:graphical-phase-evolution} remains valid. In this subsection we treat the curvature-blow-up continuation and its singular matching, and we leave the finite-angle continuation for Section \ref{sec:catenoidal-translators}, where we apply the power barriers of Subsection \ref{subsec:power-barriers-finite-angle} on its terminal interval.

\begin{proposition}\label{prop:axis-crossing-positive-matching}
	Suppose that the branch $(I_-,g_-)$ is axis-crossing and satisfies the regular positive-matching condition of Definition \ref{def:regular-positive-matching}. Let $v$ be a solution to the lower graphical equation on its maximal branch-admissible interval $[r_0,r_+)$. Then $v$ extends uniquely past $r_+$ as a positive solution, whose phase passes through $(x_0,0)$ and reaches the diagonal point $(y_*,y_*)$ at a finite radius. From this radius onward, the phase point $(x_v,y_v)$ lies on the positive branch $g_+$.
\end{proposition}

\begin{proof}
Firstly, we continue the solution through the axis. By Proposition \ref{prop:regular-lower-endpoint-dynamics}\textup{(ii)}, we have $r_+<\infty$, $v(r_+)=0$ and $v'(r_+)=x_0>0$. Moreover, by Lemma \ref{lem:matching-c1-junctions}\textup{(i)}, the concatenation of $g_-$ and $g_{\mathrm{ax}}$ is a branch of class $\mathcal C^1$ on a two-sided neighborhood of $y=0$, and $v$ solves its graphical equation on a left neighborhood of $r_+$. Then, by Proposition \ref{prop:branch-ivp}\textup{(i)} at the datum $(r_+,0)$, whose phase level is $0$, we can continue $v$ uniquely past $r_+$. In addition, since $v'(r_+)=x_0>0$, the solution enters $\{v>0\}$ right after $r_+$, where its phase $y_v=r^{-1}v(1+v^2)^{-\beta}$ is positive and follows $g_{\mathrm{ax}}$.

Next, we prove that the phase reaches the diagonal at a finite radius. Since $g_{\mathrm{ax}}(y)>y\geq0$ on $[0,y_*)$, from the graphical equation we have $v'>0$, and from \eqref{eq:graphical-phase-evolution} we have $ry_v'\geq g_{\mathrm{ax}}(y_v)-y_v>0$. Now, we fix $r_1>r_+$ and set $v_1=v(r_1)>0$. Then, for $r\geq r_1$ we have $v\geq v_1$, and we obtain
	\begin{equation*}
		ry_v' = [g_{\mathrm{ax}}(y_v)-y_v] + \alpha^{-1}v^2g_{\mathrm{ax}}(y_v) \geq \alpha^{-1}v_1^2\,g_{\mathrm{ax}}(y_v) \geq \alpha^{-1}v_1^2\,y_v,
	\end{equation*}
and therefore
	\begin{equation*}
		y_v(r)\geq y_v(r_1)\left(\frac{r}{r_1}\right)^{\alpha^{-1}v_1^2}.
	\end{equation*}
Consequently, the continuation along $\left((0,y_*),g_{\mathrm{ax}}\right)$ has a finite maximal radius $r_*$. Moreover, since $g_{\mathrm{ax}}$ is continuous on $[0,y_*]$, by Proposition \ref{prop:branch-ivp}\textup{(ii)} the phase $y_v$ reaches an endpoint of $(0,y_*)$ as $r\to r_*$. Since $y_v$ increases, that endpoint is $y_*$. Finally, at this radius the phase point is $(y_*,y_*)$, where the continuation joins $g_+$ by Lemma \ref{lem:matching-c1-junctions}\textup{(ii)}.
\end{proof}

\begin{proposition}\label{prop:curvature-blow-up-phase-asymptotics}
	Suppose that a regularity-limited lower branch admits a curvature-blow-up continuation of order $q>0$, with $c_-$ as defined in Definition~\ref{def:curvature-blow-up-continuation}, and let $v$ be a solution to the graphical equation for $\left((-\infty,0),g_-^{\mathrm{ext}}\right)$ on its maximal branch-admissible interval $[r_0,r_+)$. Then $r_+<\infty$, $y_v(r)\to0^-$ as $r\to r_+^-$, and the phase level satisfies the asymptotic expansion
	\begin{equation}\label{eq:curvature-blow-up-phase-asymptotics}
		-y_v(r)=\left(\frac{(q+1)c_-}{r_+}(r_+-r)\right)^{\frac{1}{q+1}}(1+o(1))
	\end{equation}
	as $r\to r_+^-$. Consequently, $v(r)\to0$, $x_v(r)\to+\infty$ and $v'(r)\to+\infty$ as $r\to r_+^-$, and the height $u$ extends to $r=r_+$ from the left with regularity $\mathcal C^{1,\frac{1}{q+1}}$.
\end{proposition}

\begin{proof}
Firstly, we prove that $r_+<\infty$. By Definition \ref{def:curvature-blow-up-continuation}, $g_-^{\mathrm{ext}}$ is a positive $\mathcal C^1$ branch on $(-\infty,0)$, and $g_-^{\mathrm{ext}}\geq1$ near $0^-$ because it tends to $+\infty$. Then, by Lemma \ref{lem:continuation-phase-limits}\textup{(i)} with $b=0$, we have $r_+<\infty$ and $y_v\to0^-$, and $v\to0$ by Lemma \ref{lem:continuation-phase-limits}\textup{(iii)}.

Next, we prove the expansion of the phase. With $w=-y_v$, \eqref{eq:graphical-phase-evolution} reads
	\begin{equation*}
		-rw' = \left(1+\alpha^{-1}v^2\right)g_-^{\mathrm{ext}}(-w)+w.
	\end{equation*}
Then, since $w=o(w^{-q})$ and $1+\alpha^{-1}v^2\to1$, from the expansion $g_-^{\mathrm{ext}}(-w)=c_-w^{-q}(1+o(1))$ we obtain
	\begin{equation*}
		w' = -\frac{c_-}{r_+}w^{-q}(1+o(1)).
	\end{equation*}
Therefore, by integrating $(w^{q+1})'=(q+1)w^qw'$ up to $r_+$, where $w\to0$, we obtain \eqref{eq:curvature-blow-up-phase-asymptotics}. In particular, $x_v=c_-w^{-q}(1+o(1))\to+\infty$, and from the graphical equation we have $v'=(1+v^2)^{\beta+1}x_v\to+\infty$.

Finally, we prove the regularity of the height. By \eqref{eq:curvature-blow-up-phase-asymptotics} and this expansion of $x_v$, we have $|v'(r)|\leq C(r_+-r)^{-\frac{q}{q+1}}$ near $r_+$. Then, for $r<r'<r_+$, by the subadditivity of $\tau\mapsto\tau^{\frac{1}{q+1}}$ we obtain
	\begin{align*}
		|v(r')-v(r)| &\leq \int_r^{r'}|v'(\tau)|\,d\tau \leq C(q+1)\left[(r_+-r)^{\frac{1}{q+1}}-(r_+-r')^{\frac{1}{q+1}}\right]\\
		&\leq C(q+1)(r'-r)^{\frac{1}{q+1}},
	\end{align*}
and consequently $u$ is of class $\mathcal C^{1,\frac{1}{q+1}}$ up to $r_+$ from the left.
\end{proof}

\begin{proposition}\label{prop:singular-positive-matching-phase}
	In the setting of Proposition \ref{prop:curvature-blow-up-phase-asymptotics}, assume additionally that the continuation satisfies the singular positive-matching condition of Definition \ref{def:singular-positive-matching}, with branch $g_{\mathrm{sing}}$ and constant $c_+$. Then $v$ admits exactly one positive graphical continuation through $r_+$. This continuation is defined on a right neighborhood of $r_+$ with $v(r_+)=0$, and its phase satisfies
	\begin{equation}\label{eq:positive-singular-phase-asymptotics}
		y_v(r)=\left(\frac{(q+1)c_+}{r_+}(r-r_+)\right)^{\frac{1}{q+1}}(1+o(1))
	\end{equation}
	as $r\to r_+^+$. The height $u$ extends to $r=r_+$ from the right with regularity $\mathcal C^{1,\frac{1}{q+1}}$, and the phase reaches the diagonal point $(y_*,y_*)$ at a finite radius, where it joins the positive branch $g_+$.
\end{proposition}

\begin{proof}
We divide the proof into four steps.

\emph{Step 1.} Firstly, we prove the existence of a positive graphical continuation. Since $g_{\mathrm{sing}}(y)\to+\infty$ as $y\to0^+$, the graphical equation is singular at $(r_+,0)$. For this reason, we treat the radius as a function of the slope. Indeed, a positive graphical continuation satisfies $v'=(1+v^2)^{\beta+1}g_{\mathrm{sing}}(y_v)>0$, and therefore it is the inverse of a solution of
	\begin{equation*}
		\frac{dr}{ds}=\mathsf H(r,s)=(1+s^2)^{-\beta-1}g_{\mathrm{sing}}\left(Q(r,s)\right)^{-1},\qquad r(0)=r_+.
	\end{equation*}
By the expansion of $g_{\mathrm{sing}}$, we have $\mathsf H(r,s)=c_+^{-1}Q(r,s)^{q}(1+o(1))$ as $s\to0^+$, uniformly for $r$ near $r_+$. Consequently, $\mathsf H$ extends continuously by $\mathsf H(r,0)=0$ to a neighborhood of $(r_+,0)$ in $(0,\infty)\times[0,\infty)$, and by Peano's theorem we obtain a local solution with $r(0)=r_+$.

\emph{Step 2.} Next, we prove the uniqueness, for which it suffices to show that $\mathsf H$ is nonincreasing in $r$. By differentiating $\widetilde\gamma(g_{\mathrm{sing}}(y),y)=1$, we have $g_{\mathrm{sing}}'=-\partial_y\widetilde\gamma\,(\partial_x\widetilde\gamma)^{-1}<0$, since $g_{\mathrm{sing}}(y)>y>0$ places the phase point in the positive quadrant, where both partial derivatives are positive. Moreover, for fixed $s>0$ the map $r\mapsto Q(r,s)$ is strictly decreasing by \eqref{eq:q-partial-r-section}. Therefore, $g_{\mathrm{sing}}(Q(r,s))$ increases with $r$ and $\mathsf H(r,s)$ decreases. Now, let $r_a$ and $r_b$ be two local solutions. Then $D=r_a-r_b$ satisfies $D(0)=0$ and $D'=\mathsf H(r_a,s)-\mathsf H(r_b,s)$. We note that this derivative is nonpositive where $D>0$ and nonnegative where $D<0$, and therefore $|D|$ is nonincreasing on every interval where $D\neq0$. Let us assume by contradiction that $D(s_0)\neq0$, and let $s_1$ be the last zero of $D$ in $[0,s_0)$. Then $|D(s_0)|\leq|D(s_1)|=0$, a contradiction. Consequently, $r_a=r_b$, and since $\frac{dr}{ds}>0$ for $s>0$, the inverse of $r(s)$ is the unique positive graphical continuation from $r_+$.

\emph{Step 3.} Now, we prove the expansion of the phase and the regularity of the height. Near $r_+$ we have $y_v\to0^+$ and $v\to0$. Then, since $y_v=o(y_v^{-q})$, from \eqref{eq:graphical-phase-evolution} and the expansion of $g_{\mathrm{sing}}$ we obtain
	\begin{equation*}
	 	y_v'=\frac{c_+}{r_+}\,y_v^{-q}(1+o(1)),
	\end{equation*}
and by integrating from $r_+$, where $y_v\to0$, we obtain \eqref{eq:positive-singular-phase-asymptotics}. In particular, $v>0$ right after $r_+$. Moreover, with $x_v=c_+y_v^{-q}(1+o(1))$, the argument at the end of the proof of Proposition \ref{prop:curvature-blow-up-phase-asymptotics} implies that $|v'(r)|\leq C(r-r_+)^{-\frac{q}{q+1}}$, and therefore $u$ has the $\mathcal C^{1,\frac{1}{q+1}}$ regularity at $r_+$ from the right.

\emph{Step 4.} Finally, we prove that the phase reaches the diagonal. We fix $\bar r>r_+$, so that $y_v(\bar r)>0$. Since $g_{\mathrm{sing}}(y)>y$ on $(0,y_*)$ and $g_{\mathrm{sing}}$ extends continuously to $y_*$, the growth bound in the proof of Proposition \ref{prop:axis-crossing-positive-matching} applies on $[y_v(\bar r),y_*]$. Therefore, the maximal radius is finite, and the phase reaches $(y_*,y_*)$. At this point $g_{\mathrm{sing}}$ joins $g_+$ in a $\mathcal C^1$ manner by Lemma \ref{lem:matching-c1-junctions}\textup{(ii)}, finalizing the proof.
\end{proof}
\section{Asymptotics of positive ends}\label{sec:positive-branch-asymptotics}

In this section we prove Theorem \ref{thm:main-a}, stated in the introduction, on the asymptotic behavior of positive ends, and in particular of the bowl-type solutions. The precise statement is Theorem \ref{thm:positive-bowl-asymptotics}, where a positive end is a solution of the positive graphical equation of the branch $(I_+,g_+)$ of Proposition \ref{prop:positive-branch}. To do so, we allow every $\alpha>0$, as in Section \ref{sec:barrier-method}, and we use the phase identity \eqref{eq:graphical-phase-evolution} and the endpoint barriers of Corollary \ref{cor:positive-branch-endpoint-barriers}.

\subsection{Positive ends}\label{subsec:positive-ends}

\begin{definition}\label{def:positive-end}
	Let $r_0>0$ and $v_0$ satisfy $Q(r_0,v_0)\in I_+$, and let $v$ be the solution of the positive graphical equation with $v(r_0)=v_0$ given by Proposition \ref{prop:branch-ivp}, on its maximal branch-admissible interval $[r_0,r_+)$. We call $v$ a \emph{positive end}.
\end{definition}

\begin{lemma}\label{lem:positive-end-basic}
	Let $v$ be a positive end. Then $v>0$ and $v'>0$ on $[r_0,r_+)$. If $r_+=\infty$, then $v(r)\to\infty$ as $r\to\infty$.
\end{lemma}

\begin{proof}
Firstly, we note that $v>0$ by branch admissibility, because $I_+\subset(y_*,\infty)$ with $y_*>0$ and $Q(r,\cdot)$ has the sign of its second argument. Moreover, since $g_+>0$ on $I_+$ by Proposition \ref{prop:positive-branch}, from the positive graphical equation we obtain $v'>0$. Finally, let us assume by contradiction that $r_+=\infty$ and that $v$ is bounded. Then $y_v(r)=Q(r,v(r))\to0$ as $r\to\infty$, which contradicts $y_v>y_*>0$. Consequently, $v(r)\to\infty$.
\end{proof}

\begin{lemma}\label{lem:bowl-is-positive-end}
	Let $u:[0,R_\gamma)\to\rr$ be the profile of a bowl-type solution in the sense of Definition \ref{def:bowl-type-solution}, let $v=u'$, and let $r_0\in(0,R_\gamma)$. Then the restriction of $v$ to $[r_0,R_\gamma)$ is a positive end with $r_+=R_\gamma$.
\end{lemma}

\begin{proof}
We recall that the argument of \cite[Section 8]{rengaswami2024classification} compares $v$ with the level function at $y_*$ and, in the $1$-nondegenerate case, also with the one at the upper endpoint of $I_+$. From this comparison we have $v>0$, $v'>0$ and $y_*\leq y_v<1$ on $(0,R_\gamma)$. In the $1$-degenerate case we only need the lower bound, since $I_+=(y_*,\infty)$. Moreover, at the tip the hypersurface is umbilic with curvature $y_*$, and therefore $y_v\to y_*$ as $r\to0^+$.

Next, we prove that the lower bound is strict away from the tip. Let us assume by contradiction that $y_v(r_1)=y_*$ for some $r_1\in(0,R_\gamma)$. Then $y_v$ attains its minimum at $r_1$, and consequently $y_v'(r_1)=0$. However, at that radius the phase point is $(y_*,y_*)$, where $g_+$ extends with $g_+(y_*)=y_*$ by Proposition \ref{prop:positive-branch}. Then, by \eqref{eq:graphical-phase-evolution}, we obtain $r_1y_v'(r_1)=\alpha^{-1}v(r_1)^2y_*>0$, leading to the desired contradiction.

Consequently, $y_v(r)\in I_+$ for every $r\in(0,R_\gamma)$, and therefore $v$ is a branch-admissible solution of the positive graphical equation on $[r_0,R_\gamma)$. Finally, we have $r_+=R_\gamma$ because $[0,R_\gamma)$ is the maximal interval of $u$.
\end{proof}

For such an end, Proposition \ref{prop:positive-end-maximal-interval} recovers the alternative of \cite[Theorems 1.3 and 1.4]{rengaswami2024classification}.

\subsection{Statement of the positive-end asymptotics}\label{subsec:positive-end-asymptotics-statement}

\begin{theorem}\label{thm:positive-bowl-asymptotics}
	Let $\alpha>0$ and let $v$ be a positive end. The expansions below depend only on the positive branch, and not on the initial data $(r_0,v_0)$ of $v$. In \textup{(ii)}--\textup{(iv)}, $\gamma$ is $1$-degenerate and $g_+$ has power decay in the sense of Definition \ref{def:degenerate-positive-branch-asymptotics}\textup{(i)}, that is, $g_+(y)=c_\gamma y^{-k_\gamma}(1+o(1))$ as $y\to\infty$, with $c_\gamma,k_\gamma>0$.
	
	\begin{enumerate}[label=\textup{(\roman*)}]
		\item Assume that $\gamma$ is normalized $1$-nondegenerate and that $\widetilde\gamma$ is of class $\mathcal C^3$ near $(0,1)$, so that
		\begin{equation}\label{eq:endpoint-expansion-gplus}
			g_+(y) = a_1(1-y) + a_2(1-y)^2 + O(|1-y|^3) \qquad \text{as }y\to1,
		\end{equation}
		where $a_1=-g_+'(1)>0$ and $a_2=\frac{1}{2}g_+''(1)$. Let
		\begin{align*}
			a_\gamma &= \alpha\left(\frac{\alpha}{a_1}-\beta\right),\\
			\mathcal P &= \frac{1}{2}(1-2a_\gamma)\beta\bigl(1+\beta(1-2a_\gamma)\bigr) - a_\gamma^2\beta, \\
			b_\gamma &= \alpha\left[\frac{\alpha\bigl((\beta+1)(1-2a_\gamma)-a_\gamma\bigr)}{a_1} + \frac{a_2\alpha^2}{a_1^3} - \mathcal P\right].
		\end{align*}
		Then
		\begin{equation*}
			v(r) = r^\alpha - a_\gamma r^{-\alpha} + b_\gamma r^{-3\alpha} + o(r^{-4\alpha}) \qquad \text{as }r\to\infty.
		\end{equation*}
		If $\widetilde\gamma$ is only of class $\mathcal C^2$ near $(0,1)$, the same expansion holds with a remainder $o(r^{-3\alpha})$.
		
		\item If $k_\gamma>2\alpha-1$, then
		\begin{equation*}
			v(r) = A_\gamma r^{d_\gamma}(1+o(1)) \qquad \text{as }r\to\infty,
		\end{equation*}
		where $d_\gamma=\frac{\alpha(k_\gamma+1)}{k_\gamma+1-2\alpha}$ and $A_\gamma>0$ satisfies $A_\gamma^{\frac{k_\gamma+1-2\alpha}{\alpha}}=\frac{c_\gamma}{d_\gamma}$.
		
		\item If $k_\gamma=2\alpha-1$ (so that $\alpha>\frac{1}{2}$), then
		\begin{gather*}
			\frac{d}{dr}\log v(r)=c_\gamma r^{2\alpha-1}(1+o(1)),\\
			\log v(r) = \frac{c_\gamma}{2\alpha}r^{2\alpha} + o(r^{2\alpha}) \qquad \text{as }r\to\infty.
		\end{gather*}
		
		\item Assume in addition that the positive branch has balanced second-order power decay in the sense of Definition \ref{def:degenerate-positive-branch-asymptotics}\textup{(ii)} to \textup{(iii)}, with
		\begin{equation*}
			g_+(y) = c_\gamma y^{-k_\gamma} + c_{\gamma,1}y^{-(2k_\gamma+1)} + o\left(y^{-(2k_\gamma+1)}\right) \qquad \text{as }y\to\infty.
		\end{equation*}
		Both cases below have $k_\gamma+1\geq4\alpha$, so $k_\gamma>2\alpha-1$ and $A_\gamma$, $d_\gamma$ are those of \textup{(ii)}. Set $\mathcal C_\gamma=\bigl(1+\beta(k_\gamma+1)\bigr)A_\gamma^{-2} +(c_{\gamma,1}/c_\gamma)A_\gamma^{-(k_\gamma+1)/\alpha}$. If $k_\gamma+1>4\alpha$, then
		\begin{equation*}
			v(r) = A_\gamma r^{d_\gamma} + B_\gamma r^{-d_\gamma} + o(r^{-d_\gamma}), \qquad B_\gamma = \frac{\alpha A_\gamma\mathcal C_\gamma}{k_\gamma+1-4\alpha}.
		\end{equation*}
		If $k_\gamma+1=4\alpha$, then
		\begin{equation*}
			v(r) = A_\gamma r^{d_\gamma} + L_\gamma r^{-d_\gamma}\log r + o(r^{-d_\gamma}\log r), \qquad L_\gamma = A_\gamma d_\gamma\mathcal C_\gamma.
		\end{equation*}
	\end{enumerate}
\end{theorem}

In the literature, the leading term $r^\alpha$ is due to Rengaswami \cite[Theorem 1.1 and Proposition 7.2]{rengaswami2024classification}. For the rotational translators of $S_k$ and $S_k^{1/k}$ with $k<n$ and positive angle function, de Lira--de Farias obtained the leading decay of the angle function, with its constant, in warped products asymptotic to $\rr^{n+1}$ \cite[Lemma 15 and Proposition 16]{lira2026translating}. The angle function of the fundamental translators of de Lima--Pipoli tends to zero at infinity \cite[Remark 14]{lima2025translators}. The first correction term, involving $a_\gamma$, appears in the work of Cogo--Lynch--Vi\v{c}\'anek Mart\'inez for $1$-homogeneous curvature functions \cite[Proposition A.1]{cogo-lynch-vicanek2023rotational}, as well as in Sheng--Wu for the powers $H^p$ \cite[Remark 4.1]{sheng2010rotationally}.

The proofs use $r_+=\infty$, which Subsection \ref{subsec:positive-end-interval-and-phase} establishes under the assumptions of each item.

\begin{remark}
	The first coefficient of \eqref{eq:endpoint-expansion-gplus} is determined by $\widetilde\gamma$ at the cylindrical direction. Euler's identity for the $\alpha$-homogeneous function $\widetilde\gamma$ at $(0,1)$, where $\widetilde\gamma(0,1)=1$, gives $\partial_y\widetilde\gamma(0,1)=\alpha$. Differentiating $\widetilde\gamma(g_+(y),y)=1$ at $y=1$ then gives $a_1=\alpha\bigl(\partial_x\widetilde\gamma(0,1)\bigr)^{-1}$, which is positive by ellipticity, and
	\begin{equation*}
		a_\gamma = \alpha\left(\partial_x\widetilde\gamma(0,1)-\beta\right).
	\end{equation*}
	When $\alpha=1$ this is the coefficient $\dot\gamma^1(0,1,\ldots,1)$ of \cite[Proposition A.1]{cogo-lynch-vicanek2023rotational}, which uses the velocity $\frac{1}{2}\vec e_{n+1}$ instead of $\vec e_{n+1}$. The two conventions differ by the dilation $\rho\mapsto\frac{\rho}{2}$, which does not change this coefficient. The unpublished preprint \cite{RengaswamiTorresSantaella2023Bowl} also computes $a_\gamma$ for $\alpha=1$.
\end{remark}

\begin{remark}\label{rem:powers-mean-curvature-check}
	Let $p>0$ and $\gamma=\left(\frac{H}{n-1}\right)^p$, so that $\alpha=p$, $\beta=\frac{p-1}{2p}$, $g_+(y)=(n-1)(1-y)$, $a_1=n-1$ and $a_2=0$. Theorem \ref{thm:positive-bowl-asymptotics}\textup{(i)} gives
	\begin{equation*}
		a_\gamma = \frac{p^2}{n-1} - \frac{p-1}{2}, \qquad b_\gamma = -\frac{(p-1)(3p-1)}{8} + \frac{p^2(5p-3)}{2(n-1)} - \frac{p^3(7p-1)}{2(n-1)^2}.
	\end{equation*}
	The dilation $v(r)=v_{H^p}((n-1)r)$ maps the positive ends of $H^p$ to those of $\gamma$. Under it, $a_\gamma=-(n-1)^{-p}c(n,p)$, where $c(n,p)=(n-1)^{p-1}\left(\frac{(p-1)(n-1)}{2}-p^2\right)$ is the first-correction coefficient of Sheng--Wu \cite[Remark 4.1]{sheng2010rotationally}. The next coefficient is $b_\gamma(n-1)^{3p}$, whose existence is conjectured in \cite[Section 2 and Remark 4.2]{sheng2010rotationally}. The computer-algebra values listed there for $p=2$ and $p=3$ differ from $b_\gamma(n-1)^{3p}$, which for these exponents equals $-\frac{1}{8}(n-1)^4(5n^2-122n+533)$ and $-2(n-1)^7(n^2-29n+163)$, respectively. For $p=1$ the coefficients reduce to $a_\gamma=(n-1)^{-1}$ and $b_\gamma=(n-4)(n-1)^{-2}$, in agreement with the asymptotics of \cite[Lemma 2.1]{clutterbuck2007stability}. For the unnormalized $H$, their argument also gives the second coefficient $(n-4)(n-1)$, which is the rescaled value $b_\gamma(n-1)^3$.
\end{remark}

\subsection{Maximal interval and phase at infinity}\label{subsec:positive-end-interval-and-phase}

\begin{proposition}\label{prop:positive-end-maximal-interval}
	Let $v$ be a positive end. If $\gamma$ is $1$-nondegenerate, then $r_+=\infty$. If $\gamma$ is $1$-degenerate, then either $r_+=\infty$, or $r_+<\infty$ and $v(r)\to\infty$ as $r\to r_+$.
\end{proposition}

\begin{proof}
Let us assume that $r_+<\infty$. Firstly, we prove that $y_v$ is unbounded on $[r_0,r_+)$. If $y_v$ were bounded on $[r_0,r_+)$, its range would lie in a compact interval $K\subset[y_*,\infty)$ to which $g_+$ extends continuously, as in the proof of Corollary \ref{cor:positive-branch-endpoint-barriers}. Then, by Proposition \ref{prop:branch-ivp}\textup{(ii)}, $y_v(r_+)$ would be an endpoint of $I_+$, which is excluded by Corollary \ref{cor:positive-branch-endpoint-barriers}. On the one hand, in the $1$-nondegenerate case this contradicts $y_v<1$, and consequently $r_+=\infty$. On the other hand, in the $1$-degenerate case, from \eqref{eq:q-large-slope-bound} and $r\geq r_0$ we have $y_v\leq2^{|\beta|}r_0^{-1}v^{\frac{1}{\alpha}}$ whenever $v\geq1$. Therefore, $v$ is unbounded on $[r_0,r_+)$. Since $v$ is increasing by Lemma \ref{lem:positive-end-basic}, we conclude that $v(r)\to\infty$ as $r\to r_+$.
\end{proof}

The same trapping argument appears in \cite[Lemma 2.1]{clutterbuck2007stability} for $\gamma=H$.

\begin{lemma}\label{lem:stable-scalar-bootstrap}
	Let $\phi$ be a positive continuous function on $[R_0,\infty)$ with $\int_{R_0}^\infty\phi(r)\,dr=\infty$, and let $Z\in\mathcal C^1([R_0,\infty))$.
	\begin{enumerate}[label=\textup{(\roman*)}]
		\item If there are $m\in\rr$ and $\delta>0$ such that $Z'(r)\geq\delta\phi(r)$ whenever $Z(r)\leq m$ for $r\geq R_0$, then $Z(r)>m$ for all sufficiently large $r$.
		\item If $Z'(r)=c\,\phi(r)\bigl(Z(r)-L\bigr)+E(r)$ for some $c<0$ and $L\in\rr$, where $E(r)=o\!\left((1+|Z(r)-L|)\phi(r)\right)$ as $r\to\infty$, then $Z(r)\to L$.
	\end{enumerate}
\end{lemma}

\begin{proof}
\textup{(i)} Firstly, we show that $Z(r_1)>m$ for some $r_1\geq R_0$. If $Z\leq m$ on a tail $[\bar r,\infty)$, then by integrating we have $Z(r)\geq Z(\bar r)+\delta\int_{\bar r}^{r}\phi(s)\,ds$, which tends to $\infty$, a contradiction. Next, let us assume by contradiction that $Z(r_2)\leq m$ for some $r_2>r_1$, and let $r_3$ be the smallest $r\in[r_1,r_2]$ with $Z(r)\leq m$. Then $Z(r_3)=m$ and $Z>m$ on $[r_1,r_3)$, and therefore $Z'(r_3)\leq0$, which contradicts $Z'(r_3)\geq\delta\phi(r_3)>0$.

\textup{(ii)} Let $Y=Z-L$ and $\eps>0$. By the assumption on $E$, there is $R_\eps\geq R_0$ with $|E(r)|\leq-\frac{c}{2}(1+\eps^{-1})^{-1}(1+|Y(r)|)\phi(r)$ for $r\geq R_\eps$. We note that where $|Y|\geq\eps$ we have $1+|Y|\leq(1+\eps^{-1})|Y|$. Therefore, for $r\geq R_\eps$ we obtain
	\begin{align*}
		Y(r)\geq\eps &\ \Longrightarrow\ Y'(r)\leq\frac{c}{2}\phi(r)Y(r)\leq\frac{c}{2}\eps\phi(r), \\
		Y(r)\leq-\eps &\ \Longrightarrow\ Y'(r)\geq\frac{c}{2}\phi(r)Y(r)\geq-\frac{c}{2}\eps\phi(r).
	\end{align*}
Then, by applying \textup{(i)} on $[R_\eps,\infty)$ to $Y$ and to $-Y$, with $m=-\eps$ and $\delta=-\frac{c}{2}\eps>0$, we have $|Y(r)|<\eps$ for all large $r$. Since $\eps$ is arbitrary, we conclude that $Z(r)\to L$.
\end{proof}

Both items appear, without being isolated, in \cite[Lemma 2.1]{clutterbuck2007stability} and \cite[Proposition A.1]{cogo-lynch-vicanek2023rotational}.

Rengaswami obtains the limit of the phase of the bowl-type solution in \cite[Proposition 7.2]{rengaswami2024classification} in the nondegenerate case, and in \cite[Section 6.2 and Lemma 6.1]{rengaswami2024classification} in the degenerate case. The lower bound on the phase that those arguments draw from the axis is part of branch admissibility for a positive end, so a single argument covers both cases.

\begin{lemma}\label{lem:phase-limit-positive-end}
	Assume that $\gamma$ is normalized $1$-nondegenerate or $1$-degenerate, and let $v$ be a positive end with $r_+=\infty$. Then $y_v(r)$ tends to the right endpoint of $I_+$ as $r\to\infty$, which is $1$ in the nondegenerate case and $\infty$ in the degenerate one.
\end{lemma}

\begin{proof}
We fix $m\in I_+$, so that $g_+(m)>0$. Since $v(r)\to\infty$ by Lemma \ref{lem:positive-end-basic}, there is $R_m\geq r_0$ with $\alpha^{-1}g_+(m)v(r)^2\geq2m$ for $r\geq R_m$. Then, at every such $r$ with $y_v(r)\leq m$, by \eqref{eq:graphical-phase-evolution} and the monotonicity of $g_+$ we have
	\begin{equation*}
		r y_v' = \left(1+\alpha^{-1}v^2\right)g_+(y_v)-y_v \geq \alpha^{-1}g_+(m)v^2-m \geq m.
	\end{equation*}
This means that $y_v'\geq\phi$ with $\phi(r)=mr^{-1}$, whose integral diverges. Therefore, by Lemma \ref{lem:stable-scalar-bootstrap}\textup{(i)} with $Z=y_v$ and $\delta=1$, we obtain $y_v(r)>m$ for all large $r$. Since $m\in I_+$ is arbitrary, it follows that $\liminf y_v$ is at least the right endpoint of $I_+$. Finally, we note that this endpoint is $\infty$ in the degenerate case and $1$ in the nondegenerate case, where $y_v<1$ by Corollary \ref{cor:positive-branch-endpoint-barriers}, concluding the proof.
\end{proof}

Next, we derive the form of the positive graphical equation that we use in the proofs of items \textup{(ii)} to \textup{(iv)} of Theorem \ref{thm:positive-bowl-asymptotics}. Let $\gamma$ be $1$-degenerate with power decay in the sense of Definition \ref{def:degenerate-positive-branch-asymptotics}, and let $v$ be a positive end with $r_+=\infty$. Since $y_v\to\infty$ by Lemma \ref{lem:phase-limit-positive-end}, we have $g_+(y_v)=c_\gamma y_v^{-k_\gamma}(1+o(1))$. Moreover, by the large-slope expansion \eqref{eq:eta-large-slope-asymptotic}, we have $y_v=r^{-1}v^{\frac{1}{\alpha}}(1+o(1))$, and therefore
\begin{equation*}
	g_+(y_v) = c_\gamma r^{k_\gamma}v^{-\frac{k_\gamma}{\alpha}}(1+o(1)).
\end{equation*}
Then, since $v\to\infty$, we have $(1+v^2)^{\beta+1}=v^{\frac{3\alpha-1}{\alpha}}(1+o(1))$, and the positive graphical equation becomes
\begin{equation}\label{eq:degenerate-power-main-ode}
	v' = c_\gamma r^{k_\gamma} v^{\frac{3\alpha-1-k_\gamma}{\alpha}}(1+o(1)).
\end{equation}

In addition, we show in the next lemma that items \textup{(ii)} and \textup{(iii)} of Theorem \ref{thm:positive-bowl-asymptotics} cover every power-decay branch whose positive ends are defined on a half-line.

\begin{lemma}\label{lem:subcritical-degenerate-blow-up}
	Assume that $\gamma$ is $1$-degenerate and that $g_+(y)\geq cy^{-k}$ for all large $y$, with $c>0$ and $0<k<2\alpha-1$. Then every positive end has $r_+<\infty$.
\end{lemma}

\begin{proof}
We recall that this is \cite[Theorem 1.4\,(ii)]{rengaswami2024classification} for the bowl-type solution, where it follows from a differential inequality $v'\geq c_1v^{\frac{3\alpha-1-k}{\alpha}}$ with $c_1>0$. Since $k<2\alpha-1$, the exponent exceeds $1$, and therefore $v$ blows up at a finite radius. We note that this argument uses the equation only on a tail $r\geq r_0>0$, where the phase satisfies \eqref{eq:q-large-slope-bound}. Consequently, the data at the axis play no role, and the argument applies to every positive end.
\end{proof}

For the bowl-type solution, item \textup{(ii)} of the next proposition is \cite[Theorems 1.2 and 1.4]{rengaswami2024classification}, and item \textup{(i)} is \cite[Lemma 2.1]{clutterbuck2007stability} for $\gamma=H$.

\begin{proposition}\label{prop:positive-end-at-infinity}
	Assume that $\gamma$ is normalized $1$-nondegenerate or $1$-degenerate, and let $v$ be a positive end.
	\begin{enumerate}[label=\textup{(\roman*)}]
		\item If $\gamma$ is $1$-nondegenerate, then $v(r)=r^\alpha(1+o(1))$ as $r\to\infty$.
		\item If $\gamma$ is $1$-degenerate, the decay of $g_+$ determines the alternative of Proposition \ref{prop:positive-end-maximal-interval}. If $g_+(y)\leq Cy^{-k}$ for large $y$, with $k\geq2\alpha-1$, then $r_+=\infty$, and if $g_+(y)\geq cy^{-k}$ for large $y$, with $0<k<2\alpha-1$, then $r_+<\infty$.
	\end{enumerate}
\end{proposition}

\begin{proof}
\textup{(i)} Firstly, we have $r_+=\infty$ by Proposition \ref{prop:positive-end-maximal-interval}. Then $v(r)\to\infty$ by Lemma \ref{lem:positive-end-basic}, and $y_v(r)\to1$ by Lemma \ref{lem:phase-limit-positive-end}. Therefore, from \eqref{eq:eta-large-slope-asymptotic} we obtain $v(r)^{\frac{1}{\alpha}}=ry_v(r)(1+o(1))$, that is, $v(r)=r^\alpha(1+o(1))$.

\textup{(ii)} The second assertion is Lemma \ref{lem:subcritical-degenerate-blow-up}. For the first assertion, we recall that the argument of \cite[Theorem 1.4]{rengaswami2024classification} bounds $v'$ by a multiple of $v^{\frac{3\alpha-1-k}{\alpha}}$. Since $k\geq2\alpha-1$, the exponent is at most $1$, and therefore $v$ stays bounded on bounded intervals. As in the proof of Lemma \ref{lem:subcritical-degenerate-blow-up}, this argument uses the equation only on a tail. Consequently, by Proposition \ref{prop:positive-end-maximal-interval}, we have $r_+=\infty$ for every positive end. We note that \cite[Theorem 1.4]{rengaswami2024classification} is stated for $\alpha>\frac{1}{2}$, since for $\alpha\leq\frac{1}{2}$ the bowl-type solution is entire by \cite[Theorem 1.2]{rengaswami2024classification}. The bound above holds for every $\alpha>0$, and so does the first assertion.
\end{proof}

When $r_+<\infty$, the decay of $g_+$ also controls whether the height $\int_{r_0}^{r}v(\rho)\,d\rho$ diverges as $r\to r_+$, now with threshold $k=\alpha-1$. We prove this in the next proposition.

\begin{proposition}\label{prop:finite-height-ends}
	Let $\gamma$ be $1$-degenerate and let $v$ be a positive end with $r_+<\infty$.
	\begin{enumerate}[label=\textup{(\roman*)}]
		\item If $g_+(y)\leq Cy^{-k}$ for large $y$, with $k\geq\max\{0,\alpha-1\}$, then the height diverges as $r\to r_+$.
		\item If $\alpha>1$ and $g_+(y)\geq cy^{-k}$ for large $y$, with $0<k<\alpha-1$, then the height has a finite limit as $r\to r_+$.
	\end{enumerate}
	In particular, a bowl-type solution that is not entire is cylindrical when $\alpha\leq1$.
\end{proposition}

\begin{proof}
Recall that $v(r)\to\infty$ as $r\to r_+$ by Proposition \ref{prop:positive-end-maximal-interval}. In particular, near $r_+$ we have $v\geq1$ and $r_0\leq r<r_+$, and therefore $y_v$ is comparable to $v^{\frac{1}{\alpha}}$ by \eqref{eq:q-large-slope-bound}. We use this comparison in both cases.
\begin{itemize}
	\item In case \textup{(i)}, we have $k<2\alpha-1$ by Proposition \ref{prop:positive-end-at-infinity}\textup{(ii)}, and $v'=(1+v^2)^{\beta+1}g_+(y_v)\leq C_1v^{\frac{3\alpha-1-k}{\alpha}}$ near $r_+$. Then, by integrating from $r$ to $r_+$, we obtain $v(r)\geq c_1(r_+-r)^{-\frac{\alpha}{2\alpha-1-k}}$. Since the exponent is at least $1$, the height diverges.
	\item In case \textup{(ii)}, from the reversed inequalities we obtain $v(r)\leq C_2(r_+-r)^{-\frac{\alpha}{2\alpha-1-k}}$. Since the exponent is less than $1$, the height has a finite limit.
\end{itemize}
Finally, we prove the last assertion. We have $g_+<y_*$ on $I_+$ by Proposition \ref{prop:positive-branch}, and therefore item \textup{(i)} applies with $k=0$ whenever $\alpha\leq1$. Moreover, by Lemma \ref{lem:bowl-is-positive-end}, the slope of a bowl-type solution that is not entire is a positive end with $r_+<\infty$. Consequently, its height diverges.
\end{proof}

For $\alpha>1$, a bowl-type solution that is not entire may end at a finite height, as the translating curves of Remark \ref{rem:langford-rengaswami-condition} do. In particular, assume that $\lim g_+>0$. Then the bowl-type solution of \cite[Theorem 1.3]{rengaswami2024classification} is defined on a ball, and Proposition \ref{prop:finite-height-ends} shows that it is cylindrical for $\frac{1}{2}<\alpha\leq1$. In contrast, since $g_+$ is bounded below by a positive constant, item \textup{(ii)} shows that its height is finite for $\alpha>1$. The two items do not cover every branch, since $g_+$ need not satisfy either power bound. For $\gamma=K^p$ we have $\alpha=np$ and $g_+(y)=y^{-(n-1)}$. The thresholds $k=2\alpha-1$ and $k=\alpha-1$ read $p=\frac{1}{2}$ and $p=1$, in agreement with \cite[Theorems 2 and 3]{urbas1998complete}.

\subsection{The nondegenerate endpoint}\label{subsec:one-nondegenerate-asymptotics}

We note that the coefficients of Theorem \ref{thm:positive-bowl-asymptotics}\textup{(i)} come from the branch data at $y=1$. To compute them, we use the field $\mathsf F$ of Section \ref{sec:barrier-method} with $g$ replaced by the level-one chart $g(\cdot,1)$ of Lemma \ref{lem:positive-endpoint-regularity}, which continues $g_+$ through $y=1$ into $y>1$. We need this continuation because in the estimates we evaluate $\mathsf F(r,s)$ at slopes near $r^\alpha$, where
\begin{equation*}
	Q(r,r^\alpha) = \left(1+r^{-2\alpha}\right)^{-\beta},
\end{equation*}
since $2\alpha\beta=\alpha-1$, and this value exceeds $1$ when $\alpha<1$. However, along a positive end the phase stays below $1$ by Corollary \ref{cor:positive-branch-endpoint-barriers}. There $g(\cdot,1)$ agrees with $g_+$, and therefore the coefficients in \eqref{eq:endpoint-expansion-gplus} are unchanged.

\begin{lemma}\label{lem:uniform-expansions-nondegenerate}
	Assume the hypotheses of Theorem \ref{thm:positive-bowl-asymptotics}\textup{(i)}, with $a_1$, $a_\gamma$ and $b_\gamma$ as defined there. For $k=0,\ldots,4$ let
	\begin{align*}
		&P_0=P_1=r^\alpha, \qquad P_2=P_3=r^\alpha-a_\gamma r^{-\alpha}, \qquad P_4=r^\alpha-a_\gamma r^{-\alpha}+b_\gamma r^{-3\alpha}, \\
		&L_0=L_2=L_4=0, \qquad L_1=-a_\gamma, \qquad L_3=b_\gamma,
	\end{align*}
	so that $P_k$ collects the terms of $P_4$ of order larger than $r^{-k\alpha}$ and $L_k$ is the coefficient of $r^{-k\alpha}$. Let $\Theta(r)=o(r^\alpha)$. Then, as $r\to\infty$,
	\begin{align}
		r^{k\alpha}\left[\mathsf F\left(r,P_k+\Theta\,r^{-k\alpha}\right)-P_k'\right] + k\alpha\Theta\,r^{-1} \nonumber \\
		\label{eq:normalized-field-estimate}
		\qquad = -\frac{a_1}{\alpha}(\Theta-L_k)r^{2\alpha-1} + o\!\left((1+|\Theta-L_k|)r^{2\alpha-1}\right),
	\end{align}
	with the following uniformity. For every function $\eps$ with $\eps(r)\to0$, the quotient of the remainder by $(1+|\Theta-L_k|)r^{2\alpha-1}$ tends to zero uniformly over all $\Theta$ with $|\Theta(r)|\leq\eps(r)r^\alpha$. We denote this estimate by \textup{(E$_k$)}.
\end{lemma}

\begin{proof}
We divide the proof into four steps.

\emph{Step 1.} Firstly, we rescale the slope by its leading order. Let $\xi=r^{-2\alpha}$ and write the slope as $s=r^\alpha V$. Then, by the definitions of $Q$ and $\mathsf F$ and the identity $2\alpha(\beta+1)=3\alpha-1$, we have
	\begin{equation*}
		\mathsf F\left(r,r^\alpha V\right) = r^{\alpha-1}\xi^{-1}\mathsf G(V,\xi), \qquad \mathsf G(V,\xi) = \left(V^2+\xi\right)^{\beta+1} g\!\left(V\left(V^2+\xi\right)^{-\beta},1\right),
	\end{equation*}
and, for a function $V$ of $\xi$,
	\begin{equation}\label{eq:master-residual}
		\mathsf F\left(r,r^\alpha V\right) - \frac{d}{dr}\left(r^\alpha V\right) = r^{\alpha-1}\left[\xi^{-1}\mathsf G(V,\xi) - \alpha\left(V-2\xi\dot V\right)\right],
	\end{equation}
where the dot denotes differentiation with respect to $\xi$. Recall that the chart $g(\cdot,1)$ is of class $\mathcal C^3$ on a two-sided neighborhood of $y=1$, by Lemma \ref{lem:positive-endpoint-regularity} and the hypothesis $\widetilde\gamma\in\mathcal C^3$ near $(0,1)$. Therefore, $\mathsf G$ is of class $\mathcal C^3$ near $(1,0)$, with
	\begin{equation*}
		\mathsf G(1,0)=0, \qquad \partial_V\mathsf G(1,0)=-\frac{a_1}{\alpha}, \qquad \partial_\xi\mathsf G(1,0)=\beta a_1.
	\end{equation*}

\emph{Step 2.} Next, we obtain the coefficients $a_\gamma$ and $b_\gamma$ from a formal solution of the reduced equation. By \eqref{eq:master-residual}, $v=r^\alpha V(\xi)$ solves $v'=\mathsf F(r,v)$ if and only if the bracket vanishes. For this reason, we look for a formal solution $V=1+V_1\xi+V_2\xi^2+o(\xi^2)$. Since $V-2\xi\dot V=1-V_1\xi-3V_2\xi^2+o(\xi^2)$, the constant term of the bracket is $\partial_V\mathsf G\,V_1+\partial_\xi\mathsf G-\alpha$, with derivatives at $(1,0)$. Then, by imposing that this term vanishes, we obtain
	\begin{equation*}
		-\frac{a_1}{\alpha}V_1+\beta a_1-\alpha=0, \qquad\text{so}\qquad V_1 = -\alpha\left(\frac{\alpha}{a_1}-\beta\right) = -a_\gamma.
	\end{equation*}
Next, we compute the coefficient of $\xi$ in the bracket. It comes from the quadratic part of the Taylor expansion of $\mathsf G$ at $(1,0)$ and from the term $\alpha V_1\xi$ in $-\alpha\left(V-2\xi\dot V\right)$. With $V_1=-a_\gamma$ it equals
	\begin{equation*}
		-\frac{a_1}{\alpha}V_2 + \frac{1}{2}\partial_V^2\mathsf G\,a_\gamma^2 - \partial_{V\xi}\mathsf G\,a_\gamma + \frac{1}{2}\partial_\xi^2\mathsf G - \alpha a_\gamma,
	\end{equation*}
where the second derivatives of $\mathsf G$ at $(1,0)$ are
	\begin{align*}
		\partial_V^2\mathsf G &= -\frac{2}{\alpha}\left(a_1(\beta+2)-\frac{a_2}{\alpha}\right), \qquad \partial_{V\xi}\mathsf G = a_1\left(2\beta^2+2\beta-1\right)-\frac{2a_2\beta}{\alpha}, \\
		\partial_\xi^2\mathsf G &= \beta\left(a_1(\beta+1)+2a_2\beta\right).
	\end{align*}
Then, by imposing that this coefficient vanishes, we obtain $V_2=b_\gamma$. Consequently, the formal solution is
	\begin{equation*}
		V^\ast(\xi) = 1-a_\gamma \xi+b_\gamma \xi^2+o(\xi^2).
	\end{equation*}

\emph{Step 3.} Now, we return to the normalization of the statement, which in the variable $V$ reads $V=V_k^\ast+\Theta\,\xi^{\frac{k+1}{2}}$. Here $V_k^\ast=r^{-\alpha}P_k$ is the truncation of the formal solution $V^\ast$ below the order $\xi^{\frac{k+1}{2}}$ and $L_k$ is the coefficient of that order. Recall that the reduced problem depends on $r$ only through $\xi$. Therefore, $V^\ast$ contains only integer powers of $\xi$, that is, even powers of $r^{-\alpha}$, and we have $L_k=0$ for even $k$. This means that the normalizations with even $k$ gain an order without a new coefficient, and we use the last one to improve the remainder of Theorem \ref{thm:positive-bowl-asymptotics}\textup{(i)} from $o(r^{-3\alpha})$ to $o(r^{-4\alpha})$.

Moreover, by multiplying \eqref{eq:master-residual} by $r^{k\alpha}$, we obtain the left-hand side of \eqref{eq:normalized-field-estimate} with $\Theta$ fixed, since the derivative of the normalized term contributes $-k\alpha\Theta\,r^{-k\alpha-1}$. We notice that the summand $k\alpha\Theta\,r^{-1}$ restores this contribution. Consequently, in order to prove \eqref{eq:normalized-field-estimate} it suffices to expand the bracket of \eqref{eq:master-residual} in the normalized form. Indeed, since $\eps(r)\to0$ and $\xi\to0$, we have that
	\begin{equation*}
		|V-1| \leq \left|V_k^\ast-1\right|+|\Theta|\,\xi^{\frac{k+1}{2}} = O\!\left(\xi+\eps(r)\right).
	\end{equation*}
This estimate justifies the Taylor expansions of $\mathsf G$ at $(1,0)$ that we use in Step 4. All the estimates below are uniform over $|\Theta(r)|\leq\eps(r)r^\alpha$.

\emph{Step 4.} Finally, we show that the terms of order $\xi^{\frac{k-1}{2}}$ in the bracket of \eqref{eq:master-residual} produce the main term of \eqref{eq:normalized-field-estimate}, and that all the other terms belong to the remainder of \eqref{eq:normalized-field-estimate}. In fact, the bracket has two contributions of that order. On the one hand, the contribution that is linear in $\Theta$ comes from $\partial_V\mathsf G(1,0)$, and it equals
	\begin{equation*}
		\xi^{-1}\partial_V\mathsf G(1,0)\,\Theta\,\xi^{\frac{k+1}{2}} = -\frac{a_1}{\alpha}\,\Theta\,\xi^{\frac{k-1}{2}}.
	\end{equation*}
On the other hand, the contribution that is independent of $\Theta$ is the value of the bracket at $V_k^\ast$, whose expansion contains only integer powers of $\xi$. In this expansion the power $\xi^{-1}$ is absent because $\mathsf G(1,0)=0$, and the constant term and, for $k=4$, the coefficient of $\xi$ vanish by the choice of $a_\gamma$ and $b_\gamma$. For even $k$ the exponent $\frac{k-1}{2}$ is a half-integer, and therefore this value is $o\!\left(\xi^{\frac{k-1}{2}}\right)$. For $k=1,3$ it equals $\alpha^{-1}a_1L_k\,\xi^{\frac{k-1}{2}}+o\!\left(\xi^{\frac{k-1}{2}}\right)$, because the coefficient of $\xi^{\frac{k-1}{2}}$ is affine in the omitted coefficient of $V_k^\ast$, with slope $-\alpha^{-1}a_1$, and vanishes when that coefficient equals $L_k$. Then, since multiplying by $r^{(k+1)\alpha-1}$ turns $\xi^{\frac{k-1}{2}}$ into $r^{2\alpha-1}$, the two contributions combine into
	\begin{equation*}
		r^{(k+1)\alpha-1}\left(-\frac{a_1}{\alpha}\Theta + \frac{a_1}{\alpha}L_k\right)\xi^{\frac{k-1}{2}} = -\frac{a_1}{\alpha}\left(\Theta-L_k\right)r^{2\alpha-1},
	\end{equation*}
which is the main term of \eqref{eq:normalized-field-estimate} for every $k=0,\ldots,4$.

Now, we turn to the other terms of the bracket, and we prove that each of them belongs to the remainder of \eqref{eq:normalized-field-estimate}. Firstly, the term $\alpha\left(V-2\xi\dot V\right)$ depends on $\Theta$ only one order later, through
	\begin{equation*}
		\Theta\,\xi^{\frac{k+1}{2}} - 2\xi\cdot\frac{k+1}{2}\Theta\,\xi^{\frac{k-1}{2}} = -k\,\Theta\,\xi^{\frac{k+1}{2}},
	\end{equation*}
and contributes $\alpha k\,\Theta\,\xi^{\frac{k+1}{2}}$ to the bracket. Since $r^{(k+1)\alpha-1}\xi^{\frac{k+1}{2}}=r^{-1}$, this contribution is $k\alpha\Theta\,r^{-1}$ after the normalization, and its quotient by $(1+|\Theta-L_k|)r^{2\alpha-1}$ is at most $k\alpha(1+|L_k|)r^{-2\alpha}$. Next, the $\Theta$-dependent part of the quadratic remainder of $\mathsf G$ is $O\!\left(\Theta^2\xi^{k+1}\right)+O\!\left(|\Theta|\,\xi^{\frac{k+3}{2}}\right)$. Then, by dividing by $\xi$ and multiplying by $r^{(k+1)\alpha-1}$, we obtain
	\begin{align*}
		r^{(k+1)\alpha-1}\xi^{-1}\,\Theta^2\xi^{k+1} &= \Theta^2r^{-(k+1)\alpha}\,r^{2\alpha-1}, \\
		r^{(k+1)\alpha-1}\xi^{-1}\,|\Theta|\,\xi^{\frac{k+3}{2}} &= |\Theta|\,r^{-2\alpha}\,r^{2\alpha-1}.
	\end{align*}
Since $|\Theta|\leq\left(1+|L_k|\right)\left(1+|\Theta-L_k|\right)$, the quotient of both terms by $(1+|\Theta-L_k|)r^{2\alpha-1}$ is
	\begin{equation*}
		O\!\left(\left(1+|\Theta|\right)r^{-(k+1)\alpha}+r^{-2\alpha}\right) = O\!\left(r^{-(k+1)\alpha}+\eps(r)r^{-k\alpha}+r^{-2\alpha}\right),
	\end{equation*}
which tends to zero. Finally, the $\Theta$-dependent part of the cubic remainder in \eqref{eq:endpoint-expansion-gplus} is smaller than the quadratic terms by a factor $|V-1|=O(\xi+\eps(r))$, and the remaining part contributes $O(\xi^3)\,r^{3\alpha-1}=O\!\left(r^{-3\alpha-1}\right)$ to $\mathsf F$, which becomes $O\!\left(r^{\alpha-1}\right)=o\!\left(r^{2\alpha-1}\right)$ after the strongest normalization by $r^{4\alpha}$. We note that this is the only place where we use the $\mathcal C^3$ hypothesis, finalizing the proof.
\end{proof}

We use the estimate \textup{(E$_k$)} as follows. A solution of $v'=\mathsf F(r,v)$ written as $v=P_k+\Theta\,r^{-k\alpha}$ satisfies $\mathsf F(r,v)-P_k'=\frac{d}{dr}\left(\Theta\,r^{-k\alpha}\right)$, and $r^{k\alpha}\frac{d}{dr}\left(\Theta\,r^{-k\alpha}\right)=\Theta'-k\alpha\Theta\,r^{-1}$. Therefore, the left-hand side of \textup{(E$_k$)} is $\Theta'$, which is the form required by Lemma \ref{lem:stable-scalar-bootstrap}\textup{(ii)}.

\begin{proof}[Proof of Theorem~\ref{thm:positive-bowl-asymptotics}\textup{(i)}]
Firstly, we have $\Theta_0=v-r^\alpha=o(r^\alpha)$ by Proposition \ref{prop:positive-end-at-infinity}\textup{(i)}. Since $v'=\mathsf F(r,v)$, the function $\Theta_0$ is of class $\mathcal C^1$, and by \textup{(E$_0$)} we have
	\begin{equation*}
		\Theta_0' = -\frac{a_1}{\alpha}\Theta_0 r^{2\alpha-1} + o\!\left((1+|\Theta_0|)r^{2\alpha-1}\right).
	\end{equation*}
We note that the weight $\phi=r^{2\alpha-1}$ has a divergent integral. Then, by Lemma \ref{lem:stable-scalar-bootstrap}\textup{(ii)} with $c=-\alpha^{-1}a_1<0$ and $L=0$, we obtain $\Theta_0(r)\to0$.

Next, we repeat this argument four times. At step $k=1,\ldots,4$ we write $v=P_k+\Theta_kr^{-k\alpha}$. Since $P_k=P_{k-1}+L_{k-1}r^{-(k-1)\alpha}$, we have $\Theta_k=r^\alpha(\Theta_{k-1}-L_{k-1})$, and from the previous step we obtain $\Theta_k=o(r^\alpha)$, as Lemma \ref{lem:uniform-expansions-nondegenerate} requires. Then, by \textup{(E$_k$)} and Lemma \ref{lem:stable-scalar-bootstrap}\textup{(ii)}, with the same weight and equilibrium $L_k$, we have $\Theta_k(r)\to L_k$, that is, $v=P_k+L_kr^{-k\alpha}+o(r^{-k\alpha})$. Consequently, since $L_1=-a_\gamma$, $L_2=0$, $L_3=b_\gamma$ and $L_4=0$, from the four steps we obtain successively
	\begin{equation*}
		v(r) = r^\alpha-a_\gamma r^{-\alpha}+o(r^{-\alpha}),
	\end{equation*}
the remainders $o(r^{-2\alpha})$ and $o(r^{-3\alpha})$, and the expansion of Theorem \ref{thm:positive-bowl-asymptotics}\textup{(i)} with a remainder $o(r^{-4\alpha})$.

Finally, let us assume that $\widetilde\gamma$ is only of class $\mathcal C^2$ near $(0,1)$. Then the cubic remainder in \eqref{eq:endpoint-expansion-gplus} becomes $o(|1-y|^2)$. In this case the estimates \textup{(E$_0$)} through \textup{(E$_3$)} still hold, while the last stage is left with an error $o(r^{3\alpha-1})$. Therefore, the chain ends at \textup{(E$_3$)}, and we obtain the expansion with a remainder $o(r^{-3\alpha})$ and the same coefficients $a_\gamma$ and $b_\gamma$.
\end{proof}

\subsection{Degenerate positive branches}\label{subsec:degenerate-positive-branches}
Recall that a degenerate branch has no finite endpoint, and that its phase tends to infinity by Lemma \ref{lem:phase-limit-positive-end}. In this case the decay of $g_+$ determines the growth of the end through the exponent $e_\gamma=\frac{3\alpha-1-k_\gamma}{\alpha}$ of $v$ in \eqref{eq:degenerate-power-main-ode}. We note that this exponent is smaller than, equal to or larger than $1$ according as $k_\gamma$ is larger than, equal to or smaller than $2\alpha-1$. In addition, the second-order expansion has a resonance at $k_\gamma+1=4\alpha$.

\begin{proof}[Proof of Theorem~\ref{thm:positive-bowl-asymptotics}\textup{(ii)}]
Since $k_\gamma>2\alpha-1$, we have $1-e_\gamma=\frac{k_\gamma+1-2\alpha}{\alpha}>0$. Then, by multiplying \eqref{eq:degenerate-power-main-ode} by $(1-e_\gamma)v^{-e_\gamma}$, we obtain
	\begin{equation*}
		\frac{d}{dr}\left(v^{1-e_\gamma}\right) = (1-e_\gamma)c_\gamma r^{k_\gamma}(1+o(1)).
	\end{equation*}
We note that, since $k_\gamma>0$, the error integrates to $o(r^{k_\gamma+1})$, and so does the value at the lower limit. Therefore, we have
	\begin{equation*}
		v(r)^{1-e_\gamma} = \frac{(1-e_\gamma)c_\gamma}{k_\gamma+1}r^{k_\gamma+1} + o(r^{k_\gamma+1}),
	\end{equation*}
and by raising to the power $(1-e_\gamma)^{-1}$ we conclude that $v(r)=A_\gamma r^{d_\gamma}(1+o(1))$, with $d_\gamma=\frac{k_\gamma+1}{1-e_\gamma}=\frac{\alpha(k_\gamma+1)}{k_\gamma+1-2\alpha}$ and $A_\gamma^{1-e_\gamma}=\frac{(1-e_\gamma)c_\gamma}{k_\gamma+1}=\frac{c_\gamma}{d_\gamma}$.
\end{proof}

\begin{proof}[Proof of Theorem~\ref{thm:positive-bowl-asymptotics}\textup{(iii)}]
When $k_\gamma=2\alpha-1$ we have $e_\gamma=1$. Therefore, \eqref{eq:degenerate-power-main-ode} reads $v'=c_\gamma r^{2\alpha-1}v(1+o(1))$, that is,
	\begin{equation*}
		\frac{d}{dr}\log v(r) = c_\gamma r^{2\alpha-1}(1+o(1)),
	\end{equation*}
which is the first assertion. Next, we integrate from a fixed large radius $\bar r$. As in \textup{(ii)}, the error contributes $o(r^{2\alpha})$, and $\int_{\bar r}^{r}s^{2\alpha-1}\,ds=\frac{1}{2\alpha}r^{2\alpha}+O(1)$. Consequently, we obtain $\log v(r)=\frac{c_\gamma}{2\alpha}r^{2\alpha}+o(r^{2\alpha})$.
\end{proof}

\begin{proof}[Proof of Theorem~\ref{thm:positive-bowl-asymptotics}\textup{(iv)}]\label{proof:bowl-asymptotics-iv}
Recall that $v(r)=A_\gamma r^{d_\gamma}(1+o(1))$ by item \textup{(ii)}. We divide the proof into three steps.

\emph{Step 1.} Firstly, we refine the positive graphical equation \eqref{eq:degenerate-power-main-ode} to second order. By the second-order large-slope expansion \eqref{eq:eta-large-slope-second-order}, we have $y_v=r^{-1}v^{1/\alpha}(1-\beta v^{-2}+O(v^{-4}))$, and therefore $y_v^{-(2k_\gamma+1)}=r^{2k_\gamma+1}v^{-\frac{2k_\gamma+1}{\alpha}}(1+o(1))$. Then, from the second-order expansion of $g_+$ we obtain
	\begin{equation*}
		g_+(y_v) = c_\gamma\,r^{k_\gamma} v^{-k_\gamma/\alpha} \left[ 1 + k_\gamma\beta v^{-2} + c_{\gamma,1}c_\gamma^{-1}r^{k_\gamma+1}v^{-(k_\gamma+1)/\alpha} + o(r^{-2d_\gamma}) \right],
	\end{equation*}
where $r^{k_\gamma+1}v^{-\frac{k_\gamma+1}{\alpha}}=O(r^{-2d_\gamma})$ because $(k_\gamma+1)\left(1-\frac{d_\gamma}{\alpha}\right)=-2d_\gamma$. Next, by multiplying by $(1+v^2)^{\beta+1}=v^{(3\alpha-1)/\alpha}\left(1+(\beta+1)v^{-2}+O(v^{-4})\right)$ and by substituting $v=A_\gamma r^{d_\gamma}(1+o(1))$ in the correction terms, we obtain
	\begin{equation}\label{eq:second-order-degenerate-reduced-ode}
		v' = c_\gamma\,r^{k_\gamma} v^{\frac{3\alpha-1-k_\gamma}{\alpha}} \left[ 1+\mathcal C_\gamma r^{-2d_\gamma}+o(r^{-2d_\gamma}) \right],
	\end{equation}
with $\mathcal C_\gamma$ as in the statement, since $(1+\beta(k_\gamma+1))v^{-2}=\left(1+\beta(k_\gamma+1)\right)A_\gamma^{-2}r^{-2d_\gamma}(1+o(1))$ and $c_{\gamma,1}c_\gamma^{-1}r^{k_\gamma+1}v^{-(k_\gamma+1)/\alpha}=c_{\gamma,1}c_\gamma^{-1}A_\gamma^{-\frac{k_\gamma+1}{\alpha}}r^{-2d_\gamma}(1+o(1))$.

\emph{Step 2.} Next, we normalize around the leading term of item \textup{(ii)}. Let $v=A_\gamma r^{d_\gamma}(1+z)$ with $z\to0$. We use the identities $k_\gamma+d_\gamma\,e_\gamma=d_\gamma-1$ and $c_\gamma\,A_\gamma^{e_\gamma}=d_\gamma\,A_\gamma$, the latter being the relation for $A_\gamma$ from \textup{(ii)}, and the substitution $v'=A_\gamma r^{d_\gamma-1}\left[d_\gamma(1+z)+rz'\right]$. Then \eqref{eq:second-order-degenerate-reduced-ode} becomes $d_\gamma(1+z)+rz'=d_\gamma(1+z)^{e_\gamma}\left[1+\mathcal C_\gamma r^{-2d_\gamma}+o(r^{-2d_\gamma})\right]$. Therefore, by the expansion $(1+z)^{e_\gamma}=1+e_\gamma z+O(z^2)$ and $d_\gamma(e_\gamma-1)=-(k_\gamma+1)$, we obtain
	\begin{equation}\label{eq:second-order-z-equation}
		z' = -\frac{k_\gamma+1}{r}z + d_\gamma\,\mathcal C_\gamma r^{-1-2d_\gamma} + O\!\left(\frac{z^2}{r}\right) + o(r^{-1-2d_\gamma}).
	\end{equation}

\emph{Step 3.} Finally, we integrate \eqref{eq:second-order-z-equation}, where the linear term $-\frac{k_\gamma+1}{r}z$ competes with the source of order $r^{-1-2d_\gamma}$. For this reason, we distinguish two cases.

Firstly, we consider the nonresonant case $k_\gamma+1>4\alpha$, that is, $k_\gamma+1>2d_\gamma$. Let $C_1$ bound the quadratic error, $\left|O\!\left(z^2r^{-1}\right)\right|\leq C_1z^2r^{-1}$ on a tail, and for $\delta>0$ let $R_\delta$ be such that $|z|\leq\delta$ and the last error is at most $\delta r^{-1-2d_\gamma}$ on $[R_\delta,\infty)$. We take $\delta$ so small that $\Lambda=k_\gamma+1-C_1\delta>2d_\gamma$. Then, since $z^2\leq\delta|z|$ there, the upper right Dini derivative of $|z|$ satisfies
	\begin{equation*}
		D^+|z| \leq -\frac{\Lambda}{r}|z| + C_2r^{-1-2d_\gamma}, \qquad C_2=|d_\gamma\,\mathcal C_\gamma|+\delta,
	\end{equation*}
on $[R_\delta,\infty)$, at the zeros of $z$ as well, where $D^+|z|=|z'|$. Therefore, we have $D^+\!\left(r^\Lambda|z|\right)\leq C_2r^{\Lambda-1-2d_\gamma}$, and by integrating from $R_\delta$ we obtain $r^\Lambda|z(r)|\leq R_\delta^\Lambda|z(R_\delta)|+\frac{C_2}{\Lambda-2d_\gamma}r^{\Lambda-2d_\gamma}$. Consequently, $z=O(r^{-2d_\gamma})$ and $Y=r^{2d_\gamma}z$ is bounded. Now, by the product rule, \eqref{eq:second-order-z-equation} becomes
	\begin{equation*}
		Y' = -\frac{k_\gamma+1-2d_\gamma}{r}\left(Y-\frac{d_\gamma\,\mathcal C_\gamma}{k_\gamma+1-2d_\gamma}\right) + O\!\left(\frac{Y^2}{r^{2d_\gamma+1}}\right) + o\left(\frac{1}{r}\right),
	\end{equation*}
where the quadratic term is $o(r^{-1})$ because $Y$ is bounded. Then, by Lemma \ref{lem:stable-scalar-bootstrap}\textup{(ii)}, with weight $\phi=r^{-1}$ and $c=-(k_\gamma+1-2d_\gamma)<0$, we have $Y(r)\to Y_*=\frac{d_\gamma\,\mathcal C_\gamma}{k_\gamma+1-2d_\gamma}$, and therefore $v(r)=A_\gamma r^{d_\gamma}+A_\gamma Y_*r^{-d_\gamma}+o(r^{-d_\gamma})$. Finally, since $\frac{d_\gamma}{k_\gamma+1-2d_\gamma}=\frac{\alpha}{k_\gamma+1-4\alpha}$, we obtain $B_\gamma=A_\gamma Y_*=\frac{\alpha A_\gamma\mathcal C_\gamma}{k_\gamma+1-4\alpha}$.

Next, we consider the resonant case $k_\gamma+1=4\alpha$, that is, $k_\gamma+1=2d_\gamma$. We again set $Y=r^{2d_\gamma}z$. In this case the term $2d_\gamma r^{-1}Y$ from the product rule cancels the linear term of \eqref{eq:second-order-z-equation}, and we are left with
	\begin{equation}\label{eq:second-order-resonant-Y}
		Y' = d_\gamma\,\mathcal C_\gamma r^{-1} + O\!\left(Y^2r^{-(2d_\gamma+1)}\right) + o\left(r^{-1}\right).
	\end{equation}
Since $Y^2r^{-(2d_\gamma+1)}=|Y|\,|z|\,r^{-1}$ and $z\to0$, for every $\delta>0$ there is $R_\delta$ with $|Y'|\leq(|d_\gamma\,\mathcal C_\gamma|+\delta)r^{-1}+C_0\delta|Y|r^{-1}$ on $[R_\delta,\infty)$, for a constant $C_0$ bounding the quadratic remainder. Then, by Gr\"onwall's inequality, we have $|Y(r)|\leq\left(|Y(R_\delta)|+(|d_\gamma\,\mathcal C_\gamma|+\delta)\log(r/R_\delta)\right)(r/R_\delta)^{C_0\delta}$. Next, we choose $\delta$ with $C_0\delta<d_\gamma$, so that $(r/R_\delta)^{C_0\delta}$ grows more slowly than $r^{d_\gamma}$. We notice that this means $Y=O(r^{\eps})$ for some $\eps<d_\gamma$, and hence the quadratic term in \eqref{eq:second-order-resonant-Y} is $O(r^{2\eps-2d_\gamma-1})=o(r^{-1})$. Thus, \eqref{eq:second-order-resonant-Y} becomes $Y'=d_\gamma\,\mathcal C_\gamma r^{-1}+o(r^{-1})$, and integrating it from a fixed radius we observe that $Y(r)=d_\gamma\,\mathcal C_\gamma\log r+o(\log r)$. Finally, returning to $z=r^{-2d_\gamma}Y$, we conclude that $z(r)=d_\gamma\,\mathcal C_\gamma r^{-2d_\gamma}\log r+o(r^{-2d_\gamma}\log r)$, which is the expansion of the statement with $L_\gamma=A_\gamma d_\gamma\mathcal C_\gamma$.
\end{proof}

\subsection{The Gauss-curvature root}\label{subsubsec:gauss-root-bowl}

Now, we consider a model that realizes every degenerate regime of Theorem \ref{thm:positive-bowl-asymptotics} in closed form, the family $\gamma=S_n^{1/n}=\left(\lambda_1\cdots\lambda_n\right)^{1/n}$ on $\Gamma=\Gamma_+$. We observe that $\gamma$ is a $1$-positive curvature function, since it is symmetric and positively $1$-homogeneous, with $\partial_i\gamma=\gamma(n\lambda_i)^{-1}>0$ on $\Gamma_+$. Moreover, its rotational restriction $\widetilde\gamma(x,y)=\left(xy^{n-1}\right)^{1/n}$ satisfies $\widetilde\gamma(t,1)=t^{1/n}\to0$ as $t\to0^+$, and therefore $\gamma$ is $1$-degenerate. Since only the positive cone and the positive branch enter the bowl analysis, we do not need a signed extension, and $n$ may be even or odd.

Firstly, we compute the positive branch. The level-one equation has the explicit solution $g_+(y)=y^{-(n-1)}$ on $I_+=(1,\infty)$, with $y_*=\widetilde\gamma(1,1)^{-1}=1$. In particular, the positive branch has balanced second-order power decay in the sense of Definition \ref{def:degenerate-positive-branch-asymptotics}, with $c_\gamma=1$, $k_\gamma=n-1$, $c_{\gamma,1}=0$ and vanishing remainders. The general solution of the first integral below differs from the bowl-type solution by a constant of integration. The leading growth of its positive end for $n=2$ and $n\geq3$ is stated in the unpublished preprint \cite[Theorem 1.3]{RengaswamiTorresSantaella2023Bowl}, and the expansions below recover it together with the lower-order terms.

Next, we integrate the graphical equation. Since $\alpha=1$, we have $\beta=0$ and $Q(r,v)=\frac{v}{r}$, and the graphical equation \eqref{eq:selected-branch-graphical-equation} reduces to the separable equation $v'v^{n-1}=r^{n-1}(1+v^2)$, with no large-slope approximation. Let
\begin{equation*}
	F_n(v) = \int_0^v\frac{t^{n-1}}{1+t^2}\,dt, \qquad F_1(v)=\arctan v, \qquad F_2(v)=\frac{1}{2}\log\left(1+v^2\right),
\end{equation*}
with $F_n(v)=\frac{v^{n-2}}{n-2}-F_{n-2}(v)$ for $n\geq3$. Then the bowl-type solution, with $v(0)=0$, satisfies $F_n(v(r))=\frac{r^n}{n}$ for $r\geq0$. We note that, when $n\geq2$, the function $F_n$ increases from $0$ to $+\infty$. Therefore, this identity determines $v(r)$ at every radius, and the solution is entire.

Now, we compare this solution with Theorem \ref{thm:positive-bowl-asymptotics} in each regime. When $n=2$, the decay is critical, $k_\gamma=1=2\alpha-1$, and $v(r)=\sqrt{e^{r^2}-1}$. Since $\log v=\frac{r^2}{2}+\frac{1}{2}\log\left(1-e^{-r^2}\right)$, we obtain $\log v(r)=\frac{r^2}{2}-\frac{1}{2}e^{-r^2}(1+o(1))$. This agrees with item \textup{(iii)} for $c_\gamma=1$, with an error of order $e^{-r^2}$ instead of $o(r^{2\alpha})$. In addition, the slope integrates to an entire strictly convex $S_2^{1/2}$-translator in $\rr^3$, which is the case $p=\frac{1}{2}$ of the classification of rotational $K^p$-translators of $\rr^3$ by Aydin and L\'opez \cite[Theorem 1]{aydin2023translators}.

When $n\geq3$, the decay exceeds the critical rate, and by item \textup{(ii)} we have $v(r)=A_\gamma r^{\frac{n}{n-2}}(1+o(1))$ with $d_\gamma=\frac{n}{n-2}$ and $A_\gamma=\left(\frac{n-2}{n}\right)^{\frac{1}{n-2}}$, so that $A_\gamma^{n-2}=\frac{c_\gamma}{d_\gamma}$. Moreover, when $n\geq5$, the nonresonant case $k_\gamma+1=n>4\alpha$ of item \textup{(iv)} applies with $\mathcal C_\gamma=A_\gamma^{-2}$, since $\beta=0$ and $c_{\gamma,1}=0$, and we obtain
\begin{equation*}
	v(r) = A_\gamma r^{\frac{n}{n-2}} + \frac{1}{(n-4)A_\gamma}\,r^{-\frac{n}{n-2}} + o\left(r^{-\frac{n}{n-2}}\right).
\end{equation*}
We note that the second term is the term $-\frac{v^{n-4}}{n-4}$ of the first integral, since $A_\gamma^{n-2}=\frac{n-2}{n}$.

Finally, when $n=4$, the resonance $k_\gamma+1=4\alpha$ holds. In this case, from the first integral $\frac{v^2}{2}-\frac{1}{2}\log(1+v^2)=\frac{r^4}{4}$ we have $v^2=\frac{r^4}{2}+4\log r-\log2+O(r^{-4}\log r)$, and therefore
\begin{equation*}
	v(r) = \frac{r^2}{\sqrt2} + 2\sqrt2\,\frac{\log r}{r^{2}} - \frac{\log2}{\sqrt2}\,r^{-2} + o\left(r^{-2}\right).
\end{equation*}
This agrees with the resonant case of item \textup{(iv)}, where $A_\gamma=2^{-\frac{1}{2}}$, $d_\gamma=2$, $\mathcal C_\gamma=2$ and $L_\gamma=A_\gamma d_\gamma\mathcal C_\gamma=2\sqrt2$. In particular, the closed form shows that the error $o(r^{-d_\gamma}\log r)$ of the general statement is in fact $O(r^{-2})$.
\section{Catenoidal translators}\label{sec:catenoidal-translators}

In this section we prove Theorem \ref{thm:main-c}, stated in the introduction, on the existence and classification of catenoidal continuations according to their lower branch. The precise statement is Theorem \ref{thm:rotational-translators}. For this purpose, we use the tools developed in Sections \ref{sec:curvature-functions-branches}, \ref{sec:barrier-method} and \ref{sec:positive-branch-asymptotics}. More precisely, we construct the neck through the zero level of the curvature function, and we match it with the level-one continuations on the upper side and on the reoriented lower side. Then, we show that the upper continuation is asymptotic to the bowl-type solution, and that the lower continuation depends on the type of the lower branch.

Throughout this section $\gamma$ is an $\alpha$-positive curvature function with $\alpha>\frac{1}{2}$ (see Remark \ref{rem:height-constant-threshold}), with a signed rotational level profile $\widehat\gamma$ and selected angular neck data $(A_g,\rho_0)$ in the sense of Definition \ref{def:selected-signed-neck-datum}. We use the selected chart $g$ on $\mathcal D$ of Proposition \ref{prop:conical-sheet-maximal-chart} and the maximal regular continuation $(I_-,g_-)$ of Definition \ref{def:neck-accessible-branch}. The regularity of the data enters only through the following three conditions.
\begin{enumerate}[label=\textup{(\alph*)}, ref=\textup{(\alph*)}]
	\item \label{cat:all-radii}
	Every radius $R>0$ is neck-admissible in the sense of Definition \ref{def:neck-regular-radius}.
	\item \label{cat:upper-chart}
	The set $A_g$ contains $0$, and $\Psi_g$ is of class $\mathcal C^1$ with $\Psi_g'>0$ on $(-\rho_0,0]$.
	\item \label{cat:lower-branch}
	The neck-accessible reflected branch admits a maximal regular continuation $(I_-,g_-)$.
\end{enumerate}

We use condition \ref{cat:upper-chart} on the upper side of the neck, between the zero direction and the cylindrical direction, where it makes the upper level-one chart of class $\mathcal C^1$ (Lemma \ref{lem:upper-entry-sheet-from-zero-sector}). Moreover, we note that when $\gamma$ is normalized $1$-nondegenerate, as in Theorem \ref{thm:rotational-translators}, the inequality $\Psi_g'(0)>0$ is automatic once $0\in A_g$. Indeed, in this case $(0,1)$ lies in the open set $\widetilde\Gamma$, where $\widehat\gamma=\widetilde\gamma$, and by ellipticity we have $\Psi_g'(0)=\partial_x\widetilde\gamma(0,1)>0$.

We note that regular selected data with $0\in A_g$ satisfy the three conditions. Firstly, condition \ref{cat:all-radii} follows from \eqref{eq:zero-trace-conical} and the observation after Definition \ref{def:neck-regular-radius}. Moreover, condition \ref{cat:upper-chart} follows from Lemma \ref{lem:reduced-ellipticity} and \ref{hyp:regular-ellipticity}, and condition \ref{cat:lower-branch} follows from the remark after Definition \ref{def:neck-accessible-branch}. In addition, that remark covers data that fail to be regular only on the zero ray. This is the case of two model families of Subsection \ref{subsec:model-examples-branches}, and consequently the classification covers their catenoidal translators. The first family consists of the powers $\left(\frac{H}{n-1}\right)^p$ with $p>1$, whose angular profile vanishes to order $p$ at the zero ray (Subsubsection \ref{subsubsec:axis-crossing-models}). The second family consists of the odd-root models, whose angular profile is not of class $\mathcal C^1$ there (Subsubsection \ref{subsubsec:ellipticity-limited-models}).

\subsection{The local neck}\label{subsec:local-catenoidal-neck}

Let $\mathtt c(s)=(r(s),u(s))$ be a generating curve parametrized by arclength, with $\dot{\mathtt c}=(\cos\theta,\sin\theta)$ as in Subsection \ref{subsec:rotational-generating-curves}. On the selected sheet, the arclength form \eqref{eq:rotational-translator-arclength} of the translator equation for $\widehat\gamma$ reduces, by Proposition \ref{prop:conical-sheet-maximal-chart}, to the system
\begin{equation}\label{eq:local-neck-phase-system}
		\dot r=\cos\theta,\quad \dot u=\sin\theta,\quad \dot\theta=g\left(r^{-1}\sin\theta,\cos\theta\right),\quad 
		\left(r^{-1}\sin\theta,\cos\theta\right)\in\mathcal D.
\end{equation}

\begin{proposition}\label{prop:local-catenoidal-neck}
	Let $R>0$. There exist $\eps_R>0$ and a solution $\left(r_R,u_R,\theta_R\right)$ of \eqref{eq:local-neck-phase-system} on $(-\eps_R,\eps_R)$ satisfying $\left(r_R,u_R\right)(0)=(R,0)$ and $\theta_R(0)=\frac{\pi}{2}$, with $\theta_R(s)\in(0,\pi)\setminus\{\frac{\pi}{2}\}$ for all $0<|s|<\eps_R$. Any two solutions with these initial values coincide wherever both are defined in $(-\eps_R,\eps_R)$. The generating curve $\mathtt c_R=\left(r_R,u_R\right)$ is of class $\mathcal C^2$, with $\dot\theta_R(0)=-\rho_0R^{-1}$ and $\ddot r_R(0)=\rho_0R^{-1}>0$. After decreasing $\eps_R$, there is $\delta_R>0$ such that the branches $s>0$ and $s<0$ form vertical graphs $u_R^+$ and $u_R^-$ over $(R,R+\delta_R)$. If $R$ is $\varkappa$-regular at the neck in the sense of Definition \ref{def:neck-regular-radius}, then $\dot\theta_R\in\mathcal C^{0,\varkappa}$ near $s=0$ and $\mathtt c_R$ is of class $\mathcal C^{2,\varkappa}$ there.
\end{proposition}

\begin{proof}
Firstly, we fix $R>0$. Recall that $R$ is neck-admissible by condition \ref{cat:all-radii}. Then, by Definition \ref{def:neck-regular-radius}, we have $(R^{-1},0)\in\mathcal D$, and there exist an open rectangle $\mathcal V=J\times K\subset\mathcal D$ containing $(R^{-1},0)$ and $L>0$ such that $g(\cdot,z)$ is $L$-Lipschitz for each $z\in K$. Moreover, the chart is continuous on $\mathcal D$ by Proposition \ref{prop:conical-sheet-maximal-chart}. Then, since $g(R^{-1},0)=-\rho_0R^{-1}<0$ by \eqref{eq:zero-trace-conical}, after shrinking $J$ and $K$ we may assume $g<-\frac{\rho_0}{2R}$ on $\mathcal V$ and $|z|\leq\frac{1}{2}$ on $K$.

Next, we consider a solution of \eqref{eq:local-neck-phase-system} whose phase stays in $\mathcal V$ with $\theta$ near $\frac{\pi}{2}$. We note that $z=\cos\theta\in[-\frac{1}{2},\frac{1}{2}]$ and $\sin\theta=\sqrt{1-z^2}\geq\frac{\sqrt{3}}{2}$. Then, since $\dot\theta<-\frac{\rho_0}{2R}$ on $\mathcal V$, we have
	\begin{equation}\label{eq:neck-angle-speed}
		\dot z = -\sin\theta\,\dot\theta = -\sqrt{1-z^2}\; g\!\left(r^{-1}\sqrt{1-z^2},z\right) > 0.
	\end{equation}
We note that at $z=0$ this equals $\rho_0R^{-1}$. Therefore, $z$ is strictly increasing, vanishes only at $s=0$, and serves as a local parameter. Then, by dividing the first two equations of \eqref{eq:local-neck-phase-system} by \eqref{eq:neck-angle-speed}, we obtain
	\begin{equation}\label{eq:neck-inverted-system}
		\frac{dr}{dz} = \frac{-z}{\sqrt{1-z^2}\;g\!\left(r^{-1}\sqrt{1-z^2},z\right)}, \qquad \frac{du}{dz} = \frac{-1}{g\!\left(r^{-1}\sqrt{1-z^2},z\right)},
	\end{equation}
with $r(0)=R$ and $u(0)=0$.

Now, we solve \eqref{eq:neck-inverted-system}. We note that the right-hand side of the first equation in \eqref{eq:neck-inverted-system} is continuous in $(r,z)$ where $\left(r^{-1}\sqrt{1-z^2},z\right)\in\mathcal V$. Moreover, it is Lipschitz in $r$ for $r\geq\frac{R}{2}$, because $g(\cdot,z)$ is $L$-Lipschitz, the map $r\mapsto r^{-1}\sqrt{1-z^2}$ has a derivative bounded by $4R^{-2}$, and the denominator satisfies $\sqrt{1-z^2}\,|g|\geq\frac{\sqrt{3}\rho_0}{4R}$. Then, by the Picard--Lindel\"of theorem \cite[Theorem 2.2]{teschl2012ordinary}, we have a unique solution $r(z)$ near $z=0$, and we obtain $u(z)$ by integration.

Subsequently, we recover the arclength from
	\begin{equation*}
		s(z) = -\int_0^z\left(\sqrt{1-t^2}\, g\!\left(r(t)^{-1}\sqrt{1-t^2},t\right)\right)^{-1}dt,
	\end{equation*}
which is of class $\mathcal C^1$ with $s'>0$ near $z=0$. Then, we invert it to $z(s)$ on an interval $(-\eps_R,\eps_R)$ along which $\left(r(z(s))^{-1}\sqrt{1-z(s)^2},z(s)\right)\in\mathcal V$. In particular, we have $|z(s)|\leq\frac{1}{2}$, and $z(s)$ has the sign of $s$.

Next, let $r_R(s)=r(z(s))$, $u_R(s)=u(z(s))$ and $\theta_R(s)=\arccos z(s)$. Then $\theta_R(0)=\frac{\pi}{2}$ and $\theta_R(s)\in(0,\pi)\setminus\{\frac{\pi}{2}\}$ for $0<|s|<\eps_R$. Moreover, $\theta_R$ is of class $\mathcal C^1$ because $|z|<1$. Then, since $z'(s)=-\sqrt{1-z^2}\,g\!\left(r^{-1}\sqrt{1-z^2},z\right)$, by the chain rule we obtain
	\begin{align*}
		\dot r_R &= \frac{dr}{dz}\,z'(s) = z = \cos\theta_R, \qquad \dot u_R = \frac{du}{dz}\,z'(s) = \sqrt{1-z^2} = \sin\theta_R, \\
		\dot\theta_R &= \frac{-z'(s)}{\sqrt{1-z^2}} = g\!\left(r_R^{-1}\sin\theta_R,\cos\theta_R\right),
	\end{align*}
and therefore $(r_R,u_R,\theta_R)$ solves \eqref{eq:local-neck-phase-system}.

Now, we prove the uniqueness. For this purpose, we consider any solution of \eqref{eq:local-neck-phase-system} with these initial values. We note that its phase starts at $(R^{-1},0)$ and stays in $\mathcal V$ for small $|s|$. Then \eqref{eq:neck-angle-speed} holds there and the radius solves \eqref{eq:neck-inverted-system}. Consequently, by uniqueness in the variable $z$, this solution agrees with $\mathtt c_R$ for small $|s|$. Moreover, at every point of agreement the phase is again in $\mathcal V$, and therefore the set of agreement is open. Then, since this set is closed by continuity, it follows that the two solutions coincide on their common interval of definition in $(-\eps_R,\eps_R)$.

Next, we prove the regularity of the curve. Since the field of \eqref{eq:local-neck-phase-system} is continuous, we have $\theta_R\in\mathcal C^1$ and $\mathtt c_R\in\mathcal C^2$. Moreover, $\dot\theta_R(0)=g(R^{-1},0)=-\rho_0R^{-1}$, and from $\ddot r_R=-\sin\theta_R\,\dot\theta_R$ we obtain $\ddot r_R(0)=\rho_0R^{-1}>0$. In particular, we have
	\begin{equation*}
		r_R(s)=R+\frac{\rho_0}{2R}s^2+o(s^2), \qquad \dot r_R(s)=\frac{\rho_0}{R}s+o(s), \qquad u_R(s)=s+o(s).
	\end{equation*}
Recall that $\dot r_R=z$ has the sign of $s$. Then the radius decreases strictly on $(-\eps_R,0)$ and increases strictly on $(0,\eps_R)$, with $r_R>R$ for $s\neq0$. Therefore, each half of the curve is a vertical graph $u_R^\pm$ over $(R,R+\delta_R)$.

Finally, let us assume that $R$ is $\varkappa$-regular at the neck. Then, by \eqref{eq:neck-angle-speed}, we have $z\in\mathcal C^1$ with $z'(0)=\rho_0R^{-1}$, and $y_R=r_R^{-1}\sin\theta_R$ is of class $\mathcal C^1$. Then, by splitting the increment of $\dot\theta_R=g(y_R,z)$ in $\mathcal V$ and by using the Lipschitz bound in $y$ and the $\varkappa$-H\"older bound in $z$, we obtain
	\begin{equation*}
		\left|\dot\theta_R(s_1)-\dot\theta_R(s_2)\right| \leq L\left|y_R(s_1)-y_R(s_2)\right| + C\left|z(s_1)-z(s_2)\right|^\varkappa \leq C'\left|s_1-s_2\right|^\varkappa
	\end{equation*}
for nearby $s_1,s_2$, since $\varkappa\leq1$. Consequently, $\dot\theta_R\in\mathcal C^{0,\varkappa}$ and $\mathtt c_R\in\mathcal C^{2,\varkappa}$, finalizing the proof.
\end{proof}

\begin{remark}\label{rem:neck-horizontal-graph}
We note that near the neck the generating curve is a horizontal graph $r=r(u)$, and that in this parametrization \eqref{eq:local-neck-phase-system} is a second-order equation whose right-hand side is $g$ composed with $(r,r')$. On the one hand, when the selected data are regular, $g$ is of class $\mathcal C^1$ by Proposition \ref{prop:conical-sheet-maximal-chart}. Then, by the Picard--Lindel\"of theorem, we obtain the neck directly, with $r\in\mathcal C^3$. On the other hand, this route is not available for the two model families listed at the beginning of the section, since there $g$ is only H\"older in $z=\cos\theta$ at $z=0$. For this reason, in the proof of Proposition \ref{prop:local-catenoidal-neck} we use $z$ as the independent variable, which \eqref{eq:neck-angle-speed} allows. In this way, the regularity of $g$ in $z$ enters only through continuity, and for the uniqueness we use the Lipschitz bound in the rotational variable of Definition \ref{def:neck-regular-radius}.
\end{remark}

When $p=1$, the chart of $\frac{H}{n-1}$ is smooth, and so is the generating curve. In this case we recover the graph geometry of \cite[Lemma 2.3]{clutterbuck2007stability}. In contrast, when $p>1$, the exponent $\varkappa=\frac{1}{p}$ is optimal. To see this, recall that $\rho_0=n-1$ and $g(y,z)=(n-1)\left(\operatorname{sgn}(z)|z|^{\frac{1}{p}}-y\right)$ by Subsubsection \ref{subsubsec:axis-crossing-models}. Then, by \eqref{eq:neck-angle-speed}, we have $z(s)=\rho_0R^{-1}s(1+o(1))$, while $y_R(s)=R^{-1}+O(s^2)$. Moreover, since $\frac{1}{p}<2$, the $O(s^2)$ term is of lower order, and we obtain
\begin{equation}\label{eq:powers-neck-curvature-expansion}
	\dot\theta_R(s) = -\frac{n-1}{R} + (n-1)\left(\frac{n-1}{R}\right)^{\frac{1}{p}} \operatorname{sgn}(s)\left|s\right|^{\frac{1}{p}}(1+o(1)) \qquad \text{as }s\to0.
\end{equation}
Now, we note that the unit normal $(-\sin\theta_R,\cos\theta_R)$ is of class $\mathcal C^{1,\frac{1}{p}}$ across the neck circle, because $\theta_R\in\mathcal C^{1,\frac{1}{p}}$ by Proposition \ref{prop:local-catenoidal-neck}. In addition, the rotational entry $r_R^{-1}\sin\theta_R$ of the second fundamental form is of class $\mathcal C^1$. However, the meridional entry \eqref{eq:powers-neck-curvature-expansion} is $\frac{1}{p}$-H\"older and no better, since the coefficient of $\operatorname{sgn}(s)|s|^{\frac{1}{p}}$ is nonzero. Therefore, the regularity $\mathcal C^{2,\frac{1}{p}}$ of the generating curve is sharp.

\begin{lemma}\label{lem:lower-side-neck-phase}
	Let $R>0$. Represent the lower side by the outward parametrization $\overline{\mathtt c}_R(s)=\mathtt c_R(-s)$ and $\overline\theta_R(s)=\theta_R(-s)+\pi$ of Subsection \ref{subsec:reflected-lower-components}, and let $v_R^-=(u_R^-)'$ be the reoriented lower graphical slope. Then $v_R^-(r)\to-\infty$ as $r\downarrow R$, and
	\begin{equation}\label{eq:lower-phase-at-neck}
		Q(r,v_R^-(r)) \to -\infty \qquad \text{as }r\downarrow R.
	\end{equation}
	Moreover, there is $r_R^->R$ such that on $(R,r_R^-]$ the reoriented lower side is branch-admissible for $(I_-,g_-)$ and solves the lower graphical equation \eqref{eq:selected-branch-graphical-equation}.
\end{lemma}

\begin{proof}
Firstly, we note that along the lower side the slope is $\dot u_R/\dot r_R$, which does not depend on the orientation of the curve, and whose denominator tends to $0$ and whose numerator tends to $1$ as $r\downarrow R$. Moreover, in the outward parametrization $u$ decreases as $r$ increases. Therefore, we have $v_R^-<0$ and $v_R^-(r)\to-\infty$. Recall that $Q(1,\cdot)$ is an increasing bijection of $\rr$ by Lemma \ref{lem:level-set-control}\textup{(i)}, and that $Q(r,v)=r^{-1}Q(1,v)$ with $r^{-1}$ bounded near $R$. Consequently, we obtain \eqref{eq:lower-phase-at-neck}.

Next, we compute the phase of the lower graph. On the lower side, $s<0$, the phase of the solution of \eqref{eq:local-neck-phase-system} is $(x,y,z)=(\dot\theta_R,r_R^{-1}\sin\theta_R,\cos\theta_R)$, with $y>0$, $z<0$ and $x=g(y,z)$. Then, in the outward parametrization, by Subsection \ref{subsec:reflected-lower-components} we have $(\overline x,\overline y,\overline z)=(-x,-y,-z)$, with $\overline y\to-R^{-1}$ and $\overline z\to0^{+}$ at the neck. Moreover, by the computation of Section \ref{sec:preliminaries} applied to $v=v_R^-$, the level-one phase of the lower graph is $(x_v,y_v)=\overline z^{-1/\alpha}(\overline x,\overline y)$. In particular, we have $(-y_v)^{-\alpha}=\overline z\,y^{-\alpha}$, that is, $zy^{-\alpha}=-(-y_v)^{-\alpha}$. Then, by the chart of Proposition \ref{prop:conical-sheet-maximal-chart}, we obtain
	\begin{equation*}
		x_v = -\overline z^{-1/\alpha}\,g(y,z) = -\overline z^{-1/\alpha}\,y\,\Psi_g^{-1}\!\left(zy^{-\alpha}\right) = y_v\,\Psi_g^{-1}\!\left(-(-y_v)^{-\alpha}\right) = g^\sharp(y_v,1),
	\end{equation*}
by \eqref{eq:reflected-lower-angular-relation}. Therefore, the lower side solves the graphical equation of the neck-accessible reflected branch, as observed in Subsection \ref{subsec:reflected-lower-components}. Finally, since $y_v\to-\infty$ by \eqref{eq:lower-phase-at-neck}, there is $r_R^->R$ with $Q(r,v_R^-(r))<y_{\mathrm{reg}}$ on $(R,r_R^-]$. We note that on this interval the phase lies in $I_-$, where $g^\sharp(\cdot,1)=g_-$. Consequently, $v_R^-$ is branch-admissible for $(I_-,g_-)$ and solves \eqref{eq:selected-branch-graphical-equation} on that interval.
\end{proof}

\begin{definition}\label{def:catenoidal-profile}
	Let $R>0$. We say that an arclength-parametrized generating curve $\mathtt c_R:(s_-,s_+)\to\rr^2$, with $s_-<0<s_+$ and with its neck at $s=0$, is a \emph{catenoidal profile of neck radius} $R$ if it is a maximal continuation of the local neck of Proposition \ref{prop:local-catenoidal-neck} through the graphical equations of Section \ref{sec:barrier-method}. In addition, we denote by $W_R^{\max}$ the rotational hypersurface swept by the profile. If $(s_-,s_+)=\rr$, we say that it is a \emph{complete catenoidal translator}, and we denote it by $W_R$.
\end{definition}

More precisely, the upper side $s>0$ and the lower side $s<0$ of a catenoidal profile are vertical graphs over radial intervals starting at $R$. The upper side solves the graphical equation of a branch of the level-one equation. The lower side solves an equation of the same form after it is reoriented by the outward parametrization of Lemma \ref{lem:lower-side-neck-phase}.

\subsection{The upper side}\label{subsec:upper-side-entry}

Firstly, we study the upper side of the neck. We note that $y_v\to+\infty$ at the neck on this side. Then, the phase descends along the level-one restriction of the selected chart until it reaches $I_+$, as we show in Proposition \ref{prop:upper-side-global-continuation}.

\begin{lemma}\label{lem:upper-entry-sheet-from-zero-sector}
	Assume that $\gamma$ is normalized $1$-nondegenerate and that the selected data satisfy condition \ref{cat:upper-chart}. Then $\Psi_g(0)=1$, $(1,1)\in\mathcal D$ with $g(1,1)=0$, and the upper level-one chart $h_+=g(\cdot,1)$ is of class $\mathcal C^1$ on $(1,\infty)$ and satisfies
	\begin{equation*}
		h_+(y)=y\,a(y), \qquad a(y)=\Psi_g^{-1}\left(y^{-\alpha}\right)\in(-\rho_0,0), \qquad y>1,
	\end{equation*}
	which gives $-\rho_0y<h_+(y)<0$ and $\lim_{y\to1^+}h_+(y)=0$.
\end{lemma}

\begin{proof}
	Firstly, we note that $0\in A_g$ by condition \ref{cat:upper-chart}. Moreover, since $(0,1)\in\widetilde\Gamma$ by Definition \ref{def:one-nondegeneracy}, where $\widehat\gamma=\widetilde\gamma$, the normalization implies that $\Psi_g(0)=1$. In particular, we have $1\in\Psi_g(A_g)$. Therefore, $(1,1)\in\mathcal D$ and $g(1,1)=\Psi_g^{-1}(1)=0$.
	
	Next, we consider $y>1$. Then, by the chart of Proposition \ref{prop:conical-sheet-maximal-chart}, we have $h_+(y)=g(y,1)=y\,\Psi_g^{-1}\left(y^{-\alpha}\right)=y\,a(y)$. Now, recall that $\Psi_g(-\rho_0)=0$ by Definition \ref{def:selected-signed-neck-datum} and that $\Psi_g(0)=1$. Then, since $y^{-\alpha}\in(0,1)$ and $\Psi_g$ is strictly increasing on $A_g$, we obtain $-\rho_0<a(y)<0$. Consequently, $-\rho_0y<h_+(y)<0$.
	
	In addition, we note that $\Psi_g$ is of class $\mathcal C^1$ with $\Psi_g'>0$ on $(-\rho_0,0]$. Therefore, by the inverse function theorem, $\Psi_g^{-1}$ is of class $\mathcal C^1$ on $(0,1]$, and it follows that $a$ and $h_+$ are of class $\mathcal C^1$ on $(1,\infty)$. Finally, by the continuity of $\Psi_g^{-1}$, we have $a(y)\to\Psi_g^{-1}(1)=0$ as $y\to1^+$, and we conclude that $h_+(y)\to0$.
\end{proof}

\begin{proposition}\label{prop:upper-side-global-continuation}
	Assume the hypotheses of Lemma \ref{lem:upper-entry-sheet-from-zero-sector}, and let $R>0$ with $v_R^+=(u_R^+)'$. Then there exists $r_R^{\mathrm{ent}}>R$ such that $v_R^+$ extends to $(R,\infty)$ as a solution of \eqref{eq:selected-branch-graphical-equation}, following the branch $((1,\infty),h_+)$ for $R<r<r_R^{\mathrm{ent}}$ and the positive branch $(I_+,g_+)$ for $r>r_R^{\mathrm{ent}}$. The rotational phase $Q(r,v_R^+(r))$ decreases from $+\infty$ to $1$ on $(R,r_R^{\mathrm{ent}}]$, and $v_R^+(r)=r^\alpha(1+o(1))$ as $r\to\infty$. Under the hypotheses of Theorem \ref{thm:positive-bowl-asymptotics}\textup{(i)}, the slope $v_R^+$ has the expansion of that theorem, and the height $u_R^+$ differs from that of a bowl-type solution in the sense of Definition \ref{def:bowl-type-solution} by $C_R^++o(1)$ as $r\to\infty$, for a constant $C_R^+$.
\end{proposition}

\begin{proof}
	We start the proof on the upper side of the local neck. By \eqref{eq:local-neck-phase-system}, we have $v_R^+(r_R(s))=\dot u_R(s)/\dot r_R(s)=\tan\theta_R(s)$ for small $s>0$. Recall that $\theta_R(0)=\frac{\pi}{2}$ and $\dot\theta_R(0)<0$ by Proposition \ref{prop:local-catenoidal-neck}. This means that $\theta_R(s)\to\frac{\pi}{2}$ from below as $s\to0^+$, and therefore $v_R^+(r)\to+\infty$ as $r\downarrow R$. Moreover, since $y_v=r^{-1}Q(1,v_R^+)$ and $Q(1,\cdot)$ is an increasing bijection of $\rr$ by Lemma \ref{lem:level-set-control}\textup{(i)}, we also have $y_v\to+\infty$. In particular, $y_v>1$ on some interval $(R,R+\eps)$.
	
	Next, we write the phase variables of \eqref{eq:local-neck-phase-system} in terms of the slope. With $\tan\theta=v$, they are $x=v'(1+v^2)^{-\frac{3}{2}}$, $y=r^{-1}v(1+v^2)^{-\frac{1}{2}}$ and $z=(1+v^2)^{-\frac{1}{2}}$, as in Section \ref{sec:preliminaries}. Then $(x_v,y_v)=z^{-\frac{1}{\alpha}}(x,y)$, and by the chart of Proposition \ref{prop:conical-sheet-maximal-chart} we obtain
	
		\begin{equation*}
		x_v = z^{-\frac{1}{\alpha}}g(y,z) = y_v\,\Psi_g^{-1}\left(y_v^{-\alpha}\right) = g(y_v,1) = h_+(y_v) \qquad \text{on }(R,R+\eps).
	\end{equation*}
	
	Therefore, $v_R^+$ solves \eqref{eq:selected-branch-graphical-equation} for $((1,\infty),h_+)$ on $(R,R+\eps)$. We note that this pair is a branch of the level-one equation. Indeed, we have $\widehat\gamma(h_+(y),y)=1$ by Proposition \ref{prop:conical-sheet-maximal-chart} and $h_+\in\mathcal C^1(1,\infty)$ by Lemma \ref{lem:upper-entry-sheet-from-zero-sector}. Now, we fix $r_0\in(R,R+\eps)$, and we let $[r_0,r_+)$ be the maximal branch-admissible interval of the solution of Proposition \ref{prop:branch-ivp} with datum $(r_0,v_R^+(r_0))$. Then, by uniqueness, this solution continues $v_R^+$, and $r_+$ does not depend on $r_0$. We keep the notation $v_R^+$ for the continued function.
	
	Subsequently, we follow the phase up to the cylindrical direction. Since $h_+<0$ by Lemma \ref{lem:upper-entry-sheet-from-zero-sector}, Lemma \ref{lem:continuation-phase-limits}\textup{(ii)} implies that $y_v$ decreases on $[r_0,r_+)$, and therefore on $(R,r_+)$. Moreover, it implies that $r_+<\infty$ and that $y_v(r)\to1$ as $r\to r_+$. In particular, we have $y_v([r_0,r_+))\subset[1,y_v(r_0)]$. Recall that $h_+$ extends continuously to $[1,\infty)$ with $h_+(1)=0$ by Lemma \ref{lem:upper-entry-sheet-from-zero-sector}. Then, by Proposition \ref{prop:branch-ivp}\textup{(ii)} with $K=[1,y_v(r_0)]$, we can extend $v_R^+$ to a $\mathcal C^1$ function on $[r_0,r_+]$ with $y_v(r_+)=1$, and the graphical equation implies that $v'(r_+)=0$. Consequently, the phase reaches $(0,1)$ at the finite radius $r_R^{\mathrm{ent}}=r_+$, and $y_v$ decreases from $+\infty$ to $1$ on $(R,r_R^{\mathrm{ent}}]$.
	
	Now, we continue the solution past the entry radius. By Lemma \ref{lem:positive-endpoint-regularity}, there is a $\mathcal C^1$ chart $\phi$ of the level-one set on an open interval $J\ni1$ contained in $(y_*,\infty)$, with $\phi(1)=0$. We note that $g_+$ and $h_+$ solve the level-one equation on the two sides of $y=1$, and that both tend to $0$ there. Then, by local uniqueness in the implicit function theorem, we have $\phi=g_+$ on $J\cap(y_*,1)$ and $\phi=h_+$ on $J\cap(1,\infty)$. Let $G=\phi$ on $J$. Then $(J,G)$ is a $\mathcal C^1$ branch of the level-one equation, and $v_R^+$ solves its graphical equation on a left neighborhood of $r_R^{\mathrm{ent}}$. Therefore, by Proposition \ref{prop:branch-ivp}\textup{(i)} at the datum $(r_R^{\mathrm{ent}},v_R^+(r_R^{\mathrm{ent}}))$, whose phase level $1$ lies in $J$, we can continue $v_R^+$ uniquely past $r_R^{\mathrm{ent}}$. Moreover, at $r_R^{\mathrm{ent}}$ the identity \eqref{eq:graphical-phase-evolution} implies that $ry_v'=\left(1+\alpha^{-1}v^2\right)G(1)-1=-1<0$. This means that we can find $r_1>r_R^{\mathrm{ent}}$ such that $y_v\in J\cap(y_*,1)$ on $(r_R^{\mathrm{ent}},r_1]$. On this interval we have $G=g_+$, and consequently $v_R^+$ is a branch-admissible solution of the positive graphical equation.
	
	Next, we consider the solution of the positive graphical equation with datum $(r_1,v_R^+(r_1))$, which is a positive end in the sense of Definition \ref{def:positive-end}. We note that it continues $v_R^+$. Moreover, since $\gamma$ is $1$-nondegenerate, its maximal interval is $[r_1,\infty)$ by Proposition \ref{prop:positive-end-maximal-interval}. Therefore, $v_R^+$ and $u_R^+$ extend to $(R,\infty)$. In addition, by Proposition \ref{prop:positive-end-at-infinity}\textup{(i)}, we have $v_R^+(r)=r^\alpha(1+o(1))$, and Theorem \ref{thm:positive-bowl-asymptotics}\textup{(i)} refines this expansion under its hypotheses.
	
	Finally, we compare the height with that of a bowl-type solution. Let $u_\gamma$ be the profile of a bowl-type solution, with $v_\gamma=u_\gamma'$. Then, by Lemma \ref{lem:bowl-is-positive-end}, the restriction of $v_\gamma$ to $[r_1,\infty)$ is also a positive end, and the bowl-type solution is entire by Proposition \ref{prop:positive-end-maximal-interval}. Recall that, under the hypotheses of Theorem \ref{thm:positive-bowl-asymptotics}\textup{(i)}, the expansion stated there depends only on the branch. Therefore, both slopes have this expansion, and we obtain $v_R^+-v_\gamma=o(r^{-4\alpha})$. Then, since $4\alpha>1$, this difference is integrable at infinity, and we have
	
		\begin{align*}
		u_R^+(r)-u_\gamma(r) &= \left(u_R^+-u_\gamma\right)\left(r_R^{\mathrm{ent}}\right) + \int_{r_R^{\mathrm{ent}}}^{\infty}\left(v_R^+-v_\gamma\right)dt - \int_{r}^{\infty}\left(v_R^+-v_\gamma\right)dt \\
		&= C_R^++o(1),
	\end{align*}
	
	where $C_R^+$ is the sum of the two constant terms, finalizing the proof.
\end{proof}

\begin{remark}\label{rem:height-constant-threshold}
	We note that the height constant $C_R^+$ is the only place in this section where we need an explicit lower bound on $\alpha$. In fact, this bound enters through the integrability of $v_R^+-v_\gamma$ at infinity. On the one hand, with the remainder $o(r^{-4\alpha})$ of Theorem \ref{thm:positive-bowl-asymptotics}\textup{(i)}, we need $4\alpha>1$. On the other hand, with the remainder $o(r^{-3\alpha})$, which is available when $\widetilde\gamma$ is only of class $\mathcal C^2$ near $(0,1)$, we need $3\alpha>1$. In particular, both conditions hold under the standing hypothesis $\alpha>\frac{1}{2}$ of this section.
\end{remark}

\subsection{Small necks}\label{subsec:small-necks}

Now, we study the entry radius for small necks. Recall that, between the neck and the entry radius, $y_v$ decreases from $+\infty$ to $1$, and that $a=\Psi_g^{-1}(y_v^{-\alpha})$ by Lemma \ref{lem:upper-entry-sheet-from-zero-sector}. Therefore, the angular variable $a=xy^{-1}$ of the upper side increases from the zero direction $-\rho_0$ to the cylindrical direction $0$. In addition, we will show that for small $R$ the entry occurs at the scale $R^{\rho_0/(1+\rho_0)}$, in the logarithmic sense of Proposition \ref{prop:entry-radius-estimate}. This means that the nonconvex collar of the upper side shrinks with the neck, while its radius exceeds the neck radius by the factor $R^{-1/(1+\rho_0)}$.

\begin{proposition}\label{prop:entry-radius-estimate}
	Assume the hypotheses of Lemma \ref{lem:upper-entry-sheet-from-zero-sector}, and let $R>0$. For $\delta\in(0,1)$, let $a_\delta\in(-\rho_0,0)$ denote the angular value satisfying $\Psi_g(a_\delta)=\delta$. Let $c_\delta=-a_\delta$ and $C_\delta=\bigl((4/3)^{\frac{1}{2\alpha}}\delta^{\frac{1}{\alpha}}\bigr)^{\frac{1}{1+c_\delta}}$. Then
	\begin{equation*}
		R^{\frac{\rho_0}{1+\rho_0}} < r_R^{\mathrm{ent}} \leq \delta^{-\frac{1}{\alpha}} \max\left\{ 2^{\frac{1}{c_\delta}}R,\; C_\delta R^{\frac{c_\delta}{1+c_\delta}} \right\}.
	\end{equation*}
	In particular, $r_R^{\mathrm{ent}}\to0$ and $\frac{\log r_R^{\mathrm{ent}}}{\log R}\to\frac{\rho_0}{1+\rho_0}$ as $R\to0^+$.
\end{proposition}

\begin{proof}
	Let $v=v_R^+$ and $y_v=Q(r,v)$ on $(R,r_R^{\mathrm{ent}}]$. Recall that, by Proposition \ref{prop:upper-side-global-continuation}, $y_v$ decreases from $+\infty$ to $1$ while the phase follows the branch $((1,\infty),h_+)$. Therefore, we have $x_v=h_+(y_v)=y_va(y_v)$, with $a(y)=\Psi_g^{-1}\left(y^{-\alpha}\right)>-\rho_0$ as in Lemma \ref{lem:upper-entry-sheet-from-zero-sector}. Moreover, we note that $v>0$, because $y_v\geq1$ and $Q(r,s)$ has the sign of $s$.
	
	Now, let $\sin\theta=v(1+v^2)^{-\frac{1}{2}}\in(0,1)$, where $\theta$ is the tangent angle of Section \ref{sec:preliminaries}. Then, by combining the graphical equation $v'=(1+v^2)^{1+\beta}h_+(y_v)$ with $y_v=r^{-1}v(1+v^2)^{-\beta}$, we obtain $(\log\sin\theta)'=v'v^{-1}(1+v^2)^{-1}=r^{-1}a(y_v)$. In addition, since $y_v\to+\infty$ forces $v\to+\infty$, we have $\sin\theta\to1$ as $r\downarrow R$. Consequently, by integrating, we get
	
		\begin{equation}\label{eq:sigma-angular-integral}
		\sin\theta(r) = \exp\left(\int_R^r a(y_v(t))\,\frac{dt}{t}\right), \qquad R<r\leq r_R^{\mathrm{ent}}.
	\end{equation}
	
	Next, we rewrite the phase in terms of $\sin\theta$. Since $\frac{1}{2}-\beta=\frac{1}{2\alpha}$ and $1+v^2=(1-\sin^2\theta)^{-1}$, the definition of $y_v$ becomes
	
		\begin{equation}\label{eq:phase-from-sigma}
		y_v = r^{-1}\sin\theta\,(1-\sin^2\theta)^{-\frac{1}{2\alpha}}.
	\end{equation}
	
	Firstly, we prove the lower bound. Since $a>-\rho_0$, from \eqref{eq:sigma-angular-integral} we have $\int_R^r a(y_v(t))\frac{dt}{t} > -\rho_0\log(r/R)$, and therefore $\sin\theta(r)>(R/r)^{\rho_0}$ on $(R,r_R^{\mathrm{ent}}]$. Moreover, at $r=r_R^{\mathrm{ent}}$ we have $y_v=1$, and \eqref{eq:phase-from-sigma} reduces to $\sin\theta=r_R^{\mathrm{ent}}(1-\sin^2\theta)^{\frac{1}{2\alpha}}<r_R^{\mathrm{ent}}$. Then, by combining these two inequalities, we obtain $(R/r_R^{\mathrm{ent}})^{\rho_0}<r_R^{\mathrm{ent}}$, that is, $\left(r_R^{\mathrm{ent}}\right)^{1+\rho_0}>R^{\rho_0}$, which is the lower bound.
	
	Next, we prove the upper bound. For this purpose, let $r_\delta\in(R,r_R^{\mathrm{ent}})$ be the unique radius with $y_v(r_\delta)=\delta^{-\frac{1}{\alpha}}$. We note that on $(R,r_\delta]$ we have $y_v^{-\alpha}\leq\delta$. Then, since $\Psi_g^{-1}$ is increasing, we obtain $a(y_v)\leq a_\delta=-c_\delta$, and \eqref{eq:sigma-angular-integral} implies that $\sin\theta(r)\leq(R/r)^{c_\delta}$ on this subinterval. Moreover, since $t\mapsto t(1-t^2)^{-\frac{1}{2\alpha}}$ is increasing on $(0,1)$, this bound and \eqref{eq:phase-from-sigma} at $r_\delta$ imply that
	
		\begin{equation*}
		\delta^{-\frac{1}{\alpha}} = y_v(r_\delta) \leq r_\delta^{-1} \left(\frac{R}{r_\delta}\right)^{c_\delta} \left(1-\left(\frac{R}{r_\delta}\right)^{2c_\delta}\right)^{-\frac{1}{2\alpha}}.
	\end{equation*}
	
	Now, if $r_\delta\geq2^{\frac{1}{c_\delta}}R$, the last factor is at most $(4/3)^{\frac{1}{2\alpha}}$. In this case the inequality becomes $r_\delta^{1+c_\delta}\leq(4/3)^{\frac{1}{2\alpha}}\delta^{\frac{1}{\alpha}}R^{c_\delta}$, that is, $r_\delta\leq C_\delta R^{c_\delta/(1+c_\delta)}$. Therefore, in all cases we have $r_\delta\leq\max\{2^{1/c_\delta}R,\,C_\delta R^{c_\delta/(1+c_\delta)}\}$.
	
	Subsequently, we apply Lemma \ref{lem:continuation-phase-limits}\textup{(ii)} with $r_0=r_\delta$ on the maximal branch-admissible interval $[r_\delta,r_R^{\mathrm{ent}})$, and we obtain $y_v(r)<\delta^{-\frac{1}{\alpha}}r_\delta\,r^{-1}$ for $r_\delta<r<r_R^{\mathrm{ent}}$. Then, since $y_v>1$ on this interval, we have $r_R^{\mathrm{ent}}\leq\delta^{-\frac{1}{\alpha}}r_\delta$. Consequently, by using the bound on $r_\delta$, we obtain the upper bound.
	
	Finally, we prove the two limits. From the upper bound, we have $r_R^{\mathrm{ent}}\to0$ as $R\to0^+$. Moreover, since $\log R<0$ for $R<1$, the lower bound implies that $\log r_R^{\mathrm{ent}}/\log R<\rho_0(1+\rho_0)^{-1}$. In addition, we note that in the upper bound the second term of the maximum dominates as $R\to0^+$, because its exponent is less than $1$. Therefore, we have $\log r_R^{\mathrm{ent}}\leq\log\left(\delta^{-\frac{1}{\alpha}}C_\delta\right)+c_\delta(1+c_\delta)^{-1}\log R$ for small $R$, and it follows that $\liminf_{R\to0^+}\log r_R^{\mathrm{ent}}/\log R\geq c_\delta(1+c_\delta)^{-1}$. Then, by taking $\delta\to0^+$, we have $c_\delta=-\Psi_g^{-1}(\delta)\to\rho_0$ by the continuity of $\Psi_g^{-1}$, concluding the proof of the limit.
\end{proof}

On the lower side, we obtain in Proposition \ref{prop:lower-side-collar-scale} the same scale for the maximal radius, or for the radius where the phase crosses a fixed level, and we bound the height of the side by the neck radius. For this purpose, in the outward parametrization of Lemma \ref{lem:lower-side-neck-phase} we regard $\overline\theta_R$ as a function of the radius along the lower graph $(r,u_R^-(r))$. We note that it is the tangent angle of that graph traversed outward. Therefore, we have $\tan\overline\theta_R=v_R^-$ with $\overline\theta_R\in(\frac{3\pi}{2},2\pi)$, and $-\sin\overline\theta_R=|v_R^-|\,(1+(v_R^-)^2)^{-\frac{1}{2}}\in(0,1)$. In addition, for $y\in I_-$ let $a(y)=y^{-1}g_-(y)$ be the angular variable along the lower branch. On the disjoint domain $I_-$, it plays the same role as the function $a(y)=\Psi_g^{-1}(y^{-\alpha})$ along the upper chart. Recall that $a<-\rho_0$ by \eqref{eq:lower-branch-positivity}, and that $a(y)\to-\rho_0$ as $y\to-\infty$ by \eqref{eq:reflected-lower-angular-relation}. Consequently, from the analogue of \eqref{eq:sigma-angular-integral} we obtain an upper bound on $-\sin\overline\theta_R$, opposite to the lower bound on the upper side.

\begin{proposition}\label{prop:lower-side-collar-scale}
	Let $R>0$, let $v=v_R^-$ on its maximal branch-admissible interval $(R,r_+)$ for $(I_-,g_-)$, and let $u_R^-$ denote the lower side height satisfying $(u_R^-)'=v$ and $u_R^-\to0$ at the neck. Let $\Phi_{\rho_0}(T)=\int_1^T\tau^{-\rho_0}\bigl(1-\tau^{-2\rho_0}\bigr)^{-1/2}d\tau$ for $T\geq1$.
	\begin{enumerate}[label=\textup{(\roman*)}]
		\item For $R<r<r_+$,
		\begin{equation}\label{eq:lower-sigma-angular-integral}
			-\sin\overline\theta_R(r) = \exp\left(\int_R^r a(y_v(t))\,\frac{dt}{t}\right) < \left(\frac{R}{r}\right)^{\rho_0},
		\end{equation}
		and $|u_R^-(r)|\leq R\,\Phi_{\rho_0}(r/R)$. In particular, $|u_R^-|\leq C_{\rho_0}R$ on $(R,r_+)$ when $\rho_0>1$, where $C_{\rho_0}=\Phi_{\rho_0}(\infty)<\infty$.
		\item If $(I_-,g_-)$ is axis-crossing, regularity-limited or has fractional zero output, then $r_+<\infty$ by Proposition \ref{prop:regular-lower-endpoint-dynamics}, and $r_+\to0$ with $\frac{\log r_+}{\log R}\to\frac{\rho_0}{1+\rho_0}$ as $R\to0^+$. If $(I_-,g_-)$ has finite-angle zero output, the same holds, for every $\delta>0$, for the radius where $y_v$ crosses $-\delta$. The same holds for a finite-angle continuation, along which the angular variable stays below $-\rho_0$ as in the proof of Lemma \ref{lem:finite-angle-continuation-angle}. At each point, either the phase lies in $-U_g$, where $\Psi_g(a)<0$, or $a\notin A_g$.
	\end{enumerate}
\end{proposition}

\begin{proof}
\textup{(i)} Firstly, we note that on $(R,r_+)$ the phase satisfies $y_v\in I_-\subset(-\infty,0)$. Therefore, we have $v<0$, because $Q(r,s)$ has the sign of $s$. Moreover, since $\tan\overline\theta_R=v$, we have $\overline\theta_R'=v'(1+v^2)^{-1}$ and $\cot\overline\theta_R=v^{-1}$. Then, by combining the lower graphical equation with $g_-(y_v)=y_va(y_v)$ and $y_v=r^{-1}v(1+v^2)^{-\beta}$, we obtain
	\begin{equation*}
		\frac{d}{dr}\log\left(-\sin\overline\theta_R\right) = \frac{\cos\overline\theta_R}{\sin\overline\theta_R}\,\overline\theta_R' = \frac{v'}{v(1+v^2)} = \frac{(1+v^2)^{\beta}g_-(y_v)}{v} = \frac{a(y_v)}{r}.
	\end{equation*}
Recall that $v\to-\infty$ and $y_v\to-\infty$ as $r\downarrow R$ by Lemma \ref{lem:lower-side-neck-phase}. In particular, we have $-\sin\overline\theta_R=\frac{|v|}{\sqrt{1+v^2}}\to1$ and $a(y_v)\to-\rho_0$. Then, by integrating from $R$, we obtain the equality in \eqref{eq:lower-sigma-angular-integral}, and the strict inequality follows from $a<-\rho_0$.

Next, we bound the height. By isolating $|v|$, we have $|v|=(-\sin\overline\theta_R)(1-\sin^2\overline\theta_R)^{-\frac{1}{2}}$. Then, since the function $\xi\mapsto\xi(1-\xi^2)^{-\frac{1}{2}}$ strictly increases on $(0,1)$ and $(R/r)^{\rho_0}<1$ for $r>R$, the strict inequality in \eqref{eq:lower-sigma-angular-integral} implies that $|v(r)|<(R/r)^{\rho_0}\bigl(1-(R/r)^{2\rho_0}\bigr)^{-\frac{1}{2}}$ for $R<r<r_+$. We note that the right-hand side is of order $(r-R)^{-\frac{1}{2}}$ as $r\downarrow R$, and therefore it is integrable near $R$. Consequently, by substituting $t=R\tau$ into $|u_R^-(r)|\leq\int_R^r|v(t)|\,dt$, we obtain the stated bound on the height. Finally, we note that the integrand defining $\Phi_{\rho_0}$ is integrable at $\tau=1$, where $1-\tau^{-2\rho_0}=2\rho_0(\tau-1)(1+o(1))$, and that it decays as $\tau^{-\rho_0}(1+o(1))$ at infinity. Therefore, $\Phi_{\rho_0}$ is bounded if and only if $\rho_0>1$.

\textup{(ii)} Firstly, we prove the lower bound on $r_+$. We note that the function $a$ is decreasing on $I_-$, since in \eqref{eq:reflected-lower-angular-relation} the map $y\mapsto-(-y)^{-\alpha}$ is decreasing and $\Psi_g^{-1}$ is increasing. Now, we fix $M>|y_{\mathrm{reg}}|$, so that $-M\in I_-$, and we set $A_M=-a(-M)$, which exceeds $\rho_0$ by \eqref{eq:lower-branch-positivity}. Then, by the monotonicity of $a$, we have $a\geq a(-M)=-A_M$ on $(-\infty,-M]$. Next, let $r_M\in(R,r_+)$ be the radius with $y_v(r_M)=-M$. This radius exists because $y_v$ increases strictly from $-\infty$ (by \eqref{eq:lower-phase-at-neck} and Proposition \ref{prop:graphical-phase-evolution}\textup{(ii)}) toward $y_{\mathrm{reg}}$ (by Lemma \ref{lem:continuation-phase-limits}). We note that on $(R,r_M]$ the phase satisfies $y_v\leq-M$, and therefore $a(y_v)\geq-A_M$ there. Then, by the equality in \eqref{eq:lower-sigma-angular-integral}, we have $-\sin\overline\theta_R(r)\geq(R/r)^{A_M}$. Moreover, the lower-side analogue of \eqref{eq:phase-from-sigma} reads $|y_v|=r^{-1}(-\sin\overline\theta_R)(1-\sin^2\overline\theta_R)^{-\frac{1}{2\alpha}}$. Since the last factor is at least $1$, we have $|y_v|\geq-r^{-1}\sin\overline\theta_R$. Then, by evaluating at $r_M$, we obtain $M\geq r_M^{-1}(R/r_M)^{A_M}$, that is, $r_M^{1+A_M}\geq M^{-1}R^{A_M}$, and therefore $r_M\geq M^{-\frac{1}{1+A_M}}R^{\frac{A_M}{1+A_M}}$. Now, for $R<1$, where $\log R<0$, taking logarithms and dividing by $\log R$ reverses the inequalities. Then, since $r_+>r_M$, we have
	\begin{equation*}
		\frac{\log r_+}{\log R} < \frac{\log r_M}{\log R} \leq \frac{A_M}{1+A_M} - \frac{\log M}{(1+A_M)\log R}, \qquad \text{for every }M>|y_{\mathrm{reg}}|.
	\end{equation*}
Since the last term tends to $0$ as $R\to0^+$, we have $\limsup_{R\to0^+}\log r_+/\log R\leq A_M(1+A_M)^{-1}$. Moreover, since $a(y)\to-\rho_0$ as $y\to-\infty$ by \eqref{eq:reflected-lower-angular-relation}, we have $A_M\to\rho_0$ as $M\to\infty$. Consequently, the limit superior is at most $\rho_0(1+\rho_0)^{-1}$.

Next, we prove the upper bound. For this purpose, we fix $\delta>0$, with $\delta<|y_{\mathrm{reg}}|$ if $y_{\mathrm{reg}}<0$, and so small that $g_-(y)\geq\frac{x_0}{2}$ on $(-\delta,0)$ in the axis-crossing case and $g_-(y)\geq\frac{d}{2}(-y)^\sigma$ on $(-\delta,0)$ in the fractional case, where $x_0$, $d$ and $\sigma$ are the constants of Definition \ref{def:lower-continuation-outcomes}. We note that for $2^{1/\rho_0}R\leq r<r_+$ the inequality in \eqref{eq:lower-sigma-angular-integral} implies that $-\sin\overline\theta_R<(R/r)^{\rho_0}\leq\frac{1}{2}$. Therefore, the factor $(1-\sin^2\overline\theta_R)^{-\frac{1}{2\alpha}}$ in the expression of $|y_v|$ above is at most $(4/3)^{\frac{1}{2\alpha}}$, and we have $|y_v(r)|<(4/3)^{\frac{1}{2\alpha}}r^{-1}(R/r)^{\rho_0}$. Moreover, the right-hand side is less than $\delta$ when $r^{1+\rho_0}>(4/3)^{\frac{1}{2\alpha}}\delta^{-1}R^{\rho_0}$. Consequently, we have $|y_v|<\delta$ for every $r>r_\delta$ with $r<r_+$, where
	\begin{equation*}
		r_\delta = \max\left\{ 2^{\frac{1}{\rho_0}}R,\; \bigl((4/3)^{\frac{1}{2\alpha}}\delta^{-1}R^{\rho_0}\bigr)^{\frac{1}{1+\rho_0}} \right\}.
	\end{equation*}
Now, we bound $r_+$ in terms of $r_\delta$ in each case.
\begin{itemize}
	\item \emph{Regularity-limited case.} Let us assume by contradiction that $r_+>r_\delta$. Then $|y_v|<\delta<|y_{\mathrm{reg}}|$ on $(r_\delta,r_+)$, which contradicts the branch admissibility $y_v<y_{\mathrm{reg}}$. Therefore, we have $r_+\leq r_\delta$, and we set $K_\delta=1$.
	\item \emph{Axis-crossing case.} If $r_+>r_\delta$, we set $w=-y_v$ on $(r_\delta,r_+)$, where $0<w<\delta$. Then, by the phase identity \eqref{eq:graphical-phase-evolution} and $1+\alpha^{-1}v^2\geq1$, we have $rw'=-\left(1+\alpha^{-1}v^2\right)g_-(-w)-w\leq-g_-(-w)\leq-\frac{x_0}{2}$. Then, by integrating from $r_\delta$, where $w(r_\delta)\leq\delta$, we obtain $w(r)\leq\delta-\frac{x_0}{2}\log(r/r_\delta)$ on $(r_\delta,r_+)$. Moreover, since $w>0$ there, we have $\log(r/r_\delta)<2\delta/x_0$ for every $r<r_+$. Consequently, $r_+\leq K_\delta r_\delta$ with $K_\delta=\exp(2\delta/x_0)$.
	\item \emph{Fractional case.} With the same $w$, by the phase identity we have $rw'\leq-g_-(-w)\leq-\frac{d}{2}w^\sigma$. Then, since $w>0$, we obtain $\left(w^{1-\sigma}\right)'=(1-\sigma)w^{-\sigma}w'\leq-\frac{(1-\sigma)d}{2r}$. Therefore, by integrating from $r_\delta$, we have $w(r)^{1-\sigma}\leq\delta^{1-\sigma}-\frac{(1-\sigma)d}{2}\log(r/r_\delta)$, and the positivity of $w$ implies that $r_+\leq K_\delta r_\delta$ with $K_\delta=\exp\bigl(2\delta^{1-\sigma}/((1-\sigma)d)\bigr)$.
\end{itemize}
In all three cases $r_+\leq K_\delta r_\delta$, where $K_\delta$ does not depend on $R$. We note that, as $R\to0^+$, the second term in the definition of $r_\delta$ dominates the first, because its exponent $\rho_0/(1+\rho_0)$ is less than $1$. Therefore, $r_\delta$ is a multiple of $R^{\rho_0/(1+\rho_0)}$ for small $R$, with a factor that depends on $\delta$ but not on $R$. Consequently, we have $r_+\to0$, and $\log r_+\leq c+\frac{\rho_0}{1+\rho_0}\log R$ for small $R$, with a constant $c$ independent of $R$. Then, by dividing by $\log R<0$, we obtain
	\begin{equation*}
		\liminf_{R\to0^+}\frac{\log r_+}{\log R} \geq \frac{\rho_0}{1+\rho_0}.
	\end{equation*}
Together with the limit superior above, this proves the limit.

Finally, we consider the finite-angle case. Recall that $y_v$ increases to $0^-$ by Lemma \ref{lem:continuation-phase-limits}, and therefore the radius at which $y_v$ reaches $-\delta$ exists. Moreover, it is larger than $r_M$ for every $M>\delta$, since $y_v$ is increasing, and at most $r_\delta$, since $|y_v|<\delta$ beyond $r_\delta$. Consequently, it satisfies the same upper and lower bounds, finalizing the proof.
\end{proof}

\begin{figure}[tb]
	\centering
	\captionsetup[sub]{labelfont=up}
	\begin{subfigure}[t]{0.22\textwidth}
		\centering
		\includegraphics[width=\textwidth,trim=13.2bp 10.9bp 184.56bp 12.7bp,clip]{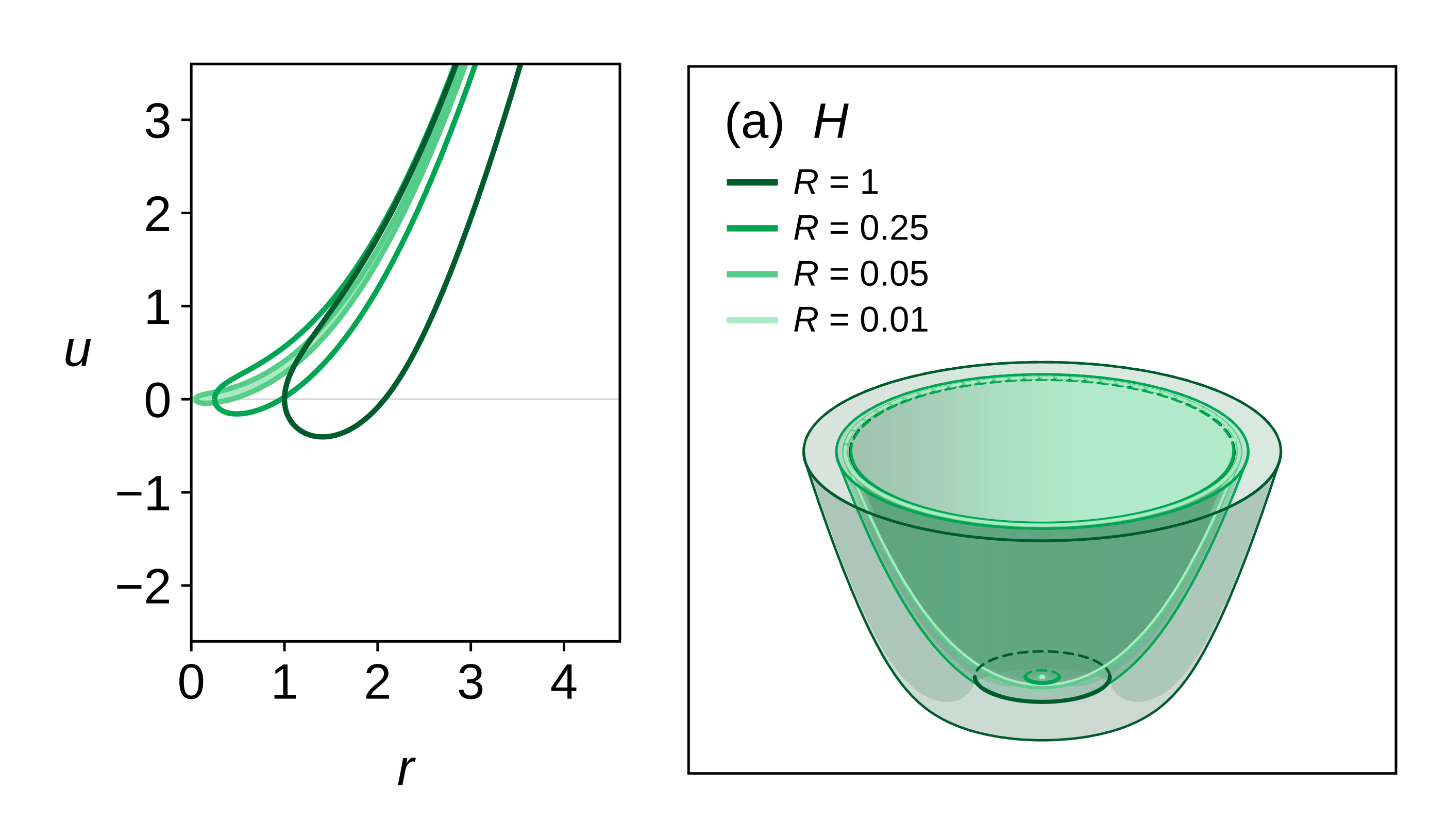}
		\caption{\scriptsize $\frac{H}{n-1}$}
	\end{subfigure}\hfill
	\begin{subfigure}[t]{0.22\textwidth}
		\centering
		\includegraphics[width=\textwidth,trim=13.2bp 10.9bp 184.56bp 12.7bp,clip]{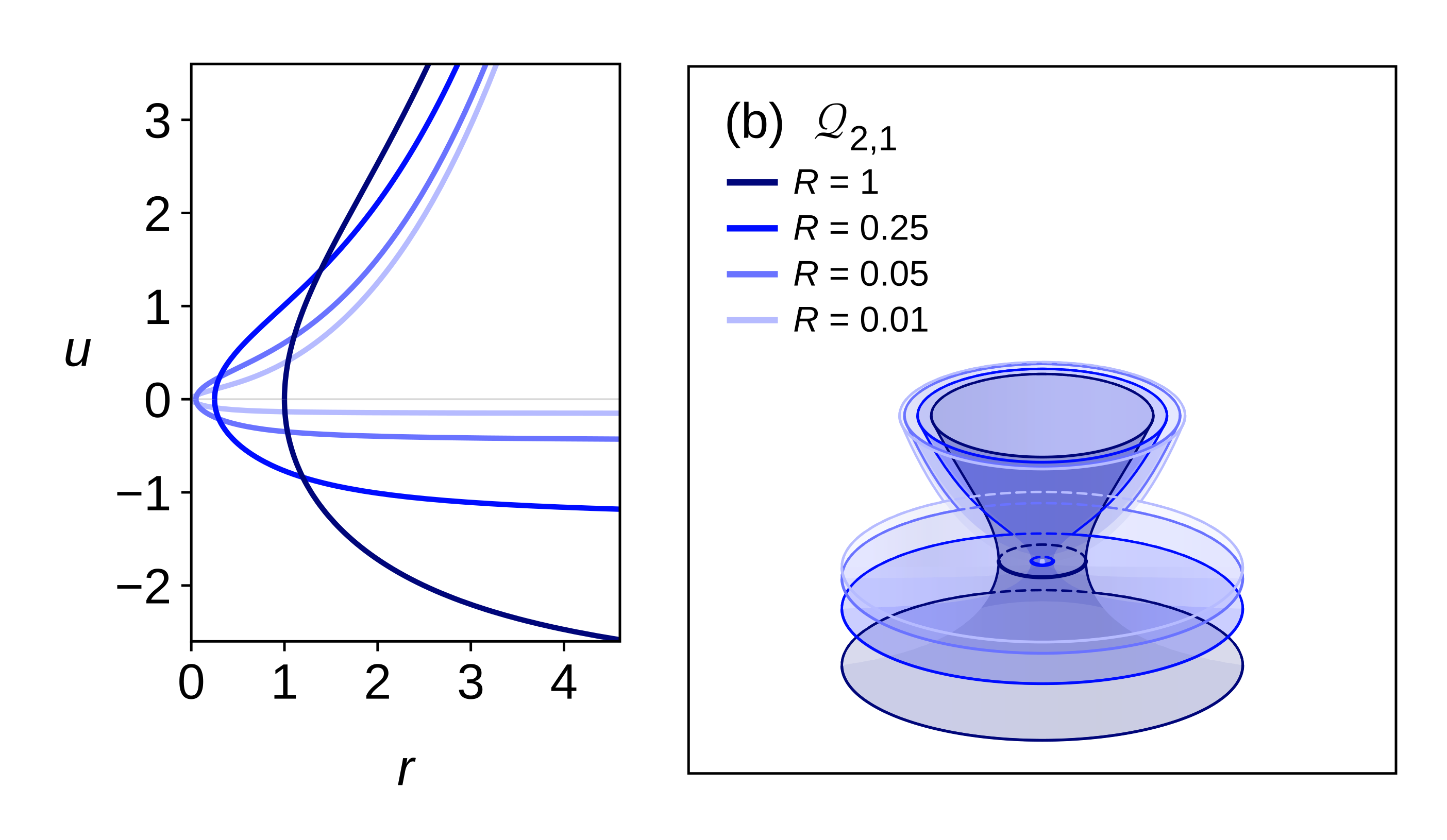}
		\caption{\scriptsize $\mathcal Q_{2,1}$}
	\end{subfigure}\hfill
	\begin{subfigure}[t]{0.22\textwidth}
		\centering
		\includegraphics[width=\textwidth,trim=13.2bp 10.9bp 184.56bp 12.7bp,clip]{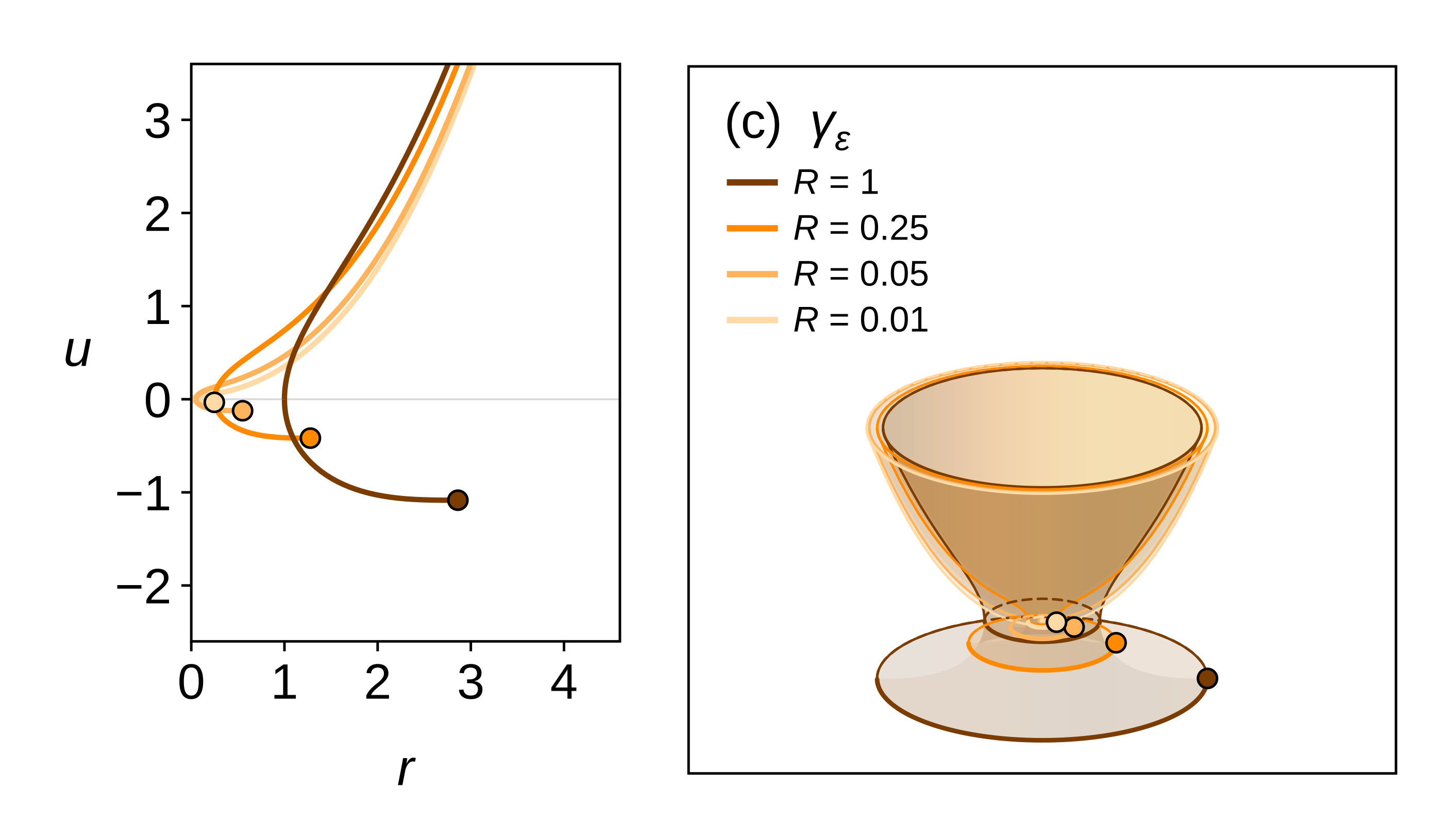}
		\caption{\scriptsize $\gamma_\eps$}
	\end{subfigure}\hfill
	\begin{subfigure}[t]{0.31\textwidth}
		\centering
		\includegraphics[width=\textwidth,trim=8.4bp 9.6bp 165.96bp 11.2bp,clip]{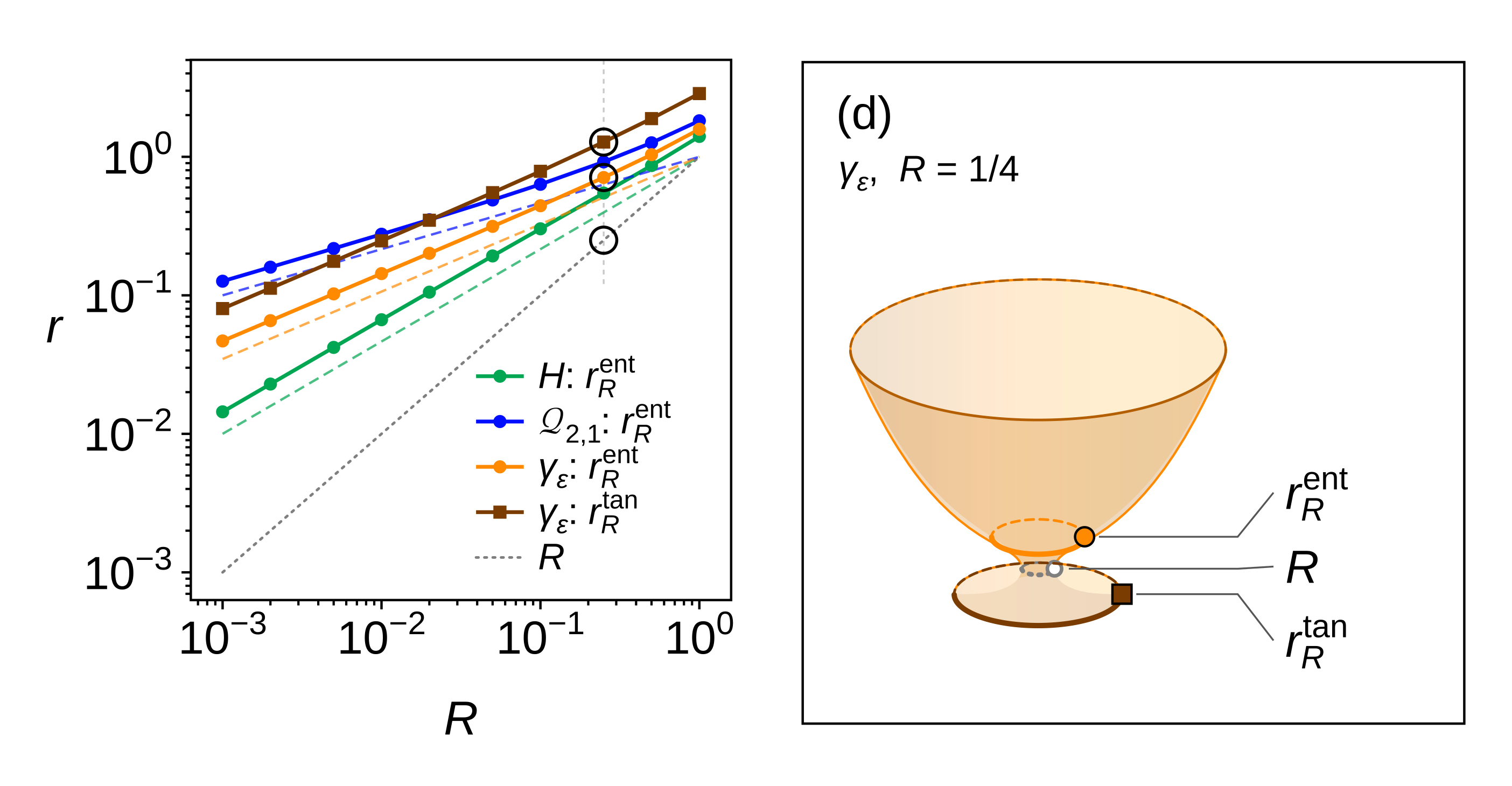}
		\caption{\scriptsize Radii against $R$}
	\end{subfigure}
	\caption{\scriptsize Catenoidal profiles for $R\in\{1,\frac{1}{4},\frac1{20},\frac1{100}\}$ and $n=3$, computed from the graphical equations of Section \ref{sec:barrier-method}. In \textup{(a)}--\textup{(c)}, the shade lightens as $R$ decreases. \textup{(a)} $H$, normalized as $\frac{H}{n-1}$, case \textup{(i)}: both positive ends converge to the bowl as $R$ decreases. \textup{(b)} $\mathcal Q_{2,1}$, case \textup{(ii)}: the lower end flattens toward the horizontal plane through the neck. \textup{(c)} $\gamma_\eps$ with $\eps=1$, case \textup{(iii)}: the maximal piece ends at the horizontal-tangency circle (dot), which collapses to the axis. \textup{(d)} The entry radius $r_R^{\mathrm{ent}}$ of the three families and the tangency radius $r_R^{\mathrm{tan}}$ of $\gamma_\eps$ against $R$, with the limits $R^{\rho_0/(1+\rho_0)}$ (dashed) and $R$ (dotted). The circles mark the values of $r_R^{\mathrm{tan}}$, $r_R^{\mathrm{ent}}$ and $R$ for $\gamma_\eps$ with $R=\frac{1}{4}$. Images courtesy of Ignacio McManus.}
	\label{fig:small-neck-limit}
\end{figure}

We note that the height $|u_R^-|$ of the lower side at the radius of Proposition \ref{prop:lower-side-collar-scale}\textup{(ii)}, which is of order $R^{\frac{\rho_0}{1+\rho_0}+o(1)}$, tends to zero with $R$. In fact, it is bounded by $C_{\rho_0}R$ when $\rho_0>1$ and by $R\,\Phi_{\rho_0}(r/R)$ otherwise. Moreover, since the integrand of $\Phi_{\rho_0}$ is $\tau^{-\rho_0}(1+o(1))$ as $\tau\to\infty$, we have $\Phi_{\rho_0}(T)=O\bigl(T^{1-\rho_0}\bigr)$ for $\rho_0<1$ and $\Phi_1(T)=O(\log T)$. Then, since at $r=R^{\frac{\rho_0}{1+\rho_0}+o(1)}$ the ratio is $T=r/R=R^{-\frac{1}{1+\rho_0}+o(1)}$, we obtain
\begin{equation}\label{eq:lower-height-small-neck}
	|u_R^-(r)| =
	\begin{cases}
		O(R), & \rho_0>1,\\
		O\bigl(R\log R^{-1}\bigr), & \rho_0=1,\\
		O\bigl(R^{\frac{2\rho_0}{1+\rho_0}+o(1)}\bigr), & \rho_0<1.
	\end{cases}
\end{equation}
Consequently, by combining \eqref{eq:lower-height-small-neck} with Proposition \ref{prop:entry-radius-estimate}, we conclude that the nonconvex collar of the upper side and the corresponding part of the lower side collapse toward the axis at the same scale as $R\to0^+$. Figure \ref{fig:small-neck-limit} shows the collapse for cases \textup{(i)}--\textup{(iii)}.

\subsection{The lower continuations}\label{subsec:lower-continuations}
Recall that, by Lemma \ref{lem:lower-side-neck-phase}, the reoriented lower side solves the lower graphical equation and issues from the infinite end of $(I_-,g_-)$. Moreover, by Proposition \ref{prop:graphical-phase-evolution}\textup{(ii)}, its phase increases toward the right endpoint $y_{\mathrm{reg}}$. We described this endpoint in Proposition \ref{prop:regular-lower-endpoint-dynamics} for the axis-crossing and regularity-limited types and for fractional zero output, all of which reach it at a finite radius.

Next, we study the lower side when the branch has finite-angle zero output or admits a finite-angle continuation. To do so, we denote by $\overline g$ the branch $g_-$ when it has finite-angle zero output in the sense of Definition \ref{def:lower-continuation-outcomes}, and the extension $g_-^{\mathrm{ext}}$ of Definition \ref{def:finite-angle-continuation} when it admits a finite-angle continuation. In both cases $\overline g$ is nonincreasing on an interval $(-\delta,0)$, by \eqref{eq:lower-branch-monotonicity} and by Definition \ref{def:finite-angle-continuation}, respectively. We use this property in the proof of Proposition \ref{prop:finite-angle-lower-completeness}, where we compare the slope with the power barriers of Proposition \ref{prop:power-barrier-flat-component} by means of the comparison principle of Proposition \ref{prop:comparison-principle-ode}. In addition, since $g_-^{\mathrm{ext}}$ is of class $\mathcal C^1$ across $y_{\mathrm{reg}}$, in the second case a solution issuing from the neck follows $\overline g$ through $y_{\mathrm{reg}}$ by Proposition \ref{prop:branch-ivp}, without a matching step.

\begin{proposition}\label{prop:finite-angle-lower-completeness}
	Assume that $(I_-,g_-)$ has finite-angle zero output with angular limit $\ell$, or that it admits a finite-angle continuation with angular limit $\ell$, and let $v$ solve the lower graphical equation on its maximal branch-admissible interval $[r_1,r_+)$ for the branch $\overline g$ on $(-\infty,0)$. Then the following hold.
	\begin{enumerate}[label=\textup{(\roman*)}]
		\item $r_+=\infty$, and $y_v(r)\to0^-$ and $v(r)\to0^-$ as $r\to\infty$.
		\item For all $\chi_1<\ell<\chi_2<0$ there are constants $c_1,c_2>0$ and a radius $\bar r\geq r_1$ with
		\begin{equation*}
			c_1r^{\chi_1} \leq -v(r) \leq c_2r^{\chi_2}, \qquad r\geq\bar r,
		\end{equation*}
		so that $\log(-v(r))/\log r\to\ell$.
		\item If the branch is flat, with $\overline g(y)=\ell y+O(|y|^{1+\mu})$ as $y\to0^-$ for some $\mu>0$, then there is a constant $c_*>0$, depending on the solution, such that
		\begin{equation*}
			v(r) = -c_*r^{\ell}(1+o(1)) \qquad \text{as }r\to\infty.
		\end{equation*}
	\end{enumerate}
\end{proposition}

\begin{proof}
Firstly, we note that item \textup{(i)} is Proposition \ref{prop:zero-output-radius-criterion} for a regular branch with finite-angle zero output and Corollary \ref{cor:radius-criterion-on-continuations} for a finite-angle continuation. In fact, from both results we have $r_+=\infty$ and $v\to0^-$, and $y_v\to0^-$ follows from Lemma \ref{lem:continuation-phase-limits}\textup{(i)}.

Next, we prove \textup{(ii)}. For this purpose, we fix $\chi_1<\ell<\chi_2<0$, we set $\eps=2^{-1}\min\{\ell-\chi_1,\,\chi_2-\ell\}$, and we choose $\delta>0$ so small that $\overline g$ is nonincreasing on $(-\delta,0)$ and
	\begin{equation*}
		\left|y^{-1}\overline g(y)-\ell\right| < \eps \qquad \text{for }-\delta<y<0.
	\end{equation*}
Then $\chi_1<\ell-\eps$ and $\ell+\eps<\chi_2<0$, as Proposition \ref{prop:power-barrier-flat-component}\textup{(ii)} requires. Now, for $i=1,2$ and $c_i>0$ to be chosen, let $\overline w_i(r)=-c_ir^{\chi_i}$, and let $t_i$ be the corresponding function $t$ of that proposition, so that $Q(r,\overline w_i(r))=-t_i(r)\in(-\infty,0)$. Then, by that proposition, $\overline w_2$ is a strict subsolution wherever $t_2(r)<\delta$, and $\overline w_1$ is a strict supersolution wherever $t_1(r)<\delta$ and $\overline w_1(r)^2\leq\eps_0$, where $\eps_0=\frac{\chi_1}{\ell-\eps}-1>0$.

Subsequently, we choose the constants. By \textup{(i)}, we can choose $\bar r\geq r_1$ with $v(\bar r)^2\leq\eps_0$ and $|y_v(\bar r)|<\delta$, and we let
	\begin{equation*}
		c_1 = -v(\bar r)\,\bar r^{-\chi_1}, \qquad c_2 = -v(\bar r)\,\bar r^{-\chi_2}.
	\end{equation*}
We note that both barriers pass through $v(\bar r)$ at $\bar r$, and that their levels there are $t_i(\bar r)=|v(\bar r)|\,\bar r^{-1}(1+v(\bar r)^2)^{-\beta}=|y_v(\bar r)|<\delta$. Moreover, each $t_i$ decreases by Proposition \ref{prop:power-barrier-flat-component}\textup{(i)}, and therefore $t_i<\delta$ on $[\bar r,\infty)$. In addition, since $\chi_1<0$, we have $|\overline w_1(r)|=c_1r^{\chi_1}\leq c_1\bar r^{\chi_1}=|v(\bar r)|$ for $r\geq\bar r$, which implies that $\overline w_1^2\leq\eps_0$ on the whole tail. Consequently, $\overline w_2$ is a strict subsolution and $\overline w_1$ a strict supersolution on $[\bar r,\infty)$.

Now, we compare $v$ with the barriers. We note that on $[\bar r,\infty)$ the phases of $v$ and of both barriers lie in $(-\delta,0)$, since the levels $-t_i$ increase from $-|y_v(\bar r)|>-\delta$ toward $0^-$ and $y_v$ increases to $0^-$ by \textup{(i)}. Then, by Proposition \ref{prop:comparison-principle-ode}, applied to the branch $((-\delta,0),\overline g)$ with $\overline w_i(\bar r)=v(\bar r)$, we have
	\begin{equation*}
		-c_2r^{\chi_2} \leq v(r) \leq -c_1r^{\chi_1}, \qquad r\geq\bar r.
	\end{equation*}
This is the bound in \textup{(ii)}. Then, by taking logarithms, dividing by $\log r>0$ and letting $\chi_1,\chi_2\to\ell$, we obtain $\log(-v(r))/\log r\to\ell$.

Finally, we prove \textup{(iii)}. For this purpose, we write the lower graphical equation as
	\begin{equation*}
		\frac{d}{dr}\log(-v) = v'v^{-1} = (1+v^2)^{\beta+1}\,v^{-1}\overline g(y_v).
	\end{equation*}
Then, by the flat expansion and the definition of $y_v$, we have
	\begin{equation*}
		v^{-1}\overline g(y_v) = \ell\,r^{-1}(1+v^2)^{-\beta} + O\!\left(|y_v|^{1+\mu}|v|^{-1}\right),
	\end{equation*}
and, since $|y_v|^{1+\mu}|v|^{-1}=O\!\left(|v|^{\mu}r^{-1-\mu}\right)$, we obtain
	\begin{equation*}
		\frac{d}{dr}\log(-v) = \frac{\ell}{r} + O\!\left(\frac{v^2}{r}\right) + O\!\left(|v|^{\mu}r^{-1-\mu}\right).
	\end{equation*}
Recall that $|v|\leq c_2r^{\chi_2}$ with $\chi_2<0$ by \textup{(ii)}. Therefore, the error terms are $O(r^{2\chi_2-1})$ and $O(r^{\mu\chi_2-1-\mu})$. We note that both exponents are less than $-1$, and consequently the errors are integrable on $[\bar r,\infty)$. Then, by integrating, we obtain $\log(-v)-\ell\log r\to\log c_*$ for some constant $c_*>0$, finalizing the proof.
\end{proof}

\subsection{Embeddedness}\label{subsec:catenoidal-embeddedness}

\begin{lemma}\label{lem:embeddedness-from-slope-comparison}
	Let $\mathtt c_R$ be a catenoidal profile (Definition \ref{def:catenoidal-profile}) whose two sides are vertical graphs $u_R^+$ and $u_R^-$ over the radial variable. Each graph is defined on a maximal interval with left endpoint $R$ and satisfies $\lim_{r\downarrow R}u_R^\pm(r)=0$, and we write $v_R^\pm=(u_R^\pm)'$. Let $(R,r_+)$ be the intersection of the two intervals of definition, and assume that at every $r\in(R,r_+)$ one of the following conditions holds.
	\begin{enumerate}[label=\textup{(\alph*)}]
		\item \label{emb:phase-separated}
		$Q(r,v_R^-(r))<Q(r,v_R^+(r))$.
		\item \label{emb:common-branch}
		There exists $r_0\in(R,r)$ such that, on $[r_0,r]$, both slopes are branch-admissible solutions to the graphical equation of one common branch $(I,g)$, where $g\in\mathcal C^1(I)$ is nonincreasing, and $v_R^-(r_0)<v_R^+(r_0)$.
	\end{enumerate}
	Then $v_R^-<v_R^+$ and $u_R^-<u_R^+$ on $(R,r_+)$. Consequently, the generating curve is embedded, and so is $W_R^{\max}$.
\end{lemma}

\begin{proof}
Firstly, we prove that $v_R^-<v_R^+$ on $(R,r_+)$. In case \ref{emb:phase-separated}, by the strict monotonicity of $Q(r,\cdot)$ we have $v_R^-(r)<v_R^+(r)$. In case \ref{emb:common-branch}, both slopes solve the graphical equation of $(I,g)$ on $[r_0,r]$, and therefore their residuals vanish. Then, by Proposition \ref{prop:comparison-principle-ode}, we have $v_R^-\leq v_R^+$ on $[r_0,r]$. Let us assume by contradiction that $v_R^-(r_1)=v_R^+(r_1)$ for some $r_1\in(r_0,r]$. Then both slopes solve the same initial value problem at $r_1$, and by Proposition \ref{prop:branch-ivp}\textup{(i)} they coincide on $[r_0,r_1]$, which contradicts $v_R^-(r_0)<v_R^+(r_0)$. Therefore, we have $v_R^-<v_R^+$ on $(R,r_+)$. Moreover, since both heights tend to $0$ at $R$, by integrating we obtain
	\begin{equation*}
		u_R^+(r)-u_R^-(r) = \int_R^r\left(v_R^+-v_R^-\right)(t)\,dt > 0
	\end{equation*}
for all $r\in(R,r_+)$.

Next, we prove that the curve is embedded. We note that each side is a vertical graph, and therefore an injective curve. Moreover, the two sides have different heights at every common radius, and outside the common domain only the side defined on the longer interval is present. In addition, the neck is the only point of the curve at radius $R$, since $r_R(s)>R$ for $s\neq0$ by Proposition \ref{prop:local-catenoidal-neck}. Consequently, the generating curve is injective. Finally, its inverse is continuous because each side is a graph over the radial variable, and it follows that the curve is embedded. We note that rotation about the vertical axis preserves this property because $r>0$ along the curve, concluding the proof.
\end{proof}

\subsection{The classification theorem}\label{subsec:catenoidal-classification}

\begin{theorem}\label{thm:rotational-translators}
	Let $\gamma$ be a normalized $1$-nondegenerate $\alpha$-positive curvature function with $\alpha>\frac{1}{2}$, endowed with selected angular neck data $(A_g,\rho_0)$ that satisfy conditions \ref{cat:all-radii} to \ref{cat:lower-branch} from the opening of this section, and let $R>0$. Then the local neck $\mathtt c_R$ of Proposition \ref{prop:local-catenoidal-neck} exists and is unique among solutions of \eqref{eq:local-neck-phase-system}, that is, among profiles whose phase lies on the selected sheet. Its upper side reaches the positive branch at a finite entry radius and continues as a complete positive end, satisfying $v_R^+=r^\alpha(1+o(1))$. If $\widetilde\gamma$ is of class $\mathcal C^3$ near $(0,1)$, then
	\begin{equation*}
		v_R^+(r)=r^\alpha-a_\gamma r^{-\alpha}+b_\gamma r^{-3\alpha}+o(r^{-4\alpha}) \qquad \text{as }r\to\infty,
	\end{equation*}
	with the coefficients $a_\gamma$ and $b_\gamma$ of Theorem \ref{thm:positive-bowl-asymptotics}\textup{(i)}, which depend on $\gamma$ but not on $R$. The height $u_R^+$ differs from that of the bowl-type solution by $C_R^++o(1)$, with the constant $C_R^+$ of Proposition \ref{prop:upper-side-global-continuation}. If $\widetilde\gamma$ is only of class $\mathcal C^2$ near $(0,1)$, the same expansion holds with a remainder $o(r^{-3\alpha})$. Across the neck circle, $\mathtt c_R$ is of class $\mathcal C^{2,\varkappa}$ when $R$ is $\varkappa$-regular at the neck, of class $\mathcal C^3$ when $g$ is of class $\mathcal C^1$ near $(R^{-1},0)$, and smooth when $g$ is smooth near that point.
	
	The continuation of the reoriented lower side is determined by the type of $(I_-,g_-)$ in Definition \ref{def:lower-continuation-outcomes}, combined with the matching and continuation conditions of Definitions \ref{def:regular-positive-matching} to \ref{def:singular-positive-matching}, structured as follows.
	\begin{enumerate}[label=\textup{(\roman*)}]
		\item \emph{Axis-crossing with regular positive matching.} In this case, $W_R$ is a complete embedded two-ended catenoidal translator, whose lower side is a second complete positive end. More precisely, the lower phase traverses $I_-$, passes through the phase points $(x_0,0)$ and $(y_*,y_*)$ at finite radii, and then follows the positive branch. Moreover, if $\widetilde\gamma$ is of class $\mathcal C^3$ near $(0,1)$, then $v_R^-$ has the same expansion as $v_R^+$ above, with the same coefficients $a_\gamma$ and $b_\gamma$, and its height differs from that of the bowl-type solution by $C_R^-+o(1)$ for a constant $C_R^-$. In particular, the two ends are asymptotic to vertical translates of the bowl, in general by different heights. In addition, $W_R$ is of class $\mathcal C^3$ away from the neck circle, and smooth away from it wherever $\widehat\gamma$ is smooth along the traversed branches.
		
		\item \emph{Finite-angle zero output with angular limit $\ell$, or a regularity-limited branch with a finite-angle continuation of angular limit $\ell$ \textup{(}Definition \ref{def:finite-angle-continuation}\textup{)}.} The hypersurface $W_R$ is a complete embedded catenoidal translator with one positive end, and its lower side is a graph over $(R,\infty)$ with a negative slope. Its slope satisfies $v_R^-(r)=-r^{\ell+o(1)}$, which refines to $v_R^-(r)=-c_*r^{\ell}(1+o(1))$ when the underlying branch is flat. Moreover, the lower height $u_R^-$ is bounded when $\ell<-1$, which holds whenever $\rho_0\geq1$, and it tends to $-\infty$ when $\ell>-1$. When the branch is flat, we have $u_R^-(r)=-\frac{c_*}{\ell+1}r^{\ell+1}(1+o(1))$ if $\ell>-1$, $u_R^-(r)=-c_*\log r\,(1+o(1))$ if $\ell=-1$, and $u_R^-(r)\to u_\infty^-$ for a finite limit $u_\infty^-<0$ if $\ell<-1$. For a regularity-limited branch with a finite-angle continuation, $W_R$ is of class $\mathcal C^3$ across the circle where the phase passes $y_{\mathrm{reg}}$, and smooth there provided $\widehat\gamma$ is smooth near $(x_{\mathrm{reg}},y_{\mathrm{reg}})$.
		
		\item \emph{Fractional zero output with exponent $\sigma\in(0,1)$.} The maximal continuation of the lower side terminates at a finite radius $r_R^{\mathrm{tan}}$, at which $v_R^-\to0$ and the phase of $v_R^-$ approaches the origin $(0,0)$. Its rotational phase decays at the rate \eqref{eq:fractional-phase-endpoint-asymptotics}. The hypersurface $W_R^{\max}$ is embedded and bounded by the horizontal-tangency circle $r=r_R^{\mathrm{tan}}$, beyond which no continuation is asserted.
		
		\item \emph{Regularity-limited branches.} The lower phase reaches $y_{\mathrm{reg}}$ at a finite radius, where $v_R^-$ and its derivative have finite one-sided limits. Then we have the following alternatives.
		\begin{enumerate}[label=\textup{(\arabic*)}]
			\item Without a further continuation, the selected continuation terminates at this radius, and we obtain the maximal hypersurface $W_R^{\max}$.
			\item If the branch admits a finite-angle continuation, then the lower side continues through this radius and case \textup{(ii)} applies.
			\item If the branch admits a curvature-blow-up continuation of order $q>0$ and the singular positive-matching condition of Definition \ref{def:singular-positive-matching} holds, then we obtain a complete embedded two-ended catenoidal translator of global class $\mathcal C^{1,\frac{1}{q+1}}$, whose lower side is a second complete positive end. The continuation along these branches is unique and reaches a circle $r=r_R^{\mathrm{sing}}$, past which the phase joins the positive branch. At this circle, the two pieces of the lower side match with regularity $\mathcal C^{1,\frac{1}{q+1}}$ and no better, the unit normal extends continuously, and the meridional curvature diverges to $+\infty$ from both sides with the same order $q$. In addition, the hypersurface is of class $\mathcal C^3$ away from the singular circle and the neck circle, and smooth away from those circles wherever $\widehat\gamma$ is smooth along the traversed branches.
		\end{enumerate}
		\item \label{case6:small-necks}\emph{Small necks.} As $R\to0^+$, the nonconvex collars of both sides of the neck collapse toward the axis at the scale $R^{\frac{\rho_0}{1+\rho_0}}$. To be more precise, let $r_R$ denote either the entry radius $r_R^{\mathrm{ent}}$ of the upper side or the radius at which the lower side leaves $(I_-,g_-)$, where in case \textup{(ii)} we take instead the radius at which the lower phase crosses a fixed level $-\delta<0$. Then, in both cases, $r_R\to0$ with $\log r_R/\log R\to\frac{\rho_0}{1+\rho_0}$, and the height of the lower side at its radius tends to zero with $R$. In particular, the ratio $r_R^{\mathrm{ent}}/R$ tends to infinity.
	\end{enumerate}
	
	If the selected branch is axis-crossing without regular positive matching, or has infinite-angle zero output that is not fractional, the selected continuation yields $W_R^{\max}$ and we make no further assertion.
\end{theorem}

\begin{proof}
Firstly, we note that the local neck is constructed in Proposition \ref{prop:local-catenoidal-neck}, and that the statement about the upper side, namely the entry radius, the leading asymptotic order and the height constant, is Proposition \ref{prop:upper-side-global-continuation}. Moreover, the regularity at the neck circle is that of Proposition \ref{prop:local-catenoidal-neck}, namely $\mathcal C^{2,\varkappa}$ when $R$ is $\varkappa$-regular at the neck. In addition, when $g$ is of class $\mathcal C^1$ near $(R^{-1},0)$, the field of \eqref{eq:local-neck-phase-system} is of class $\mathcal C^1$. Therefore, $\theta_R$ is of class $\mathcal C^2$ and $\mathtt c_R$ of class $\mathcal C^3$ across the neck, and smooth when $g$ is smooth. Recall that, by Lemma \ref{lem:lower-side-neck-phase}, the reoriented lower side lies on $(I_-,g_-)$ near the neck. Now, we prove each case separately.

\emph{Case \textup{(i)}.} By Proposition \ref{prop:regular-lower-endpoint-dynamics}\textup{(ii)}, there is a finite radius at which $v_R^-=0$, with $(v_R^-)'\to x_0>0$. Then, by Proposition \ref{prop:axis-crossing-positive-matching}, we can continue the solution uniquely across $(x_0,0)$, with a positive slope and an increasing phase, up to the diagonal point $(y_*,y_*)$ at a finite radius $r_*$. Moreover, by Lemma \ref{lem:matching-c1-junctions}\textup{(ii)}, the concatenation of $g_{\mathrm{ax}}$ and $g_+$ is a $\mathcal C^1$ graph through the diagonal, and therefore the continuation through $r_*$ is unique.

Next, we continue the solution past $r_*$. By evaluating the phase identity \eqref{eq:graphical-phase-evolution} at $r_*$ with $y_v=y_*$ and $g_+(y_*)=y_*$, we obtain $ry_v'=\left(1+\alpha^{-1}(v_R^-)^2\right)y_*-y_*=\alpha^{-1}(v_R^-)^2y_*>0$. Therefore, the phase enters the interior of $I_+$ immediately after $r_*$. Then, by Proposition \ref{prop:positive-end-maximal-interval}, we have the global continuation, and by Proposition \ref{prop:positive-end-at-infinity} the leading order $v_R^-=r^\alpha(1+o(1))$. Moreover, Theorem \ref{thm:positive-bowl-asymptotics} refines it under its hypotheses. In addition, by repeating for the lower side the integration step in the proof of Proposition \ref{prop:upper-side-global-continuation}, we obtain the constant $C_R^-$.

Now, we show that $W_R$ is complete. Indeed, both sides are defined on $(R,\infty)$, and therefore the generating curve consists of two graphs attached at the neck. Since the arclength of each graph from the neck exceeds $r-R$, the curve is complete, and so is $W_R$.

Moreover, $W_R$ is of class $\mathcal C^3$ away from the neck circle, and smooth wherever the traversed branches are smooth. Indeed, a slope is of class $\mathcal C^2$ whenever it solves a graphical equation whose right-hand side is of class $\mathcal C^1$. On the one hand, the lower slope solves the graphical equation of the single branch obtained by concatenating $g_-$, $g_{\mathrm{ax}}$ and $g_+$, which is of class $\mathcal C^1$ across $y=0$ and across $y_*$ by the regular positive-matching condition and Lemma \ref{lem:matching-c1-junctions}. On the other hand, the upper slope solves the graphical equation of the branch $(J,G)$ of Proposition \ref{prop:upper-side-global-continuation}, which is of class $\mathcal C^1$ across $y=1$. Consequently, both slopes are of class $\mathcal C^2$, which proves the claim.

Finally, we prove that $W_R$ is embedded. We note that on $(R,r_*]$ the lower phase satisfies $Q(r,v_R^-)\leq y_*$, since $y_v\leq0$ while $v_R^-\leq0$ and $y_v\in(0,y_*]$ on the matching branch $g_{\mathrm{ax}}$ by Proposition \ref{prop:axis-crossing-positive-matching}. Moreover, beyond $r_*$ it lies in $I_+$, by the departure from the diagonal and Proposition \ref{prop:positive-end-maximal-interval}. In particular, it is less than $1$ at every radius. In contrast, the upper phase is at least $1$ on $(R,r_R^{\mathrm{ent}}]$ and lies in $I_+\subset(y_*,1)$ beyond $r_R^{\mathrm{ent}}$, by Proposition \ref{prop:upper-side-global-continuation}. Therefore, we have $Q(r,v_R^-)<Q(r,v_R^+)$ up to the larger of the radii $r_*$ and $r_R^{\mathrm{ent}}$, which is alternative \ref{emb:phase-separated} of Lemma \ref{lem:embeddedness-from-slope-comparison} there. Then, by continuity, we have $v_R^-<v_R^+$ on a short interval beyond that radius. Moreover, from that radius on both phases lie in $I_+$, and therefore alternative \ref{emb:common-branch} holds at every larger $r$ with the positive branch and any $r_0<r$ in that short interval. Consequently, we have $u_R^+>u_R^-$ on the common domain, and $W_R$ is embedded.

\emph{Case \textup{(ii)}.} We note that the argument covers both alternatives. Recall that, by Lemma \ref{lem:lower-side-neck-phase}, the reoriented lower side lies on $\overline g$ near the neck. In the second alternative, the solution follows $\overline g=g_-^{\mathrm{ext}}$ through $y_{\mathrm{reg}}$ by Proposition \ref{prop:branch-ivp}\textup{(i)}, because the extended branch is $\mathcal C^1$ on $(-\infty,0)$. Moreover, the phase increases through $y_{\mathrm{reg}}$ by Proposition \ref{prop:graphical-phase-evolution}\textup{(ii)}.

Now, in both alternatives, Proposition \ref{prop:finite-angle-lower-completeness} applies to $v_R^-$. Then, $r_+=\infty$ and $v_R^-<0$, and the slope decays as in the statement, with its flat refinement. This means that the lower side is a graph over $(R,\infty)$ whose arclength from the neck exceeds $r-R$, and consequently it is complete.

It remains to study the height $u_R^-$ and the regularity across $y_{\mathrm{reg}}$. Recall that $(u_R^-)'=v_R^-$. When the branch is flat, the behavior of $u_R^-$ in the statement follows by integrating the expansion of Proposition \ref{prop:finite-angle-lower-completeness}\textup{(iii)} from the neck. In general, if $\ell<-1$, we fix $\chi_2\in(\ell,-1)$ in Proposition \ref{prop:finite-angle-lower-completeness}\textup{(ii)}, and we see that $v_R^-$ is integrable at infinity, so that $u_R^-$ is bounded. If instead $\ell>-1$, we fix $\chi_1\in(-1,\ell)$, and the lower bound of the same proposition shows that $u_R^-\to-\infty$. We note that the first case occurs whenever $\rho_0\geq1$, because $\ell<-\rho_0$ by Proposition \ref{prop:general-lower-phase-alternative} for a regular branch and by Lemma \ref{lem:finite-angle-continuation-angle} for a continuation.

Next, across the circle where the phase passes $y_{\mathrm{reg}}$, the graphical equation has a $\mathcal C^1$ right-hand side. Therefore, $v_R^-$ is of class $\mathcal C^2$ and $W_R$ of class $\mathcal C^3$ there, and smooth if $\widehat\gamma$ is smooth near $(x_{\mathrm{reg}},y_{\mathrm{reg}})$. Finally, the embeddedness follows from alternative \ref{emb:phase-separated} of Lemma \ref{lem:embeddedness-from-slope-comparison}, since $Q(r,v_R^-)<0<Q(r,v_R^+)$ at every common radius.

\emph{Case \textup{(iii)}.} By Proposition \ref{prop:regular-lower-endpoint-dynamics}\textup{(iii)}, we have the finite radius $r_R^{\mathrm{tan}}$, the phase decay, and the limits $v_R^-\to0$ and $(v_R^-)'\to0$. Therefore, the lower graph extends to $r=r_R^{\mathrm{tan}}$ with a horizontal tangent, and $W_R^{\max}$ is the maximal continuation along the selected branch. We note that at that circle the right-hand side of the graphical equation behaves like $|y_v|^\sigma$ with $\sigma<1$, which is not Lipschitz. For this reason, a continuation beyond $r_R^{\mathrm{tan}}$ need not be unique. Finally, since $Q(r,v_R^-)<0<Q(r,v_R^+)$ at every common radius, alternative \ref{emb:phase-separated} of Lemma \ref{lem:embeddedness-from-slope-comparison} holds, and we conclude that $W_R^{\max}$ is embedded.

\emph{Case \textup{(iv)}.} By Proposition \ref{prop:regular-lower-endpoint-dynamics}\textup{(i)}, we have the finite exit radius and finite one-sided limits of $v_R^-$ and $(v_R^-)'$, which bound the maximal continuation $W_R^{\max}$ along the selected branch. This is alternative \textup{(1)}. Recall that a finite-angle continuation, which is alternative \textup{(2)}, is treated in case \textup{(ii)}. Now, we prove alternative \textup{(3)}. Let us assume that the branch admits a curvature-blow-up continuation. Then the extended branch is of class $\mathcal C^1$ through $(x_{\mathrm{reg}},y_{\mathrm{reg}})$, and by Proposition \ref{prop:branch-ivp}\textup{(i)} applied to $\left((-\infty,0),g_-^{\mathrm{ext}}\right)$ we can continue the solution uniquely through that point. Subsequently, by Proposition \ref{prop:curvature-blow-up-phase-asymptotics}, the phase reaches $y_v=0$ at a finite radius $r_R^{\mathrm{sing}}$, with $v_R^-\to0^-$, $(v_R^-)'\to+\infty$ and one-sided $\mathcal C^{1,\frac{1}{q+1}}$ regularity of the height. Moreover, under the singular positive-matching condition, by Proposition \ref{prop:singular-positive-matching-phase} we have the unique positive continuation with $v\to0^+$ at the same circle, whose phase reaches the diagonal, where it meets $g_+$. We note that the departure from the diagonal is as in case \textup{(i)}, because Lemma \ref{lem:matching-c1-junctions}\textup{(ii)} provides a $\mathcal C^1$ junction of $g_{\mathrm{sing}}$ and $g_+$ at $y_*$. Then, by Propositions \ref{prop:positive-end-maximal-interval} and \ref{prop:positive-end-at-infinity}, we obtain the second complete positive end.

Next, we study the regularity at the circle $r=r_R^{\mathrm{sing}}$. We note that both one-sided slopes vanish there. Therefore, the height is $\mathcal C^1$ there, with a horizontal tangent and a vertical, continuous unit normal. Moreover, the one-sided $\mathcal C^{1,\frac{1}{q+1}}$ bounds of Propositions \ref{prop:curvature-blow-up-phase-asymptotics} and \ref{prop:singular-positive-matching-phase} combine by the triangle inequality. Indeed, for $r_1<r_R^{\mathrm{sing}}<r_2$, we have
	\begin{align*}
		|v(r_2)-v(r_1)| &\leq |v(r_2)-v(r_R^{\mathrm{sing}})| + |v(r_R^{\mathrm{sing}})-v(r_1)| \\
		&\leq C(r_2-r_R^{\mathrm{sing}})^{\frac{1}{q+1}} + C(r_R^{\mathrm{sing}}-r_1)^{\frac{1}{q+1}} \leq 2C(r_2-r_1)^{\frac{1}{q+1}},
	\end{align*}
which establishes the matching regularity $\mathcal C^{1,\frac{1}{q+1}}$. In addition, this regularity is sharp. To see this, we note that by \eqref{eq:curvature-blow-up-phase-asymptotics} and $v=r_R^{\mathrm{sing}}\,y_v(1+o(1))$ on the left, $|v(r)|$ is bounded below by a positive multiple of $(r_R^{\mathrm{sing}}-r)^{\frac{1}{q+1}}$. Consequently, $v$ is not $\mathcal C^{0,\varkappa}$ at the circle for any $\varkappa>\frac{1}{q+1}$. Moreover, both one-sided meridional curvatures diverge to $+\infty$, with the same order $q$ and the constants $c_\pm$ of Definitions \ref{def:curvature-blow-up-continuation} and \ref{def:singular-positive-matching}. This means that the singularity is a failure of the second fundamental form to stay bounded, not of the immersion or of its normal. Now, away from the singular circle and the neck circle, the traversed branches $g_-^{\mathrm{ext}}$, $g_{\mathrm{sing}}$ and $g_+$ are $\mathcal C^1$, and the concatenation of $g_{\mathrm{sing}}$ and $g_+$ is $\mathcal C^1$ across $y_*$ by Lemma \ref{lem:matching-c1-junctions}\textup{(ii)}. Therefore, the hypersurface is of class $\mathcal C^3$ away from these two circles, and smooth wherever $\widehat\gamma$ is smooth along the traversed branches. The embeddedness follows as in case \textup{(i)}, with $g_{\mathrm{sing}}$ in place of $g_{\mathrm{ax}}$ in the phase comparison.

\emph{Case \textup{(v)}.}\label{proof:rotational-translators-v} Finally, we connect the small-neck estimates of Propositions \ref{prop:entry-radius-estimate} and \ref{prop:lower-side-collar-scale} and of \eqref{eq:lower-height-small-neck} with each alternative of the theorem. We notice that the upper side behaves in the same way in every alternative, whereas the lower side depends on the case, as we study below.
\begin{itemize}
	\item For the upper side, Proposition \ref{prop:entry-radius-estimate} shows that the entry radius $r_R^{\mathrm{ent}}$ tends to zero as $R\to0^+$, with $\log r_R^{\mathrm{ent}}/\log R\to\frac{\rho_0}{1+\rho_0}$.
	\item For the lower side in cases \textup{(i)}, \textup{(iii)} and \textup{(iv)}, the radius at which it leaves $(I_-,g_-)$ is the turning radius, the radius $r_R^{\mathrm{tan}}$ and the exit radius, respectively. In these cases the branch is axis-crossing, fractional or regularity-limited, and Proposition \ref{prop:lower-side-collar-scale}\textup{(ii)} shows that this radius has the same limit as the entry radius.
	\item In case \textup{(ii)}, the lower side is complete, and Proposition \ref{prop:lower-side-collar-scale}\textup{(ii)} provides the same limit, in both alternatives, for the radius at which the lower phase crosses a fixed level $-\delta<0$.
\end{itemize}
In every case, the height of the lower side at these radii tends to zero with $R$ by \eqref{eq:lower-height-small-neck}, concluding the proof.
\end{proof}

\begin{remark}\label{rem:catenoidal-radius-dependence}
	We remark that the type of the lower side in Theorem \ref{thm:rotational-translators} does not depend on $R$. The dependence on $R$ of the catenoidal translator lies only in its neck, in the radii of the circles where the profile changes branch, and in the height constants. The first of these radii is the entry radius, estimated in Proposition \ref{prop:entry-radius-estimate}, and the height constants are $C_R^+$ for the upper side and $C_R^-$ for the lower side in case \textup{(i)}. We refer the reader to Table \ref{tab:model-summary} for the case realized by each model of Subsection \ref{subsec:model-examples-branches}, and to Figures \ref{fig:catenoid-cases} and \ref{fig:catenoid-profiles} for the five behaviors.
\end{remark}

\begin{figure}[tb]
	\centering
	\captionsetup[sub]{labelfont=up,justification=raggedright,singlelinecheck=false}
	\begin{subfigure}[t]{0.19\textwidth}
		\centering
		\includegraphics[width=\textwidth,trim=0 3.4bp 269.16bp 4.2bp,clip]{fig_intro_catenoid_H.png}
		\caption{\scriptsize $\frac{H}{n-1}$, case \textup{(i)}.}
	\end{subfigure}\hfill
	\begin{subfigure}[t]{0.19\textwidth}
		\centering
		\includegraphics[width=\textwidth,trim=0 3.4bp 269.16bp 4.2bp,clip]{fig_intro_catenoid_Q21.png}
		\caption{\scriptsize $\mathcal Q_{2,1}$, case \textup{(ii)}.}
	\end{subfigure}\hfill
	\begin{subfigure}[t]{0.19\textwidth}
		\centering
		\includegraphics[width=\textwidth,trim=0 3.4bp 269.16bp 4.2bp,clip]{fig_intro_catenoid_frac.png}
		\caption{\scriptsize $\gamma_\eps$, case \textup{(iii)}.}
	\end{subfigure}\hfill
	\begin{subfigure}[t]{0.19\textwidth}
		\centering
		\includegraphics[width=\textwidth,trim=0 3.4bp 269.16bp 4.2bp,clip]{fig_intro_catenoid_tri.png}
		\caption{\scriptsize $\gamma_{\mathrm{tri}}$, case \textup{(iv)}.}
	\end{subfigure}\hfill
	\begin{subfigure}[t]{0.19\textwidth}
		\centering
		\includegraphics[width=\textwidth,trim=0 3.4bp 269.16bp 4.2bp,clip]{fig_intro_catenoid_root.png}
		\caption{\scriptsize normalized $S_3^{1/3}$, case \textup{(iv)}.}
	\end{subfigure}
	\caption{Generating curves of the catenoidal translators of Figure \ref{fig:catenoid-cases}, with the same parameters. The neck radius is dashed and, in \textup{(b)}, so is the asymptotic height of the lower end. The dots mark the circle where the lower side has a horizontal tangent plane in \textup{(a)}, and the circles $r=r_R^{\mathrm{tan}}$ in \textup{(c)}, $r=r_+$ in \textup{(d)} and $r=r_R^{\mathrm{sing}}$ in \textup{(e)}. Images courtesy of Ignacio McManus.}
	\label{fig:catenoid-profiles}
\end{figure}

\begin{remark}\label{rem:catenoidal-admissible-region}
	In contrast with the powers $H^p$, for which the translator equation \eqref{eq:translator} holds on the whole catenoidal translator, only the upper side satisfies it when the curvature function is defined on a G{\aa}rding cone $\Gamma_k$ with $k\geq2$. Indeed, the lower side lies outside $\Gamma_k$ for either normal by Lemma \ref{lem:lower-side-outside-garding}, and the same holds for the lower side of the annular solutions of de Lira--de Farias for $S_k$ with $k$ even, which has a negative angle function \cite[Remark 6]{lira2026translating}. The position of the neck circle with respect to $\overline\Gamma$ is decided in Lemma \ref{lem:neck-cone-boundary}, and when it lies on $\partial\Gamma$, as for $H^p$, the equation for $\gamma$ extended by zero holds at the neck circle, while otherwise only the signed equation holds near it. For $H^p$ we read the lower side with the reversed normal of Subsection \ref{subsec:reflected-lower-components}, and for $H$ the equation holds on the whole catenoid with a continuous normal. In this sense, P\'erez-Garc\'ia used the translating catenoids of the mean curvature flow as barriers to prove that a noncompact embedded translator with compact boundary cannot be contained in a round cylinder whose axis is parallel to the direction of translation \cite[Theorem 2.2]{perez-garcia2016translating}. For other curvature functions, we expect that the catenoidal translators can play this role in the viscosity sense, through the $k$-plurisubharmonicity of certain solutions and the Dirichlet duality of nonlinear potential theory \cite[Example 1.13 and Chapter 3]{cirant2023comparison}.
\end{remark}

\begin{remark}\label{rem:catenoidal-novelty}
	With respect to the catenoidal examples recalled in Subsection \ref{subsec:intro-background}, which Theorem \ref{thm:rotational-translators} contains in case \textup{(i)} and, for $k$ odd, in case \textup{(iv)}, the novelty of this theorem lies in the sharp regularity $\mathcal C^{1,\frac{1}{q+1}}$ at the singular sphere $\{r=r_R^{\mathrm{sing}}\}$ and $\mathcal C^{2,\frac{1}{p}}$ at the neck for the powers, in the uniqueness of the neck when the chart is not differentiable there, and in the independence of the type from $R$.
\end{remark}
\section{Uniqueness of the bowl and translators in vertical cylinders}\label{sec:uniqueness-catenoidal-barriers}

In this section we prove Theorems \ref{thm:main-b} and \ref{thm:main-d}, stated in the introduction, whose precise statements are Theorems \ref{thm:bowl-asymptotic-uniqueness}, \ref{thm:finite-height-bowl-uniqueness} and \ref{thm:vertical-cylinder-graphs}. To do so, we use the bowl-type solution $B$ of Definition \ref{def:bowl-type-solution} as a barrier, together with the tangency principle of \cite[Theorem 1.2]{Yo}. In part \textup{(a)} of Theorem \ref{thm:bowl-asymptotic-uniqueness} and in Theorem \ref{thm:vertical-cylinder-graphs} we move $B$ vertically until it touches the translator, whereas in part \textup{(b)} of Theorem \ref{thm:bowl-asymptotic-uniqueness} and in Theorem \ref{thm:finite-height-bowl-uniqueness} we use the reflection method over $\rr^n$ and over a ball. In addition, in Corollary \ref{cor:all-round-cylinders} we combine Theorem \ref{thm:vertical-cylinder-graphs} with the grim reaper cylinders of Remark \ref{rem:langford-rengaswami-condition}.

Throughout this section, $\gamma$ is an $\alpha$-positive curvature function that is smooth on $\Gamma_+$, and a $\gamma$-translator is an immersed hypersurface of class $\mathcal C^2$ with a continuous unit normal $\vec\nu$ satisfying $\gamma(\lambda)=\langle\vec\nu,\vec e_{n+1}\rangle$ with $\lambda\in\Gamma$ at every point. We note that, since $\gamma>0$ on $\Gamma$, the angle function is positive. Therefore, $\vec\nu$ is the upward normal and, near each point, the translator is the vertical graph of a function of class $\mathcal C^2$. Moreover, the isometries of $\rr^{n+1}$ that fix $\vec e_{n+1}$ preserve the angle function and the principal curvatures with respect to the upward normal, so they map $\gamma$-translators to $\gamma$-translators. In addition, a $\gamma$-translator is never totally geodesic, since $\lambda\in\Gamma$ and $0\notin\Gamma$, because otherwise homogeneity would give $\gamma(0)=0$, which contradicts $\gamma>0$ on $\Gamma$.

In contrast to Section \ref{sec:catenoidal-translators}, we do not assume $\alpha>\frac{1}{2}$. In fact, this bound enters only through the critical regime of Theorem \ref{thm:positive-bowl-asymptotics}, which we use in case \ref{case7:critical} of Theorem \ref{thm:bowl-asymptotic-uniqueness}\textup{(b)} and in Corollary \ref{cor:growth-entire-graphs}\textup{(iv)}. Moreover, case \ref{case7:nondegenerate} of Theorem \ref{thm:bowl-asymptotic-uniqueness}\textup{(b)} requires $\alpha\geq1$, and Theorem \ref{thm:bowl-asymptotic-uniqueness}\textup{(a)} requires only that the bowl-type solution be entire.

\begin{theorem}[Tangency principle, {\cite[Theorem 1.2]{Yo}}]\label{thm:tangency-principle}
	Let $\gamma:\Gamma\to[0,\infty)$ be smooth, symmetric, $\alpha$-homogeneous with $\alpha>0$, and elliptic on $\Gamma_+$. Let $\Sigma_1,\Sigma_2\subset\rr^{n+1}$ be connected $\gamma$-translators such that $\Sigma_1$ is strictly convex, that is, $\lambda\in\Gamma_+$ on $\Sigma_1$, and $\Sigma_2$ is convex, that is, $\lambda\in\overline{\Gamma}_+$ on $\Sigma_2$, and is not totally geodesic. Then the following hold.
	\begin{enumerate}[label=\textup{(\alph*)}]
		\item \textup{(}Interior tangency principle\textup{)} Assume that there exists an interior point $p\in\Sigma_1\cap\Sigma_2$ such that the tangent spaces coincide at $p$. If $\Sigma_1$ lies on one side of $\Sigma_2$, then $\Sigma_1=\Sigma_2$.
		\item \textup{(}Boundary tangency principle\textup{)} Assume that the boundaries $\partial\Sigma_i$ lie in the same hyperplane $\Pi$ and that each $\Sigma_i$ meets $\Pi$ transversally. If $\Sigma_1$ lies on one side of $\Sigma_2$ and there exists $p\in\partial\Sigma_1\cap\partial\Sigma_2$ at which $\Sigma_1$ and $\Sigma_2$ have the same tangent space, and so do $\partial\Sigma_1$ and $\partial\Sigma_2$, then $\Sigma_1=\Sigma_2$.
	\end{enumerate}
\end{theorem}

We notice that the second fundamental form of the graph of a function $u$ with the upward normal is $D^2u\,(1+|Du|^2)^{-\frac{1}{2}}$, and therefore convexity and strict convexity in Theorem \ref{thm:tangency-principle} mean $D^2u\geq0$ and $D^2u>0$, respectively. The proof of \cite[Theorem 1.2]{Yo} is inspired by the one of M{\o}ller for the mean curvature \cite{moller2014nonexistence}. More precisely, near $p$ the hypersurfaces $\Sigma_1$ and $\Sigma_2$ are graphs of two functions over their common tangent space, and the statement that $\Sigma_1$ lies on one side of $\Sigma_2$ means that these functions are ordered. Since $\Sigma_1$ is strictly convex and $\Sigma_2$ is convex, every convex combination of the two functions with positive weights has a strictly convex graph, whose principal curvatures lie in $\Gamma_+$, where $\gamma$ is elliptic. Consequently, their difference satisfies a linear elliptic equation, and the strong maximum principle and the Hopf lemma give parts \textup{(a)} and \textup{(b)}, respectively. This argument uses only two derivatives of the functions, and therefore it applies to the $\gamma$-translators of class $\mathcal C^2$ of this section.

\subsection{Uniqueness for bowl-type asymptotics}\label{subsec:bowl-uniqueness}

Mart\'in--Savas-Halilaj--Smoczyk characterized the bowl among the complete embedded translators of the mean curvature flow with a single end smoothly asymptotic to it \cite[Theorem A]{Paco_2014}. Their reflection argument was extended by Mart\'inez--Mart\'inez-Trivi\~no to $[\varphi,e_3]$-minimal surfaces in $\rr^3$ with a single end asymptotic to the corresponding bowl \cite[Theorem B]{martinez2022equilibrium}. In \cite[Theorem 1.4]{Yo} we adapted it to $\gamma$-translators that are graphs over a ball with an end asymptotic to a vertical cylinder, where the tangency principle deals with the possible intersections of the hypersurface with its reflections. Theorem \ref{thm:bowl-asymptotic-uniqueness} is the counterpart of \cite[Theorem A]{Paco_2014} for $\gamma$-translators. In this case the model is the entire bowl, Lemma \ref{lem:one-end-entire-graph} shows that the translator is an entire graph, and we compare the reflected caps with respect to the vertical order. For this purpose, we first fix the class of translators that we consider.

\begin{definition}\label{def:graphical-end}
	Let $M$ be a connected $n$-dimensional manifold without boundary, let $F:M\to\rr^{n+1}$ be a proper immersion of class $\mathcal C^2$ that is a $\gamma$-translator, and let $\mathrm{pr}:\rr^{n+1}\to\rr^n$ be the projection $\mathrm{pr}(x,x_{n+1})=x$. We say that $F$ has a graphical end if there is $R_0>0$ such that the following hold.
	\begin{enumerate}[label=\textup{(\alph*)}]
		\item $F(M)\cap\mathrm{pr}^{-1}(\{|x|\geq R_0\})$ is contained in the graph of a continuous function on $\{|x|\geq R_0\}$.
		\item $F(M)\cap\mathrm{pr}^{-1}(\{|x|<R_0\})$ is bounded.
	\end{enumerate}
\end{definition}

\begin{lemma}\label{lem:one-end-entire-graph}
	Let $F:M\to\rr^{n+1}$ have a graphical end. Then $\mathrm{pr}\circ F$ is a diffeomorphism of class $\mathcal C^2$ from $M$ onto $\rr^n$. Consequently, $F$ is an embedding, and $F(M)$ is the graph of a function $u\in\mathcal C^2(\rr^n)$.
\end{lemma}

\begin{proof}
	Firstly, we prove that $\mathrm{pr}\circ F$ is a local diffeomorphism. Let $p\in M$ and $w\in T_pM$ with $d(\mathrm{pr}\circ F)_p(w)=0$. Then $dF_p(w)=s\,\vec e_{n+1}$ for some $s\in\rr$, and since $dF_p(w)\in dF_p(T_pM)$, we have $0=\langle dF_p(w),\vec\nu\rangle=s\langle\vec e_{n+1},\vec\nu\rangle$. Since the angle function is positive, we obtain $s=0$, and therefore $w=0$ because $F$ is an immersion. Consequently, $d(\mathrm{pr}\circ F)_p$ is an isomorphism, and by the inverse function theorem $\mathrm{pr}\circ F$ is a local diffeomorphism of class $\mathcal C^2$.

	Next, we prove that $\mathrm{pr}\circ F$ is proper. For this purpose, let $C\subset\rr^n$ be compact. On the one hand, by condition \textup{(a)} of Definition \ref{def:graphical-end}, $F$ maps $(\mathrm{pr}\circ F)^{-1}(C\cap\{|x|\geq R_0\})$ into the graph of a continuous function over the compact set $C\cap\{|x|\geq R_0\}$, which is bounded. On the other hand, by condition \textup{(b)}, the image of $(\mathrm{pr}\circ F)^{-1}(C\cap\{|x|<R_0\})$ is bounded. Therefore, $F((\mathrm{pr}\circ F)^{-1}(C))$ is contained in a closed ball of $\rr^{n+1}$, and $(\mathrm{pr}\circ F)^{-1}(C)$ is a closed subset of the preimage of that ball under $F$, which is compact since $F$ is proper.

	Finally, $\mathrm{pr}\circ F$ is a proper local diffeomorphism of class $\mathcal C^2$ between nonempty connected manifolds, and therefore it is a covering map \cite[Proposition 4.46]{lee2013smooth}, whose proof uses only that the map is a proper local homeomorphism. Since $\rr^n$ is simply connected, the covering has a single sheet, and so $\mathrm{pr}\circ F$ is a diffeomorphism of class $\mathcal C^2$. Consequently, the function $u=x_{n+1}\circ F\circ(\mathrm{pr}\circ F)^{-1}$ is of class $\mathcal C^2$, $F(M)$ is its graph, and $F$ is an embedding, concluding the proof.
\end{proof}

We write $u_B:[0,\infty)\to\rr$ for the profile of the bowl-type solution when it is entire, $v_B=u_B'$, and $B=\{(x,u_B(|x|)):x\in\rr^n\}$. By Lemma \ref{lem:bowl-is-positive-end}, $v_B$ is a positive end. In part \textup{(b)} of Theorem \ref{thm:bowl-asymptotic-uniqueness} we consider the following cases.
\begin{enumerate}[label=\textup{(B\arabic*)},ref=\textup{(B\arabic*)}]
	\item \label{case7:nondegenerate} $\gamma$ is normalized $1$-nondegenerate and $\alpha\geq1$.
	\item \label{case7:power-growth} The positive end of the bowl is in the power-growth regime of Theorem \ref{thm:positive-bowl-asymptotics}\textup{(ii)}, with $d_\gamma\geq1$.
	\item \label{case7:critical} The positive end of the bowl is in the critical regime of Theorem \ref{thm:positive-bowl-asymptotics}\textup{(iii)}.
\end{enumerate}
In each case, the bowl is entire by Propositions \ref{prop:positive-end-maximal-interval} and \ref{prop:positive-end-at-infinity}.

\begin{theorem}\label{thm:bowl-asymptotic-uniqueness}
	Assume that the bowl-type solution for the $\gamma$-flow is entire. Let $F:M\to\rr^{n+1}$ have a graphical end, let $u\in\mathcal C^2(\rr^n)$ be the function of Lemma \ref{lem:one-end-entire-graph}, let $c\in\rr$, and let $\varphi(x)=u(x)-u_B(|x|)-c$. Then the following hold.
	\begin{enumerate}[label=\textup{(\alph*)}]
		\item If $F(M)$ is convex and $\varphi(x)\to0$ as $|x|\to\infty$, then $F(M)=B+c\,\vec e_{n+1}$.
		\item In cases \ref{case7:nondegenerate} through \ref{case7:critical}, if $F(M)$ is strictly convex and
		\begin{equation}\label{eq:weighted-asymptotics}
			|\varphi(x)|+|D\varphi(x)|=o\bigl(|x|^{-1}v_B(|x|)\bigr) \qquad \text{as }|x|\to\infty,
		\end{equation}
		with $|x|^{-1}$ replaced by $|x|^{1-2\alpha}$ in case \ref{case7:critical} when $\alpha>1$, then $F(M)$ is a vertical translate of $B$.
	\end{enumerate}
\end{theorem}

We note that condition \eqref{eq:weighted-asymptotics} is the analogue, for the bowl, of the $\mathcal C^k$-asymptotic condition to a vertical cylinder of \cite[Definition 4.1]{Yo}, and that its rate is the scale of the separation that the reflection produces at infinity (Lemma \ref{lem:start-moving-plane-at-infinity}). By Theorem \ref{thm:positive-bowl-asymptotics}, in cases \ref{case7:nondegenerate} and \ref{case7:power-growth} this rate is comparable to $|x|^{\alpha-1}$ and $|x|^{d_\gamma-1}$, respectively. Therefore, when $\alpha=1$ or $d_\gamma=1$, condition \eqref{eq:weighted-asymptotics} means that $\varphi$ and $D\varphi$ tend to zero, and part \textup{(b)} follows from part \textup{(a)}. Otherwise, the rate tends to infinity, and \eqref{eq:weighted-asymptotics} allows $\varphi$ to be unbounded and, in case \ref{case7:critical}, to grow exponentially.

The asymptotic condition of part \textup{(a)} does not single out $B$ by itself. Indeed, when Theorem \ref{thm:positive-bowl-asymptotics}\textup{(i)} applies, the slopes of both sides of a two-ended catenoidal translator of Theorem \ref{thm:rotational-translators}\textup{(i)} have the expansion of Theorem \ref{thm:positive-bowl-asymptotics}\textup{(i)}, which is also that of $v_B$. Then the difference of the slopes is $o(r^{-4\alpha})$, which is integrable, and therefore each side satisfies the asymptotic condition of Theorem \ref{thm:bowl-asymptotic-uniqueness}\textup{(a)}, with its own constant $c$. The translating catenoids of Clutterbuck--Schn\"urer--Schulze \cite[Lemma 2.3]{clutterbuck2007stability} behave in the same way for the mean curvature flow. However, these hypersurfaces are not convex, and they consist of two graphs over the complement of a ball. Moreover, they satisfy the translator equation with $\lambda\in\Gamma$ only sidewise, each side with respect to its own upward normal, and consequently the hypotheses of Theorem \ref{thm:bowl-asymptotic-uniqueness} exclude them.

\begin{proof}[Proof of Theorem \ref{thm:bowl-asymptotic-uniqueness}\textup{(a)}]
	By Lemma \ref{lem:one-end-entire-graph}, $F(M)$ is the graph of $u$, and $\varphi$ is continuous on $\rr^n$ and tends to zero at infinity. Let us assume by contradiction that $m=\sup_{\rr^n}\varphi$ is positive. Since $\varphi$ tends to zero at infinity, there is $R>0$ such that $\varphi<\frac{m}{2}$ on $\{|x|\geq R\}$, and therefore $m$ is the maximum of $\varphi$ on $\overline B_R$, attained at some $x_0$ with $|x_0|<R$. Then $u\leq u_B(|\cdot|)+c+m$ on $\rr^n$, with equality at $x_0$, and $D\varphi(x_0)=0$. In particular, $\Sigma_1=B+(c+m)\vec e_{n+1}$, which is strictly convex as recalled after Definition \ref{def:bowl-type-solution}, and $\Sigma_2=F(M)$ pass through $(x_0,u(x_0))$ with the same tangent space, and $\Sigma_1$ lies above $\Sigma_2$. Consequently, by Theorem \ref{thm:tangency-principle}\textup{(a)} we have $u=u_B(|\cdot|)+c+m$, that is, $\varphi\equiv m$, which contradicts $\varphi\to0$. Similarly, if $m'=\inf_{\rr^n}\varphi$ is negative, the same argument with $B+(c+m')\vec e_{n+1}$, which lies below $F(M)$, implies $\varphi\equiv m'$, leading again to a contradiction. Therefore, $\varphi\equiv0$, and $F(M)=B+c\,\vec e_{n+1}$.
\end{proof}

For the proof of part \textup{(b)}, we identify $M$ with $F(M)=\operatorname{graph}(u)$ by Lemma \ref{lem:one-end-entire-graph} and, after replacing $B$ by $B+c\,\vec e_{n+1}$, we assume that $c=0$. Moreover, we write $x=(x_1,x')\in\rr\times\rr^{n-1}$ for the points of $\rr^n$, and we use the notation of \cite[Section 3]{Paco_2014} and \cite[Section 4]{Yo}. That is, for $t\in\rr$, the reflecting hyperplane is $\Pi_t=\{x_1=t\}$, which we regard as a hyperplane of $\rr^n$ or as a vertical hyperplane of $\rr^{n+1}$, according to the context. The parts of $M$ on each side of it are $M_+(t)=M\cap\{x_1\geq t\}$ and $M_-(t)=M\cap\{x_1\leq t\}$, and $M_+^*(t)$ is the image of $M_+(t)$ under the reflection $(x_1,x',x_{n+1})\mapsto(2t-x_1,x',x_{n+1})$ across $\Pi_t$. Since the reflection fixes the vertical direction, $M_+^*(t)$ is the vertical graph over $\{x_1\leq t\}$ of $u_t^*(x)=u(x^t)$, where $x^t=(2t-x_1,x')$ is the reflection of $x\in\rr^n$ across $\Pi_t$. In the reflection argument we move $\Pi_t$ from far away toward $\Pi_0$ while $u_t^*\geq u$ holds on $\{x_1<t\}$. The model at infinity enters in the first step, which is the content of the following lemma, the counterpart of \cite[Claim 1]{Paco_2014}.

\begin{lemma}\label{lem:start-moving-plane-at-infinity}
	In cases \ref{case7:nondegenerate} through \ref{case7:critical}, let $u\in\mathcal C^2(\rr^n)$ be such that $u-u_B(|\cdot|)$ satisfies \eqref{eq:weighted-asymptotics} in place of $\varphi$. Then the following hold.
	\begin{enumerate}[label=\textup{(\alph*)}]
		\item There is $t_0>0$ such that $u_t^*>u$ on $\{x_1<t\}$ for every $t\geq t_0$.
		\item For every $\sigma>0$ and every compact $J\subset(0,\infty)$, there are $R>0$ and $\mu>0$ such that $u_t^*-u\geq\mu$ on $\{|x|\geq R,\ x_1\leq t-\sigma\}$ for every $t\in J$.
		\item For every $\tau_0>0$, there is $R'>0$ such that $\partial_1u>0$ on $\{x_1\geq\tau_0,\ |x'|\geq R'\}$.
		\item For every $a>0$, the function $u_a^*-u$ does not vanish identically on $\{x_1<a\}$.
	\end{enumerate}
\end{lemma}

\begin{proof}
	Let $\varphi=u-u_B(|\cdot|)$ and $\varpi(r)=r^{-1}v_B(r)$, with $r^{-1}$ replaced by $r^{1-2\alpha}$ in case \ref{case7:critical} when $\alpha>1$, so that the hypothesis reads $|\varphi(x)|+|D\varphi(x)|=o(\varpi(|x|))$. Firstly, we prove that there are $\kappa_0\in(0,1]$ and $r_1\geq1$ such that, for $r_1\leq s<s'$,
	\begin{subequations}\label{eq:bowl-separation}
		\begin{gather}
			\kappa_0\leq\varpi(s)\leq\kappa_0^{-1}\varpi(s'),\label{eq:bowl-separation-a}\\
			v_B(s)\geq\kappa_0\,s\,\varpi(s),\label{eq:bowl-separation-b}\\
			u_B(s')-u_B(s)\geq\kappa_0\min\{1,s'^2-s^2\}\,\varpi(s').\label{eq:bowl-separation-c}
		\end{gather}
	\end{subequations}

	In cases \ref{case7:nondegenerate} and \ref{case7:power-growth}, by Proposition \ref{prop:positive-end-at-infinity}\textup{(i)} and Theorem \ref{thm:positive-bowl-asymptotics}\textup{(ii)} we have $v_B(s)=As^d(1+o(1))$, with $A=1$ and $d=\alpha$ in case \ref{case7:nondegenerate}, and $A=A_\gamma$ and $d=d_\gamma$ in case \ref{case7:power-growth}. Then, after choosing $r_1\geq1$ sufficiently large, we have $\frac{A}{2}s^{d-1}\leq\varpi(s)\leq2As^{d-1}$ for $s\geq r_1$, and \eqref{eq:bowl-separation-a} follows from $d\geq1$, while \eqref{eq:bowl-separation-b} holds because $s\varpi(s)=v_B(s)$ and $\kappa_0\leq1$. For \eqref{eq:bowl-separation-c}, if $s'\leq2s$, then
	\begin{equation*}
		u_B(s')-u_B(s)=\int_s^{s'}\sigma\varpi(\sigma)\,d\sigma\geq\frac{A}{4}s^{d-1}\left(s'^2-s^2\right)\geq2^{-d-2}\varpi(s')\left(s'^2-s^2\right),
	\end{equation*}
	and if $s'>2s$, then $u_B(s')-u_B(s)\geq\int_{s'/2}^{s'}\sigma\varpi(\sigma)\,d\sigma\geq2^{-d-3}\varpi(s')$, because $s'\geq1$. Consequently, \eqref{eq:bowl-separation} holds in these cases with $\kappa_0=\min\{\frac{A}{2},2^{-d-3}\}$.

	In case \ref{case7:critical}, by Theorem \ref{thm:positive-bowl-asymptotics}\textup{(iii)}, after choosing $r_1\geq1$ sufficiently large, we have $\frac{c_\gamma}{2}s^{2\alpha-1}v_B\leq v_B'\leq2c_\gamma s^{2\alpha-1}v_B$ for $s\geq r_1$. Let $m=\max\{1,2\alpha-1\}$, so that $\varpi(s)=s^{-m}v_B(s)$. Then $(\log\varpi)'(s)\geq\frac{c_\gamma}{2}s^{2\alpha-1}-ms^{-1}$, which is positive for large $s$ because $2\alpha-1>0$. Therefore, after enlarging $r_1$, the function $\varpi$ is increasing on $[r_1,\infty)$, which implies \eqref{eq:bowl-separation-a} for every $\kappa_0\leq\varpi(r_1)$, while \eqref{eq:bowl-separation-b} follows from $s^m\geq s$. For \eqref{eq:bowl-separation-c}, since $v_B$ is increasing, we have
	\begin{equation*}
		u_B(s')-u_B(s)\geq v_B(s)(s'-s)\geq v_B(s)\frac{s'^2-s^2}{2s'},
	\end{equation*}
	and since $v_B(\sigma)\geq(2c_\gamma s'^m)^{-1}v_B'(\sigma)$ for $\sigma\in[s,s']$, we also have
	\begin{equation*}
		u_B(s')-u_B(s)\geq\frac{v_B(s')-v_B(s)}{2c_\gamma s'^m}.
	\end{equation*}
	On the one hand, if $v_B(s')\geq2v_B(s)$, the second bound implies $u_B(s')-u_B(s)\geq(4c_\gamma)^{-1}\varpi(s')$. On the other hand, if $v_B(s)>\frac{1}{2}v_B(s')$, the first bound implies $u_B(s')-u_B(s)\geq\frac{1}{4}(s'^2-s^2)\varpi(s')$. This proves \eqref{eq:bowl-separation} in case \ref{case7:critical} with $\kappa_0=\min\{\varpi(r_1),(4c_\gamma)^{-1},\frac{1}{4}\}$.

	Next, we fix the notation for the four items. By \eqref{eq:weighted-asymptotics}, for each $\varepsilon>0$ there is $R_\varepsilon\geq r_1$ such that $|\varphi(x)|+|D\varphi(x)|\leq\varepsilon\varpi(|x|)$ for $|x|\geq R_\varepsilon$. For $t>0$, we set $\rho=|x|$ and $\rho_t=|x^t|$, so that $u_t^*(x)-u(x)=u_B(\rho_t)-u_B(\rho)+\varphi(x^t)-\varphi(x)$ and $\rho_t^2-\rho^2=4t(t-x_1)$.

	We prove item \textup{(c)} first, since we use it in the proof of \textup{(a)}. Let $\tau_0>0$ and $R'=R_\varepsilon$ with $\varepsilon=\frac{\kappa_0\tau_0}{2}$. If $\tau\geq\tau_0$ and $|x'|\geq R'$, then $\rho=|(\tau,x')|\geq R'$, and by \eqref{eq:bowl-separation} we obtain
	\begin{equation*}
		\partial_1u(\tau,x')=v_B(\rho)\frac{\tau}{\rho}+\partial_1\varphi(\tau,x')\geq\kappa_0\tau\varpi(\rho)-\frac{\kappa_0\tau_0}{2}\varpi(\rho)>0.
	\end{equation*}

	Now, let $\sigma>0$, let $J\subset(0,\infty)$ be compact, let $t_1=\min J$, and let $R=R_\varepsilon$ with $\varepsilon=\frac{\kappa_0^2}{4}\min\{1,4t_1\sigma\}$, in order to prove \textup{(b)}. For $t\in J$ and $x$ with $x_1\leq t-\sigma$ and $\rho\geq R$, we have $\rho_t>\rho\geq R$ and $\rho_t^2-\rho^2\geq4t_1\sigma$. Then, by \eqref{eq:bowl-separation}, we obtain
	\begin{equation*}
		u_t^*(x)-u(x)\geq\kappa_0\min\{1,4t_1\sigma\}\varpi(\rho_t)-2\varepsilon\kappa_0^{-1}\varpi(\rho_t)\geq\frac{\kappa_0^2}{2}\min\{1,4t_1\sigma\},
	\end{equation*}
	which proves \textup{(b)} with $\mu=\frac{\kappa_0^2}{2}\min\{1,4t_1\sigma\}$.

	To prove \textup{(a)}, let $R_2=R_\varepsilon$ with $\varepsilon=\frac{\kappa_0^2}{4}$, and let $t\geq2R_2$ and $x$ with $x_1<t$. On the one hand, if $\frac{t}{2}<x_1<t$, then every $\tau\in[x_1,2t-x_1]$ satisfies $|(\tau,x')|\geq\tau>\frac{t}{2}\geq R_2\geq1$. Therefore, as in \textup{(c)}, we have $\partial_1u(\tau,x')\geq(\kappa_0\tau-\frac{\kappa_0^2}{4})\varpi(|(\tau,x')|)>0$, and
	\begin{equation*}
		u_t^*(x)-u(x)=\int_{x_1}^{2t-x_1}\partial_1u(\tau,x')\,d\tau>0.
	\end{equation*}
	On the other hand, if $x_1\leq\frac{t}{2}$, then $\rho_t^2-\rho^2\geq2t^2\geq1$ and $\rho_t\geq\sqrt2\,t$. When $\rho\geq R_2$, by \eqref{eq:bowl-separation} we have $u_t^*(x)-u(x)\geq\frac{\kappa_0}{2}\varpi(\rho_t)>0$. When $\rho<R_2$, we apply \eqref{eq:bowl-separation} with $s=R_2$ and $s'=\rho_t$, for which $\rho_t^2-R_2^2\geq7R_2^2\geq1$. Together with the monotonicity of $u_B$, we obtain $u_B(\rho_t)-u_B(\rho)\geq u_B(\rho_t)-u_B(R_2)\geq\kappa_0\varpi(\rho_t)$, and therefore
	\begin{equation*}
		u_t^*(x)-u(x)\geq\frac{1}{2}\left(u_B(\sqrt2\,t)-u_B(R_2)\right)+\frac{\kappa_0}{2}\varpi(\rho_t)-\max_{|y|\leq R_2}|\varphi(y)|-\frac{\kappa_0^2}{4}\varpi(\rho_t),
	\end{equation*}
	which is positive for large $t$, because $u_B(r)\to\infty$ as $r\to\infty$.

	Finally, we prove \textup{(d)}. Let $a>0$, and let us assume by contradiction that $u_a^*=u$ on $\{x_1<a\}$. Let $\varepsilon=\frac{\kappa_0^2}{4}$ and let $R>2a$ be so large that $R-2a\geq R_\varepsilon$ and $4a(R-a)\geq1$. At $x=(2a-R)\vec e_1$, for which $x^a=R\vec e_1$, the assumption reads $u_B(R)-u_B(R-2a)=\varphi((2a-R)\vec e_1)-\varphi(R\vec e_1)$. However, by \eqref{eq:bowl-separation} with $s=R-2a$ and $s'=R$, the left-hand side is at least $\kappa_0\varpi(R)$, while the right-hand side is at most $2\varepsilon\kappa_0^{-1}\varpi(R)=\frac{\kappa_0}{2}\varpi(R)$, leading to the desired contradiction.
\end{proof}

\begin{proof}[Proof of Theorem \ref{thm:bowl-asymptotic-uniqueness}\textup{(b)}]
	\emph{Step 1.} Firstly, we reduce the proof to a reflection symmetry. Indeed, it is enough to prove that $M$ is rotationally symmetric about the vertical axis. Let $\mathcal R$ be a rotation of $\rr^n$. The graph of $u\circ\mathcal R$ is the image of $M$ under a rotation of $\rr^{n+1}$ fixing $\vec e_{n+1}$, and therefore it is again a strictly convex $\gamma$-translator. Moreover, we have $u\circ\mathcal R=u_B(|\cdot|)+\varphi\circ\mathcal R$, with $|\varphi\circ\mathcal R|=|\varphi|\circ\mathcal R$ and $|D(\varphi\circ\mathcal R)|=|D\varphi|\circ\mathcal R$, and consequently $u\circ\mathcal R$ also satisfies \eqref{eq:weighted-asymptotics}. In particular, it suffices to prove that $M$ is symmetric with respect to $\Pi_0$, since the same argument in every horizontal direction then provides symmetry with respect to every vertical hyperplane through the origin, which implies rotational symmetry about the vertical axis.
	
	\emph{Step 2.} Next, we set up the reflection argument of \cite[Theorem A]{Paco_2014}, as in \cite[Theorem 1.4]{Yo}. Given $t\geq0$ and $x$ with $x_1<t$, we set
	\begin{equation*}
		h_t(x)=u_t^*(x)-u(x)=u(2t-x_1,x')-u(x_1,x'),
	\end{equation*}
	so that $h_t\geq0$ on $\{x_1<t\}$ means that the reflected cap $M_+^*(t)$ lies above $M_-(t)$, and we define
	\begin{equation*}
		\mathcal A = \left\{t\in[0,\infty):h_s\geq0\text{ on }\{x_1<s\}\text{ for every }s>t\right\}.
	\end{equation*}
	By Lemma \ref{lem:start-moving-plane-at-infinity}\textup{(a)}, $\mathcal A$ contains every sufficiently large $t$. Moreover, if $t\in\mathcal A$ and $t'\geq t$, then every $s>t'$ satisfies $s>t$, and therefore $h_s\geq0$ on $\{x_1<s\}$. In particular, we have $[t,\infty)\subset\mathcal A$.
	
	We note that $\mathcal A$ is closed. Indeed, let $t_j\in\mathcal A$ converge to $t_0$, so that $t_0\geq0$. For $s>t_0$, there is $j$ with $t_j<s$, since $t_j\to t_0<s$, and $t_j\in\mathcal A$ implies $h_s\geq0$ on $\{x_1<s\}$. Since $s>t_0$ is arbitrary, it follows that $t_0\in\mathcal A$.
	
	Consequently, $\mathcal A$ is nonempty and closed in $[0,\infty)$, and therefore it contains its infimum $a\geq0$, with $\mathcal A=[a,\infty)$. It remains to prove that $a=0$, and for this purpose we assume by contradiction that $a>0$ in Steps 3 and 4.
	
	\emph{Step 3.} We show that the reflected cap $M_+^*(a)$ lies strictly above $M_-(a)$. Since $a\in\mathcal A$, we have $h_s(x)\geq0$ for every $s>a$ and every $x$ with $x_1<a$. Then, by letting $s\to a^+$ and by the continuity of $u$, we obtain
	\begin{equation}\label{eq:order-at-a}
		h_a\geq0 \qquad \text{on }\{x_1<a\}.
	\end{equation}
	
	We claim that this order is strict, that is,
	\begin{equation}\label{eq:strict-order-at-a}
		h_a>0 \qquad \text{on }\{x_1<a\}.
	\end{equation}
	Indeed, let us assume by contradiction that $h_a(x_0)=0$ for some $x_0=(x_{0,1},x_0')$ with $x_{0,1}<a$, as in Figure \ref{fig:bowl-reflection}. By \eqref{eq:order-at-a}, the function $h_a$ attains its minimum at the interior point $x_0$, and therefore $Dh_a(x_0)=0$, that is, $Du_a^*(x_0)=Du(x_0)$. Next, we apply Theorem \ref{thm:tangency-principle}\textup{(a)} to the graphs $\Sigma_1=M_+^*(a)$ and $\Sigma_2=M_-(a)$ of $u_a^*$ and $u$ over the closed half-space $\{x_1\leq a\}$. They pass through $(x_0,u(x_0))$ with the same tangent space, and $\Sigma_1$ lies on one side of $\Sigma_2$ by \eqref{eq:order-at-a}. Then we obtain $\Sigma_1=\Sigma_2$, that is, $h_a\equiv0$ on $\{x_1<a\}$. Consequently, $M$ is symmetric with respect to $\Pi_a$, which contradicts Lemma \ref{lem:start-moving-plane-at-infinity}\textup{(d)} and proves \eqref{eq:strict-order-at-a}.
	
	\begin{figure}[!htb]
		\centering
		\includegraphics[width=0.4\textwidth,trim=10bp 10bp 10bp 10bp,clip]{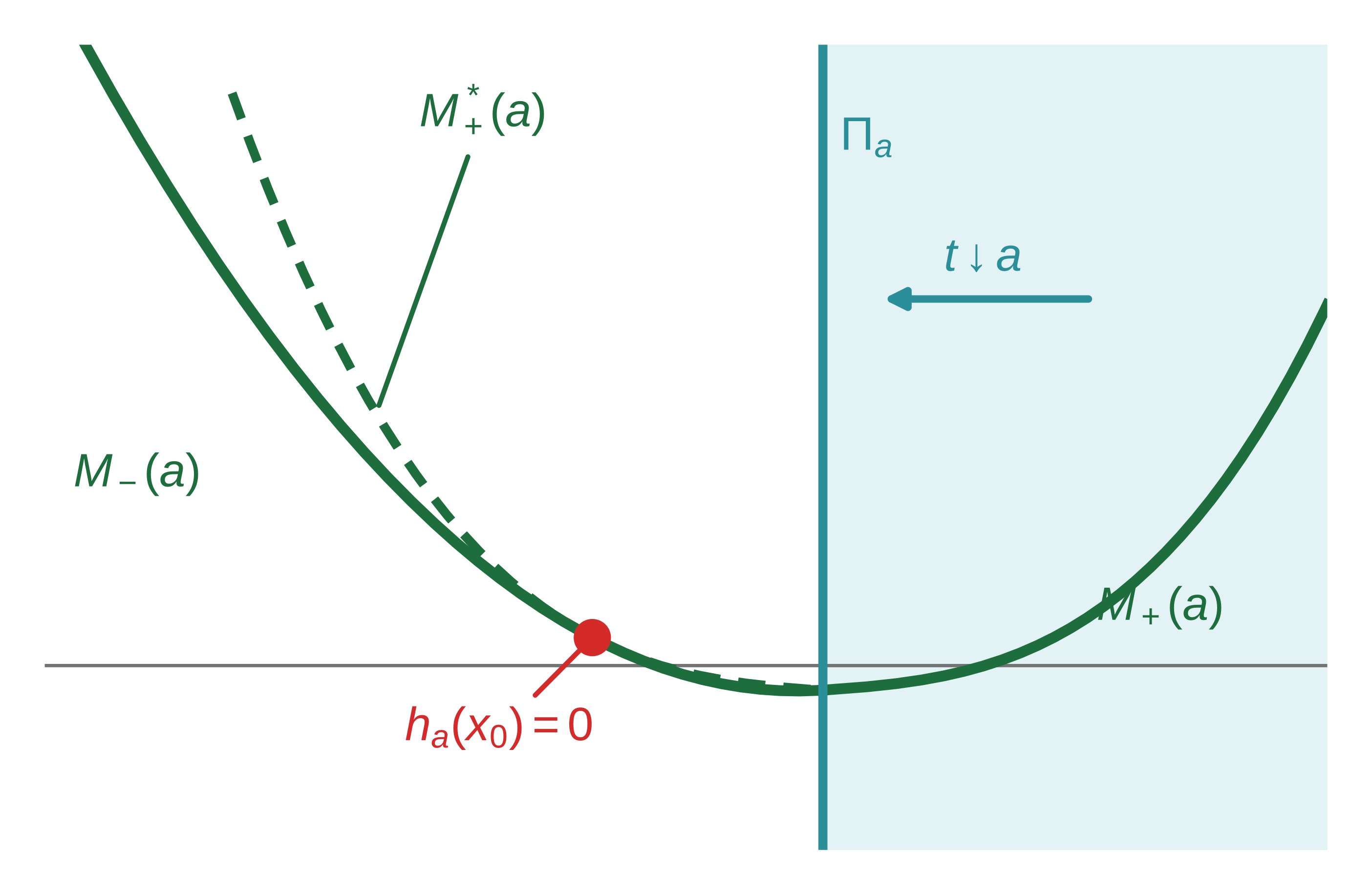}
		\caption{\scriptsize The plane $\Pi_a$ with $a>0$, the cap $M_+(a)$, its reflection $M_+^*(a)$ (dashed) above $M_-(a)$, and an interior contact point $x_0$ with $h_a(x_0)=0$. The arrow indicates that $\Pi_t$ moves as $t$ decreases to $a$. The proof of Theorem \ref{thm:bowl-asymptotic-uniqueness}\textup{(b)} excludes this configuration. Image courtesy of Ignacio McManus.}
		\label{fig:bowl-reflection}
	\end{figure}
	
	\emph{Step 4.} Now, we contradict the minimality of $a$ by showing that $a-\eps\in\mathcal A$ for some $\eps>0$. We follow the openness step of \cite[Claim 3]{Paco_2014}, where we replace the paraboloid asymptotics used there by the estimates \eqref{eq:positive-derivative-strip-near-a}--\eqref{eq:inside-ball-order-near-a} below. By \eqref{eq:order-at-a} and $a\in\mathcal A$, it suffices to prove that $h_s\geq0$ on $\{x_1<s\}$ for every $s\in(a-\eps,a)$, and we prove it separately in a strip next to $\Pi_s$, outside a large ball, and in the compact set that remains.
	
	Firstly, we study the derivative of $u$ across $\Pi_a$. We notice that the formula for $h_a$ defines a function on $\{x_1\leq a\}$ that vanishes on $\{x_1=a\}$ and, by \eqref{eq:strict-order-at-a}, is positive on $\{x_1<a\}$. Therefore, we have $\partial_1h_a(a,x')\leq0$, and since $\partial_1h_a(a,x')=-2\partial_1u(a,x')$, we obtain $\partial_1u(a,x')\geq0$ for every $x'$.
	
	Moreover, this inequality is strict by the boundary tangency principle. Indeed, the caps $\Sigma_1=M_+^*(a)$ and $\Sigma_2=M_-(a)$ of Step 3 have the common boundary $M\cap\Pi_a$, and they meet $\Pi_a$ transversally because their tangent spaces are never vertical. In addition, by \eqref{eq:order-at-a}, $\Sigma_1$ lies on one side of $\Sigma_2$, and on $\{x_1=a\}$ the functions $u_a^*$ and $u$ agree, their derivatives in the directions of $\Pi_a$ agree, and $\partial_1u_a^*(a,x')=-\partial_1u(a,x')$. Let us assume by contradiction that $\partial_1u(a,x')=0$ for some $x'$. Then $Du_a^*(a,x')=Du(a,x')$, and therefore the tangent spaces of $\Sigma_1$ and $\Sigma_2$ coincide at the boundary point $(a,x',u(a,x'))$, and so do those of $\partial\Sigma_1=\partial\Sigma_2$. Consequently, by Theorem \ref{thm:tangency-principle}\textup{(b)} we have $\Sigma_1=\Sigma_2$, that is, $u_a^*=u$ on $\{x_1<a\}$, which contradicts \eqref{eq:strict-order-at-a}. We conclude that
	\begin{equation}\label{eq:positive-normal-derivative-at-a}
		\partial_1u(a,x')>0 \qquad \text{for every }x'\in\rr^{n-1}.
	\end{equation}
	
	Next, we extend this positivity to a strip around $\Pi_a$. By Lemma \ref{lem:start-moving-plane-at-infinity}\textup{(c)} applied with $\tau_0=\frac{a}{2}$, there is $R'>0$ such that $\partial_1u(\tau,x')>0$ for $\tau\geq\frac{a}{2}$ and $|x'|\geq R'$. Moreover, by \eqref{eq:positive-normal-derivative-at-a} and continuity, $\partial_1u>0$ on a neighborhood of the compact set $\{a\}\times\{|x'|\leq R'\}$, which contains $\{|\tau-a|<2\delta,\ |x'|\leq R'\}$ for some $\delta\in(0,\frac{a}{8})$. Since $|\tau-a|<2\delta$ implies $\tau>\frac{3a}{4}$, we obtain
	\begin{equation}\label{eq:positive-derivative-strip-near-a}
		\partial_1u(\tau,x')>0 \qquad \text{for } |\tau-a|<2\delta \text{ and every }x'\in\rr^{n-1}.
	\end{equation}
	
	Next, with this $\delta$ fixed, we consider the points outside a large ball. By Lemma \ref{lem:start-moving-plane-at-infinity}\textup{(b)} applied with $\sigma=\delta$ and $J=[\frac{a}{2},\frac{3a}{2}]$, there are $R_1>0$ and $\mu_1>0$ such that
	\begin{equation}\label{eq:outside-ball-order-near-a}
		h_s(x)\geq\mu_1\quad \text{for }|s-a|\leq\frac{a}{2},\:|x|\geq R_1,\text{ and } x_1\leq s-\delta.
	\end{equation}

	It remains to consider, for $s<a$, the points of $\{x_1<s\}$ with $x_1\leq s-\delta$ and $|x|<R_1$. They lie in the compact set $K=\{x\in\rr^n: |x|\leq R_1+1,\ x_1\leq a-\frac{\delta}{2}\}$, where $h_a$ is positive by \eqref{eq:strict-order-at-a}, and therefore $\mu_2=\min_Kh_a>0$. In order to compare $h_s$ with $h_a$ on $K$, we notice that, for $|s-a|\leq\frac{a}{2}$ and $x\in K$, we have $h_s(x)-h_a(x)=u(2s-x_1,x')-u(2a-x_1,x')$, and both points lie in the closed ball of radius $3a+2R_1+2$ centered at the origin. Then, if $L$ is the maximum of $|Du|$ on that ball, the mean value theorem gives $\sup_K|h_s-h_a|\leq2L|s-a|$. Consequently, there is $\eps_1\in(0,\frac{a}{2}]$ such that $\sup_K|h_s-h_a|<\mu_2$ for $|s-a|<\eps_1$. In particular, we obtain
	\begin{equation}\label{eq:inside-ball-order-near-a}
		h_s>0 \quad \text{on }K\quad \text{for } |s-a|<\eps_1.
	\end{equation}

	Finally, we combine these three estimates. Let $\eps\in(0,\min\{\eps_1,\frac{\delta}{2}\})$, $s\in(a-\eps,a)$, and $x$ with $x_1<s$.
	\begin{itemize}
		\item If $s-\delta<x_1<s$, then $x_1$ and $2s-x_1$ lie in $(a-2\delta,a+2\delta)$, and by \eqref{eq:positive-derivative-strip-near-a} we have $h_s(x)=\int_{x_1}^{2s-x_1}\partial_1u(\tau,x')\,d\tau>0$.
		\item If $x_1\leq s-\delta$ and $|x|\geq R_1$, then $h_s(x)>0$ by \eqref{eq:outside-ball-order-near-a}.
		\item If $x_1\leq s-\delta$ and $|x|<R_1$, then $x_1<a-\frac{\delta}{2}$ and $x\in K$, and therefore $h_s(x)>0$ by \eqref{eq:inside-ball-order-near-a}.
	\end{itemize}
	Consequently, $h_s>0$ on $\{x_1<s\}$ for every $s\in(a-\eps,a)$, which implies $a-\eps\in\mathcal A$ and contradicts $a=\min\mathcal A$.
	
	\emph{Step 5.} Finally, we conclude the symmetry and identify $M$. By Steps 3 and 4, we have $a=0$, and therefore $h_s\geq0$ on $\{x_1<s\}$ for every $s>0$. Then, by letting $s\to0^+$, we obtain $u(-x_1,x')\geq u(x_1,x')$ for every $x_1<0$. The same argument in the direction $-\vec e_1$ provides the opposite inequality, and therefore $u(-x_1,x')=u(x_1,x')$, that is, $M$ is symmetric with respect to $\Pi_0$. By Step 1, $M$ is rotationally symmetric about the vertical axis, and $u(x)=\tilde u(|x|)$, where $\tilde u(r)=u(r\vec e_1)$ is of class $\mathcal C^2$ on $[0,\infty)$ and $\tilde u'(0)=\partial_1u(0)=0$, because $u$ is even in $x_1$.
	
	In particular, $M$ is the rotational graph of $\tilde u$ on $[0,\infty)$, and its principal curvatures lie in $\Gamma_+\subset\Gamma$ by strict convexity. Since $\tilde u\in\mathcal C^2([0,\infty))$, the slope $\tilde u'$ solves \eqref{eq:vertical-graph-translator} with $\tilde u'(0)=0$, and therefore $M$ is an entire bowl-type solution in the sense of Definition \ref{def:bowl-type-solution}. Consequently, by the uniqueness in \cite[Theorem 1.1 and Proposition 3.5]{rengaswami2024classification}, recalled after that definition, $\tilde u-u_B$ is constant, and $M$ is a vertical translate of $B$, concluding the proof.
\end{proof}

Theorem \ref{thm:bowl-asymptotic-uniqueness} requires the bowl to be entire. Now, we consider the bowls of finite height, which appear in part \textup{(c)} of Theorem \ref{thm:main-b}. When $R_\gamma<\infty$, we also write $B$ for the graph of $u_B$ over $\{|x|<R_\gamma\}$. By Propositions \ref{prop:positive-end-at-infinity}\textup{(ii)} and \ref{prop:finite-height-ends}\textup{(ii)}, we have that $R_\gamma<\infty$ and that $u_B$ is bounded if $\gamma$ is $1$-degenerate, $\alpha>1$ and $g_+(y)\geq cy^{-k}$ for large $y$ with $0<k<\alpha-1$. For instance, this holds for $\bigl(\frac{S_n}{S_{n-1}}\bigr)^\alpha$ with $\alpha>1$, whose positive branch $g_+(y)=\frac{y}{y-n+1}$ tends to $1$, and for the powers $K^p$ with $p>1$, whose positive branch $g_+(y)=y^{1-n}$ satisfies this condition with $k=n-1<np-1=\alpha-1$. At the boundary of such a bowl the tangent hyperplanes are vertical and $\gamma(\lambda)$ tends to zero, so we cannot apply the tangency principle there. However, the reflection method of part \textup{(b)} of Theorem \ref{thm:bowl-asymptotic-uniqueness} still works over a ball. The difference is that the reflected caps separate at the boundary of the ball, since the slope of the translator is infinite there, instead of separating at infinity by the asymptotics of Section \ref{sec:positive-branch-asymptotics}. In this way, we obtain the counterpart of \cite[Theorem 1.4]{Yo} for bowls of finite height.

\begin{theorem}\label{thm:finite-height-bowl-uniqueness}
	Let $M$ be a strictly convex $\gamma$-translator that is the graph of a function $u\in\mathcal C^2(B_R)\cap\mathcal C^0(\overline B_R)$, where $B_R=\{x\in\rr^n:|x|<R\}$, and assume that $|Du(x)|\to\infty$ as $|x|\to R$. Then $R_\gamma=R$, the profile $u_B$ is bounded, and $M$ is a vertical translate of $B$.
\end{theorem}

In the proof we keep the notation of the proof of Theorem \ref{thm:bowl-asymptotic-uniqueness}\textup{(b)}, and we replace the half-space $\{x_1<t\}$ by the reflection of $B_R\cap\{x_1>t\}$ across $\Pi_t$, that is,
\begin{equation*}
	\Omega_t=\{x\in\rr^n:\ x_1<t,\ x^t\in B_R\},\qquad t\in(0,R).
\end{equation*}
Since $|x^t|^2-|x|^2=4t(t-x_1)$, we have $\Omega_t\subset B_R$, and the boundary of $\Omega_t$ consists of $\{x_1=t\}\cap\overline B_R$ and of the reflection of $\partial B_R\cap\{x_1>t\}$, which lies in $B_R$. Moreover, $\Omega_t$ is convex, and $h_t$ is continuous on $\overline\Omega_t$ because $u\in\mathcal C^0(\overline B_R)$. The following lemma replaces Lemma \ref{lem:start-moving-plane-at-infinity}.

\begin{lemma}\label{lem:moving-plane-ball}
	Let $u$ be as in Theorem \ref{thm:finite-height-bowl-uniqueness}. Then the following hold.
	\begin{enumerate}[label=\textup{(\alph*)}]
		\item There is $t_0\in(0,R)$ such that $h_t>0$ on $\Omega_t$ for every $t\in[t_0,R)$.
		\item If $a\in(0,R)$ and $h_a\geq0$ on $\Omega_a$, then $h_a>0$ on $\partial\Omega_a\setminus\{x_1=a\}$.
		\item For every $\tau_0\in(0,R)$, there is $\rho\in(0,R)$ such that $\partial_1u>0$ on $\{x\in B_R:\ x_1\geq\tau_0,\ |x|>R-\rho\}$.
	\end{enumerate}
\end{lemma}

\begin{proof}
	Firstly, we show that the direction of $Du$ becomes radial at $\partial B_R$. By the convexity of $u$, the inequality $\langle Du(x),z-x\rangle\leq u(z)-u(x)\leq2\max_{\overline B_R}|u|$ holds for $x,z\in B_R$. Then, if $Du(x)\neq0$, by choosing $z=sR|Du(x)|^{-1}Du(x)$ and letting $s\to1^-$, we obtain $R|Du(x)|-\langle Du(x),x\rangle\leq2\max_{\overline B_R}|u|$. Since $|x|<R$, it follows that
	\begin{equation}\label{eq:ball-gradient-alignment}
		\left|\frac{Du(x)}{|Du(x)|}-\frac{x}{R}\right|^2\leq\frac{2}{R}\left(R-\frac{\langle Du(x),x\rangle}{|Du(x)|}\right)\leq\frac{4\max_{\overline B_R}|u|}{R|Du(x)|},
	\end{equation}
	where the right-hand side tends to zero as $|x|\to R$, because $|Du(x)|\to\infty$.

	Next, we prove item \textup{(c)}. Let $\tau_0\in(0,R)$. On the set of \textup{(c)}, the first coordinate of $\frac{x}{R}$ is at least $\frac{\tau_0}{R}$, and we compare the direction of $Du(x)$ with $\frac{x}{R}$. Indeed, for $x\in B_R$ with $Du(x)\neq0$, by the Cauchy--Schwarz inequality we have
	\begin{equation*}
		\partial_1u(x)=|Du(x)|\left(\frac{x_1}{R}+\left\langle\frac{Du(x)}{|Du(x)|}-\frac{x}{R},\vec e_1\right\rangle\right)\geq|Du(x)|\left(\frac{x_1}{R}-\left|\frac{Du(x)}{|Du(x)|}-\frac{x}{R}\right|\right).
	\end{equation*}
	Let $\rho\in(0,R)$ be such that $|Du(x)|\geq16R\tau_0^{-2}\max_{\overline B_R}|u|$ for $|x|>R-\rho$, which exists because $|Du(x)|\to\infty$ as $|x|\to R$. Then, for such $x$, \eqref{eq:ball-gradient-alignment} gives $\bigl|\frac{Du(x)}{|Du(x)|}-\frac{x}{R}\bigr|\leq\frac{\tau_0}{2R}$, and if in addition $x_1\geq\tau_0$, the previous inequality gives $\partial_1u(x)\geq\frac{\tau_0}{2R}|Du(x)|>0$.

	Now, we deduce \textup{(a)} from \textup{(c)}. Let $\rho$ be given by \textup{(c)} with $\tau_0=\frac{R}{2}$, and let $t_0=\max\{\frac{3R}{4},R-\frac{\rho}{2}\}$. By the choice of $t_0$, for $t\in[t_0,R)$ we have
	\begin{equation*}
		2t-R\geq2t_0-R\geq\max\left\{\frac{R}{2},R-\rho\right\}.
	\end{equation*}
	Given such $t$ and $x\in\Omega_t$, the first coordinate of $x^t\in B_R$ is $2t-x_1<R$, that is, $x_1>2t-R$. Moreover, every point $y$ of the segment joining $x$ and $x^t$ lies in $B_R$, since $B_R$ is convex, and satisfies $|y|\geq y_1\geq x_1>2t-R$. Consequently, $y_1>\frac{R}{2}$ and $|y|>R-\rho$, so that $\partial_1u(y)>0$ by \textup{(c)}, and therefore $h_t(x)=\int_{x_1}^{2t-x_1}\partial_1u(\tau,x')\,d\tau>0$.

	Finally, we prove \textup{(b)}, where the infinite slope of $u$ at $\partial B_R$ places the reflected boundary strictly above $M$. Let $p\in\partial\Omega_a\setminus\{x_1=a\}$. Then $p\in B_R$, and its reflection across $\Pi_a$ is a point $q\in\partial B_R$ with $q_1>a$. Let $\nu=R^{-1}q$ and $\nu^*=(-\nu_1,\nu')$. For small $d>0$, the point $q-d\nu$ lies in $B_R\cap\{x_1>a\}$, and its reflection across $\Pi_a$ is $p-d\nu^*$. Hence, $p-d\nu^*\in\Omega_a$ and
	\begin{equation*}
		h_a(p-d\nu^*)=h_a(p)-\bigl(u(q)-u(q-d\nu)\bigr)+\bigl(u(p)-u(p-d\nu^*)\bigr).
	\end{equation*}
	Let us assume by contradiction that $h_a(p)=0$. On the one hand, the last term is $O(d)$ as $d\to0^+$, because $p\in B_R$. On the other hand, the points $y=q-s\nu=(R-s)\nu$ satisfy
	\begin{equation*}
		\left|\frac{Du(y)}{|Du(y)|}-\nu\right|\leq\left|\frac{Du(y)}{|Du(y)|}-\frac{y}{R}\right|+\frac{s}{R},
	\end{equation*}
	where the right-hand side tends to zero as $s\to0^+$ by \eqref{eq:ball-gradient-alignment}. Consequently, $\langle Du(y),\nu\rangle=|Du(y)|\langle|Du(y)|^{-1}Du(y),\nu\rangle\to\infty$. Then, since $u$ is continuous at $q$, the identity $u(q)-u(q-d\nu)=\int_0^d\langle Du(q-s\nu),\nu\rangle\,ds$ implies that $d^{-1}\bigl(u(q)-u(q-d\nu)\bigr)\to\infty$ as $d\to0^+$. Therefore, $h_a(p-d\nu^*)<0$ for small $d>0$, which contradicts $h_a\geq0$ on $\Omega_a$.
\end{proof}

\begin{proof}[Proof of Theorem \ref{thm:finite-height-bowl-uniqueness}]
	The proof follows the five steps of the proof of Theorem \ref{thm:bowl-asymptotic-uniqueness}\textup{(b)}, with $\{x_1<t\}$ replaced by $\Omega_t$ and Lemma \ref{lem:start-moving-plane-at-infinity} replaced by Lemma \ref{lem:moving-plane-ball}. Thus, it suffices to describe the changes in each step.

	Firstly, Step 1 applies without changes, since the rotations of $\rr^{n+1}$ that fix $\vec e_{n+1}$ preserve $B_R$. Next, in Step 2 the set $\mathcal A$ consists of the $t\in[0,R)$ such that $h_s\geq0$ on $\Omega_s$ for every $s\in(t,R)$. It contains $t_0$ by Lemma \ref{lem:moving-plane-ball}\textup{(a)}, and it is closed in $[0,R)$ by the same argument. Let us assume by contradiction that $a=\min\mathcal A$ is positive.

	Now, in Step 3 we obtain $h_a\geq0$ on $\Omega_a$ as before, because every $x\in\Omega_a$ lies in $\Omega_s$ for $s>a$ close to $a$. Let us assume by contradiction that $h_a$ vanishes at a point of $\Omega_a$. Then the tangency principle, applied to the graphs of $u_a^*$ and $u$ over the convex set $\Omega_a$, gives $u_a^*=u$ on $\Omega_a$, and by continuity $h_a=0$ on $\partial\Omega_a$. This contradicts Lemma \ref{lem:moving-plane-ball}\textup{(b)}, which takes the place of Lemma \ref{lem:start-moving-plane-at-infinity}\textup{(d)}. Therefore, $h_a>0$ on $\Omega_a$, and together with Lemma \ref{lem:moving-plane-ball}\textup{(b)} we obtain
	\begin{equation}\label{eq:ball-strict-order-at-a}
		h_a>0\qquad\text{on }\overline\Omega_a\setminus\{x_1=a\}.
	\end{equation}

	Subsequently, Step 4 changes in three places.
	\begin{itemize}
		\item The proof of \eqref{eq:positive-normal-derivative-at-a} for every $x'$ with $(a,x')\in B_R$ is the same, since near each such point $\Omega_a$ coincides with $\{x_1<a\}$. In the boundary tangency principle we use the graphs of $u_a^*$ and $u$ over $\Omega_a\cup(\{x_1=a\}\cap B_R)$, so that their boundaries lie in $\Pi_a$.
		\item In the strip estimate, Lemma \ref{lem:moving-plane-ball}\textup{(c)} with $\tau_0=\frac{a}{2}$ replaces Lemma \ref{lem:start-moving-plane-at-infinity}\textup{(c)}, and the compact set $\{a\}\times\{|x'|\leq R'\}$ is replaced by $\{x_1=a,\ |x|\leq R-\rho\}$. Thus, \eqref{eq:positive-derivative-strip-near-a} holds at the points of $B_R$ for some $\delta\in(0,\min\{\frac{a}{8},R-a\})$. The bound $\delta<R-a$ ensures that the set $\overline\Omega_a\cap\{x_1\leq a-\frac{\delta}{2}\}$ below is nonempty.
		\item Since $\overline\Omega_a$ is compact, the exterior estimate \eqref{eq:outside-ball-order-near-a} is no longer needed. However, $\Omega_s$ need not be contained in $\Omega_a$, and we replace \eqref{eq:inside-ball-order-near-a} by a comparison of $h_s$ with $h_a$ at shifted points. By \eqref{eq:ball-strict-order-at-a}, $h_a$ has a positive minimum $\mu$ on the compact set $\overline\Omega_a\cap\{x_1\leq a-\frac{\delta}{2}\}$. Since $u$ is uniformly continuous on $\overline B_R$, there is $\eps\in(0,\frac{\delta}{2})$ such that $|u(y)-u(y')|<\mu$ for $y,y'\in\overline B_R$ with $|y-y'|<2\eps$. Given $s\in(a-\eps,a)$ and $x\in\Omega_s$ with $x_1\leq s-\delta$, the point $z=x+2(a-s)\vec e_1$ satisfies $z_1\leq a-\delta+(a-s)<a-\frac{\delta}{2}$ and $z^a=(2a-z_1,x')=(2s-x_1,x')=x^s$, which lies in $B_R$ because $x\in\Omega_s$. Then $z\in\Omega_a$ by the definition of $\Omega_a$, so that $h_a(z)\geq\mu$, and since $|z-x|=2(a-s)<2\eps$, we obtain $h_s(x)=u(z^a)-u(x)=h_a(z)+u(z)-u(x)>0$.
	\end{itemize}
	Then, as at the end of Step 4, let $s\in(a-\eps,a)$ and $x\in\Omega_s$. If $s-\delta<x_1<s$, the strip estimate gives $h_s(x)>0$, since the segment joining $x$ and $x^s$ lies in $B_R$. If $x_1\leq s-\delta$, the previous estimate gives $h_s(x)>0$. Consequently, $h_s>0$ on $\Omega_s$ for every $s\in(a-\eps,a)$, and $a-\eps\in\mathcal A$ contradicts the minimality of $a$.

	Finally, Step 5 gives the symmetry with respect to $\Pi_0$, since $x\in\Omega_s$ for small $s>0$ whenever $x\in B_R$ and $x_1<0$. Then, by Step 1, we have $u(x)=\tilde u(|x|)$, where $\tilde u\in\mathcal C^2([0,R))\cap\mathcal C^0([0,R])$, $\tilde u'(0)=0$, and $\tilde u'$ solves \eqref{eq:vertical-graph-translator}. Moreover, $|\tilde u'(r)|=|Du(r\vec e_1)|\to\infty$ as $r\to R$, so $[0,R)$ is the maximal interval of $\tilde u$, and the graph of $\tilde u$ is a bowl-type solution in the sense of Definition \ref{def:bowl-type-solution}. Therefore, by the uniqueness recalled after that definition, $R_\gamma=R$ and $\tilde u-u_B$ is constant. In particular, $u_B$ is bounded, and $M$ is a vertical translate of $B$.
\end{proof}

\begin{remark}\label{rem:finite-height-bowl}
	By \cite[Theorem 1.1]{rengaswami2024classification} and Proposition \ref{prop:finite-height-ends}, the bowl-type solution is entire or cylindrical when $\alpha\leq1$, and therefore no translator as in Theorem \ref{thm:finite-height-bowl-uniqueness} exists in this range. Moreover, the assumption that the domain is a ball cannot be dropped. Indeed, let $\gamma=K^p$ with $p>1$, and let $\Omega$ be a bounded convex domain with $\partial\Omega\in\mathcal C^{1,1}$ that satisfies an enclosing sphere condition at each boundary point. After a dilation of $\Omega$ that normalizes the velocity, there is a convex translating graph $u\in\mathcal C^\infty(\Omega)\cap\mathcal C^0(\overline\Omega)$ over $\Omega$, unique up to vertical translations, with $\det D^2u>0$ and $|Du|\to\infty$ at $\partial\Omega$ \cite[Theorem 2 and Section 4]{urbas1998complete}. In particular, this translator is strictly convex and has finite height, and when $\Omega$ is a nonround ellipsoid it is not rotationally symmetric. For a ball, the same uniqueness gives Theorem \ref{thm:finite-height-bowl-uniqueness} in the case of $K^p$.
\end{remark}

\subsection{Nonexistence of convex graphs in vertical cylinders}\label{subsec:catenoidal-barriers-bounded-domains}

In \cite[Corollary 1.3]{Yo} we proved that a cylindrical bowl-type solution over a ball excludes convex entire graphs. Theorem \ref{thm:vertical-cylinder-graphs} is the counterpart of this result when the bowl is entire, where the convex graph is confined to a vertical cylinder. Recall that, by \cite[Theorems 1.1, 1.2 and 1.4]{rengaswami2024classification}, the bowl is entire when $\gamma$ is $1$-nondegenerate or $\alpha\leq\frac{1}{2}$, and, for $1$-degenerate $\gamma$ with $\alpha>\frac{1}{2}$, when $g_+$ satisfies the decay condition $g_+(y)\leq Cy^{-k}$ with $k\geq2\alpha-1$ of Proposition \ref{prop:positive-end-at-infinity}\textup{(ii)}. In particular, this condition holds with $k=k_\gamma$ in the power-growth and critical regimes of Theorem \ref{thm:positive-bowl-asymptotics}, where $k_\gamma\geq2\alpha-1$.

\begin{theorem}\label{thm:vertical-cylinder-graphs}
	Assume that the bowl-type solution for the $\gamma$-flow is entire. Then there is no convex $\gamma$-translator that is the graph of a function $u\in\mathcal C^2(\Omega)$ over a bounded domain $\Omega\subset\rr^n$ with $u(x)\to+\infty$ as $x\to\partial\Omega$.
\end{theorem}

\begin{proof}
	Let us assume by contradiction that there exists such a translator $M=\operatorname{graph}(u)$, and let $w(x)=u_B(|x|)$ for $x\in\rr^n$, so that $B+t\,\vec e_{n+1}$ is the graph of $w+t$. Firstly, since $\Omega$ is bounded, $w$ is bounded on $\Omega$, and therefore $u-w\to+\infty$ at $\partial\Omega$. Then every sublevel set $\{u-w\leq c\}$ is a closed subset of $\Omega$ at a positive distance from $\partial\Omega$, and therefore compact. In particular, $u-w$ attains its minimum $t_1$ at some $x_0\in\Omega$, so that $u\geq w+t_1$ on $\Omega$, with equality at $x_0$, and $Du(x_0)=Dw(x_0)$.
	
	Next, let $\Sigma_1$ be the graph of $w+t_1$ over $\Omega$, which is strictly convex since $B$ is, and $\Sigma_2=M$. They pass through $(x_0,u(x_0))$ with the same tangent space, and $\Sigma_1$ lies on one side of $\Sigma_2$, since $u\geq w+t_1$. Consequently, by Theorem \ref{thm:tangency-principle}\textup{(a)} we have $\Sigma_1=\Sigma_2$, that is, $u=w+t_1$ on $\Omega$. However, $w+t_1$ is bounded on $\Omega$, while $u\to+\infty$ at $\partial\Omega$, leading to the desired contradiction.
\end{proof}

Mart\'in--Savas-Halilaj--Smoczyk observed that an entire translating graph of the mean curvature flow cannot grow uniformly slower or faster than $|x|^2$ \cite[Remark 3.1]{Paco_2014}. The barrier of Theorem \ref{thm:vertical-cylinder-graphs} provides the analogue for $\gamma$-translators in the three regimes of Theorem \ref{thm:positive-bowl-asymptotics}, as follows.

\begin{corollary}\label{cor:growth-entire-graphs}
	Let $u\in\mathcal C^2(\rr^n)$ be such that $\operatorname{graph}(u)$ is a convex $\gamma$-translator.
	\begin{enumerate}[label=\textup{(\roman*)}]
		\item \label{case7:growth-general} If the bowl-type solution is entire, then, with $w(x)=u_B(|x|)$, neither $u-w\to+\infty$ nor $w-u\to+\infty$ as $|x|\to\infty$.
		\item \label{case7:growth-nondegenerate} If $\gamma$ is normalized $1$-nondegenerate, there are no constants $C>0$, $r_0>0$ and $\zeta>0$ such that, for $|x|\geq r_0$, either
		\begin{equation*}
			u(x)\leq C|x|^{\zeta}, \quad\text{with } \zeta<\alpha+1, \text{ or } \zeta=\alpha+1 \text{ and } C<\frac{1}{\alpha+1},
		\end{equation*}
		or
		\begin{equation*}
			u(x)\geq C|x|^{\zeta}, \quad\text{with } \zeta>\alpha+1, \text{ or } \zeta=\alpha+1 \text{ and } C>\frac{1}{\alpha+1}.
		\end{equation*}
		\item \label{case7:growth-power} If the positive end of the bowl is in the power-growth regime of Theorem \ref{thm:positive-bowl-asymptotics}\textup{(ii)}, the conclusion of \ref{case7:growth-nondegenerate} holds with $\alpha+1$ replaced by $d_\gamma+1$ and $\frac{1}{\alpha+1}$ replaced by $\frac{A_\gamma}{d_\gamma+1}$.
		\item \label{case7:growth-critical} If the positive end of the bowl is in the critical regime of Theorem \ref{thm:positive-bowl-asymptotics}\textup{(iii)}, there are no constants $C>0$, $r_0>0$ and $\zeta>0$ such that, for $|x|\geq r_0$, either
		\begin{equation*}
			u(x)\leq e^{C|x|^{\zeta}}, \quad\text{with } \zeta<2\alpha, \text{ or } \zeta=2\alpha \text{ and } C<\frac{c_\gamma}{2\alpha},
		\end{equation*}
		or
		\begin{equation*}
			u(x)\geq e^{C|x|^{\zeta}}, \quad\text{with } \zeta>2\alpha, \text{ or } \zeta=2\alpha \text{ and } C>\frac{c_\gamma}{2\alpha}.
		\end{equation*}
	\end{enumerate}
\end{corollary}

\begin{proof}
	Firstly, we prove \ref{case7:growth-general}. Let us assume by contradiction that $u-w\to+\infty$ as $|x|\to\infty$. Then the sublevel sets of $u-w$ are closed and bounded, and therefore $u-w$ attains its minimum $t_1$ at some $x_0\in\rr^n$, so that $u\geq w+t_1$ on $\rr^n$, with equality at $x_0$, and $Du(x_0)=Dw(x_0)$. Consequently, $B+t_1\vec e_{n+1}$, which is strictly convex, and $\operatorname{graph}(u)$ pass through $(x_0,u(x_0))$ with the same tangent space, and the first lies below the second. Then, by Theorem \ref{thm:tangency-principle}\textup{(a)}, we have $u=w+t_1$, which contradicts $u-w\to+\infty$. Similarly, if $w-u\to+\infty$, then $u-w$ attains its maximum $t_2$ at some $x_1\in\rr^n$, so that $u\leq w+t_2$ with equality at $x_1$, and Theorem \ref{thm:tangency-principle}\textup{(a)} gives $u=w+t_2$, a contradiction.

	Next, we reduce \ref{case7:growth-nondegenerate} through \ref{case7:growth-critical} to \ref{case7:growth-general}. The bowl is entire in case \ref{case7:growth-nondegenerate} by \cite[Theorem 1.1]{rengaswami2024classification}. In cases \ref{case7:growth-power} and \ref{case7:growth-critical}, its slope $v_B$ is a positive end by Lemma \ref{lem:bowl-is-positive-end}, and $g_+(y)\leq2c_\gamma y^{-k_\gamma}$ for large $y$ with $k_\gamma\geq2\alpha-1$, so that the bowl is entire by Proposition \ref{prop:positive-end-at-infinity}\textup{(ii)}. Consequently, it suffices to determine the growth of $w$.
	
	In case \ref{case7:growth-nondegenerate}, since $\widetilde\gamma(0,1)=1$, \cite[Theorem 1.1]{rengaswami2024classification} also gives $w(x)=\frac{1}{\alpha+1}|x|^{\alpha+1}(1+o(1))$ as $|x|\to\infty$ for every $\alpha>0$. On the one hand, if $u\leq C|x|^\zeta$ for $|x|\geq r_0$ with $\zeta<\alpha+1$, or with $\zeta=\alpha+1$ and $C<\frac{1}{\alpha+1}$, then $w-u\geq\delta|x|^{\alpha+1}$ for some $\delta>0$ and all large $|x|$, and therefore $w-u\to+\infty$. On the other hand, if $u\geq C|x|^\zeta$ with $\zeta>\alpha+1$, or with $\zeta=\alpha+1$ and $C>\frac{1}{\alpha+1}$, then $u-w\geq\delta|x|^{\alpha+1}$ for large $|x|$, and therefore $u-w\to+\infty$. Both alternatives contradict \ref{case7:growth-general}. In case \ref{case7:growth-power}, by Theorem \ref{thm:positive-bowl-asymptotics}\textup{(ii)} we have $v_B(r)=A_\gamma r^{d_\gamma}(1+o(1))$, and therefore $w(x)=\frac{A_\gamma}{d_\gamma+1}|x|^{d_\gamma+1}(1+o(1))$, so that the previous argument applies with $d_\gamma+1$ and $\frac{A_\gamma}{d_\gamma+1}$ in place of $\alpha+1$ and $\frac{1}{\alpha+1}$.
	
	Finally, in case \ref{case7:growth-critical}, by Theorem \ref{thm:positive-bowl-asymptotics}\textup{(iii)} we have $\log v_B(r)=\frac{c_\gamma}{2\alpha}r^{2\alpha}+o(r^{2\alpha})$. Since $B$ is strictly convex, $v_B$ is increasing, and therefore $u_B(0)+v_B(r-1)\leq u_B(r)\leq u_B(0)+r\,v_B(r)$ for $r\geq1$, where both $\log v_B(r-1)$ and $\log(r\,v_B(r))$ equal $\frac{c_\gamma}{2\alpha}r^{2\alpha}+o(r^{2\alpha})$. Therefore, $\log w(x)=\frac{c_\gamma}{2\alpha}|x|^{2\alpha}+o(|x|^{2\alpha})$. On the one hand, if $u\leq e^{C|x|^\zeta}$ with $\zeta<2\alpha$, or with $\zeta=2\alpha$ and $C<\frac{c_\gamma}{2\alpha}$, then $e^{C|x|^\zeta}\leq w\,e^{-\delta|x|^{2\alpha}}$ for some $\delta>0$ and all large $|x|$, and consequently $w-u\geq w(1-e^{-\delta|x|^{2\alpha}})\to+\infty$. On the other hand, if $u\geq e^{C|x|^\zeta}$ with $\zeta>2\alpha$, or with $\zeta=2\alpha$ and $C>\frac{c_\gamma}{2\alpha}$, then $u\geq w\,e^{\delta|x|^{2\alpha}}$ for large $|x|$, and consequently $u-w\geq w(e^{\delta|x|^{2\alpha}}-1)\to+\infty$. Both alternatives contradict \ref{case7:growth-general}, concluding the proof.
\end{proof}

When $\gamma$ is $1$-homogeneous with $\gamma(1,0,\ldots,0)>0$, the grim reaper cylinders of Remark \ref{rem:langford-rengaswami-condition} play, in nonvertical round cylinders, the role that the bowl plays in vertical ones. Through \cite[Theorem 1.3]{Yo2}, they exclude strictly convex translators from nonvertical round cylinders. Moreover, if in addition $\gamma(0,1,\ldots,1)>0$ and $\gamma$ extends continuously by zero to $\partial\Gamma$, then $(0,1,\ldots,1)\in\Gamma$, and therefore $\gamma$ is $1$-nondegenerate. In this case, the bowl-type solution is entire by the opening of this subsection. Consequently, Theorem \ref{thm:vertical-cylinder-graphs} applies as well, and we obtain the following nonexistence result.

\begin{corollary}\label{cor:all-round-cylinders}
	Let $\gamma$ be an $\alpha$-positive curvature function with $\alpha=1$ that is smooth on $\Gamma$, extends continuously by zero to $\partial\Gamma$, and satisfies $\gamma(0,1,\ldots,1)>0$ and $\gamma(1,0,\ldots,0)>0$. Then the following hold.
	\begin{enumerate}[label=\textup{(\roman*)}]
		\item There is no convex $\gamma$-translator that is the graph of a function $u\in\mathcal C^2(\Omega)$ over a bounded domain $\Omega\subset\rr^n$ with $u(x)\to+\infty$ as $x\to\partial\Omega$.
		\item There is no complete, noncompact, properly embedded, strictly convex $\gamma$-translator $\Sigma$ with compact (possibly empty) boundary such that $\Sigma$ is contained in a nonvertical round cylinder.
	\end{enumerate}
\end{corollary}

\begin{proof}
	Part \textup{(i)} follows from Theorem \ref{thm:vertical-cylinder-graphs}, since $\gamma$ is $1$-nondegenerate and its bowl-type solution is entire.
	
	For part \textup{(ii)}, we rescale $\gamma$ so that $\gamma(1,0,\ldots,0)=1$, and therefore we can use \cite[Theorem 1.3]{Yo2}. More precisely, with $c=\gamma(1,0,\ldots,0)>0$, the function $c^{-1}\gamma$ meets the hypotheses of that theorem. We notice that, if $\Sigma$ is a $\gamma$-translator, then the dilated hypersurface $c^{-1}\Sigma$, whose principal curvatures are $c\lambda$ with the same unit normal, is a $(c^{-1}\gamma)$-translator by $1$-homogeneity. Since the properties listed in \textup{(ii)} are invariant under dilations, a translator as in \textup{(ii)} for $\gamma$ would provide one for $c^{-1}\gamma$, which contradicts \cite[Theorem 1.3]{Yo2}, concluding the proof.
\end{proof}

\begin{remark}\label{rem:vertical-cylinders-literature}
	M{\o}ller \cite[Corollaries 1.4 and 1.5]{moller2014nonexistence} and P\'erez-Garc\'ia \cite[Theorem 2.2]{perez-garcia2016translating} proved the results for the mean curvature flow recalled in Subsection \ref{subsec:intro-results} by comparing two $H$-translators with the maximum principle, without any convexity assumption, since the translator equation \eqref{eq:translator} for $H$ is elliptic at every $\lambda\in\rr^n$. In contrast, for a fully nonlinear curvature function the equation is elliptic only on $\Gamma$, and the convexity assumptions of Theorem \ref{thm:vertical-cylinder-graphs} and Corollary \ref{cor:all-round-cylinders} come from the tangency principle. Kim--Pyo established half-space theorems \cite{kim-pyo2021halfspace}, but their proof contains an error, and Alencar--Bessa--Silva Neto proved the corresponding statements for properly immersed translators \cite[Corollaries 1.1 and 1.2 and Remark 1.5]{alencar2026halfspace}. Without convexity, they also exclude complete $S_k$-translators from half-spaces and from complements of rotational cones under an ellipticity condition and growth conditions on $S_{k-1}$, and extend the bi-halfspace theorem of Chini--M{\o}ller \cite{chini-moller2021bihalfspace} to these translators \cite[Theorems 1.1--1.3]{alencar2026halfspace}. Finally, we remark that the assumption that the bowl is entire cannot be dropped. Indeed, by \cite[Theorem 1.3]{rengaswami2024classification} and Proposition \ref{prop:finite-height-ends}, a bowl-type solution that is not entire is cylindrical when $\alpha\leq1$, and then it is itself a strictly convex graph over a ball that diverges at its boundary.
\end{remark}

\section{Further directions}\label{sec:conclusions}

We end with five open problems suggested by the results above.

\begin{op}\label{op:bowl-entire}
	Theorem \ref{thm:bowl-asymptotic-uniqueness} characterizes an entire bowl-type solution $B$ among the convex $\gamma$-translators with a graphical end asymptotic to a vertical translate of $B$. We notice that the horizontal translates of $B$ limit the asymptotic classes in which such a characterization can hold. Indeed, let us assume that $v_B(r)=Ar^{d}(1+o(1))$ with $A,d>0$, as in the nondegenerate case and in the power-growth regime of Theorem \ref{thm:positive-bowl-asymptotics}. Then, for $a\in\rr^n\setminus\{0\}$, the horizontal translate $B+(a,0)$ is the graph of $u_B(|x-a|)$, and by the mean value theorem and the power growth of $v_B$ we have
	\begin{equation*}
		u_B(|x-a|)-u_B(|x|)=-v_B(|x|)\frac{\langle a,x\rangle}{|x|}+o\bigl(v_B(|x|)\bigr) \qquad \text{as }|x|\to\infty,
	\end{equation*}
	whose right-hand side is $-|a|\,v_B(|x|)(1+o(1))$ along the ray $\{ta:t>0\}$, and therefore is not $o(v_B(|x|))$. Thus, $o(v_B(|x|))$ is the natural class for a characterization of $B$, and we propose the following question.
	\begin{enumerate}[label=\textup{(\alph*)},ref=\textup{(\alph*)}]
		\item\label{op1:asymptotic-class} Is a convex $\gamma$-translator with a graphical end a vertical translate of $B$ whenever its function $u$ of Lemma \ref{lem:one-end-entire-graph} satisfies $u(x)-u_B(|x|)-c=o(v_B(|x|))$ for some $c\in\rr$?
	\end{enumerate}
	Theorem \ref{thm:bowl-asymptotic-uniqueness} answers this question positively when $u-u_B(|\cdot|)-c\to0$ and, for strictly convex translators with $d\geq1$, under \eqref{eq:weighted-asymptotics}, which here requires $u-u_B(|\cdot|)-c$ and its gradient to be $o(v_B(|x|)/|x|)$. Moreover, the answer is also positive for $K^{\frac{1}{n+2}}$. Indeed, its translator equation is $\det D^2u=1$, and J\"orgens, Calabi and Pogorelov proved that the entire convex solutions of this equation are quadratic polynomials \cite{jorgens1954losungen,calabi1958improper,pogorelov1972improper}. Therefore, the entire strictly convex translating graphs are the elliptic paraboloids. In particular, the bowl is $u_B=\frac{|x|^2}{2}$, so that $o(v_B(|x|))=o(|x|)$, and this class cannot be enlarged because of the horizontal translates of $B$ and the nonround paraboloids.
	
	Without an asymptotic assumption, the uniqueness of the bowl may fail. For the mean curvature flow, it holds in $\rr^3$ \cite[Theorem 1.1]{wang2011convex} and fails in higher dimensions by the theorems of Wang recalled in Subsection \ref{subsec:intro-results}, and it is recovered under noncollapsing and two-convexity conditions \cite[Theorem 1.2]{haslhofer2015uniqueness}, \cite[Theorem 1.1]{brendle-choi2019uniqueness}, \cite[Theorem 1.1]{brendle-choi2021uniqueness}, \cite[Theorem 1.1]{cogo-lynch-vicanek2023rotational}. In $\rr^3$, Hershkovits proved the uniqueness without convexity assumptions for properly embedded translators whose asymptotic shrinker at $-\infty$ is a cylinder \cite[Theorem 3]{hershkovits2020translators}. In $\rr^4$, the two-convexity cannot be removed, since Choi--Haslhofer--Hershkovits proved that the noncollapsed translators are the round bowl, the product of a line with the bowl soliton of $\rr^3$, and the oval bowls constructed by Hoffman--Ilmanen--Mart\'in--White \cite[Theorem 1.1]{choi-haslhofer-hershkovits2023translators}. The powers $K^p$ of the Gauss curvature are $1$-degenerate, and for $p\leq\frac{1}{2}$ their bowl-type solution is entire. In this case uniqueness fails at least for $p<\frac{1}{n+1}$, since Jian--Wang constructed infinitely many smooth convex translators of the $K^p$-flow without rotational symmetry, which are entire graphs in this range \cite[Theorem 1.2]{jian2014existence}.
	
	The case $n=2$ leads to two further questions.
	\begin{enumerate}[resume*]
		\item\label{op1:nondegenerate-R3} For a $1$-nondegenerate $\alpha$-positive curvature function, is every strictly convex entire graphical $\gamma$-translator in $\rr^3$ a translate of the bowl-type solution?
	\end{enumerate}
	Here strict convexity excludes the grim reaper cylinder of $\gamma$, which is an entire convex graph when $\alpha\leq\frac{1}{2}$ and $e_1\in\Gamma$ (Remark \ref{rem:langford-rengaswami-condition}).
	
	In the $1$-degenerate case, the theorem of J\"orgens, Calabi and Pogorelov \cite{jorgens1954losungen,calabi1958improper,pogorelov1972improper} shows that the answer to the analogue of \ref{op1:nondegenerate-R3} is negative for $K^{1/4}$, whose translators in $\rr^3$ include the nonround paraboloids. Moreover, Choi--Choi--Kim classified the translators of $K^p$ in $\rr^3$ for $p<\frac{1}{4}$ \cite{choi2025translating}. The degree of $K^{1/4}$ is $\alpha=\frac{1}{2}$, which is the threshold of Rengaswami's classification \cite[Theorems 1.2--1.4]{rengaswami2024classification}.
	\begin{enumerate}[resume*]
		\item Is the bowl-type solution of a $1$-degenerate curvature function with $\alpha>\frac{1}{2}$ unique in the sense of \ref{op1:nondegenerate-R3} whenever it is entire?
	\end{enumerate}
\end{op}

\begin{op}\label{op:flying-wings}
	For the mean curvature flow, the flying wings are complete strictly convex graphical translators over slabs of width greater than $\pi$, asymptotic to two tilted grim reaper cylinders \cite[Theorem 3]{bourni2020existence}, where rotational symmetry in the remaining $n-1$ variables is assumed when $n\geq3$. Their existence was conjectured by Hamilton (see \cite{bourni2022ancient}), and it is recalled together with their uniqueness in $\rr^3$ in Subsection \ref{subsec:intro-related}.
	
	We ask whether flying wings exist for other curvature functions. Let $\gamma$ be an $\alpha$-positive curvature function with $\alpha>\frac{1}{2}$ and $e_1\in\Gamma$, and let $w$ be the width of its grim reaper cylinder (Remark \ref{rem:langford-rengaswami-condition}). By rotating and dilating this cylinder, we obtain tilted grim reaper cylinders of $\gamma$ in slabs of every width greater than $w$. Under which hypotheses on $\gamma$ does each of these slabs contain a complete strictly convex $\gamma$-translator that is asymptotic to the corresponding tilted grim reaper cylinder? An answer would extend the conjecture of Hamilton to other curvature functions.
\end{op}

\begin{op}
	The main result of Clutterbuck--Schn\"urer--Schulze is the asymptotic stability of the bowl under the mean curvature flow for continuous perturbations. More precisely, a graphical solution whose continuous initial datum converges to the bowl at spatial infinity, with no decay rate imposed, converges uniformly to the translating bowl as $t\to\infty$ \cite[Theorem 1.1]{clutterbuck2007stability}. Recently, asymptotic stability has been proved in \cite{park2026stability} for other families of translators of the mean curvature flow, including the bowl under different hypotheses.
	
	Let $\gamma$ have an entire bowl-type solution $B$, and let $u_0\in\mathcal C^2(\rr^n)$ be strictly convex with $u_0-u_B(|\cdot|)$ bounded and $D(u_0-u_B(|\cdot|))\to0$ at infinity. Does the graphical $\gamma$-flow starting at $u_0$ exist for all time and converge, after subtracting the vertical translation of $B$, to a vertical translate of $B$?
	
	The proof of Clutterbuck--Schn\"urer--Schulze uses the two sides of the translating catenoids, shifted vertically, as barriers on the complement of a ball \cite[Lemma 4.1]{clutterbuck2007stability}. The natural candidates for this role are the catenoidal translators of Theorem \ref{thm:rotational-translators}, whose upper sides are asymptotic to vertical translates of $B$ under the hypotheses of Theorem \ref{thm:positive-bowl-asymptotics}\textup{(i)}. However, they are not convex and, in contrast with the mean curvature, they satisfy the translator equation in general only sidewise (Remark \ref{rem:catenoidal-admissible-region}), so that the tangency principle of Theorem \ref{thm:tangency-principle} does not apply to them. We expect that maximum and tangency principles for such barriers can be obtained in the viscosity sense, through the potential theory mentioned in Remark \ref{rem:catenoidal-admissible-region}, although the translator equation depends on $Du$ through the induced metric of the graph, and without a convexity constraint this dependence lies outside the monotonicity conditions required in \cite[Section 1.8]{cirant2023comparison}.
\end{op}

\begin{op}\label{op:boundary-nondegenerate}
	Rengaswami's classification includes the boundary $1$-nondegenerate curvature functions among the nondegenerate ones, and for them the bowl-type solution is entire, with leading term $\frac{|x|^{\alpha+1}}{\alpha+1}$ after normalization \cite[Definition 3 and Theorem 1.1]{rengaswami2024classification}. However, the expansion of Theorem \ref{thm:positive-bowl-asymptotics}\textup{(i)} need not hold in this class. For instance, for $\gamma=\frac{1}{n-1}\bigl(H+S_n^{1/n}\bigr)$ on $\Gamma_+$, the positive branch vanishes to order $n$ at $y=1$, and a formal computation, which agrees with a numerical integration, gives $v_B(r)=r-\frac{1}{n-1}r^{1-\frac{2}{n}}(1+o(1))$ as $r\to\infty$. This raises two questions.
	\begin{enumerate}[label=\textup{(\alph*)}]
		\item How does the vanishing order of $g_+$ at $y=1$ determine the asymptotics of the positive ends of a boundary $1$-nondegenerate curvature function?
		\item Do the catenoidal translators of Theorem \ref{thm:rotational-translators} and the uniqueness of Theorem \ref{thm:bowl-asymptotic-uniqueness} extend to this class?
	\end{enumerate}
\end{op}

\begin{op}\label{op:finite-height-domains}
	Let $\gamma$ be an $\alpha$-positive curvature function whose bowl-type solution has finite height, as for $\bigl(\frac{S_n}{S_{n-1}}\bigr)^\alpha$ with $\alpha>1$ (Proposition \ref{prop:finite-height-ends}). By Theorem \ref{thm:finite-height-bowl-uniqueness}, a strictly convex $\gamma$-translator given by the graph of a function over a ball, continuous up to the boundary and with $|Du|\to\infty$ there, is a translate of the bowl. For $K^p$ with $p>1$, the theorem of Urbas recalled in Remark \ref{rem:finite-height-bowl} gives more, since it provides such a translator, unique up to vertical translations, over a dilation of every domain as in that remark. This leads to the following question. Does the same hold for every curvature function as above?
\end{op}

\appendix

\section{Elementary properties of $Q(r,s)$ and $y_v(r)$}\label{app:q-eta-properties}

In this appendix we prove the properties of the function $Q(r,s)=s\,r^{-1}(1+s^2)^{-\beta}$ and of the phase $y_v(r)=Q(r,v(r))$ used along Sections \ref{sec:barrier-method}, \ref{sec:positive-branch-asymptotics} and \ref{sec:catenoidal-translators}. Recall that $\beta=\frac{\alpha-1}{2\alpha}$, which implies that $1-2\beta=\alpha^{-1}>0$. These properties are proved in \cite[Lemmas 2.7--2.8]{rengaswami2024classification} for $m>0$, and the identity \eqref{eq:wm-derivative-level} is the inverse of \cite[(6)]{rengaswami2024classification}. For the sake of completeness, we include their proofs, which also cover the case $m\leq0$, since $Q(r,-s)=-Q(r,s)$.

\begin{lemma}\label{lem:level-set-control}
	Let $r>0$.
	\begin{enumerate}[label=\textup{(\roman*)}]
		\item $Q(r,\cdot)$ is an odd, strictly increasing homeomorphism of $\rr$ onto $\rr$ with $Q(r,0)=0$, and $Q(r,s)\to\pm\infty$ as $s\to\pm\infty$. In particular, if $r$ ranges in a compact subset of $(0,\infty)$, then $y_v\to0$ if and only if $v\to0$, and $y_v\to\pm\infty$ if and only if $v\to\pm\infty$.
		\item For every $m\in\rr$ there is a unique $w_m(r)\in\rr$ with $Q(r,w_m(r))=m$. The function $w_m$ is of class $\mathcal C^1$ on $(0,\infty)$, $w_m(r)$ has the sign of $m$, $w_0=0$, and $w_a(r)<w_b(r)$ whenever $a<b$. Moreover, for every function $v$ defined at $r$,
		\begin{equation*}
			w_a(r)\leq v(r)\leq w_b(r) \quad\Longleftrightarrow\quad a\leq y_v(r)\leq b.
		\end{equation*}
	\end{enumerate}
\end{lemma}

\begin{proof}
	Firstly, by differentiating $Q$, we obtain
	\begin{equation}\label{eq:q-partial-s-section}
		\partial_sQ(r,s) = \frac{1+(1-2\beta)s^2}{r(1+s^2)^{\beta+1}} >0,
	\end{equation}
	where the sign comes from $1-2\beta>0$, and
	\begin{equation}\label{eq:q-partial-r-section}
		\partial_rQ(r,s) = -r^{-1}Q(r,s).
	\end{equation}
	In particular, $Q(r,\cdot)$ is strictly increasing for each $r>0$. Moreover, whenever $s\neq0$ we have $Q(r,s)=r^{-1}\operatorname{sgn}(s)|s|^{\frac{1}{\alpha}}(1+s^{-2})^{-\beta}$, and therefore $Q(r,s)\to\pm\infty$ as $s\to\pm\infty$. Consequently, since $Q(r,\cdot)$ is continuous, odd and strictly increasing, it is an order-preserving bijection of $\rr$ onto $\rr$ with $Q(r,0)=0$. This proves the existence and uniqueness of $w_m(r)$, its sign, $w_0=0$, the ordering of $w_a$ and $w_b$ for $a<b$, and the equivalence between $w_a(r)\leq v(r)\leq w_b(r)$ and $a\leq Q(r,v(r))\leq b$.
	
	Next, we prove the equivalences in \textup{(i)}. Since $Q(r,s)=r^{-1}Q(1,s)$ and $Q(1,\cdot)$ is a homeomorphism fixing $0$, they hold when $r$ stays in a compact subset of $(0,\infty)$, where $r^{-1}$ is bounded above and below.
	
	Finally, since $Q$ is smooth on $(0,\infty)\times\rr$ and $\partial_sQ>0$, the function $w_m$ is of class $\mathcal C^1$ by the implicit function theorem, concluding the proof.
\end{proof}

The derivative of $w_m$ is computed in the same way. Indeed, by differentiating $Q(r,w_m(r))=m$, and by using \eqref{eq:q-partial-s-section} and \eqref{eq:q-partial-r-section} together with $Q(r,w_m(r))=m$ once more, we obtain
\begin{equation}\label{eq:wm-derivative-level}
	w_m'(r) = (1+w_m(r)^2)^{\beta+1} \frac{m}{1+(1-2\beta)w_m(r)^2}.
\end{equation}

For the phase $y_v$ we have two expansions at large slopes. Indeed, when $v(r)>0$ the expression of $Q$ above reads $y_v(r)=v(r)^{\frac{1}{\alpha}}r^{-1}(1+v(r)^{-2})^{-\beta}$. In particular, as $v(r)\to+\infty$ we obtain
\begin{equation}\label{eq:eta-large-slope-asymptotic}
	y_v(r) = v(r)^{\frac{1}{\alpha}}r^{-1}(1+o(1)),
\end{equation}
and, by the expansion $(1+\tau)^{-\beta}=1-\beta\tau+O(\tau^2)$ as $\tau\to0$,
\begin{equation}\label{eq:eta-large-slope-second-order}
	y_v(r) = v(r)^{\frac{1}{\alpha}}r^{-1} \left(1-\beta v(r)^{-2}+O(v(r)^{-4})\right).
\end{equation}

Finally, we bound $Q$ from both sides for $s\geq1$. In this range we have $s^2\leq1+s^2\leq2s^2$, and therefore $2^{-|\beta|}s^{-2\beta}\leq(1+s^2)^{-\beta}\leq2^{|\beta|}s^{-2\beta}$ for either sign of $\beta$. Then, by the definition of $Q$, we obtain
\begin{equation}\label{eq:q-large-slope-bound}
	\frac{2^{-|\beta|}}{r}\,s^{\frac{1}{\alpha}} \leq Q(r,s) \leq \frac{2^{|\beta|}}{r}\,s^{\frac{1}{\alpha}}, \qquad s\geq1,\ r>0.
\end{equation}

\section*{Acknowledgements}
We thank S. Rengaswami for the joint work \cite{RengaswamiTorresSantaella2023Bowl} and for many fruitful discussions related to the questions developed in the present article. We also thank Mariel S\'aez for her support and for many helpful conversations, Silvia Vald\'es for helpful conversations, the referees for their careful reading and suggestions, and Ignacio McManus for preparing the figures of this article.

\section*{Statements and Declarations}\label{sec:declarations}
\textbf{Funding:}
No funding was received to assist with the preparation of this manuscript.

\textbf{Competing interests:}
The author has no relevant financial or non-financial interests to disclose.

\textbf{Data availability:}
Data sharing is not applicable to this article, as no datasets were generated or analyzed during the current study.

\textbf{Use of large language models:}
The author used the large language model Claude (from Anthropic) as an assistant in the preparation of this article, to review and restructure the exposition, to numerically verify the coefficients of the asymptotic expansions, and to write the code that computed the curves and the models shown in the figures, which were prepared by Ignacio McManus. It was also used to organize the structure of the responses to the referees, which the author revised in their final version. The results, their proofs, and the references were developed and verified by the author, who assumes full responsibility for their content.

	
\end{document}